\documentclass[a4paper,11pt]{article}
\usepackage{amsmath,amssymb,amsthm,graphicx}
\usepackage{hyperref,url,enumitem,bm,esint}
\usepackage{color}
\usepackage{comment}
\usepackage{algorithm}
\usepackage{algpseudocode}
\usepackage[T1]{fontenc}
\usepackage{graphicx}
\usepackage{multirow}
\usepackage{makecell} 
\usepackage{subcaption}
\usepackage{mwe}
\usepackage{overpic}
\usepackage[margin=30mm]{geometry}
\usepackage{booktabs}

\newcommand{\bb}{\mathbf b}
\newcommand{\bx}{\mathbf x}
\newcommand{\bn}{\mathbf n}
\newcommand{\br}{\mathbf r}
\newcommand{\bA}{\mathbf A}

\newcommand{\anna}[2][black]{\textcolor{#1}{#2}}

\theoremstyle{plain}

\numberwithin{equation}{section}

\def\cl {\nonumber \\}
\def\el {\nonumber }

\title{A finite element framework for simulating residential burglary in realistic urban geometries}

\author{Baoli Hao$^1$, Kamrun Mily$^2$, Annalisa Quaini$^2$, and Ming Zhong$^2$}

\date{}
\begin{document}
\maketitle

\begin{center}
$^1$ Department of Applied Mathematics, Illinois Institute of Technology, 10 West 35th St, Chicago, IL 60616 \\ {\tt bhao2@hawk.illinoistech.edu}

\noindent $^2$ Department of Mathematics, University of Houston, 3551 Cullen Blvd, Houston, TX 77204, USA \\
{\tt kmily@cougarnet.uh.edu, aquaini@central.uh.edu, mzhong3@central.uh.edu}
\end{center}

\begin{abstract}
We consider a partial differential equation (PDE) model to predict residential burglary derived from a probabilistic agent-based model through a mean-field limit operation. The PDE model is a nonlinear, coupled system of two equations in two variables (attractiveness of residential sites and density of criminals), similar to the Keller–Segel model for aggregation based on chemotaxis.  Unlike previous works, which applied  periodic boundary conditions, we enforce boundary conditions that arise naturally from the variational formulation of the PDE problem, i.e., the starting point for the application of a finite element method.  These conditions specify the value of the normal derivatives of the system variables at the boundary.  For the numerical solution of the PDE problem discretized in time and space, we propose a scheme that decouples the computation of the attractiveness from the computation of the criminal density at each time step, resulting in the solution of two linear algebraic systems per iteration. Through numerous numerical tests, we demonstrate the robustness and computational efficiency of this approach. Leveraging the flexibility allowed by the finite element method, we show results for spatially heterogeneous model parameters and a realistic geometry (city of Chicago). \anna{The paper includes a discussion of future perspectives to build multiscale, ``multi-physics'' models that can become a tool for the community.
The robust and efficient code developed for this paper, which is shared open-source, is intended as the solid base for this broader research program.}
\end{abstract}

\noindent \emph{Key words}: Discrete and continuum crime models; Finite Element approximation; Partitioned methods; Pattern formation.  


\section{Introduction}\label{sec:intro}

\begin{figure}[htb!]
\centering
\begin{overpic}[width=0.45\textwidth, grid=false]{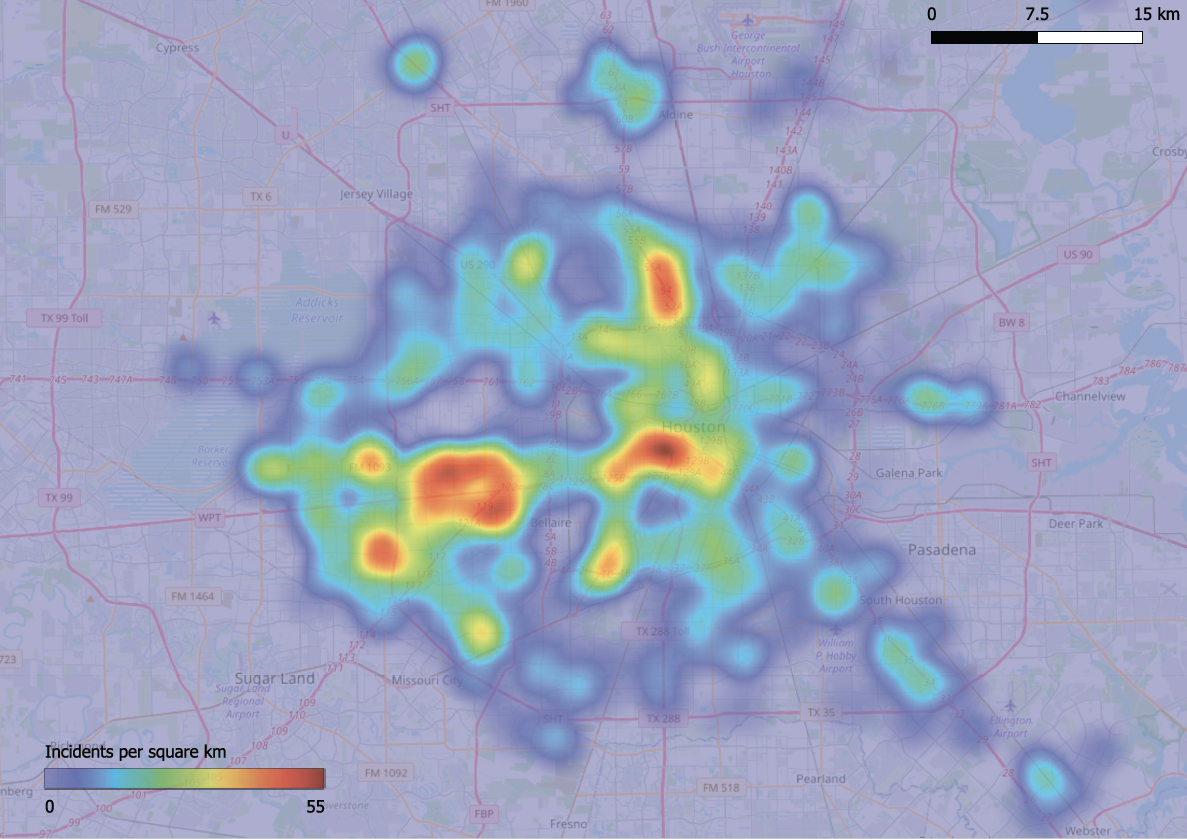}
\put(20,73){January to March 2024}
\end{overpic} \quad 
\begin{overpic}[width=0.45\textwidth, grid=false]{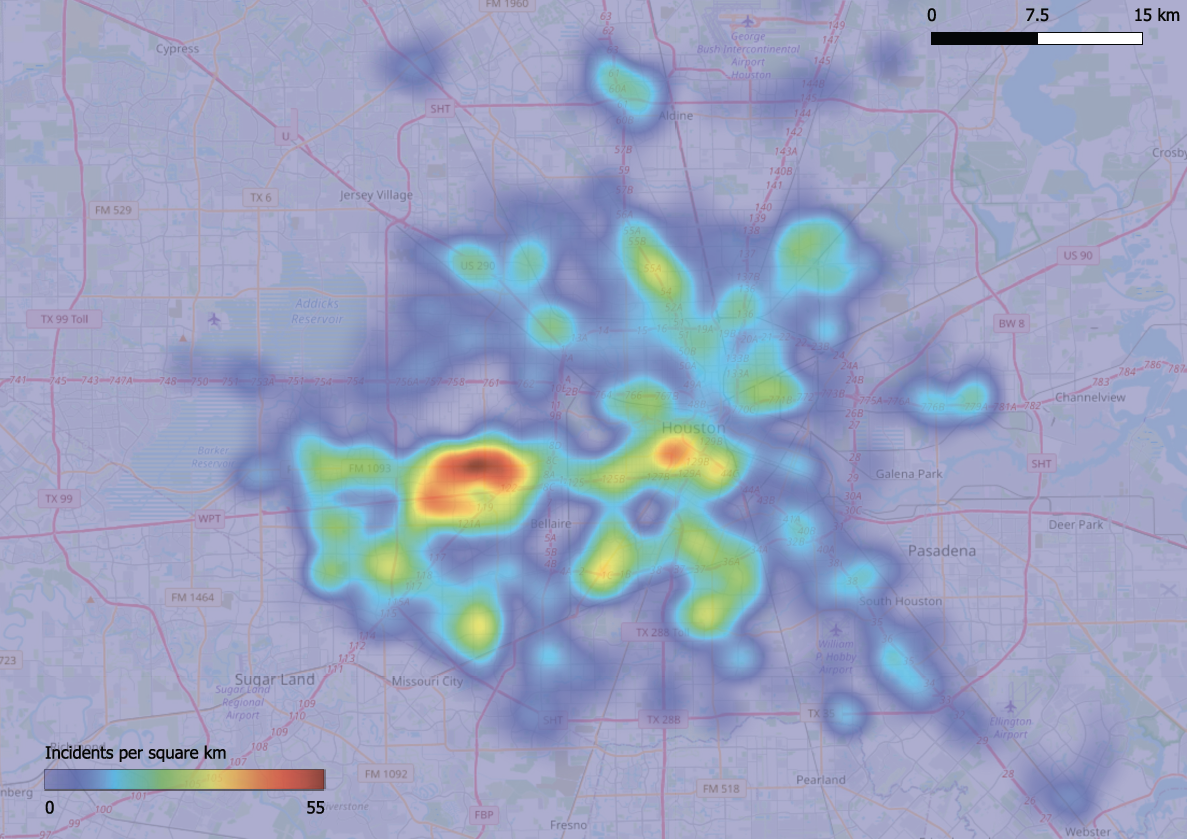}
\put(28,73){April to June 2024}
\end{overpic} \\
\vskip .8cm
\begin{overpic}[width=0.45\textwidth, grid=false]{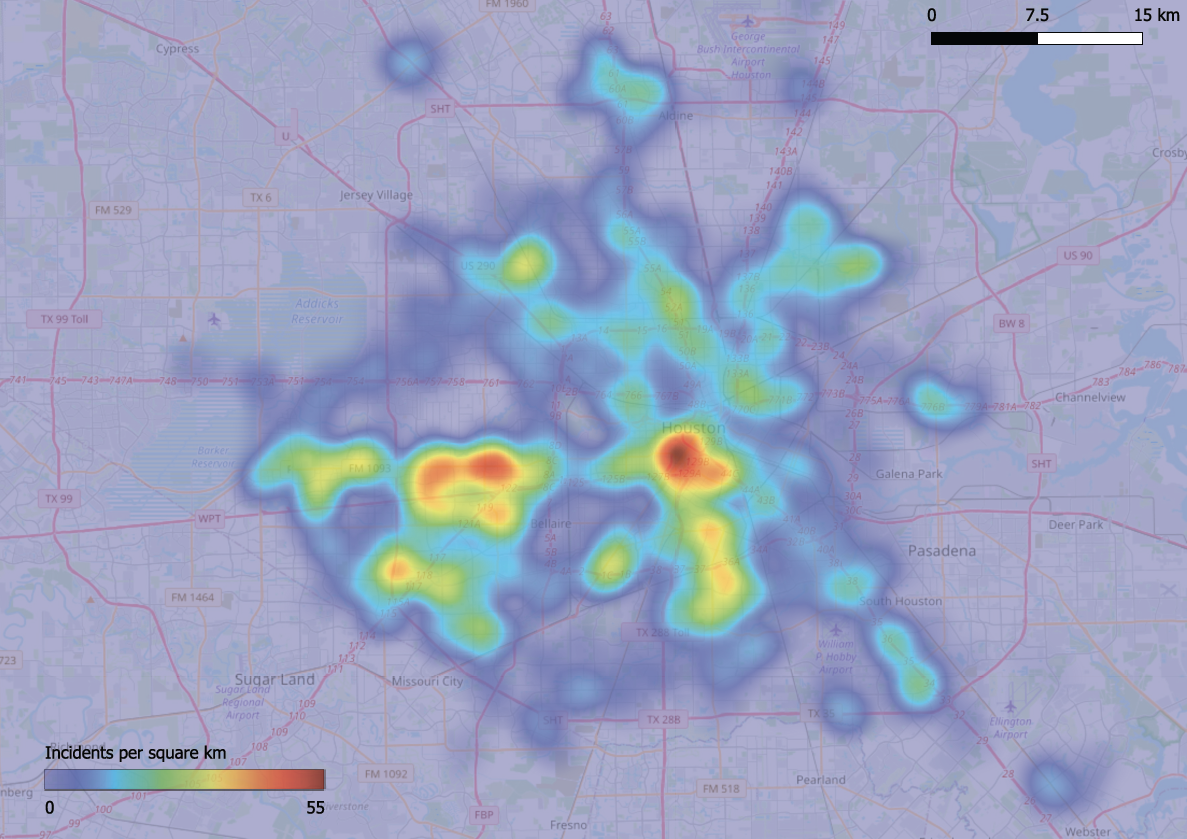} 
\put(20,73){July to September 2024}
\end{overpic} \quad
\begin{overpic}[width=0.45\textwidth, grid=false]{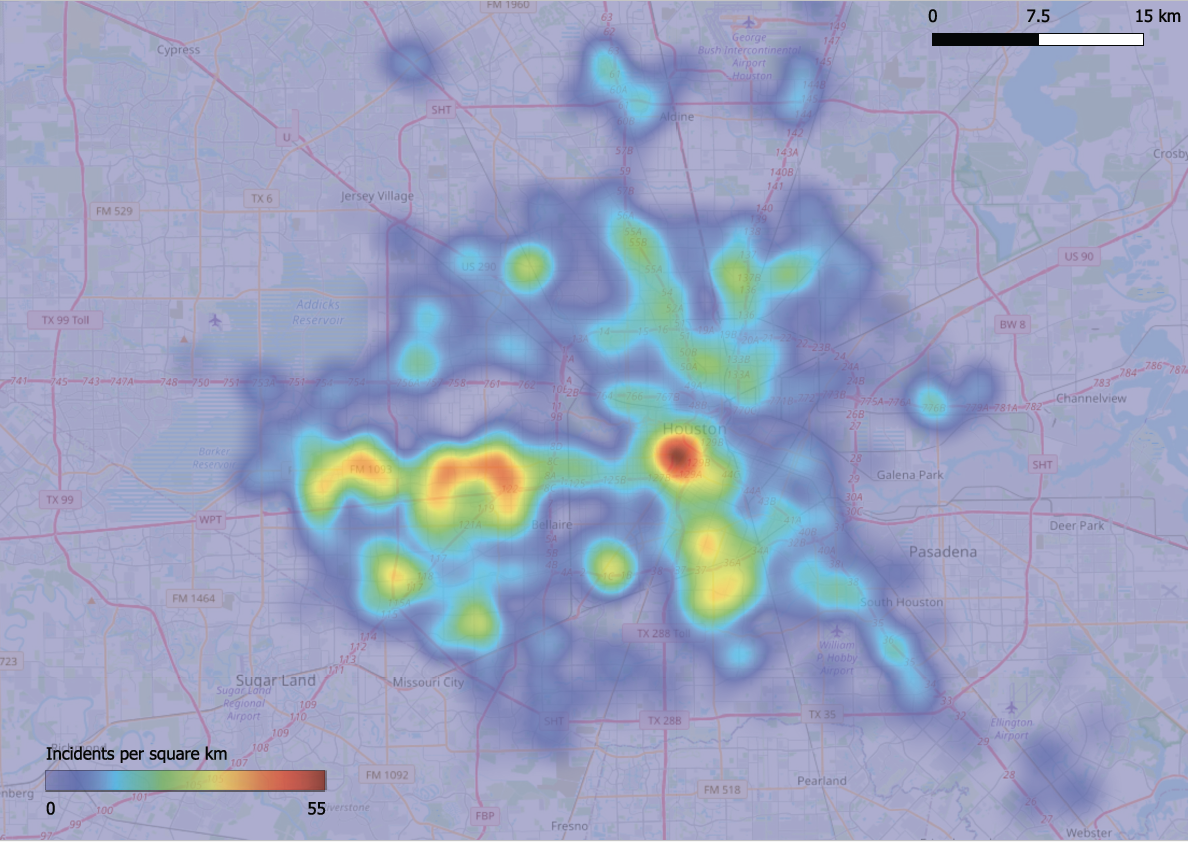}
\put(15,73){October to December 2024}
\end{overpic} 
\caption{Number of burglaries per km$^2$ in the Houston area
for four three-month periods in 2024. Police data are taken from a publicly available database \cite{Police_data} and visualized with 
open-source visualization software QGIS \cite{QGIS}, using a radius of 2.5 Km
to calculate the densities.
}
\label{fig:Houston}
\end{figure}

Predictability and control of urban crime are fundamental for public safety. Despite advances in understanding the mechanisms driving crime spread and dissipation \cite{Brantinghambook,Eck1995,Chaineybook}, it remains extremely hard to predict crime due to its complex nature. 
Prediction strategies based on 
mapped crime data have limited efficacy because crime distributions exhibit low regularity in both space and time and display significantly different patterns at different spatio-temporal scales \cite{Weisburd2015,mohler2019reducing}.
See Fig.~\ref{fig:Houston}, which shows the emergence of large burglary hotspots in the Houston area over three-month periods in 2024. 
These hotspots vary in lifetime and scale depending on geographic, economic, and seasonal factors. Also, 
different crime types exhibit distinct patterns of emergence, diffusion, and dissipation.


The need for better strategies to predict crime hotspots, and thus address or prevent them,
has prompted mathematicians to develop crime models. In this paper, we focus on a specific kind of crime: residential burglary, which is simpler than other crimes since mobile offenders target stationary sites. This crime has been extensively studied in criminology and sociology, where researchers have observed how offenders are driven by environmental cues \cite{Wilson1982,Woodworth2024} and respond to the success or failure of their 
activities \cite{Johnson1997,Wrightbook}. 
The pioneering paper that translated 
these studies into mathematical modeling 
is \cite{short2008statistical}. Therein, two types of models are proposed. 
The first type is a probabilistic agent-based model, where the interacting “agents” are offenders and housing units. 
The second type of model is derived from the first one through a mean-field limit operation  and consists in a system of PDEs that describes the average behavior at large scale, i.e., for a large number of agents. 
Both the probabilistic model and the PDE model feature complex, possibly dynamic, parameters such as the attractiveness of residences \cite{Woodworth2024} and strength of 
neighborhood effects \cite{Johnson1997,Wilson1982}. \anna{The use of these models for predictive policing, i.e., to suggest law enforcement where crime is most likely to occur, has attracted criticism \cite{Haskins2019,Castelvecchi2020}.
While the models can certainly be improved to reduce bias, we see them as the basis to develop a
tool to assist municipalities in deciding where to direct community investments (instead of policemen) to address systemic needs and reduce crime.}

A rigorous analysis of the hotspot formation in the PDE model from \cite{short2008statistical}
is presented in \cite{short2010dissipation}. 
In \cite{short2010nonlinear}, a 
weakly nonlinear analysis is used to study hotspot formation, bifurcations and suppression of crime patterns in 
a simplified version of the 
PDE model from  \cite{short2008statistical}. 
A study on the local existence and uniqueness of the solution to such simplified PDE model is reported in \cite{rodriguez2010}. The simplified PDE 
model contains a deterrence term 
for the effect of police deployment.
The modeling of police intervention is improved in \cite{zipkin2014} through
a dynamic optimal control problem, while
the work in \cite{Woodworth2024} incorporates spatial data about housing to improve
the estimates on residential crime density. 
In \cite{pan2018}, the Poisson process 
used in the agent-based model from \cite{short2008statistical} to represent the movement of criminal agents (i.e., random walk) is replaced by a truncated L\'{e}vy distribution. This leads to a 
different system PDEs with respect to the one in \cite{short2008statistical} and
results for problems in one dimension are presented in \cite{pan2018}.
A different agent-based model with police enforcement can be found in \cite{jones2010}, while other 
chemotaxis and nonlinear diffusion models can be found in \cite{mei2020,rodriguez2020,BellomoWinkler2022,Fuest2024}. \anna{In particular, reference \cite{Fuest2024} describes and analyzes a model for the 
interactions of two rivaling
gangs attempting to establish or defend territories by spraying intimidating graffiti. The analysis includes construction of  global, bounded classical
solutions and proof of their convergence to homogeneous steady states, while the numerical results show steady-state solutions (in diffusion-dominated regime)
and nontrivial asymptotic behaviors (in convection-dominated regime).
}

\anna{
As it may be apparent from the literature review
presented above, most models have been developed at the microscopic (i.e., individual-based) or macroscopic (i.e., based on continuum assumptions) scale. To the best of our knowledge, the only model
developed at the mesoscopic scale, which is 
intermediate between the small and the large
scales, is presented in \cite{Bellomo2015}. This mesoscopic model
incorporates tools from game theory 
to study how interactions between the general population, criminals, and policemen lead to patterns in criminal behavior. 
It is found that small, individual changes influenced by factors, such as wealth differences,  
can cause large, discernible differences in crime patterns. 
These factors are included in an additional variable, called {activity}.
While the activity variable in 
\cite{Bellomo2015} is constant in time for
simplicity, it could evolve in time as a result of interactions and, consequently, affect the evolution of crime. We believe that more work is needed at the mesoscopic scale and, perhaps
even more importantly, to integrate the scales in a multiscale vision. Indeed, only one
representation scale is not sufficient to 
describe the overall collective dynamics of crime.
See Sec.~\ref{sec:future} for more details on these future perspectives.}


More recently, data-driven methods
have been developed. A
deep learning technique called 
Spatio-Temporal Residual Networks (ST-ResNet) for real time crime forecasting is proposed in \cite{wang2017} and further
developed in \cite{wang2019}. 
A statistics-based model is introduced in \cite{cai2022} and statistical study on the effect of
social distancing during the COVID-19 pandemic on residential burglary and domestic violence is presented in \cite{MOHLER2020101692}.  
Finally, 
reported gang-related violent crimes in Los Angeles from 2014 to 2017 are modeled 
in \cite{Park02102021}
using spatial–temporal marked Hawkes point processes with covariates.

The reader interested in more in-depth 
reviews of criminal behavior dynamics modeling, including data-driven approaches, is referred to \cite{blumstein2002, DORSOGNA20151, 10937481}. 

In this paper, we consider the two original models from \cite{short2008statistical} and introduce several elements of novelty. First, we replace the periodic boundary conditions used in \cite{short2008statistical} with 
natural boundary conditions, i.e., 
conditions that specify the value of the normal derivatives of the system variables at the boundary. 
While periodic boundary conditions are
convenient for the analysis, they assume an infinitely large domain and thus are of little practical use if one is interested in simulating crime in a real city. The boundary conditions we enforce arise naturally from the variational formulation of the PDE problem, which is the starting point for the application of a finite element method. As such, these boundary conditions are enforced automatically, i.e., without any external constraint. 
Of the references
listed above,
numerical results for the PDE model from \cite{short2008statistical} can be found only in \cite{short2008statistical,short2010dissipation,zipkin2014}. These references adopt a spectral method to approximate the solution of the PDE model. 
Spectral methods can be computationally less expensive and easier to implement than finite element methods. 
Moreover, they can be highly accurate in simple domains with smooth solutions. However, because of their global nature, they lead to linear systems with dense matrices. Thus, computational efficiency suffers as the size of the problem increases. For larger problems, finite element methods generally work better due to sparse matrices and more flexibility in handling non-trivial geometries, including better modeling of sharp bends. 
Because of our interest in simulating crime in real cities, we opt for a finite element-based solver for the PDE model.

The PDE model is a nonlinear, coupled system of two equations in two variables (attractiveness of residential sites and density of criminals), similar to the Keller–Segel model for aggregation based on chemotaxis \cite{KELLER1970399,Byrne2004,doi:10.1137/050637923}. \anna{To learn more about chemotaxis models and cross-diffusion models in complex environments (i.e., interacting with a time-evolving external system), the reader is referred to \cite{BellomoWinkler2022}.}
For the numerical solution of the PDE model from \cite{short2008statistical}, we propose a scheme that 
decouples the computation of the 
attractiveness from the computation of the criminal density at each time step, resulting in the solution of two linear algebraic systems per iteration.
Through numerous numerical tests, 
we demonstrate the robustness and computational efficiency of this approach.
Since we cannot compare with the results from \cite{short2008statistical,short2010dissipation,zipkin2014} due to the difference in boundary conditions, we validate the results
obtained with the PDE model against the results provided by the agent-based model. Finally, leveraging the flexibility allowed by the finite element method, we show results 
for spatially heterogeneous model parameters and a realistic geometry, i.e., the city of Chicago.  

\anna{
The results in
\cite{short2008statistical,short2010dissipation,zipkin2014} assume parameter values uniform in space and time. This oversimplification does not allow to account space heterogeneity in, e.g., attractiveness of housing sites, and behavioral heterogeneity in criminals and policemen. So, we believe that the work in this paper is an important effort towards making the PDE model applicable to realistic scenarios. 
A way to introduce further heterogeneity is through the above-mentioned activity variable. 
See Sec.~\ref{sec:future} for more details.}


The rest of the paper is organized as follows. 
Sec.~\ref{sec:pb_def} discusses the agent-based model and associated PDE model from \cite{short2008statistical}, together with the new boundary conditions. In Sec.~\ref{sec:num_met}, we present the time and space discretization of the PDE model 
and introduce our iterative partitioned algorithm to solve the discretized problem. Sec.~\ref{sec:num_res}
reports the numerical results obtained with both models. \anna{Future perspectives are discussed
in Sec.~\ref{sec:future} and} conclusions are drawn
in Sec.~\ref{sec:concl}.

\section{Problem definition}\label{sec:pb_def}

\subsection{The agent-based model}\label{sec:disc_model}

This section summarizes the discrete
statistical model from \cite{short2008statistical}. It is relevant to our work for two reasons. 
First, the PDE model for which we propose a finite element method is derived from this discrete model. 
Second, Sec.~\ref{sec:validation}
compares the results obtained from the PDE model and the agent-based model.

Let us consider a smooth and bounded  domain $\Omega\subset \mathbb{R}^2$ where there are two kinds of ``agents'': the burglars and the houses where the burglaries can occur. The houses are located on a two-dimensional lattice. 
For simplicity, we will consider a rectangular domain and 
consider a Cartesian lattice with uniform
spacing $h = \Delta x = \Delta y$.
Similarly, time is discretized with a time step $\Delta t$. Houses are located at the lattice
sites $s$ and are associated with a time-dependent
level of attractiveness $A_s(t)$, which measures the burglars’
perception of the attractiveness of the home at site $s$.
In general, we will use subindex $s$ to denote the variables in the discrete model.
Variable $A_s(t)$ is written as the sum of a static (possibly
spatially varying) component $A_s^{st}$ and a dynamic component 
$B_s(t)$ associated with repeat
and near-repeat victimization:
\begin{equation}\label{eq:attr}
    A_s(t) = A_s^{st} + B_s(t).
\end{equation}  

Criminals either burglarize a house or move to a neighboring
house with a random walk biased toward areas of high attractiveness. The probability of each crime occurring at site $s$ during time interval $[t, t+\Delta t]$ is given by:
\begin{equation}\label{eq:prob}
    p_s(t)=1-e^{-A_s(t)\Delta t},
\end{equation}
in accordance with a standard Poisson process.
When a site $s$ is burglarized, the corresponding criminal agent is removed from the lattice to represent
the tendency of burglars to flee the crime
location. Obviously, removed burglars can become active again, which can be seen as a regeneration of criminals. 
To represent this ``burglar generation'' in the model, a parameter $\Gamma$ is introduced at each lattice site $s$.

\begin{figure}[htb!]
    \centering
    \begin{overpic}[width=0.5\linewidth,grid=false]{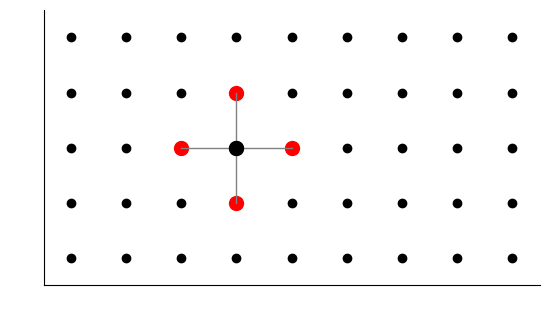}
\put(50,3){$x$}
\put(3,35){$y$}
\put(46,36){$s$}
    \end{overpic}
    \caption{Visualization of a 5-point stencil on a Cartesian lattice with 
    site $s$ at the center to its four neighboring sites (red dots) in the horizontal and vertical directions.}
    \label{fig:stencil}
\end{figure}

The repeat victimization (i.e., a house has an elevated risk of being re-victimized shortly
after a first break in), near-repeat victimization (i.e., the houses near a burglarized site are also at increased risk), and the broken-window effect (i.e., the perception of a neighborhood as crime-tolerant increases the risk of the houses in such neighborhood) are
modeled through the dynamics behavior of $B_s$. 
To capture repeat victimization, the attractiveness $B_s(t)$ at location $s$ 
is increased by $\theta$ for each crime 
committed at $s$ in the preceding time period \cite{short2008statistical,short2010dissipation}. Let parameter $\omega$ be 
a time scale over which repeat victimizations are most likely to occur
and $E_s(t)$ the number of burglaries occurred at site $s$ during the time
interval $(t - \Delta t, t)$, which we discuss how to compute later. 
Thus, to account for repeat victimization
the equation for $B_s(t + \Delta t)$ needs to include term:
\begin{equation*}
    B_s(t)(1-\omega\Delta t)+\theta E_s(t)
\end{equation*}
To model near-repeat victimization and the broken-window effect, we need to allow 
$B_s(t)$ to spread spatially from a house to its neighbors. So, $B_s(t + \Delta t)$ needs to include term of this kind:
\begin{equation*}
     \eta\bigg[\frac{1}{4}\sum_{s' \sim s} B_{s'}(t) - B_s(t) \bigg](1-\omega\Delta t),
\end{equation*}
where parameter $\eta \in [0, 1]$ characterizes the strength of the 
near-repeat victimization and the broken-window effect
and $s' \sim s$ denotes the set of houses neighboring site $s$.
{For a Cartesian lattice, the cardinality of this set is four (see Fig.~\ref{fig:stencil}), which explains the factor of 4 in the formula above.}

Putting the two terms discussed above together, we get
\begin{equation}\label{eq:discrete_bertozzi-0}
    B_{s}(t+\Delta t) = \bigg[{(1 - \eta)}B_s(t)+\frac{\eta}{4}\sum_{s' \sim s} B_{s'}(t)\bigg](1-\omega\Delta t)+\theta E_s(t).
\end{equation}
Note that eq.~\eqref{eq:discrete_bertozzi-0} can be rewritten as
\begin{equation}\label{eq:discrete_bertozzi-1}
    B_{s}(t+\Delta t) = \bigg[B_s(t)+\frac{\eta h^2}{4} \Delta_h B_{s}(t)\bigg](1-\omega\Delta t)+\theta E_s(t),
\end{equation}
where
\begin{equation}\label{eq:Delta_Bst}
    \Delta_h B_{s}(t) = \frac{1}{h^2} \left(\sum_{s' \sim s} B_{s'}(t) - 4 B_{s}(t) \right)
\end{equation}
is the discrete spatial Laplacian operator computed using the stencil in 
Fig.~\ref{fig:stencil}. One way to estimate $E_s(t)$ in \eqref{eq:discrete_bertozzi-1} is given by:
\begin{equation}\label{eq:Es}
    E_s(t) = n_s(t) p_s(t),
\end{equation}
where $n_s (t)$ is the number of criminals at site $s$ at time $t$ and $p_s(t)$ is defined in \eqref{eq:prob}.

Next, we need to describe the evolution of $n_s (t)$.
Recall that a criminal at site $s$ at a given time either burglarizes it or moves to a neighboring site following a biased
random walk, with probability of moving to site $r$:
\begin{equation*}
    q_{s \rightarrow r} = \frac{A_r (t)}{\displaystyle{\sum_{s' \sim s} A_{s'}(t)}}.
\end{equation*}
Thus, the criminals at site $s$ at the following time step
have either arrived after failing to burglarize a 
neighboring site or have been generated there at rate $\Gamma$.
This gives us the following:
\begin{equation}
    n_s(t+\Delta t) = A_s \sum_{s'\sim s} \frac{n_{s'}(t)[1-p_{s'}(t)]}{T_{s'}(t)}+\Gamma \Delta t, \quad 
    T_{s'}(t) = \sum_{s''\sim s'}A_{s''}(t).
    \label{eq:discrete_bertozzi-2}
\end{equation}
{Note that, to prevent the denominator of $q_{s \rightarrow r}$ from becoming zero,
we assume $A_s(t) \geq a_0 > 0$ at all times $t$ under consideration and for all sites $s$ in the lattice.}

\begin{figure}[htb!]
    \centering
    \begin{overpic}[width=0.25\linewidth,grid=false]{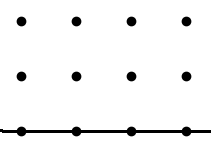}
\put(75,45){$\Omega$}
\put(110,10){$\partial \Omega$}
\put(44,38){$s$}
\put(44,15){$r$}
\put(36,9){\vector(0,-1){12}}
\put(44,-3){$\bn$}
    \end{overpic}
    \caption{{Part of a simplified domain $\Omega$ with boundary $\partial \Omega$ and normal to the boundary $\bn$. Internal site $s$ lies in the normal direction to the boundary with respect to site $r$, which lies on the boundary.}}
    \label{fig:BC}
\end{figure}

Problem \eqref{eq:discrete_bertozzi-1},\eqref{eq:discrete_bertozzi-2}, with $E_s(t)$ as in \eqref{eq:Es} and $p_s(t)$ as in \eqref{eq:prob}, represents the probabilistic agent-based model.
{As mentioned in Sec.~\ref{sec:intro}, we consider a finite domain and natural boundary conditions. In all the numerical results in Sec.~\ref{sec:num_res}, we enforce homogeneous 
natural boundary conditions, meaning that at the boundary the normal derivatives of the system variables {are} set to zero. In the agent-based model, we enforce such boundary conditions for
{for the level of attractiveness and the number of criminals. Let us consider that simple setting in Fig.~\ref{fig:BC} to clarify how we handle these boundary conditions: we impose that $A_s (t) = A_r(t)$ and $n_s (t) = n_r(t)$, so that the derivatives of $A$ and $n$ in the normal direction to the boundary $\bn$ are zero at the boundary. Note that that boundary condition $A_s (t) = A_r(t)$, which can be rewritten as $B_s (t) + A_s^{st} = B_r(t) + A_r^{st} $, has also the effect of making site $s$ no more attractive than site $r$ and thus criminals would not move from $s$ to $r$ based on the above}
random walk probability $q_{s \rightarrow r}$.}

It is shown in \cite{short2008statistical} that,
depending on the values of the parameters, 
the system described by \eqref{eq:discrete_bertozzi-1}-\eqref{eq:discrete_bertozzi-2} can exhibit different 
hotspot patterns. 
Additionally, when all sites have the same level of attractiveness $\bar{A} = A^{st}+\bar{B}$ and, on average, the same number of criminals $\bar{n}$, the system has a homogeneous equilibrium
solution with:
\begin{equation}\label{eq:equil_sol}
    \bar{B} =\frac{\theta\Gamma}{\omega}, \quad \quad \quad \quad  \bar{n} =\frac{\Gamma \Delta t}{1-e^{-\bar{A}\Delta t}}.
\end{equation}
While it might be unsettling to see the time step, i.e., a chosen discretization quantity, in the above homogeneous equilibrium solution, we note that at the leading order the denominator of $\bar{n}$
in \eqref{eq:equil_sol} is approximated by $\bar{A}\Delta t$ and so $\bar{n} \approx \Gamma/\bar{A}$.

\subsection{The PDE model}
Let $\rho = n_s(t)/h^2$ be the density
of criminals.
By manipulating eq.~\eqref{eq:discrete_bertozzi-1}, with $E_s(t)$ as in \eqref{eq:Es}, and taking its limit for $h, \Delta t \rightarrow 0$, while keeping the ratio $\frac{h^2}{\Delta t}$ fixed to $D$ and quantity $\theta \Delta t$ fixed to $\epsilon$, we obtain:
\begin{equation}\label{eq:contB_dim}
    \frac{\partial B}{\partial t} = \frac{\eta D}{4} \Delta B - \omega B + \epsilon D \rho A, \quad A(\bx, t) = A^{st}(\bx) + B(\bx, t),
\end{equation}
where $\bx$ gives the spatial coordinates. In \eqref{eq:contB_dim},
we have removed the subindex $s$ {from} the variables since the model is not discrete anymore. 
Similarly, by manipulating eq.~\eqref{eq:discrete_bertozzi-2} and taking its limit for $h, \Delta t \rightarrow 0$, while keeping the ratio $\frac{\Gamma}{h^2}$ fixed to $\gamma$, we get:
\begin{equation}\label{eq:contRho_dim}
     \frac{\partial \rho}{\partial t}= \frac{D}{4} \nabla \cdot \Big(\nabla\rho - \frac{2 \nabla A}{A} \rho \Big)-\rho A +\gamma.
\end{equation}
{Similarly to the strict positivity requirement on $A_s(t)$ stated after eq.~\eqref{eq:discrete_bertozzi-2}, we assume 
$A(\bx, t) \geq a_0 > 0$ for $\bx \in \Omega$ and $t$ in the time interval of interest.}

To simplify problem \eqref{eq:contB_dim}-\eqref{eq:contRho_dim}, we introduce the following non-dimensional quantities:
\begin{equation}\label{eq:non_dim}
    \tilde{A} = \frac{A}{\omega} = \frac{A^{st}}{\omega} + \frac{B}{\omega},\quad \tilde{\rho} = \epsilon L^2 \rho, \quad \tilde{\bx} = 2 \frac{\bx}{L}, \quad \tilde{t} = \omega t,
\end{equation}
where 
$L = \sqrt{{D}/{\omega}}$ is a characteristic length. 
{We note that the 2 in the third equation in \eqref{eq:non_dim}
comes from $\sqrt{4}$, where 4 is the number of neighbors of a given site (see Fig.~\ref{fig:stencil}).}
For simplicity of notation, we will omit the tilde on the non-dimensional variables and from now on every variable in the continuum (or PDE) model will be non-dimensional.

The non-dimensional continuum model can be written as: find $A(\bx, t)$ and $\rho(\bx, t)$ such that 
\begin{align}
     & \frac{\partial A}{\partial t} - \eta \Delta A + A - \rho A = - \eta \Delta A^{st} + A^{st}, \label{eq:continous_bertozzi-1}\\
     &\frac{\partial \rho}{\partial t} - \nabla \cdot \Big(\nabla\rho - \frac{2 \nabla A}{A} \rho\Big)+\rho A = \frac{\Gamma \theta}{\omega^2},\label{eq:continous_bertozzi-2}
\end{align}
in $\Omega \times (0,T]$, where $(0,T]$ is a time interval of interest. 
Eq.~\eqref{eq:continous_bertozzi-1}-\eqref{eq:continous_bertozzi-2} are arranged so that 
known forcing terms are at the right-hand side, while
all the terms involving variables are at the left-hand side.


Let us denote a homogeneous equilibrium solution to \eqref{eq:continous_bertozzi-1}-\eqref{eq:continous_bertozzi-2}
with $\bar{A}$ and $\bar{\rho}$. 
We note that there is a slight abuse of notation with respect to Sec.~\ref{sec:disc_model}, since here $\bar{A}$ and $\bar{\rho}$ are nodimensionalized according to \eqref{eq:non_dim}.
Assuming $A^{st}$ is homogeneous in space, solution $(\bar{A},\bar{\rho})$ satisfies:
\begin{align}
     & \bar{A} - \bar{\rho} \bar{A} = A^{st}, \label{eq:continous_bertozzi-1ss} \\
     &\bar{\rho} \bar{A} = \frac{\Gamma \theta}{\omega^2}, \label{eq:continous_bertozzi-2ss}
\end{align}
where the right-hand side of \eqref{eq:continous_bertozzi-2ss}
is $\bar{B}$ is given in 
\eqref{eq:equil_sol} and nodimensionalized according to \eqref{eq:non_dim}. We continue with the same abuse of notation and denote the right-hand side of \eqref{eq:continous_bertozzi-2ss} with $\bar{B}$. Thus, \eqref{eq:continous_bertozzi-2ss} can be written as $\bar{\rho} = \bar{B}/\bar{A}$, which can be plugged into \eqref{eq:continous_bertozzi-1ss} to get the obvious $\bar{A} = A^{st}+\bar{B}$. Hence, once the value of $\bar{B}= {\Gamma \theta}/{\omega^2}$ is selected, 
we get automatically $\bar{A} = A^{st}+\bar{B}$, since $A^{st}$ is given, and then we can find $\bar{\rho} = \bar{B}/\bar{A}$.

The linear stability analysis of problem \eqref{eq:continous_bertozzi-1}-\eqref{eq:continous_bertozzi-2} with periodic boundary conditions
determines which combinations of parameters $(\eta, A^{st},\bar{B})$ lead to the instability around the homogeneous steady-state that generates hotspots
\cite{short2008statistical}. Such instability arises when  $\bar{B}>A^{st}/2$ and the following sufficient condition is
satisfied:
\begin{equation}\label{eq:suff_cond}
    \eta <\frac{3\bar{\rho}+1-\sqrt{12\bar{\rho}}}{\bar{A}}.
\end{equation}
This means that, provided that $A^{st} < 2 \Gamma \theta/\omega^2$, hotspots form only when the neighborhood effects are weak enough. {Although we do not use periodic boundary conditions, we have reported this result because in 
Sec.~\ref{sec:validation} we will use \eqref{eq:suff_cond}
as a heuristic indication of when to expect the above-mentioned instability, i.e., below a given threshold for $\eta$. Of course, the value of the threshold would be different in our case
since we use Neumann boundary conditions.}

Obviously, problem \eqref{eq:continous_bertozzi-1}-\eqref{eq:continous_bertozzi-2} needs to be supplemented with appropriate boundary and initial conditions. 
Given that periodic boundary conditions are not very realistic, 
we opt for the following natural boundary conditions
\begin{align}
   \nabla A {\cdot} \bn &= 0,  \quad \text{on } \partial\Omega,\; t>0, \label{eq:BC1} \\
   \nabla \rho {\cdot} \bn &= 0, \quad \text{on } \partial\Omega,\; t>0, \label{eq:BC2} 
\end{align}
where $\bn$ is the outward normal unit vector. Below, we will show how these boundary conditions arise naturally from the variational formulation of problem \eqref{eq:continous_bertozzi-1}-\eqref{eq:continous_bertozzi-2}. As for the initial conditions, we set: 
\begin{align}
     A(\bm{x}, 0) &= A^{st}(\bm{x}) +  B_0(\bm{x})+\chi_\delta(\bm{x})\xi_B(\bm{x})\in L^2(\Omega),\label{eq:ic1} \\
\rho(\bm{x}, 0) &= \rho_0(\bm{x})+\chi_\delta(\bm{x})\xi_\rho(\bm{x}) \in  L^2(\Omega), \label{eq:ic2}
\end{align}
where 
$A^{st}(\bm{x})+B_0(\bm{x})$ is the prescribed initial level of attractiveness and
$\rho_0(\bm{x})$ is the initial density, while $\xi_B(\bm{x})$  and $\xi_\rho(\bm{x})$ represent noise perturbations, drawn from a Gaussian distribution $\mathcal{N}(0,\sigma^2)$. Function $\chi_\delta(\bm{x})$ encodes the sparsity of the noise: $\chi_\delta(\bm{x})=1$ with probability $\delta \in [0,1]$, and $\chi_\delta(\bm{x})=0$ otherwise, independently for each 
$\bm{x}$. Note that $\delta = 1$ corresponds to spatially uniform Gaussian noise, while $\delta = 0$ corresponds to no noise. For $0 < \delta < 1$, the noise is sparse. This formulation of the initial conditions enables a systematic investigation of the emergence of realistic patterns under both uniform and sparse noise perturbations.

For the numerical method that will be discussed in Sec.~\ref{sec:num_met}, 
we need a variational formulation of problem
\eqref{eq:continous_bertozzi-1}-\eqref{eq:continous_bertozzi-2}. To devise it, one multiplies \eqref{eq:continous_bertozzi-1} by $u\in H^1(\Omega)$
and \eqref{eq:continous_bertozzi-2} by $v\in H^1(\Omega)$, integrates over $\Omega$ and employs the integration by parts identity.  This leads to the variational formulation:
Find $A(\bx, t)\in H^1(\Omega) {\cap
L^\infty(\Omega)}$ and $\rho(\bx, t) \in H^1(\Omega)$ such that 
\begin{align}
        &\int_\Omega \frac{\partial A}{\partial t} u\;d\Omega + \int_\Omega \eta\nabla A \cdot \nabla u\;d\Omega - \int_{\partial \Omega} \eta(\nabla A {\cdot} \bn)u\;d\Gamma 
        +\int_\Omega Au \;d\Omega \cl
        & \quad \quad \quad 
        - \int_\Omega \rho A u \;d\Omega  = \int_\Omega (A^{st} - \eta \Delta A^{st}) u \;d\Omega , \cl
        & \int_\Omega \frac{\partial \rho}{\partial t} v \;d\Omega +\int_\Omega  \Big(\nabla\rho - \frac{2 \nabla A}{A} \rho \Big) \nabla v\;d\Omega 
        -\int_{\partial \Omega}  \Big(\nabla\rho - \frac{2 \nabla A}{A} \rho \Big){\cdot}\bn v\;d\Gamma \cl
        & \quad \quad \quad 
        +\int_\Omega \rho A v\;d\Omega =  \int_\Omega \frac{\Gamma \theta}{\omega^2} v\;d\Omega,  \el
\end{align}
for all $(u,v)\in H^1(\Omega)\times H^1(\Omega)$. 
Boundary integrals involving $\nabla A {\cdot} \bn$ and
$\nabla \rho {\cdot} \bn$ have naturally appeared 
in the above weak problem. Thanks to boundary 
conditions \eqref{eq:BC1}-\eqref{eq:BC2}, such integrals vanish. Thus, 
we obtain the following simpler weak problem: Find $A(\bx, t)\in H^1(\Omega){\cap
L^\infty(\Omega)}$ and $\rho(\bx, t) \in H^1(\Omega)$ such that
\begin{align}
        &\int_\Omega \frac{\partial A}{\partial t} u\;d\Omega + \int_\Omega \eta\nabla A \cdot \nabla u\;d\Omega 
        +\int_\Omega Au \;d\Omega  
        - \int_\Omega \rho A u \;d\Omega  = \int_\Omega (A^{st} - \eta \Delta A^{st}) u \;d\Omega , \label{eq:weak_bertozzi-1} \\
        & \int_\Omega \frac{\partial \rho}{\partial t} v \;d\Omega +\int_\Omega  \Big(\nabla\rho - \frac{2 \nabla A}{A} \rho \Big) \nabla v\;d\Omega 
        +\int_\Omega \rho A v\;d\Omega =  \int_\Omega \frac{\Gamma \theta}{\omega^2} v\;d\Omega,  \label{eq:weak_bertozzi-2} 
\end{align}
for all $(u,v)\in H^1(\Omega)\times H^1(\Omega)$.


In \eqref{eq:weak_bertozzi-1}, we have not {integrated} by parts the Laplacian of $A^{st}$, which means that we require $A^{st} \in H^2(\Omega)$. That is the case for all the numerical tests discussed in Sec.~\ref{sec:num_res}. If one wants to require less regularity on $A^{st}$, i.e., $A^{st} \in H^1(\Omega)$, then eq.~\eqref{eq:weak_bertozzi-1} becomes
\begin{align}
        &\int_\Omega \frac{\partial A}{\partial t} u\;d\Omega + \int_\Omega \eta\nabla A \cdot \nabla u\;d\Omega 
        +\int_\Omega Au \;d\Omega  
        - \int_\Omega \rho A u \;d\Omega  = \int_\Omega A^{st} u \;d\Omega \cl
        & \quad + \int_\Omega \eta\nabla A^{st} \cdot \nabla u\;d\Omega - \int_{\partial \Omega} \eta(\nabla A^{st} {\cdot}\bn)u\;d\Gamma. \el
\end{align}

\section{Numerical method for the PDE model}\label{sec:num_met}

PDE problem \eqref{eq:weak_bertozzi-1}-\eqref{eq:weak_bertozzi-2} poses some interesting challenges.
The first main challenge is related to its two-way coupled nature: both variables $\rho$ and $A$ appear
in both eq.~\eqref{eq:weak_bertozzi-1} and eq.~\eqref{eq:weak_bertozzi-2}. The coupling terms
create an intricate feedback mechanisms between the criminal density $\rho$ and the attractiveness field $A$. 
The second main challenge is the strong nonlinear nature of the problem. In particular, the second term in  
eq.~\eqref{eq:weak_bertozzi-2}
is degenerate and nonlinear and, as such, 
it could complicate stability and convergence in 
a numerical scheme. Thus, the design of an accurate
and efficient numerical method for problem \eqref{eq:weak_bertozzi-1}-\eqref{eq:weak_bertozzi-2}
is non-trivial. 

Let us start from the time discretization. For this, 
we introduce a time step $\Delta t = T/N_t$ for given $N_t$, set $t^n =n\Delta t $ for $n = 1, \dots, N_t$
and denote with $y^n$ the approximation of quantity 
$y$ at time $t^n$. 
For simplicity, we apply Backward Differentiation Formula
of order 1 to \eqref{eq:weak_bertozzi-1}-\eqref{eq:weak_bertozzi-2} and get the following problem at the generic time $t^{n+1}$: for $n\geq 0$, 
given $A^n$ and $\rho^n$, find $A^{n+1}\in H^1(\Omega) {\cap
L^\infty(\Omega)}$ and $\rho^{n+1}\in H^1(\Omega)$
such that
\begin{align}
           & \frac{1}{\Delta t}\int_\Omega A^{n+1} u\;d\Omega + 
           \int_\Omega \eta\nabla A^{n+1} \cdot \nabla u\;d\Omega +\int_\Omega A^{n+1}u \;d\Omega \cl
           & \quad \quad -\int_\Omega \rho^{n+1} A^{n+1} u \;d\Omega  = \int_\Omega (A^{st} - \eta \Delta A^{st} ) u \;d\Omega + \frac{1}{\Delta t}\int_\Omega A^{n} u\;d\Omega, \label{eq:implicit-1} \\
         & \frac{1}{\Delta t} \int_\Omega \rho^{n+1} v \;d\Omega +\int_\Omega  \Big(\nabla\rho^{n+1} - \frac{2\nabla A^{n+1}}{A^{n+1}} \rho^{n+1}\Big) \nabla v\;d\Omega  \cl
         & \quad \quad +\int_\Omega \rho^{n+1} A^{n+1} v\;d\Omega = \int_\Omega \frac{\Gamma \theta}{\omega^2} v\;d\Omega + \frac{1}{\Delta t} \int_\Omega \rho^{n} v \;d\Omega, \label{eq:implicit-2}
\end{align}
for all $(u,v)\in H^1(\Omega)\times H^1(\Omega)$.
System \eqref{eq:implicit-1}-\eqref{eq:implicit-2} is
implicit and fully coupled. Hence, if one adopts a fine mesh for space discretization, it could lead to 
a large nonlinear system that needs to be solved
at each time step through an
iterative solvers like, e.g., Newton’s method. The strongly 
nonlinear term $\frac{2 \nabla A^{n+1}}{A^{n+1}}\rho^{n+1}$ could make the convergence of such 
iterative method slow, especially for large-scale simulations \cite{quarteroni1994numerical,saad2003iterative}. 

One way to contain the computational cost is to decouple system \eqref{eq:implicit-1}-\eqref{eq:implicit-2}, while linearizing it, so that one solves two smaller linear systems, instead of one large nonlinear system, per time step. 
At time $t^{n+1}$, given $A^n$ and $\rho^n$, the decoupled, linear algorithm performs two sequential steps: 
\begin{itemize}
    \item[-] \emph{Step 1}:  find $A^{n+1}\in H^1(\Omega)$ such that
    \begin{align}
           &\frac{1}{\Delta t}\int_\Omega A^{n+1} u\;d\Omega + 
           \int_\Omega \eta\nabla A^{n+1} \cdot \nabla u\;d\Omega +\int_\Omega A^{n+1}u \;d\Omega \cl
           & \quad \quad -\int_\Omega \rho^{n} A^{n+1} u \;d\Omega  = \int_\Omega (A^{st} - \eta \Delta A^{st}) u \;d\Omega + \frac{1}{\Delta t}\int_\Omega A^{n} u\;d\Omega, \label{eq:explicit-1} 
\end{align}
for all $u \in H^1(\Omega)$.
 \item[-] \emph{Step 2}:  find $\rho^{n+1}\in H^1(\Omega)$ such that
 \begin{align}
         & \frac{1}{\Delta t} \int_\Omega \rho^{n+1} v \;d\Omega +\int_\Omega  \Big(\nabla\rho^{n+1} - \frac{2\nabla A^{n+1}}{A^{n+1}} \rho^{n+1}\Big) \nabla v\;d\Omega  \cl
         & \quad \quad +\int_\Omega \rho^{n+1} A^{n+1} v\;d\Omega = \int_\Omega \frac{\Gamma \theta}{\omega^2} v\;d\Omega + \frac{1}{\Delta t} \int_\Omega \rho^{n} v \;d\Omega, \label{eq:explicit-2}
\end{align}
for all $v \in H^1(\Omega)$.
\end{itemize}
While successfully containing the computational time, the above scheme may suffer
from instability for certain combination of parameter values. Algorithms like \eqref{eq:explicit-1}-\eqref{eq:explicit-2} are called loosely-coupled because
they do not enforce the strong coupling between the two subproblems at each time step.

One way to cure possible instability 
in \eqref{eq:explicit-1}-\eqref{eq:explicit-2}
is to iterate between Step 1 and Step 2 at a given time
step until a suitable convergence criterion is satisfied. This means that at every time step $t^{n+1}$, given $A^n$ and $\rho^n$, the following steps are performed at iteration $k+1$, $k \ge 0$:
\begin{itemize}
    \item[-] \emph{Step 1}:  find $A^{k+1}\in H^1(\Omega)$ such that
    \begin{align}
           &\frac{1}{\Delta t}\int_\Omega A^{k+1} u\;d\Omega + 
           \int_\Omega \eta\nabla A^{k+1} \cdot \nabla u\;d\Omega +\int_\Omega A^{k+1}u \;d\Omega \cl
           & \quad \quad -\int_\Omega \rho^{k} A^{k+1} u \;d\Omega  = \int_\Omega (A^{st} - \eta \Delta A^{st}) u \;d\Omega + \frac{1}{\Delta t}\int_\Omega A^{n} u\;d\Omega, \label{eq:iterative-1} 
\end{align}
for all $u \in H^1(\Omega)$. For $k = 0$, 
$\rho^{k} = \rho^{n}$.
 \item[-] \emph{Step 2}:  find $\rho^{k+1}\in H^1(\Omega)$ such that
 \begin{align}
          & \frac{1}{\Delta t} \int_\Omega \rho^{k+1} v \;d\Omega +\int_\Omega  \Big(\nabla\rho^{k+1} - \frac{2\nabla A^{k+1}}{A^{k+1}} \rho^{k+1}\Big) \nabla v\;d\Omega  \cl
         & \quad \quad +\int_\Omega \rho^{k+1} A^{k+1} v\;d\Omega = \int_\Omega \frac{\Gamma \theta}{\omega^2} v\;d\Omega + \frac{1}{\Delta t} \int_\Omega \rho^{n} v \;d\Omega, \label{eq:iterative-2}
\end{align}
for all $v \in H^1(\Omega)$.
\item[-] \emph{Step 3}: Check the stopping criterion, e.g.
\begin{equation}\label{eq:crit}
    \frac{\|A^{k+1}-A^{k}\|_{L^2(\Omega)}}{\|A^{k}\|_{L^2(\Omega)}}< \texttt{tol}_1 \quad \text{and} \quad \frac{\|\rho^{k+1}-\rho^{k}\|_{L^2(\Omega)}}{\|\rho^{k}\|_{L^2(\Omega)}}< \texttt{tol}_2,
\end{equation}
where $\texttt{tol}_1$ and $\texttt{tol}_2$ are given stopping tolerances. If not satisfied, repeat steps 1-3. If satisfied, set $A^{n+1} = A^{k+1}$ and $\rho^{n+1} = \rho^{k+1}$.
\end{itemize}
At a given time step, algorithm \eqref{eq:iterative-1}-\eqref{eq:crit} solves
two smaller linear problems, i.e., \eqref{eq:iterative-1}-\eqref{eq:iterative-2}, as many times as needed to satisfy
criterion \eqref{eq:crit}. It can be viewed as a fixed-point scheme for nonlinear system 
\eqref{eq:implicit-1}-\eqref{eq:implicit-2}
that allows to decouple the problem while linearizing it. 
So, algorithms like \eqref{eq:iterative-1}-\eqref{eq:crit} are called strongly-coupled, as they enforce the coupling between two subproblems up to user-prescribed tolerances.

Next, we present the space discretization of \eqref{eq:iterative-1}-\eqref{eq:iterative-2}.
We introduce a partition $\mathcal{T}_h$ of $\Omega$ into $N_e$
quadrilaterals $Q_k$ such that
\begin{equation*}
    \overline{\Omega} = \bigcup_{k = 1}^{N_e} \overline{Q}_{k},
\end{equation*}
with 
\begin{equation*}
    h = \max_{Q_k \in \mathcal{T}_h} h_k, \quad h_k = \text{diam}(Q_k), \quad k= 1, \dots, N_e.
\end{equation*}
Let $\mathbb{Q}_N (Q_k)$ be the set of algebraic polynomials, defined on $Q_k$, of degree less than or equal to $N$ in each space variable, and let $V_h (\Omega)$ be the space of global continuous functions on $\overline{\Omega}$ that are
polynomials of degree $N$ on each $Q_k \in \mathcal{T}_h $:
\begin{equation*}
    V_h (\Omega) = \{ v_h \in C^0 (\overline{\Omega})~:~v_h |_{Q_k} \in \mathbb{Q}_N, \forall Q_k \in \mathcal{T}_h  \}.
\end{equation*}
We denote with $\{ \phi_i \}_{i = 1}^{N_Q}$ a basis for $V_h$, where $N_Q$ the total number of degrees of freedom in $\Omega$.

The space-discrete counterpart of \eqref{eq:iterative-1} reads: find $A^{k+1}_h \in V_h (\Omega)$ such that
    \begin{align}
           &\frac{1}{\Delta t}\int_\Omega A^{k+1}_h u\;d\Omega + 
           \int_\Omega \eta\nabla A^{k+1}_h \cdot \nabla u\;d\Omega +\int_\Omega A^{k+1}_h u \;d\Omega \cl
           & \quad \quad -\int_\Omega \rho^{k}_h A^{k+1}_h u \;d\Omega  = \int_\Omega (A_{0,h} - \eta \Delta A_{0,h}) u \;d\Omega + \frac{1}{\Delta t}\int_\Omega A^{n}_h u\;d\Omega, \label{eq:iterative-1-sd} 
\end{align}
for all $u \in V_h (\Omega)$.
Similarly, the space-discrete counterpart of 
\eqref{eq:iterative-2} reads: 
find $\rho^{k+1}_h\in V_h (\Omega)$ such that
 \begin{align}
          & \frac{1}{\Delta t} \int_\Omega \rho^{k+1}_h v \;d\Omega +\int_\Omega  \Big(\nabla\rho^{k+1}_h - \frac{2\nabla A^{k+1}_h}{A^{k+1}_h} \rho^{k+1}_h\Big) \nabla v\;d\Omega  \cl
         & \quad \quad +\int_\Omega \rho^{k+1}_h A^{k+1}_h v\;d\Omega = \int_\Omega \frac{\Gamma \theta}{\omega^2} v\;d\Omega + \frac{1}{\Delta t} \int_\Omega \rho^{n}_h v \;d\Omega, \label{eq:iterative-2-sd}
\end{align}
for all $v \in V_h (\Omega)$.
Obviously, stopping criterion \eqref{eq:crit} becomes
\begin{equation}\label{eq:crit-sd}
    \frac{\|A^{k+1}_h - A^{k}_h\|_{L^2(\Omega)}}{\|A^{k}_h\|_{L^2(\Omega)}}< \texttt{tol}_1 \quad \text{and} \quad \frac{\|\rho^{k+1}_h -\rho^{k}_h \|_{L^2(\Omega)}}{\|\rho^{k}_h \|_{L^2(\Omega)}}< \texttt{tol}_2.
\end{equation}

In order to write the matrix formulation of \eqref{eq:iterative-1-sd}-\eqref{eq:iterative-2-sd}, 
let $M$ be the mass matrix, $K$ the stiffness matrix, $N(w)$ a weighted mass matrix 
and $D(w)$ and ``weighted'' divergence matrix, for suitable weight $w \in V_h (\Omega)$:
\begin{align}
    &M_{i,j}= \int_\Omega \phi_j \phi_i \;d\Omega, \quad
    K_{i,j} = \int_\Omega \nabla \phi_j \cdot \nabla\phi_i \;d\Omega, \cl 
    &N(w)_{i,j} = \int_\Omega w \phi_j \phi_i \;d\Omega, \quad
    D(w)_{i,j} = \int_\Omega \frac{2 \nabla w}{w} \phi_j \nabla \phi_i \;d\Omega, \el
\end{align}
for $i, j = 1, \dots, N_Q$. Moreover, let $\bA^{k+1}, \br^{k+1} \in \mathbb{R}^{N_Q}$ be the vectors of  
the nodal values of $A^{k+1}_h$ and $\rho^{k+1}_h$.

The matrix form of \eqref{eq:iterative-1-sd}
is given by: find vector $\bA^{k+1} \in \mathbb{R}^{N_Q}$ such that
\begin{equation}\label{eq:mat1}
    \frac{1 + \Delta t}{\Delta t} M \bA^{k+1} + \eta K \bA^{k+1} - N(\br^{k}) \bA^{k+1} = \bb_A^{k+1},
\end{equation}
where the vector at the right-hand side $\bb_A^{k+1}$ has $i$-th component
\begin{equation*}
    \bb_{A,i}^{k+1} = \int_\Omega (- \eta \Delta A_{0,h} - A_{0,h}) \phi_i \;d\Omega + \frac{1}{\Delta t} (M\bA^{k})_i,
\end{equation*}
for $i = 1, \dots, N_Q$. Note that the above matrix form assumes that $\eta$ is constant in space. This means that step 1 can be written as: solve linear system:
\begin{equation}\label{eq:mat1bis}
     \left[ \frac{1 + \Delta t}{\Delta t} M + \eta K - N(\br^{k}) \right] \bA^{k+1} = \bb_A^{k+1}.
\end{equation}

The matrix form of \eqref{eq:iterative-2-sd}
is given by: find vector $\br^{k+1} \in \mathbb{R}^{N_Q}$ such that
\begin{equation}\label{eq:mat2}
    \frac{1}{\Delta t} M \br^{k+1} + K \br^{k+1} - D(\bA^{k+1}) \br^{k+1} + N(\bA^{k+1}) \br^{k+1} = \bb_\rho^{k+1},
\end{equation}
where the vector at the right-hand side $\bb_\rho^{k+1}$ has $i$-th component
\begin{equation*}
    \bb_{\rho,i}^{k+1} = \int_\Omega \frac{\Gamma \theta}{\omega^2} \phi_i \;d\Omega + \frac{1}{\Delta t} (M\br^{k})_i,
\end{equation*}
for $i = 1, \dots, N_Q$. Then, step 2 can be written in the equivalent form: 
\begin{equation}\label{eq:mat2bis}
    \left[ \frac{1}{\Delta t} M + K - D(\bA^{k+1}) + N(\bA^{k+1}) \right] \br^{k+1} = \bb_\rho^{k+1},
\end{equation}
where $\bA^{k+1}$ is known from solving \eqref{eq:mat1bis}. Systems \eqref{eq:mat1bis},\eqref{eq:mat2bis}
clearly show that algorithm \eqref{eq:iterative-1}-\eqref{eq:crit}
is a fixed-point scheme for problem 
\eqref{eq:implicit-1}-\eqref{eq:implicit-2}.

\section{Numerical experiments}\label{sec:num_res}

We assess the proposed finite element-based solver for the PDE model through several numerical experiments. 
We start with the validation in Sec.~\ref{sec:validation}
and then study the robustness and computational efficiency of the solver in Sec.~\ref{sec:res_solver}. For the results in Sec.~\ref{sec:validation} and \ref{sec:res_solver}, we use a simple square geometry and constant (in time and space) model parameters. 
Results for spatially heterogeneous parameters and a nontrivial geometry are presented in Sec.~\ref{sec:space_dep}.

We implemented the PDE solver using an iterative partitioned scheme in FEniCSx \cite{BarattaEtal2023,BasixJoss,AlnaesEtal2014}, a popular open-source computing platform for solving PDE problems with the finite element method. 
For the agent-based model, we developed  an in-house solver in MATLAB. 
Both the agent-based and the PDE solvers are available on GitHub \cite{OurCode}.  
For all the simulations based on the PDE model, we used first-order Lagrange finite elements on quadrilateral cells
$\mathbb{Q}_1$.

\subsection{Validation of the solver for the PDE model}\label{sec:validation}

As mentioned in Sec.~\ref{sec:intro},
we validate the results obtained with the PDE model against the
results provided by the agent-based model. For this purpose, we choose 3 representative cases. 
Since the PDE model is appropriate
when the number of criminals is large
(i.e., the criminal population can be looked at as a continuum versus separate individuals), all three cases
have a large initial density. For each case, 
we first present the solution given by the PDE model and motivate the choice of model parameters. Then, we compare with the solution obtained with the agent-based model.

All the quantities mentioned in the following are non-dimensional. 
The computational domain is $\Omega = [0,16]\times [0,16]$ and time interval of interest is $[0, 200]$.
We consider three sets of parameters and initial conditions that yield representative solutions: 
\begin{itemize}
    \item Case 1:  $\eta = 0.9$, $ A^{st} = 1/30$, $\Gamma \theta /\omega^2 = B_0 = 1$, and $\rho_0 = 0.8$. No noise is added to the initial data, i.e., $\chi_\delta(\bm{x}) = 0$ in \eqref{eq:ic1}–\eqref{eq:ic2}.
    \item  Case 2:  $\eta = 0.3 $, $A^{st} = 1/30$, $\Gamma \theta /\omega^2 = {B}_0 = 1$,  and  ${\rho}_0 = 0.8$. Uniform Gaussian noise is added  to the
initial data, with $\chi_\delta(\bm{x}) = 1$, $\xi_B(\bm{x}) \sim \mathcal{N}(0, 0.05^2)$ and $\xi_\rho(\bm{x}) \sim \mathcal{N}(0, 0.01^2)$ in \eqref{eq:ic1}–\eqref{eq:ic2}, sampled independently at each spatial point $\bm{x}$. 
    \item Case 3: $\eta = 0.03 $, $A^{st} = 1/30$, $\Gamma \theta /\omega^2 = {B}_0 = 1$, and ${\rho}_0 = 0.8$. The same noise as in case 2 is added to the
initial data. 
\end{itemize}

In all three case, we have set ${B}_0 = \Gamma \theta /\omega^2$, which corresponds to $\bar{B}$ in \eqref{eq:equil_sol} after non-dimensionalization
as in \eqref{eq:non_dim}. 
Additionally, we note that all the cases have nondimensional $\bar{B}>A^{st}/2$. {In this section, we set $A^{st}$ constant because that is typically the case for the numerical results 
given by the PDE model from \cite{short2008statistical}
presented in the literature \cite{short2008statistical,short2010dissipation,zipkin2014}.
We will consider a space-dependent $A^{st} = A^{st} (\bm{x})$ 
in Sec.~\ref{sec:space_dep}. For realistic simulations, the setting of $A^{st}$ is an open question that needs to be object of future work. See Sec.~\ref{sec:future}.}


To simulate the crime dynamics in the above three cases, 
we consider time step to $\Delta t = 1/50$ and a mesh with size $h = 16/200$ for both
the PDE solver and the agent-based solver.
The PDE solver uses the iterative partitioned algorithm \eqref{eq:iterative-1}-\eqref{eq:crit}, with stopping
tolerances $\texttt{tol}_1 = \texttt{tol}_2 = 10^{-9}.$

\begin{figure}[htb!]
     \centering
         \begin{overpic}[percent,width=0.19\textwidth]{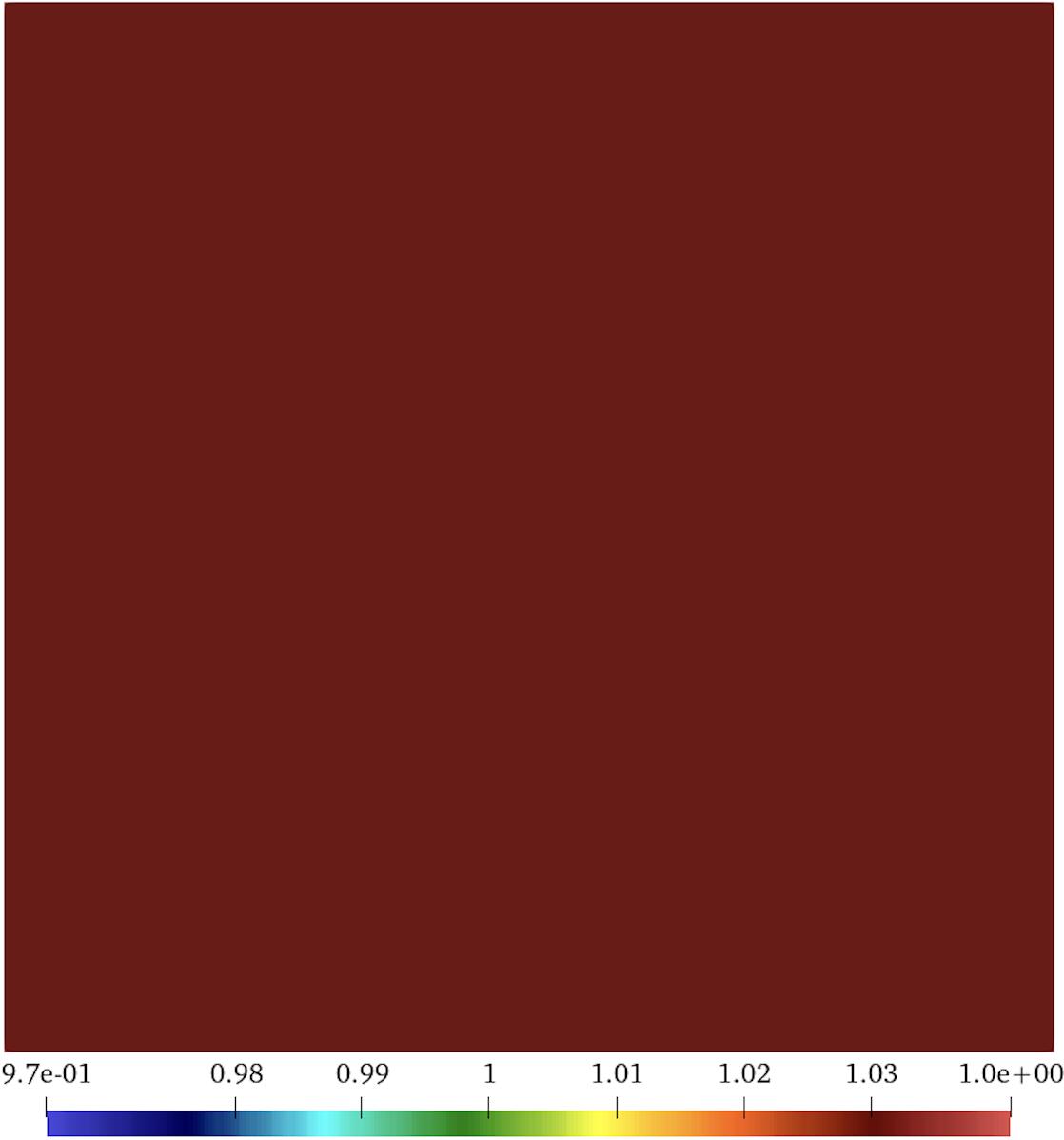}
    \put(30,103){\footnotesize{$t = 0$}}
    \put(-10,55){\footnotesize{$A$}}
    \end{overpic} 
     \begin{overpic}[percent,width=0.19\textwidth, grid=false]{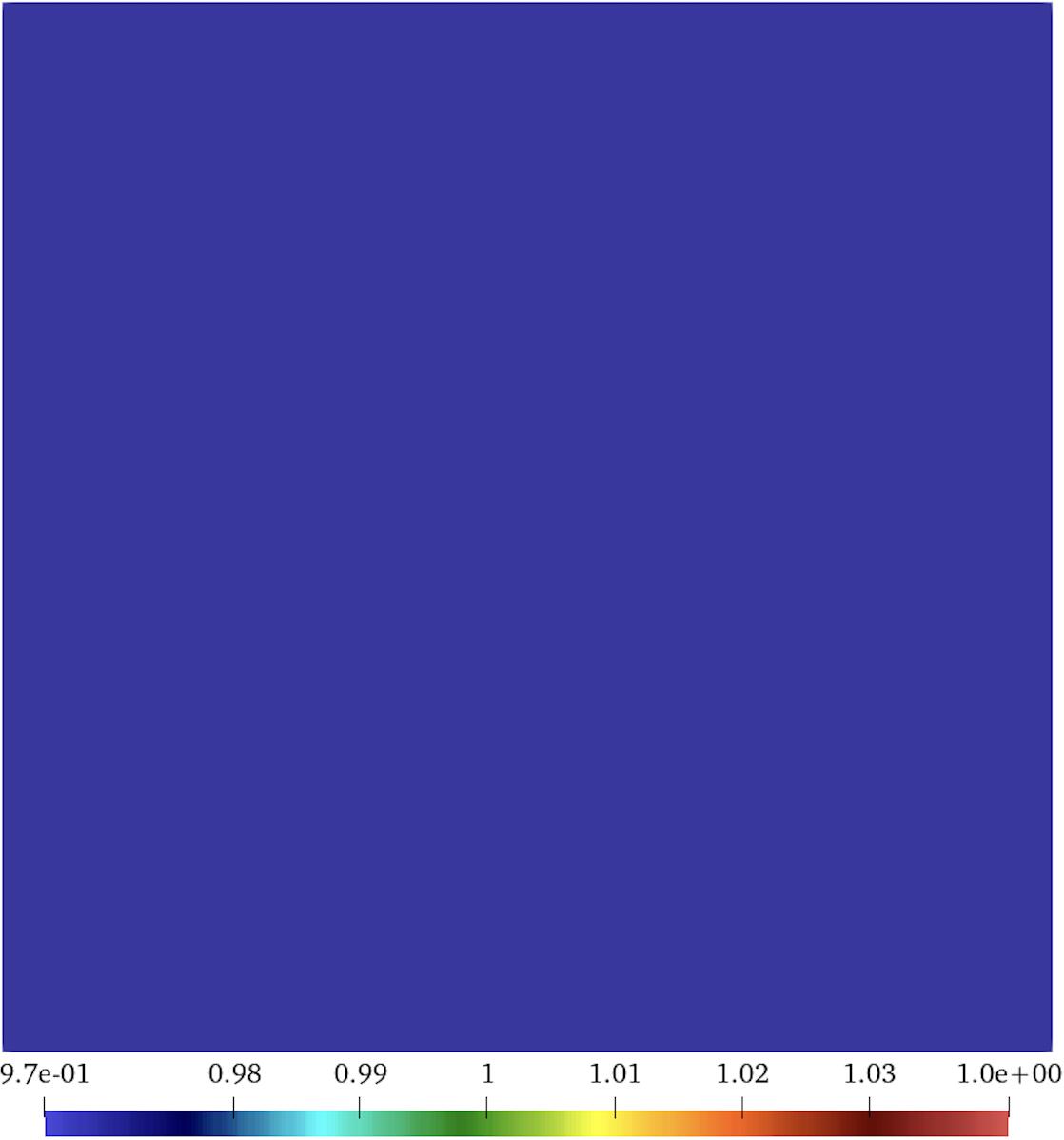}
    \put(35,103){\footnotesize{$t = 2$}}
    \end{overpic}
    \begin{overpic}[percent,width=0.19\textwidth]{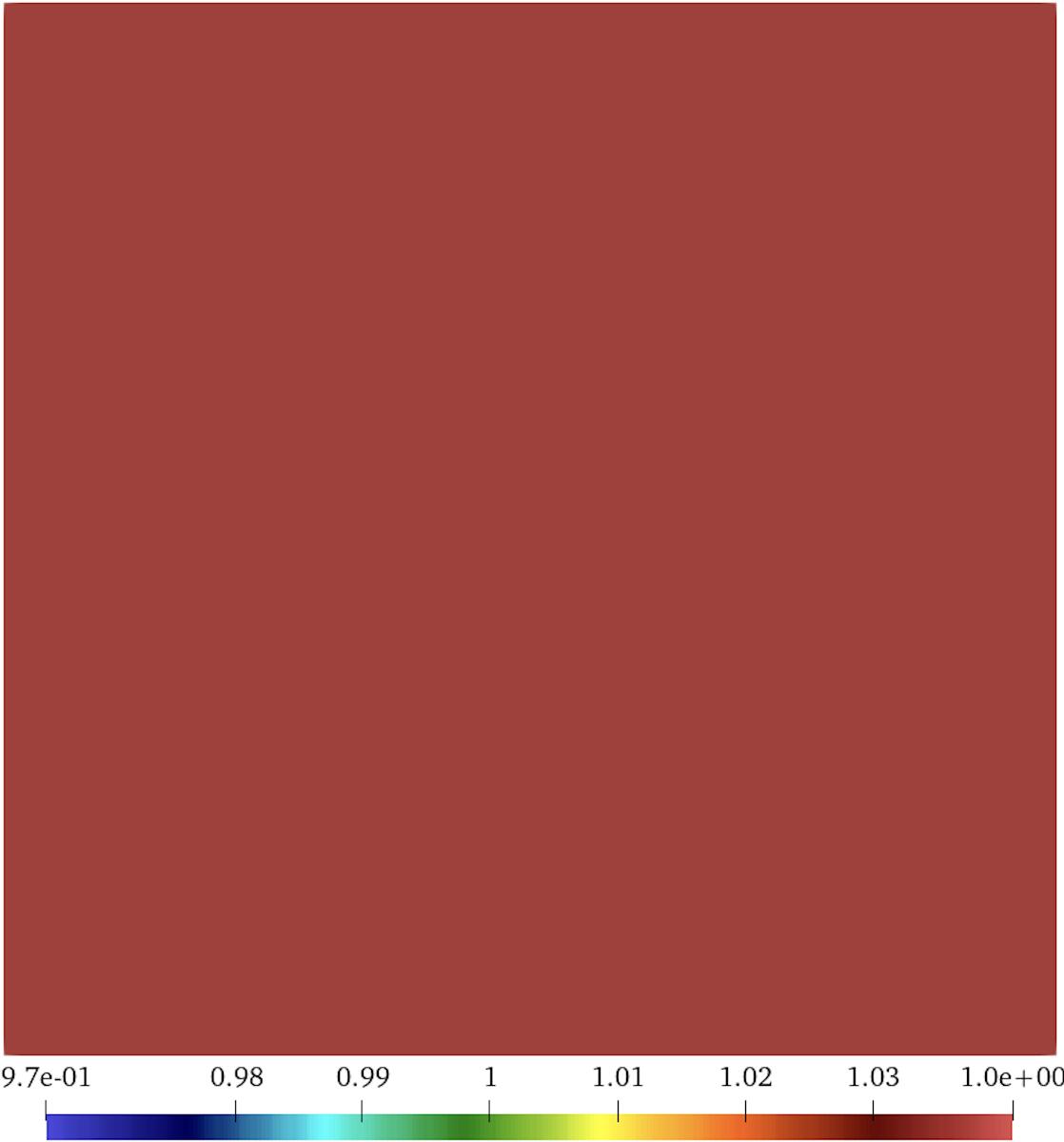}
    \put(35,103){\footnotesize{$t = 4$}}
    \end{overpic} 
        \begin{overpic}[percent,width=0.19\textwidth]{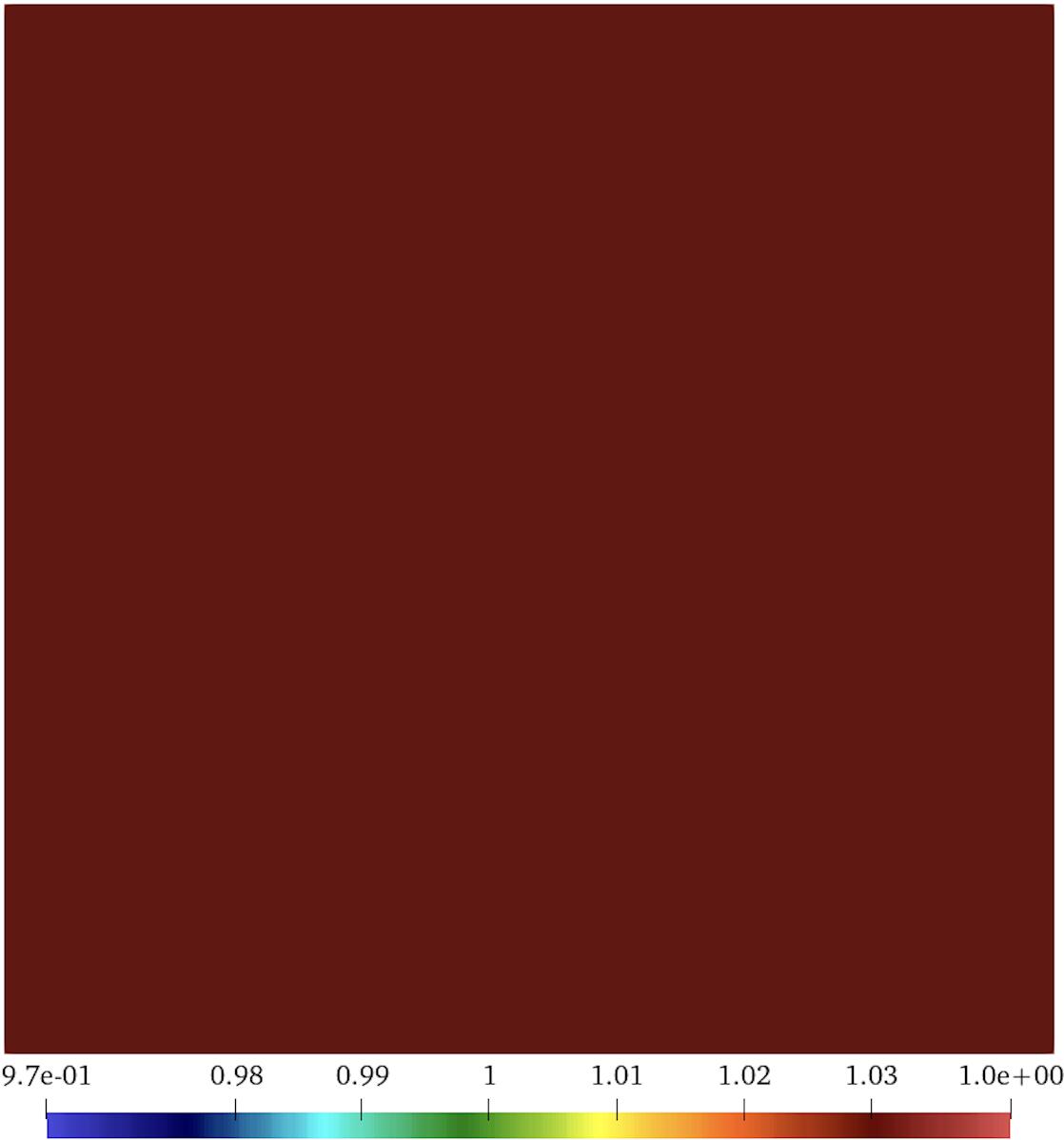}
    \put(32,103){\footnotesize{$t = 10$}}
    \end{overpic} 
    \begin{overpic}[percent,width=0.19\textwidth]{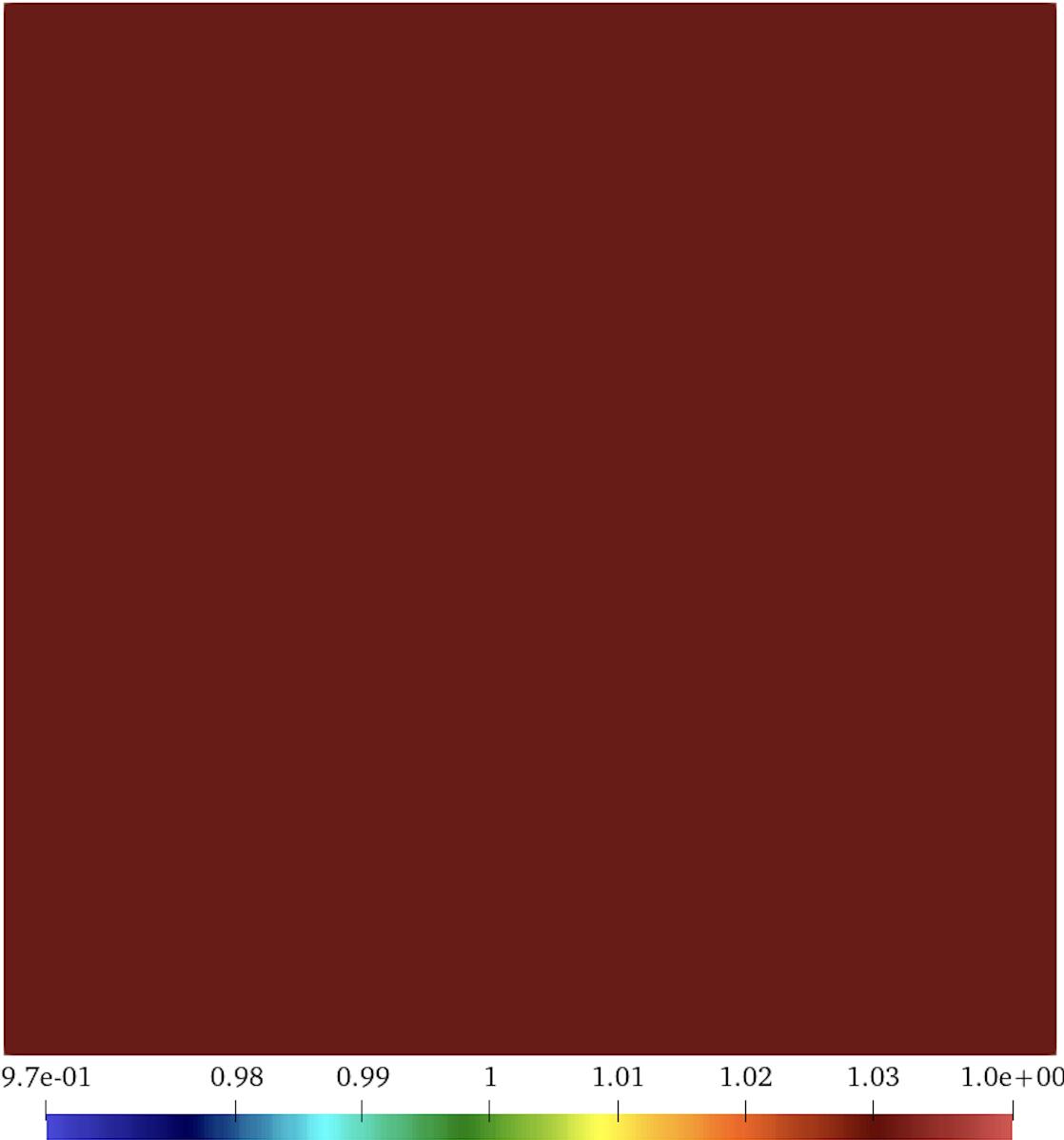}
    \put(30,103){\footnotesize{$t = 200$}}
    \end{overpic}
    \\
         \begin{overpic}[percent,width=0.19\textwidth]{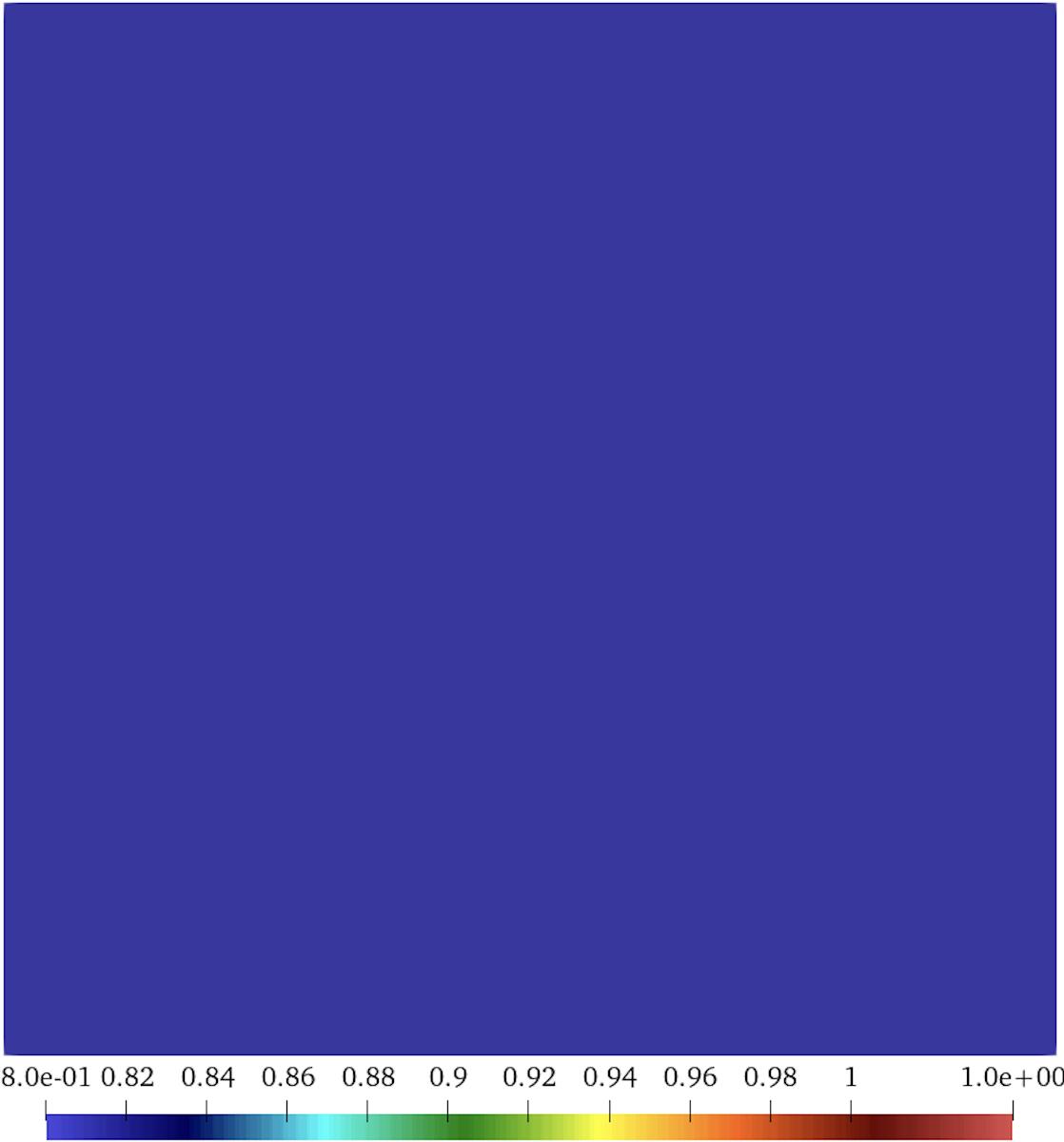}
    \put(-10,55){\footnotesize{$\rho$}}
    \end{overpic} 
     \begin{overpic}[percent,width=0.19\textwidth, grid=false]{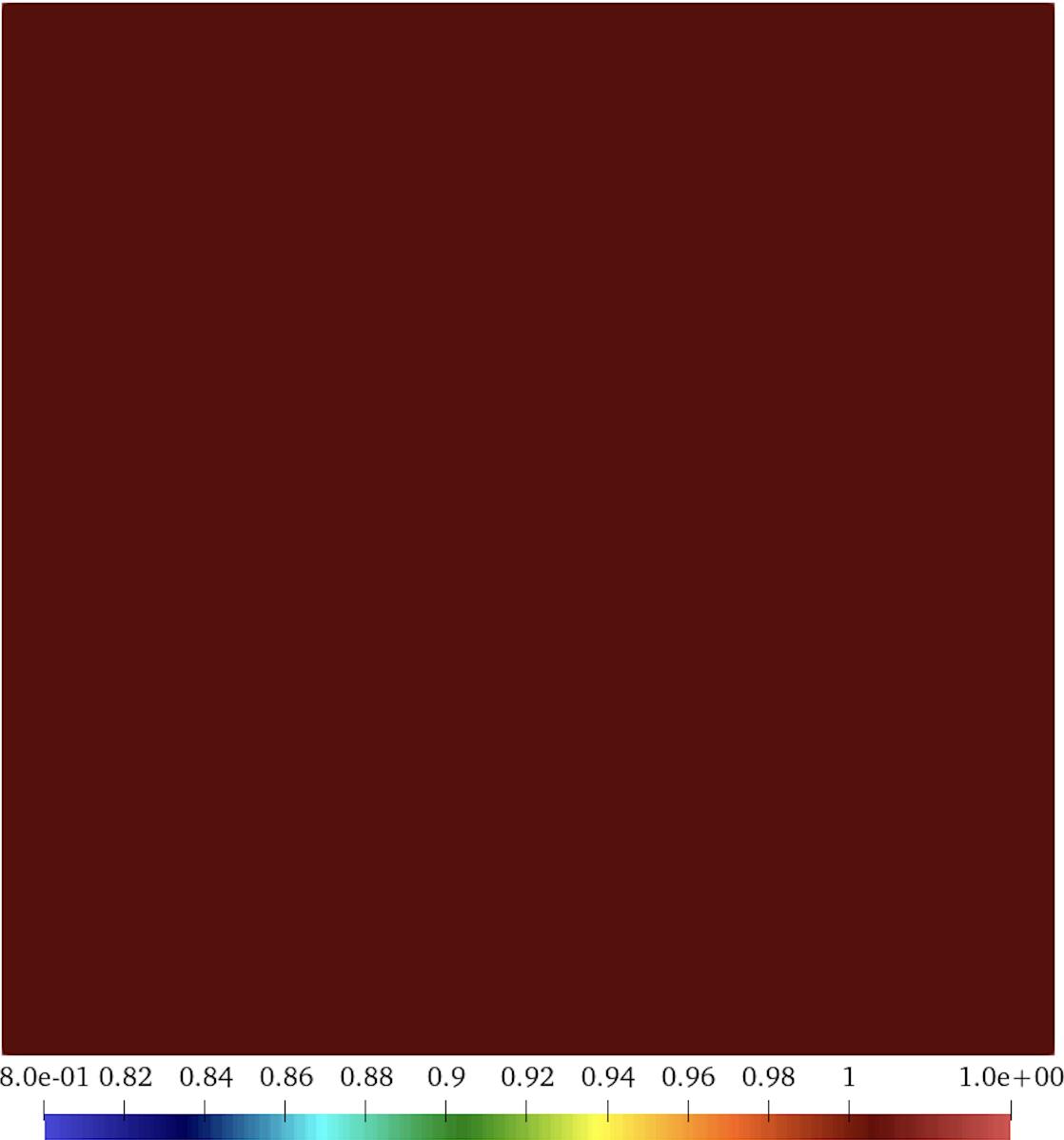}
    \end{overpic}
    \begin{overpic}[percent,width=0.19\textwidth]{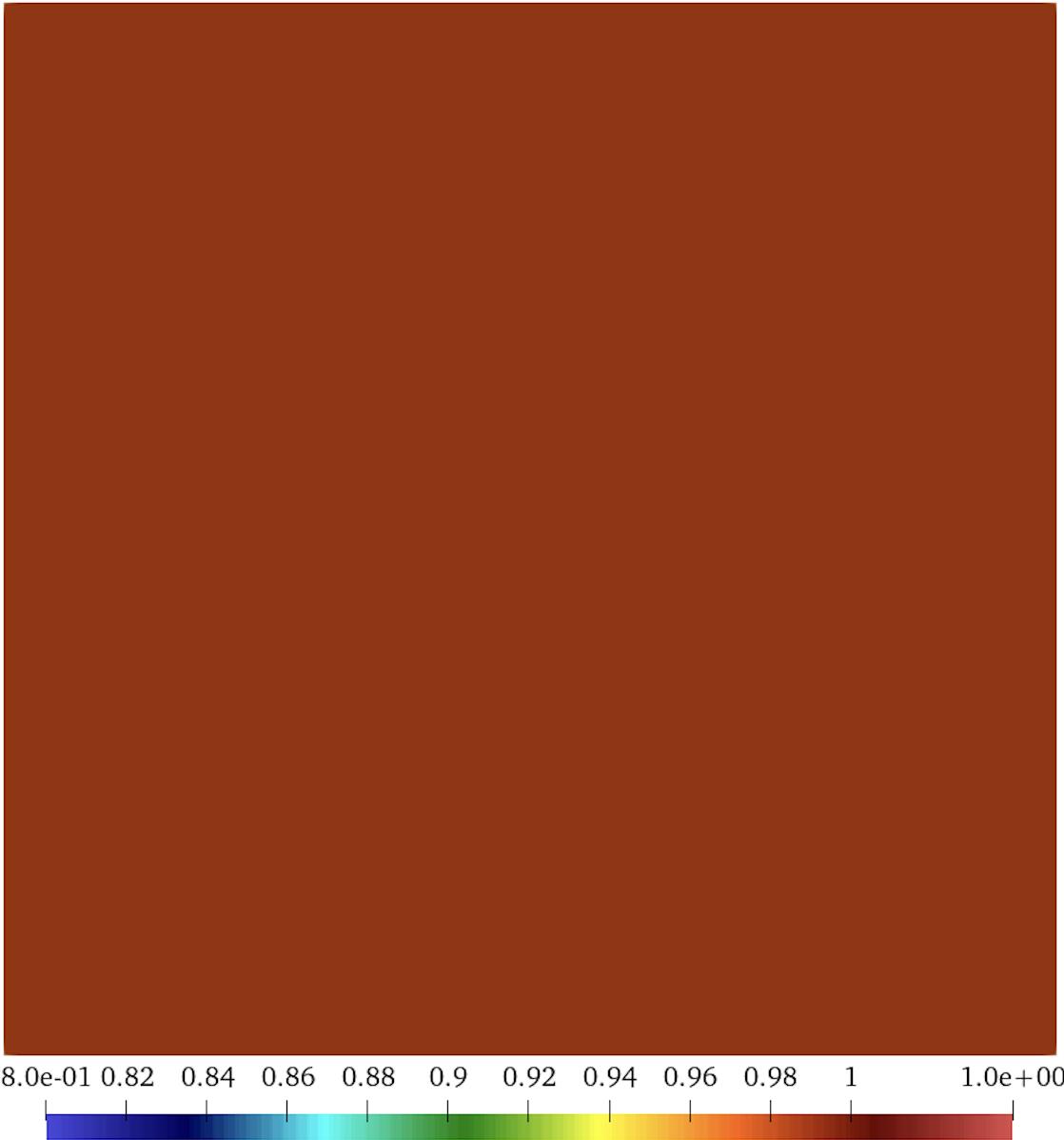}
    \end{overpic} 
        \begin{overpic}[percent,width=0.19\textwidth]{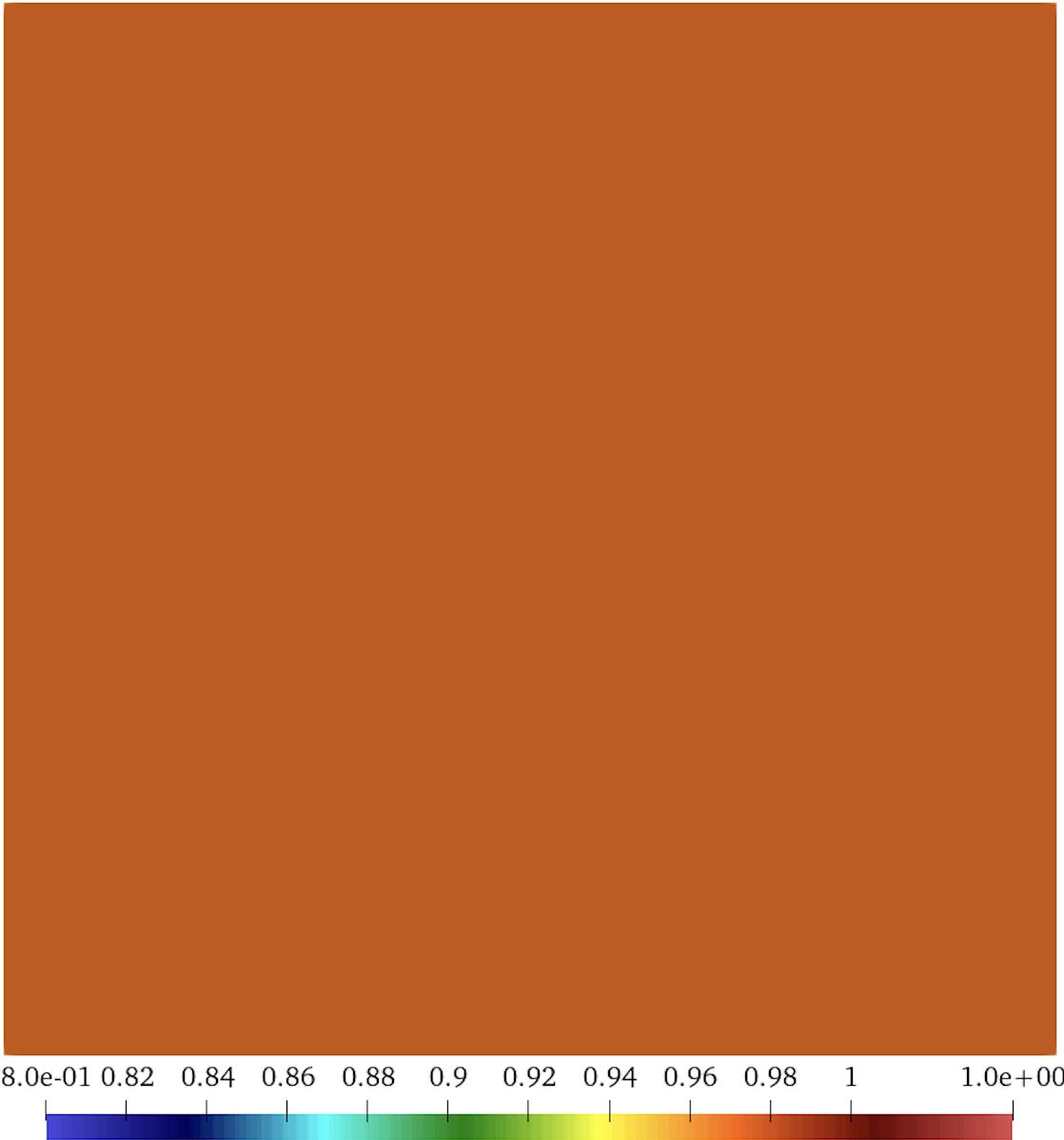}
    \end{overpic} 
    \begin{overpic}[percent,width=0.19\textwidth]{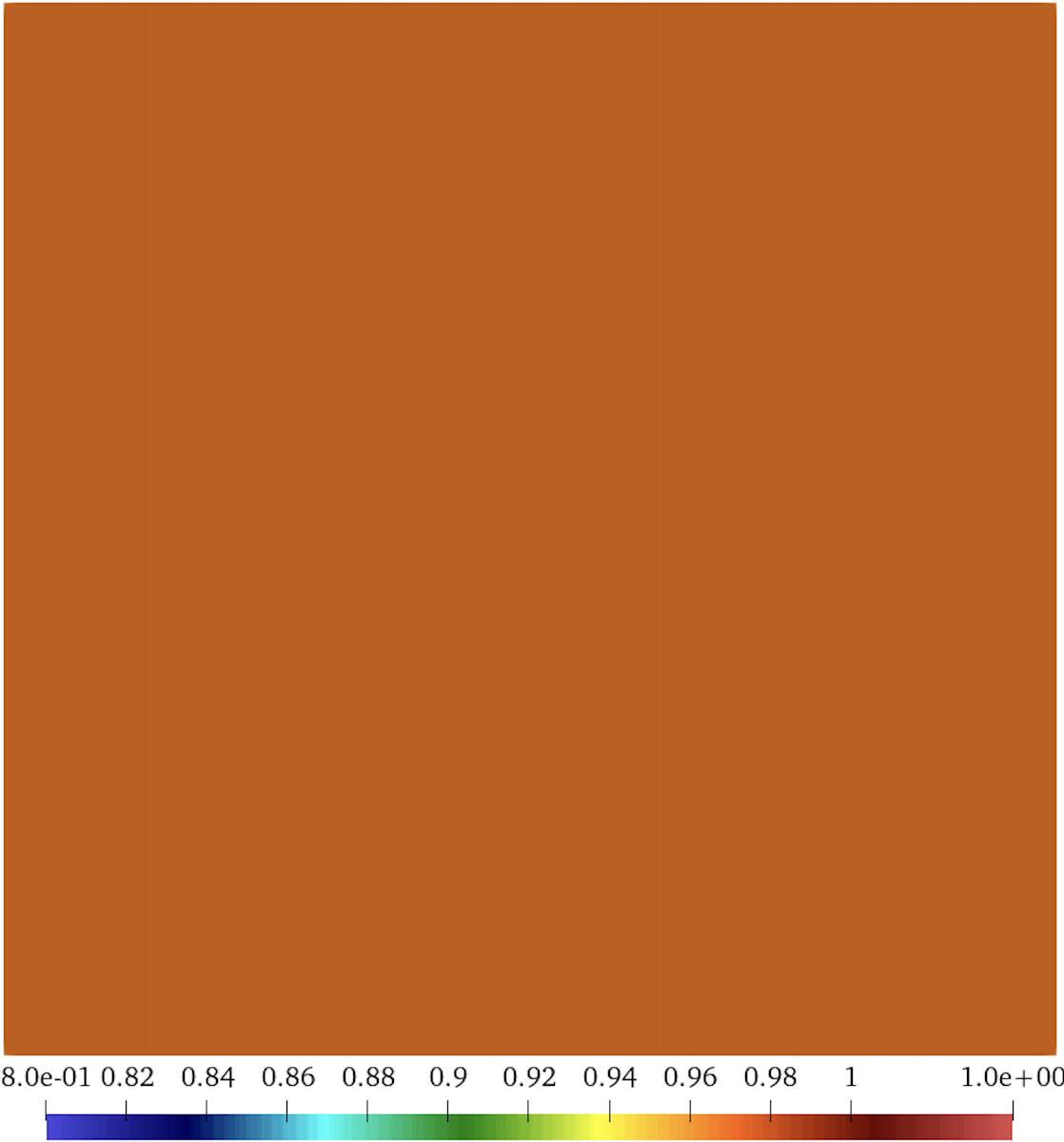}
    \end{overpic}
    \caption{Case 1 - PDE model: Evolution of the level of attractiveness (top) and density (bottom).}
    \label{fig:case1}
\end{figure}

Fig.~\ref{fig:case1} shows the evolution of the level of attractiveness and density for case 1 given by the PDE solver.
We see that the system quickly reaches a steady-state solution that is uniform in space, i.e., no hotspots form. 
The computed value of $A$ shown in Fig.~\ref{fig:case1} is roughly $1.0333$, which approximates $\bar{A} = A^{st} + \Gamma \theta /\omega^2  = 1/30 + 1$, while the computed value of $\rho$ is roughly 0.9677, which approximates $\bar{\rho} = \Gamma \theta/(\omega^2\bar{A}) = 1/(1/30 + 1)$. Recall that these $\bar{A}$ and $\bar{\rho}$ satisfy \eqref{eq:continous_bertozzi-1ss}-\eqref{eq:continous_bertozzi-2ss}. Since we are not adding noise, initial conditions \eqref{eq:ic1}-\eqref{eq:ic2} become $A = \bar{A}$ and $\rho = \rho_0 \neq \bar{\rho}$. The difference between  $\rho_0$ and $\bar{\rho}$ explains the transition to the steady state observed in Fig.~\ref{fig:case1}.

If one adds noise to the initial 
data of case 1, the computed solutions look similar to the ones reported in Fig.~\ref{fig:case1}, with the addition of a small noisy perturbation that we suspect would be dampened if the simulation time is long enough.

\begin{figure}[htb!]
     \centering
         \begin{overpic}[percent,width=0.19\textwidth]{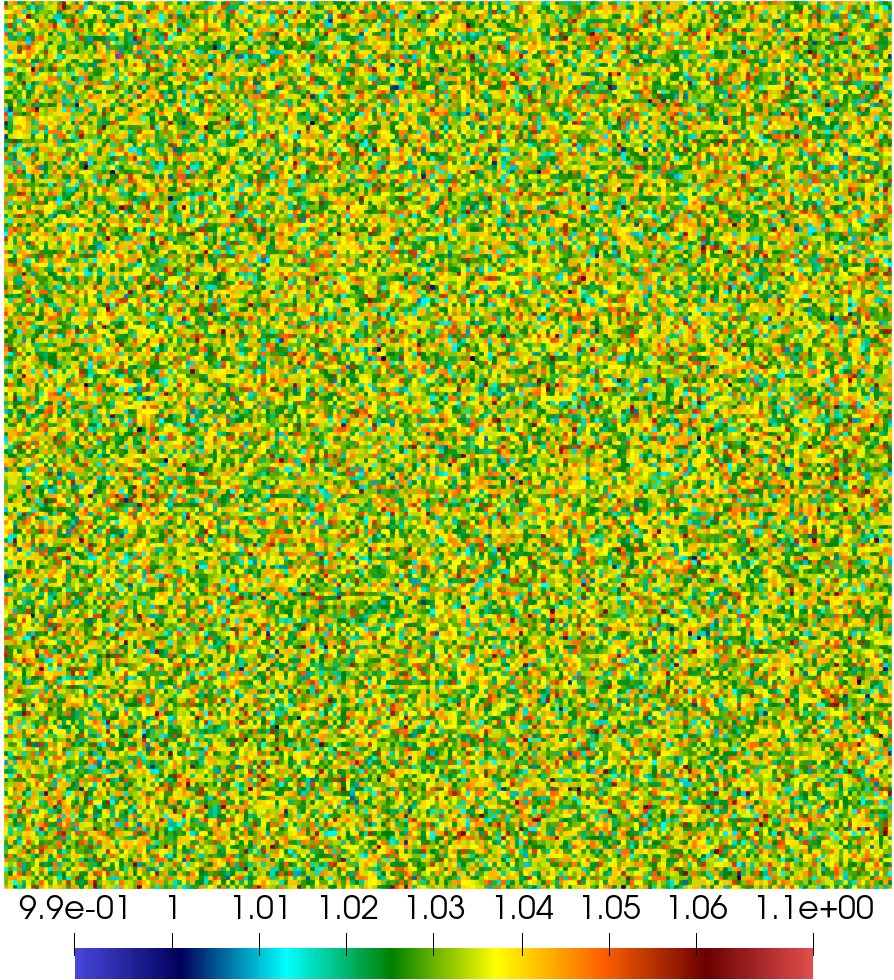}
    \put(30,103){\footnotesize{$t = 0$}}
    \put(-10,55){\footnotesize{$A$}}
    \end{overpic} 
     \begin{overpic}[percent,width=0.19\textwidth, grid=false]{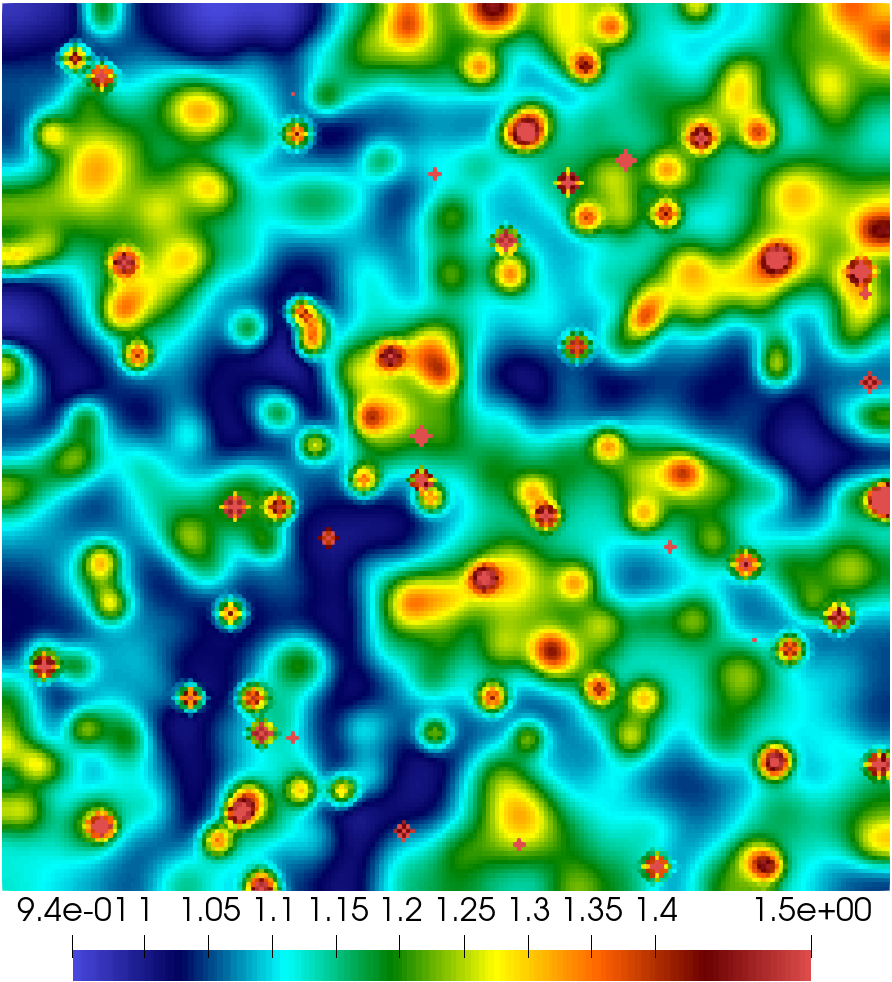}
    \put(35,103){\footnotesize{$t = 3.6$}}
    \end{overpic}
    \begin{overpic}[percent,width=0.19\textwidth]{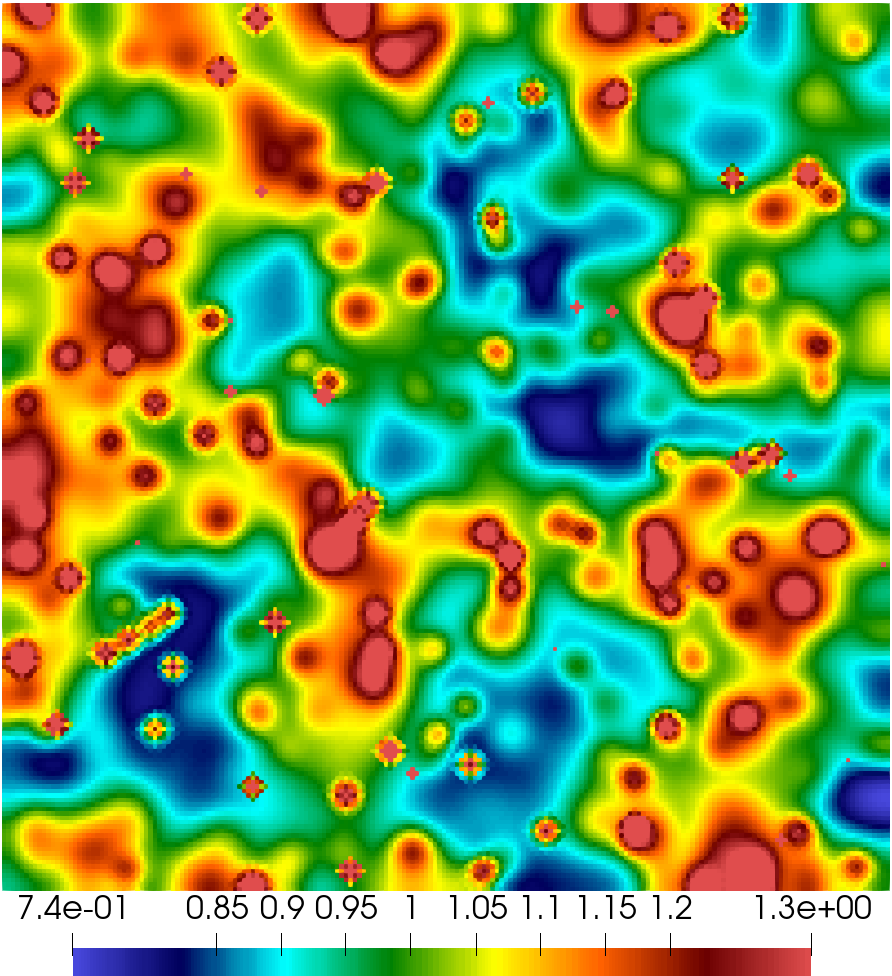}
    \put(35,103){\footnotesize{$t = 10$}}
    \end{overpic} 
        \begin{overpic}[percent,width=0.19\textwidth]{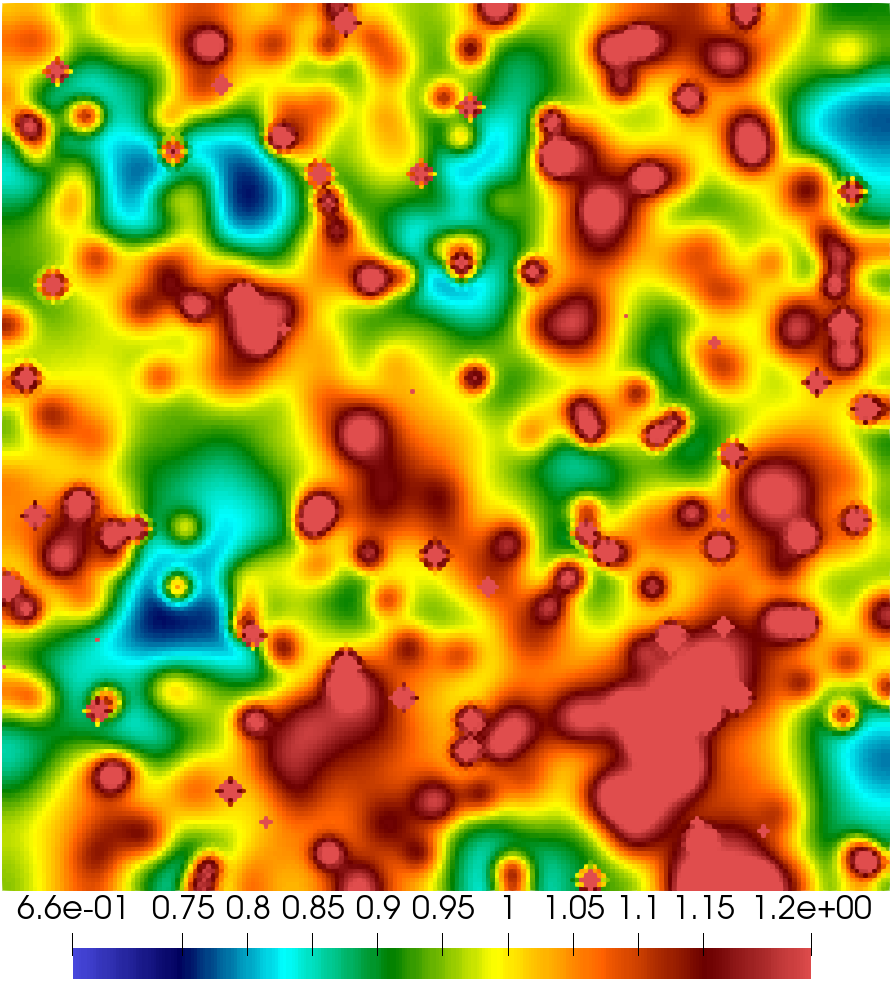}
    \put(32,103){\footnotesize{$t = 15$}}
    \end{overpic} 
    \begin{overpic}[percent,width=0.19\textwidth]{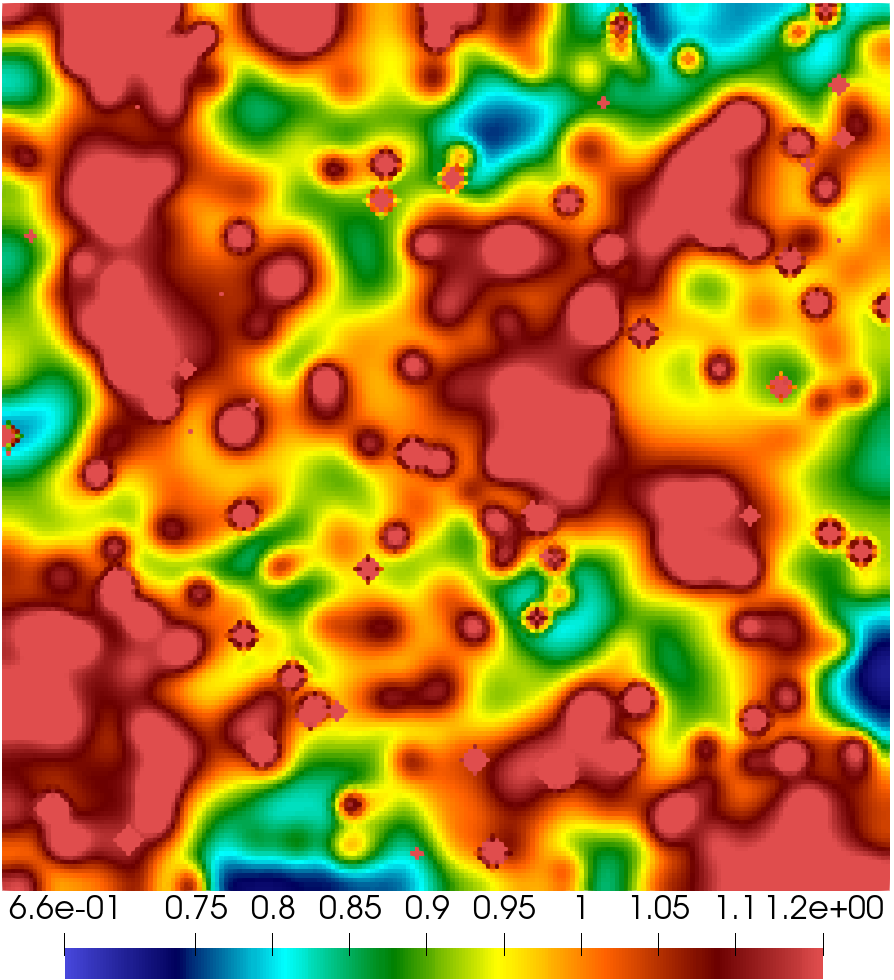}
    \put(30,103){\footnotesize{$t = 200$}}
    \end{overpic}
    \\
         \begin{overpic}[percent,width=0.19\textwidth]{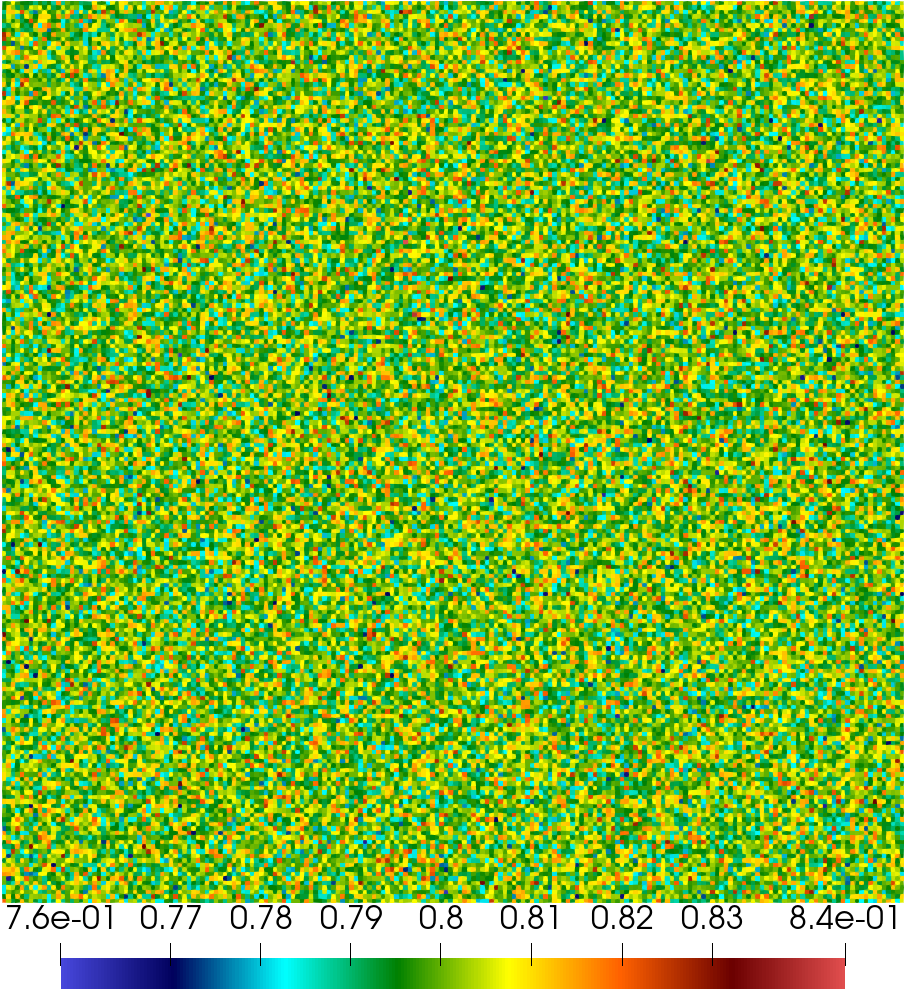}
    \put(-10,55){\footnotesize{$\rho$}}
    \end{overpic} 
     \begin{overpic}[percent,width=0.19\textwidth, grid=false]{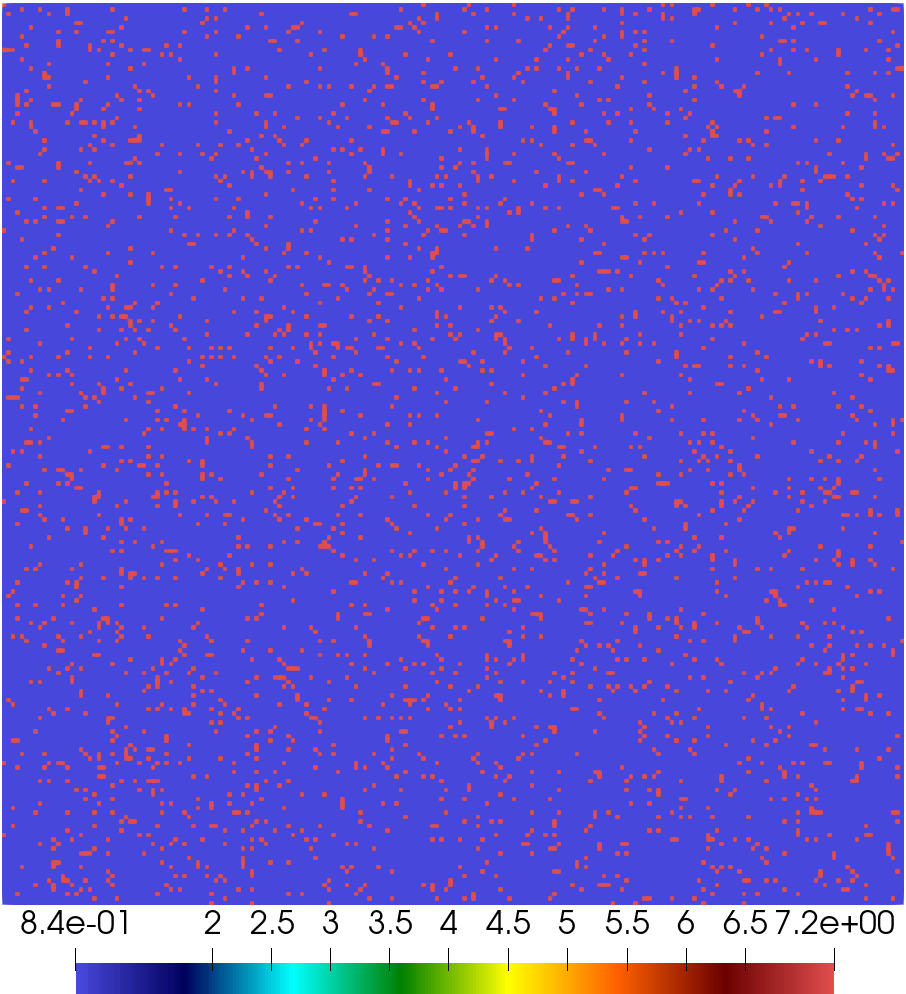}
    \end{overpic}
    \begin{overpic}[percent,width=0.19\textwidth]{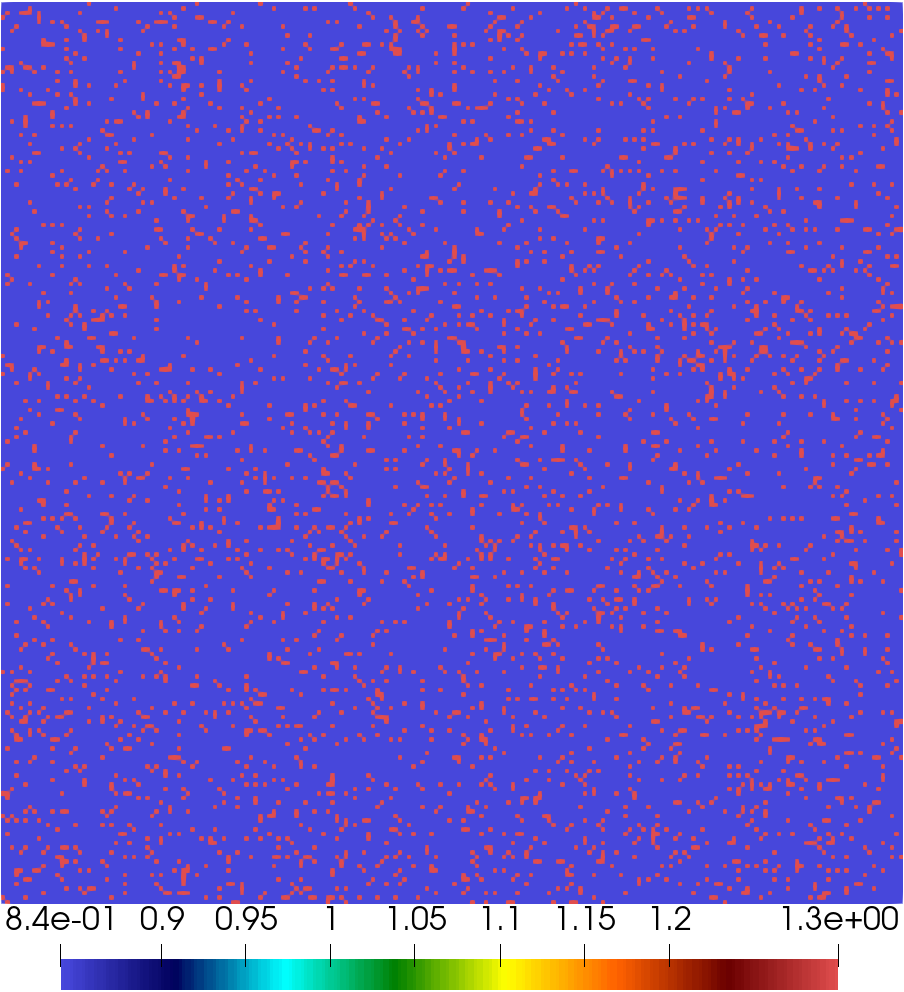}
    \end{overpic} 
        \begin{overpic}[percent,width=0.19\textwidth]{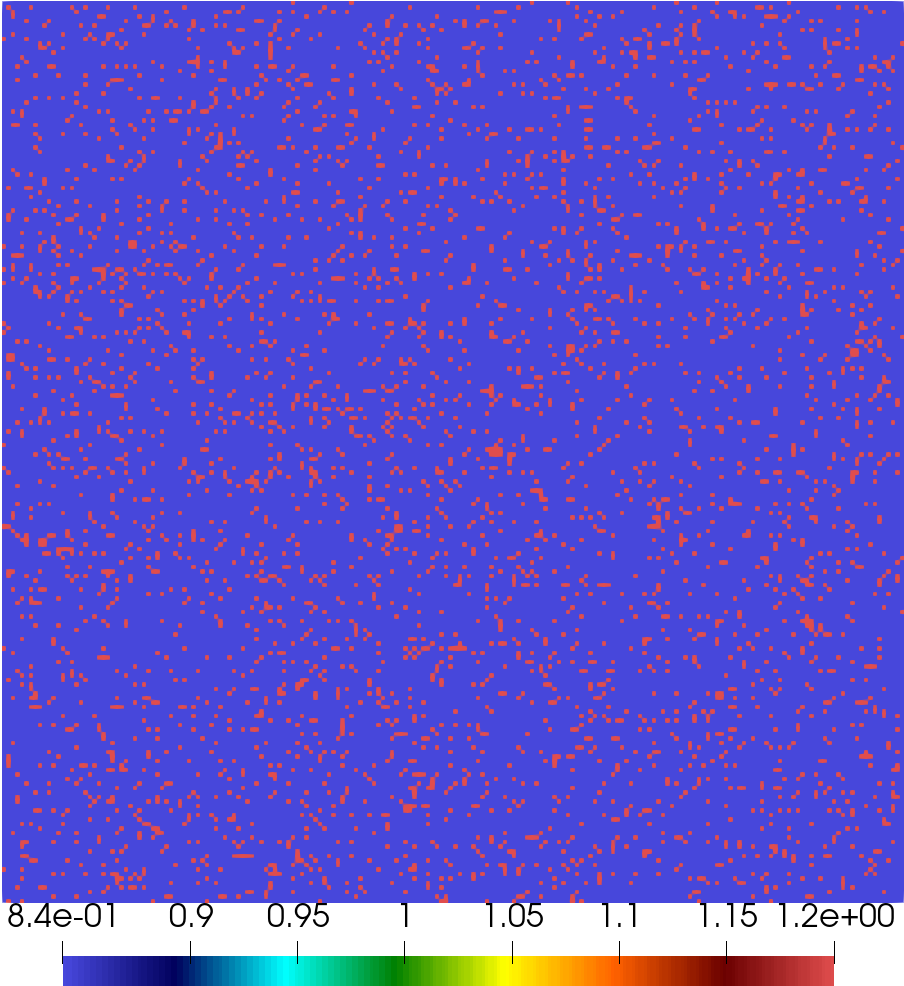}
    \end{overpic} 
    \begin{overpic}[percent,width=0.19\textwidth]{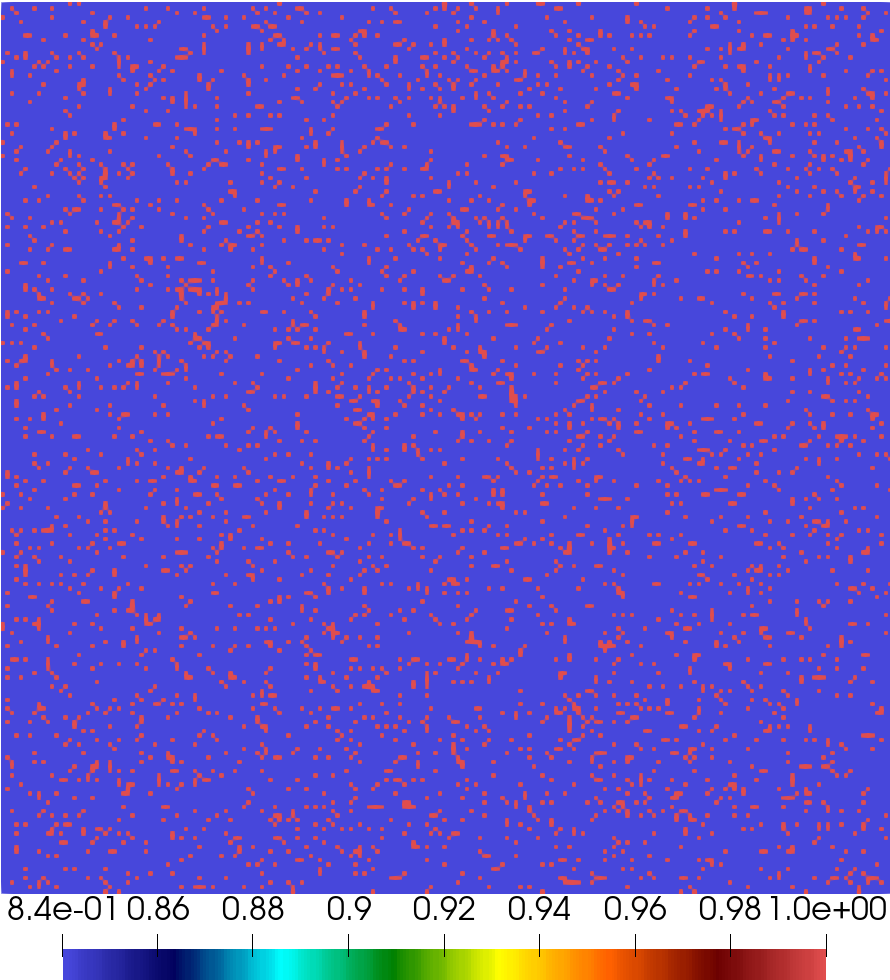}
    \end{overpic}
    \caption{Case 1 - agent-based model: Evolution of the level of attractiveness (top) and criminals (bottom).}
    \label{fig:case1DM}
\end{figure}

\begin{figure}[htb!]
     \centering
         \begin{overpic}[percent,width=0.19\textwidth]{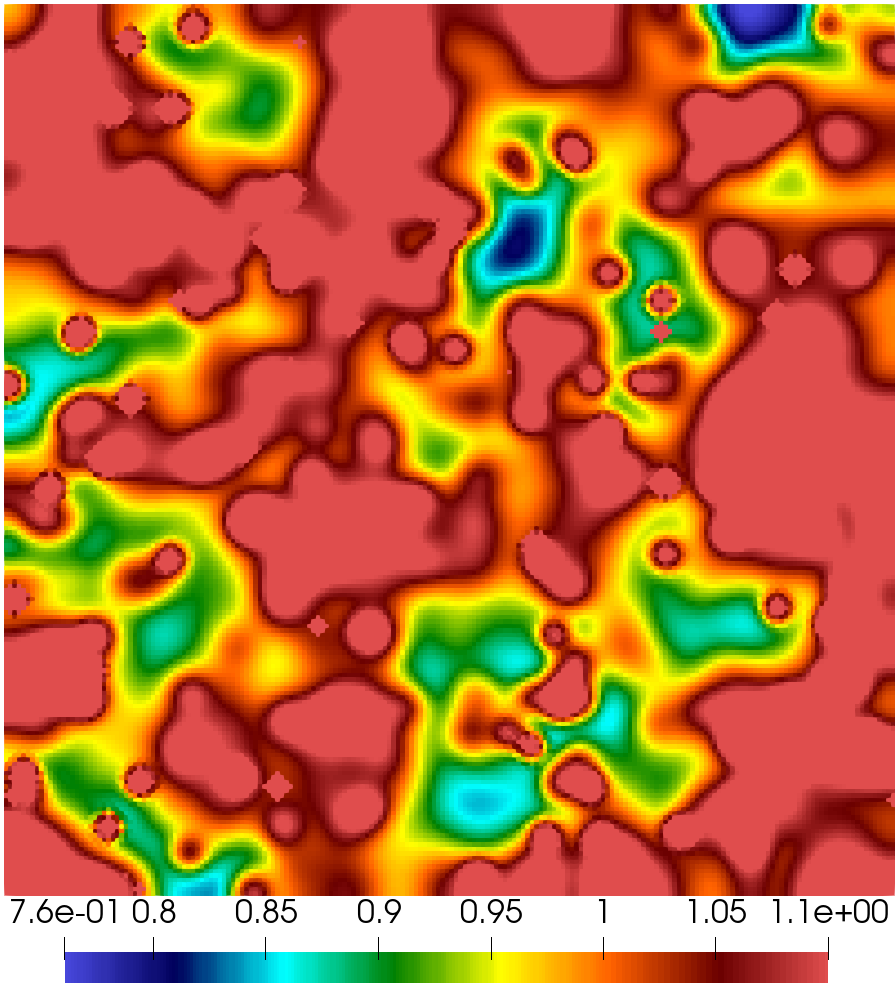}
    \put(30,103){\footnotesize{$t = 225$}}
    \end{overpic} 
     \begin{overpic}[percent,width=0.19\textwidth, grid=false]{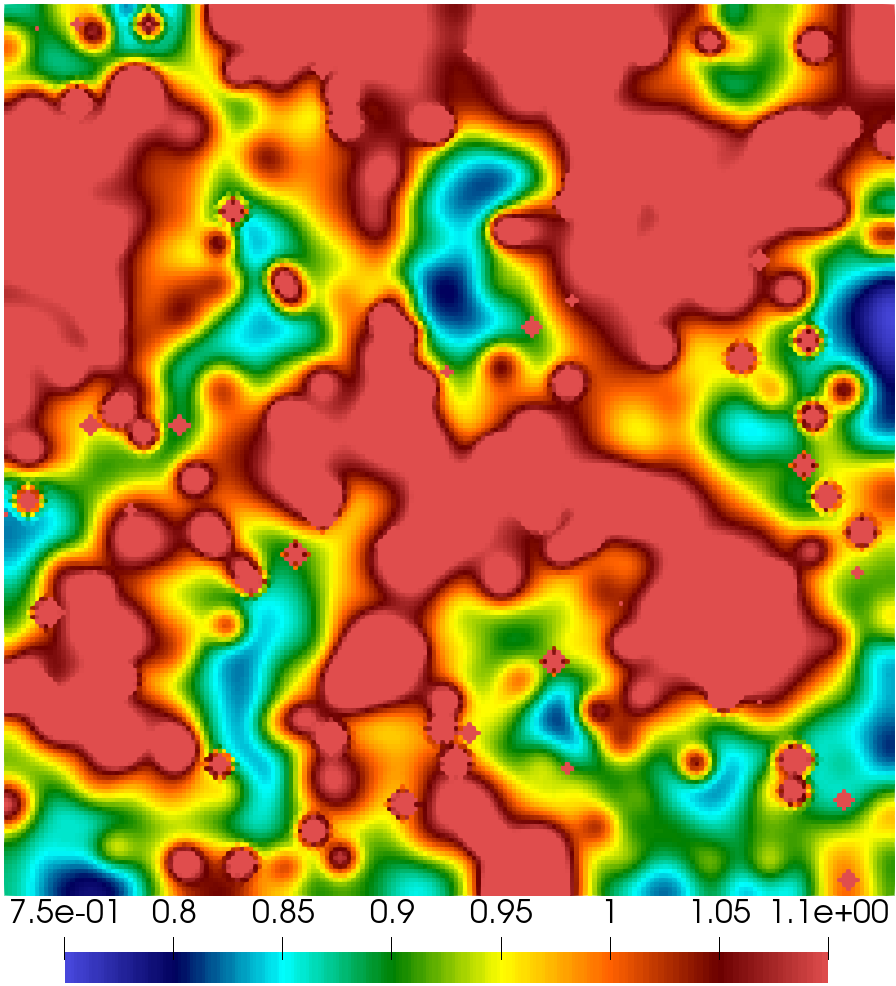}
    \put(35,103){\footnotesize{$t = 250$}}
    \end{overpic}
    \begin{overpic}[percent,width=0.19\textwidth]{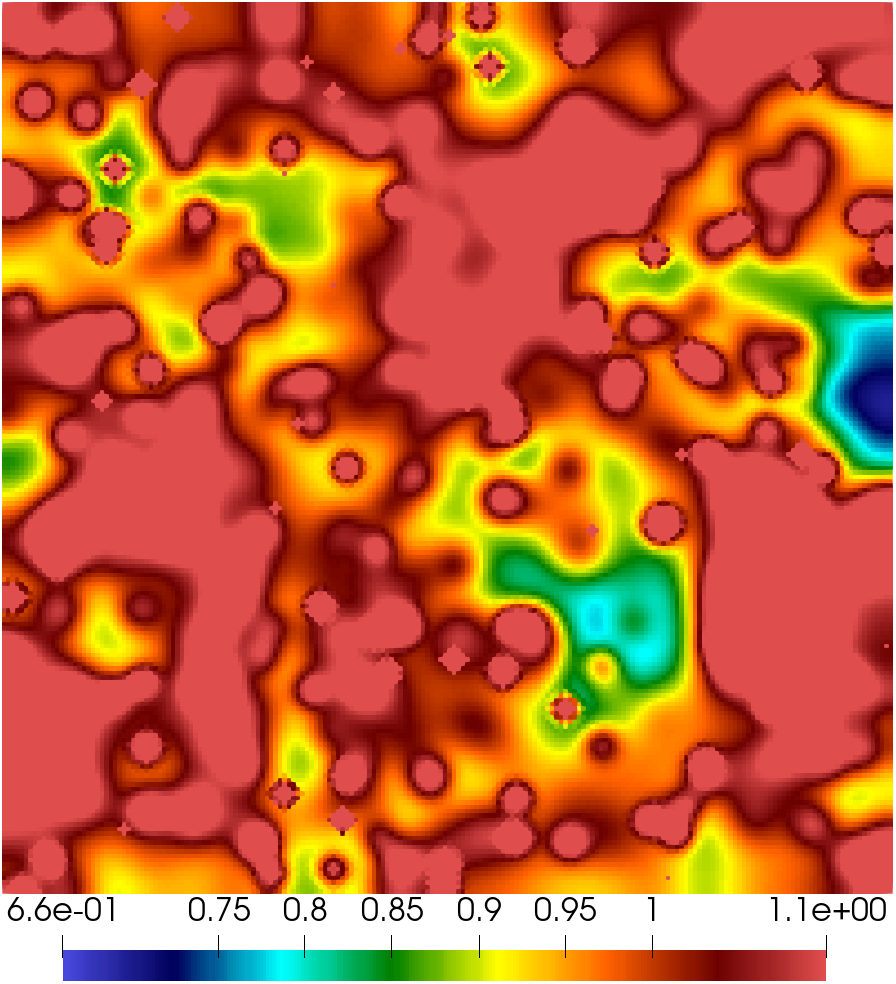}
    \put(35,103){\footnotesize{$t = 275$}}
    \end{overpic} 
        \begin{overpic}[percent,width=0.19\textwidth]{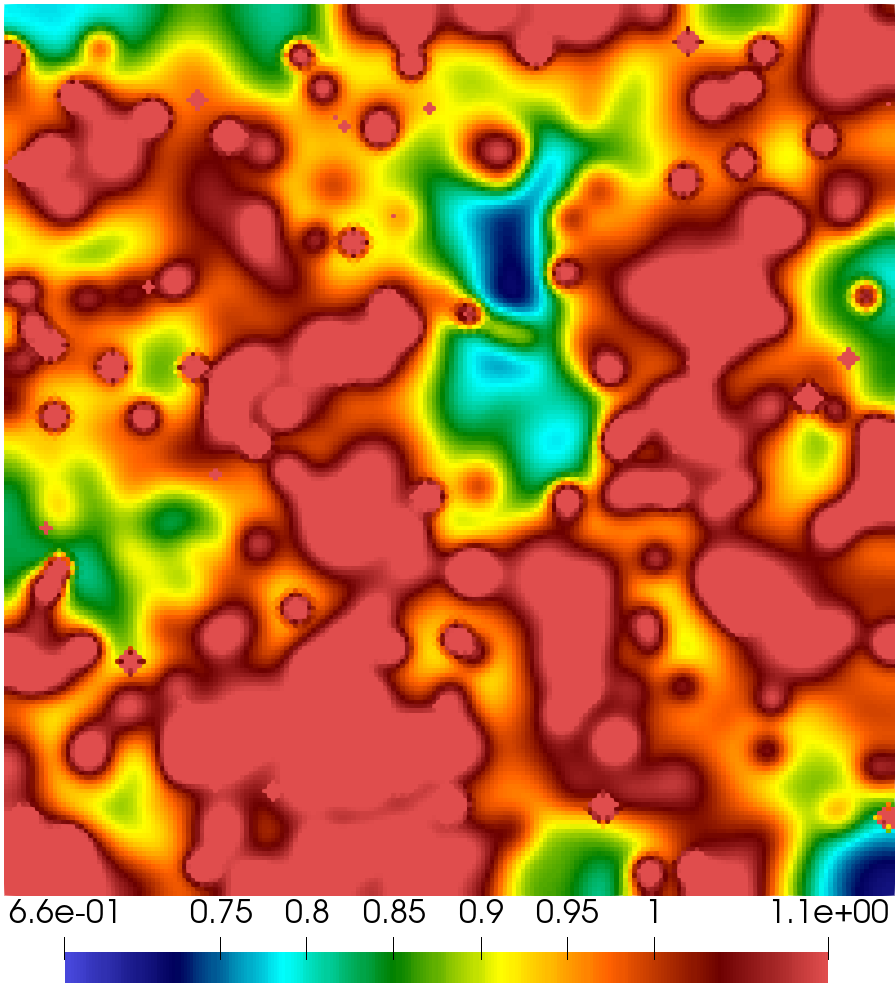}
    \put(32,103){\footnotesize{$t = 300$}}
    \end{overpic} 
    \caption{Case 1 - agent-based model: Evolution of $A$ past $t = 200$.}
    \label{fig:case1DM-bis}
\end{figure}

The initial conditions and the parameters used for the results in Fig.~\ref{fig:case1} are suitably dimensionalized and passed to the agent-based solver, which uses all dimensional quantities. For the agent-based solver, we set $\Gamma = 0.0077$, $\theta = 0.58$, and $\omega = 1/15$ to obtain $\Gamma \theta /\omega^2 \approx 1$, as required by case 1. The evolution of $A$ and $\rho$ for case 1 given by the agent-based solver after non-dimensionalization are shown in Fig.~\ref{fig:case1DM}. We see that, while
$A$ is not homogeneous in space like in Fig.~\ref{fig:case1},
the attractiveness field has at most places a value close to $\bar{A}$. Any local variation in values disappears rather quickly. 
See Fig.~\ref{fig:case1DM-bis} for $A$
past $t = 200$.
The same is true for $\rho$.
We note that also in \cite{short2008statistical}
a spatially homogeneous PDE solution in Fig.~4(a) corresponds to parameter values that give agent-based solutions
not exactly spatially homogeneous, shown
in Fig.~3(a) and 3(b).

The stability analysis from 
\cite{short2008statistical} suggests 
that to see the formation of hotspots, one should lower the value of $\eta$ to satisfy condition \eqref{eq:suff_cond}. While that analysis assumes periodic boundary conditions, which we do not impose, we used condition \eqref{eq:suff_cond} as an indication to design case 2. Note that case 2 has the same parameters as case 1, with the exception of $\eta$, which is lower than in case 1. This means that case 2 is characterized by 
weaker neighborhood effects.

\begin{figure}[htb!]
     \centering
         \begin{overpic}[percent,width=0.19\textwidth]{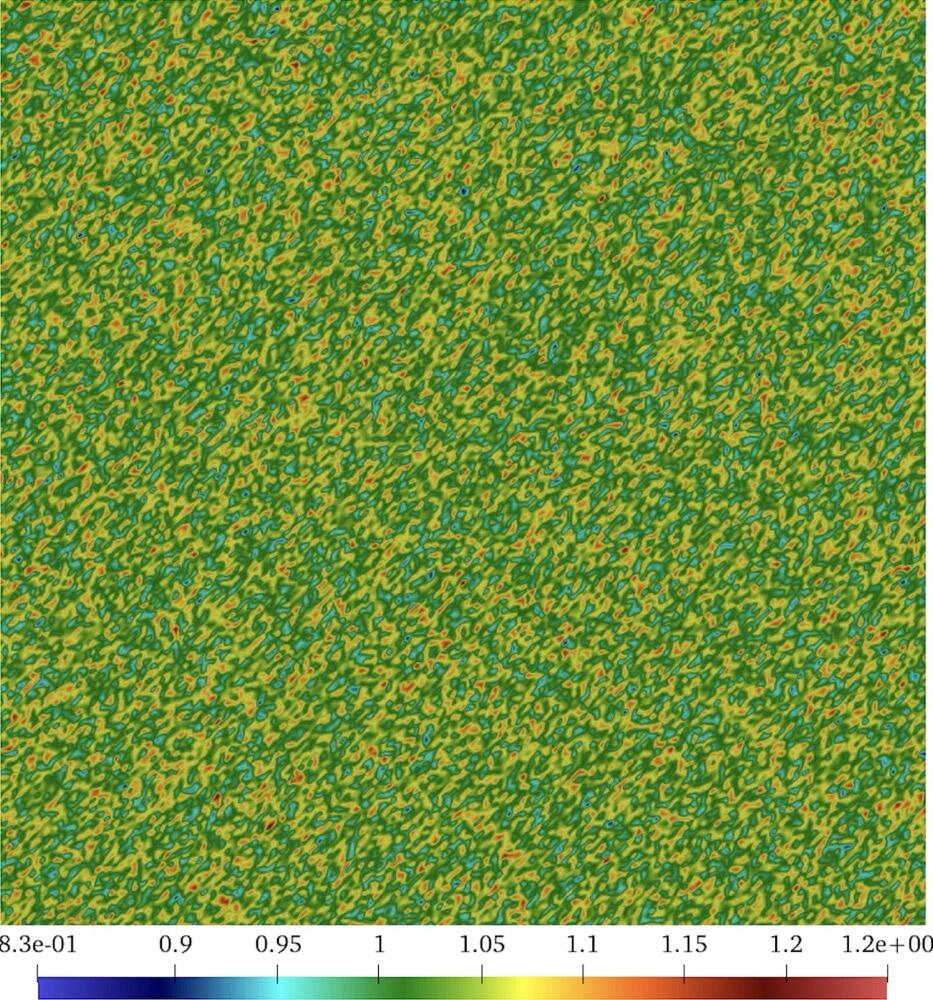}
    \put(30,103){\footnotesize{$t = 0$}}
    \put(-10,55){\footnotesize{$A$}}
    \end{overpic} 
     \begin{overpic}[percent,width=0.19\textwidth, grid=false]{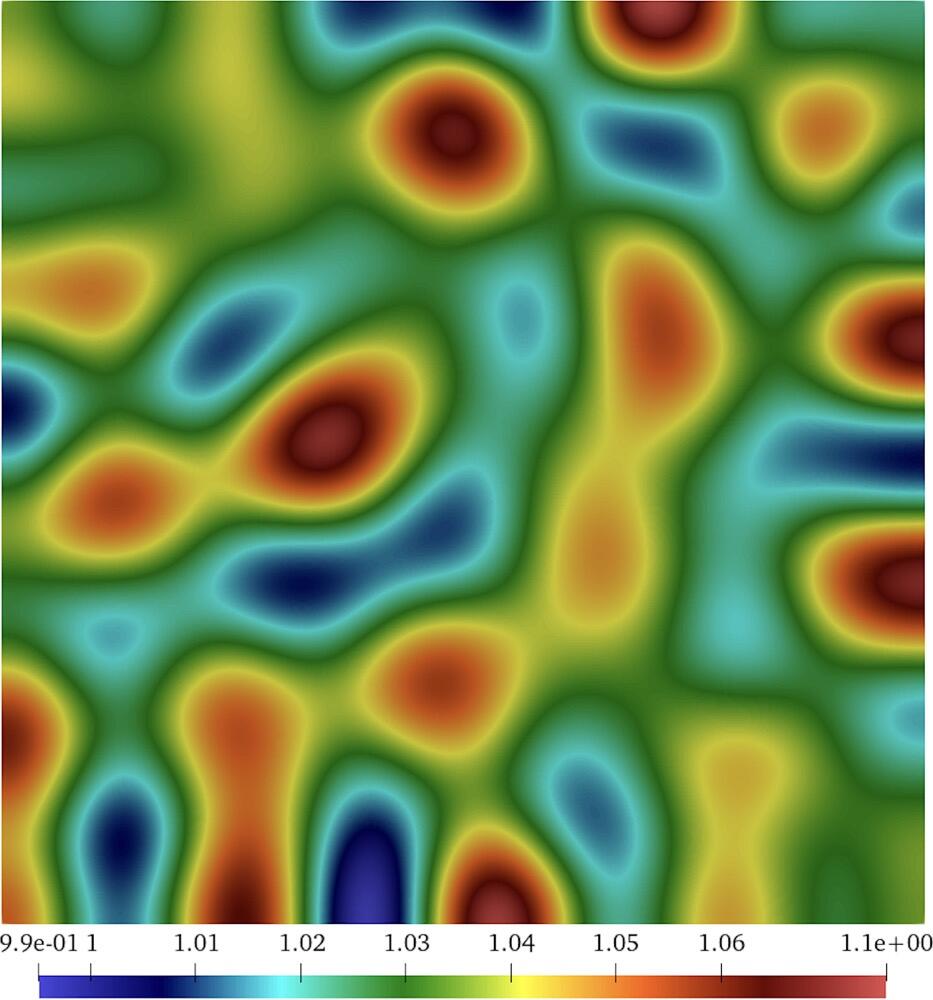}
    \put(35,103){\footnotesize{$t = 10$}}
    \end{overpic}
    \begin{overpic}[percent,width=0.19\textwidth]{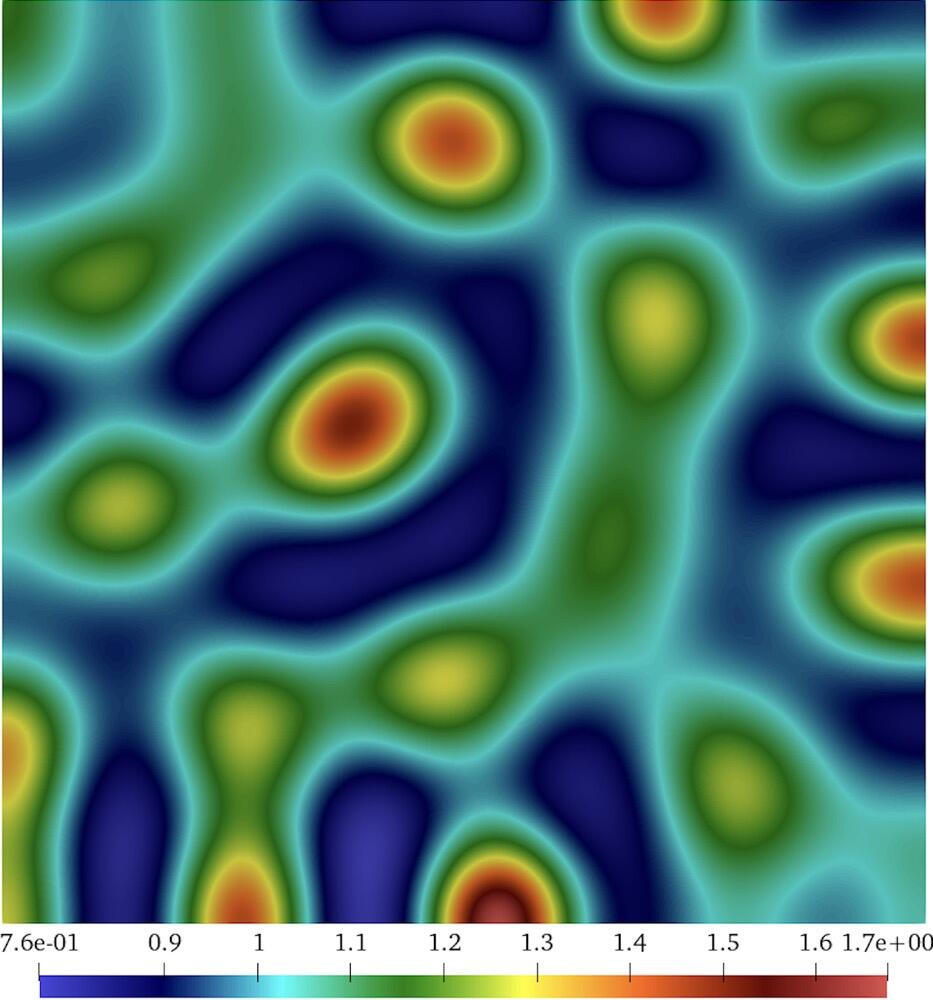}
    \put(35,103){\footnotesize{$t = 20$}}
    \end{overpic} 
        \begin{overpic}[percent,width=0.19\textwidth]{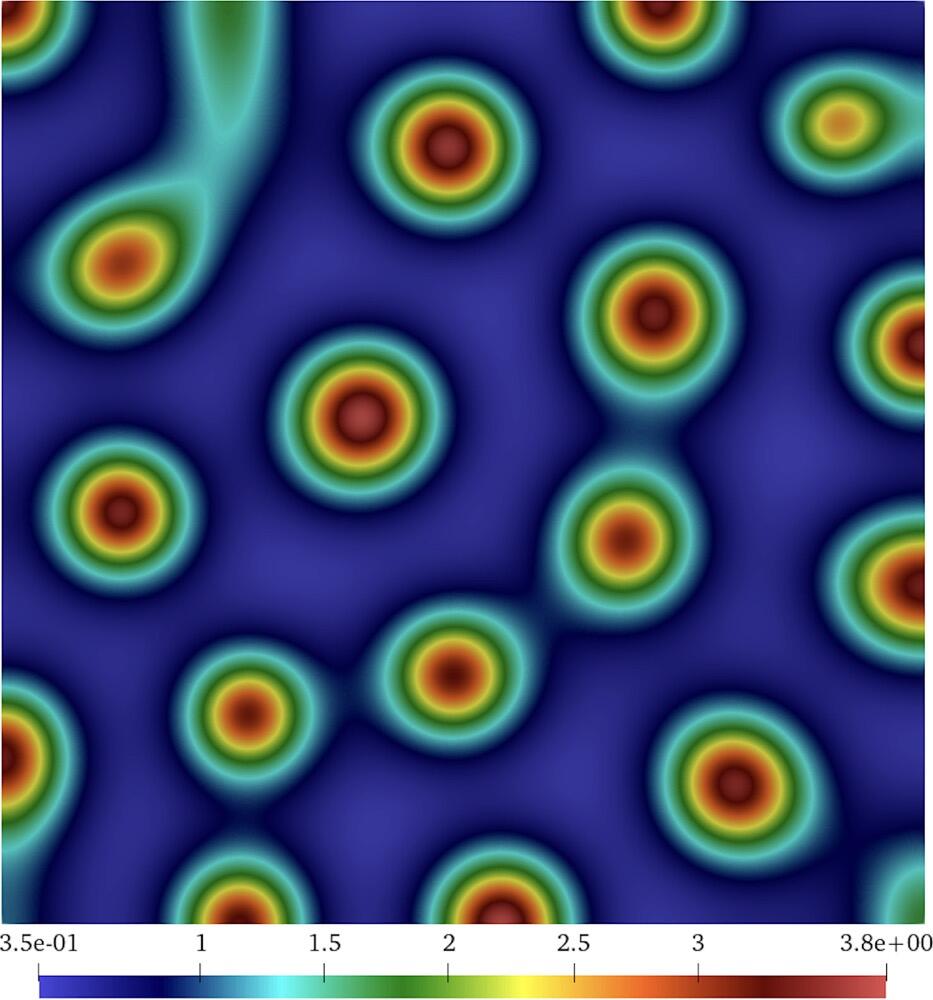}
    \put(32,103){\footnotesize{$t = 30$}}
    \end{overpic} 
    \begin{overpic}[percent,width=0.19\textwidth]{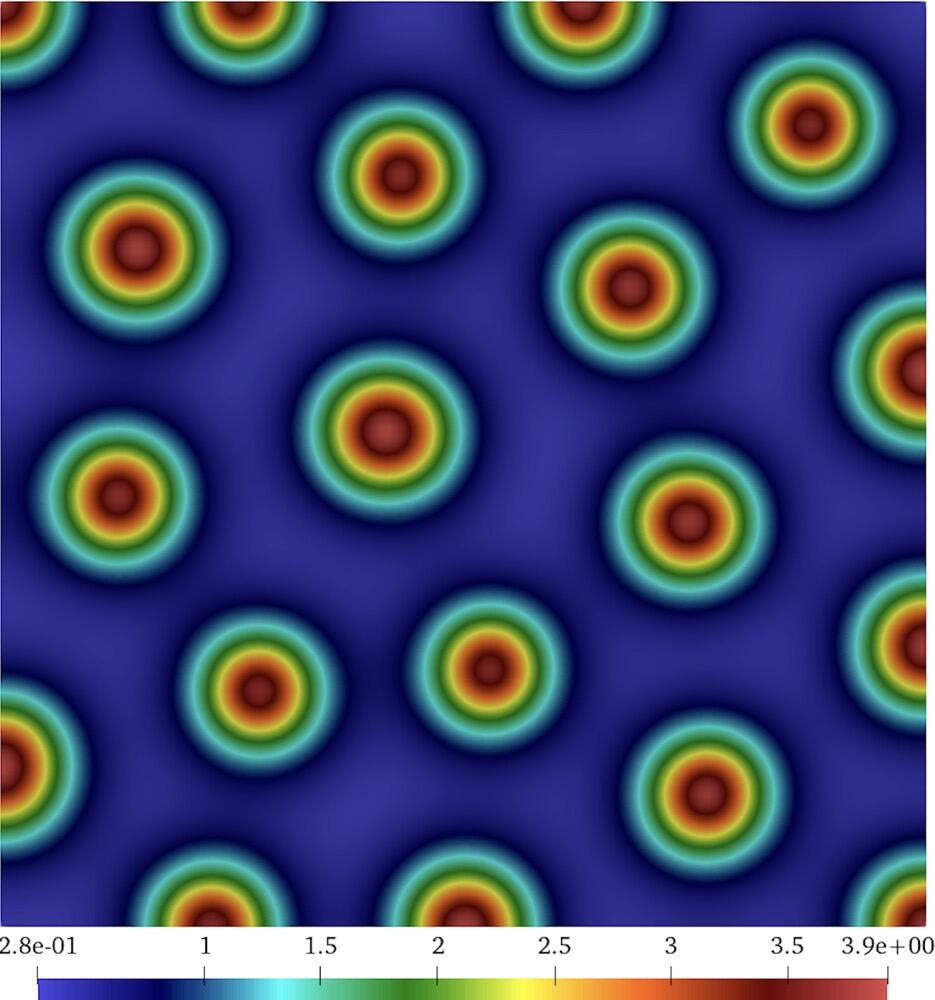}
    \put(30,103){\footnotesize{$t = 200$}}
    \end{overpic}
    \\
         \begin{overpic}[percent,width=0.19\textwidth]{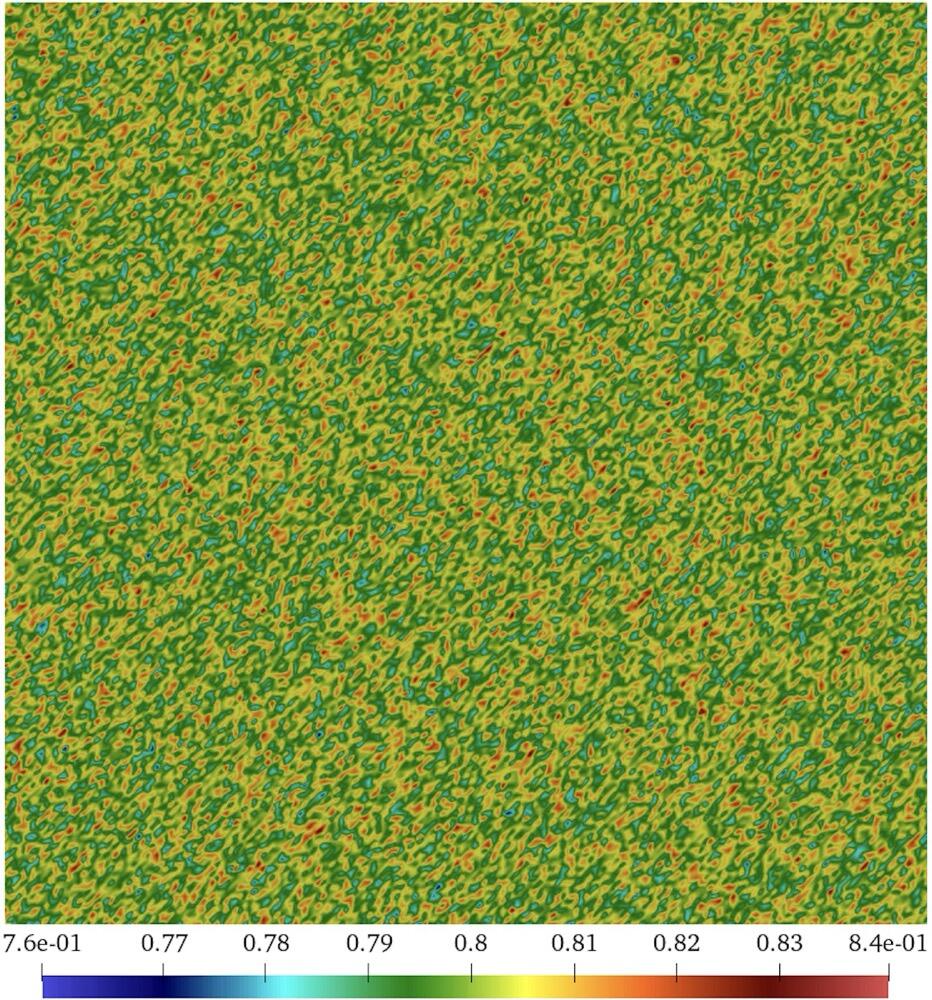}
    \put(-10,55){\footnotesize{$\rho$}}
    \end{overpic} 
     \begin{overpic}[percent,width=0.19\textwidth, grid=false]{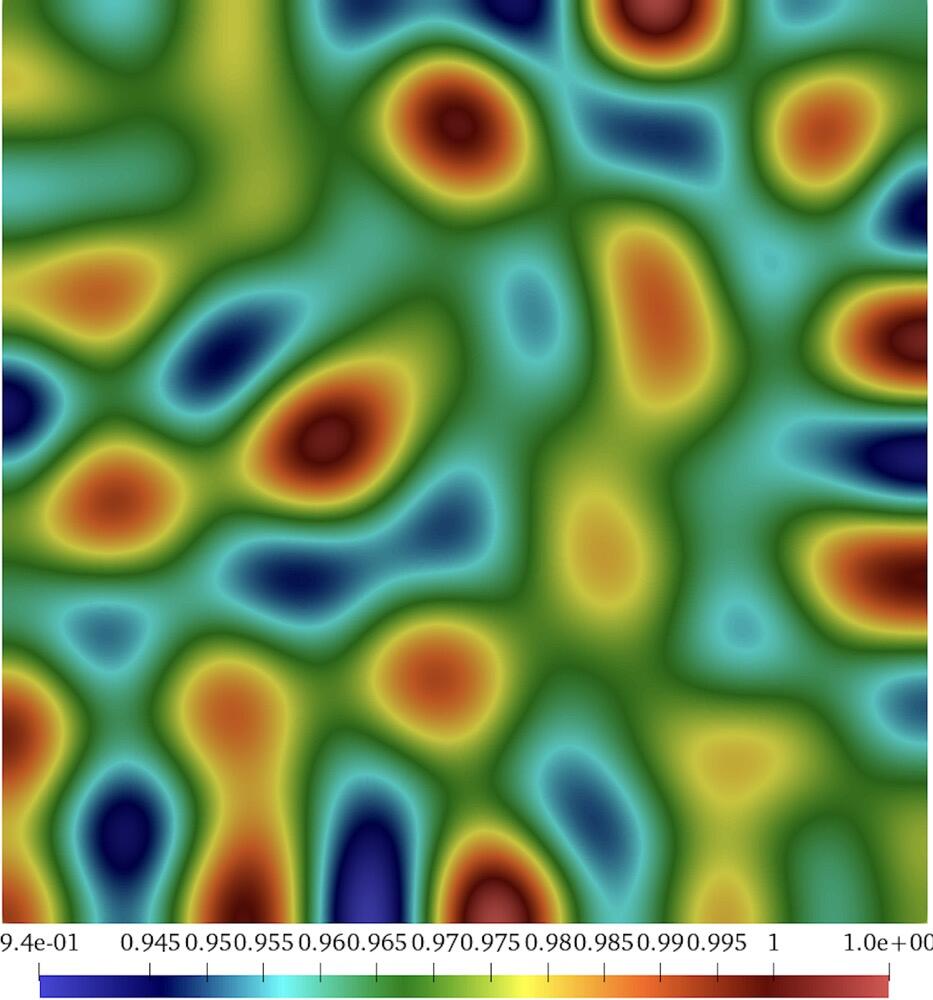}
    \end{overpic}
    \begin{overpic}[percent,width=0.19\textwidth]{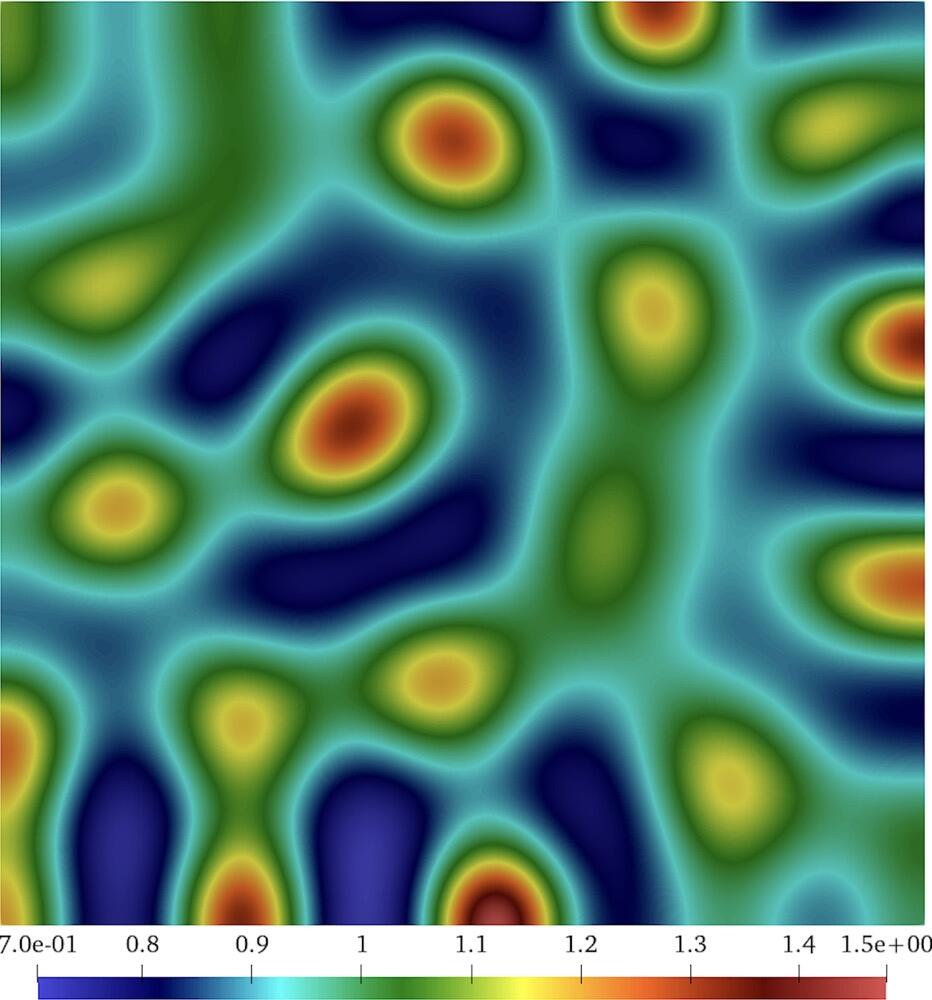}
    \end{overpic} 
        \begin{overpic}[percent,width=0.19\textwidth]{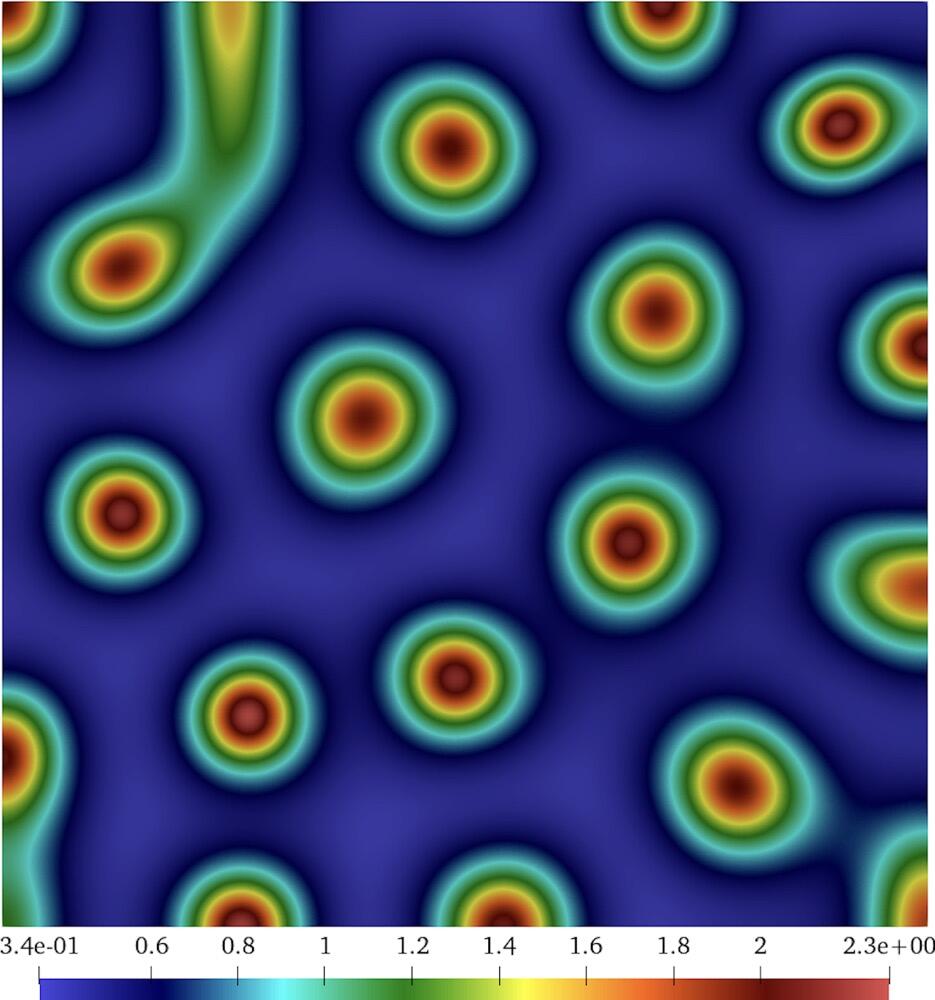}
    \end{overpic} 
    \begin{overpic}[percent,width=0.19\textwidth]{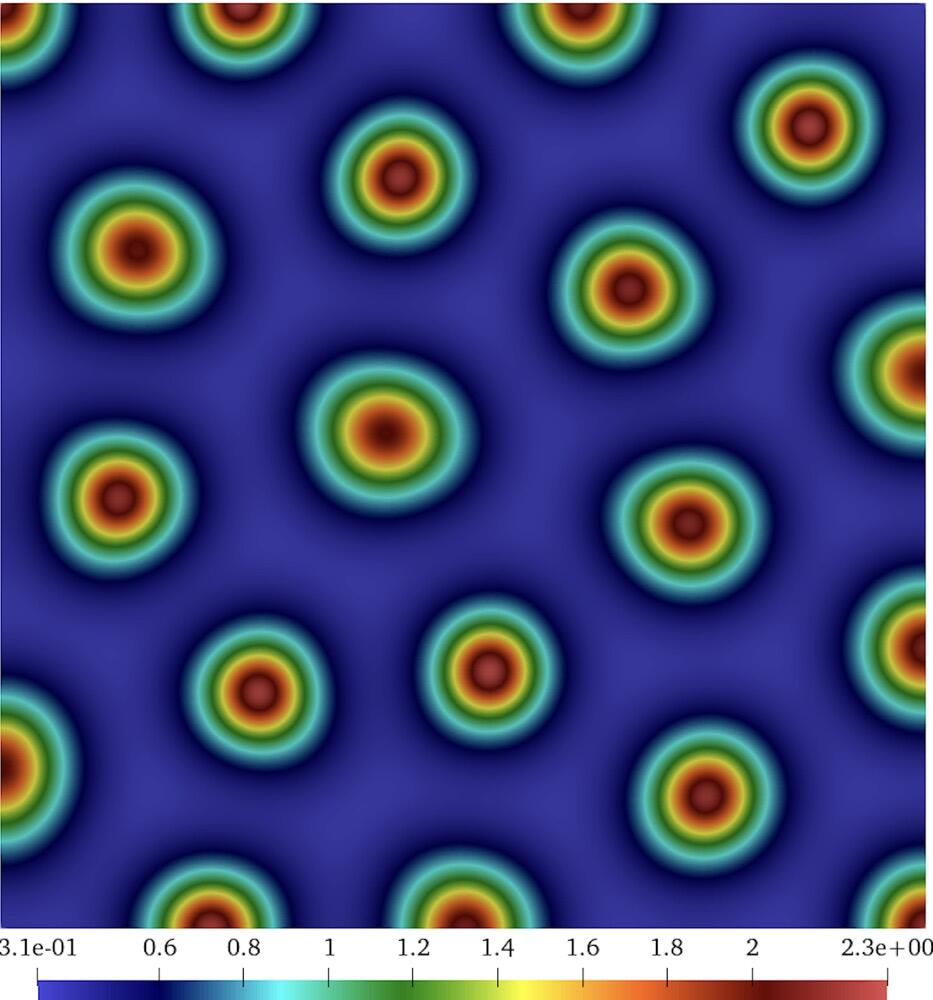}
    \end{overpic}
    \caption{Case 2 - PDE model: Evolution of the level of attractiveness (top) and density (bottom).}
    \label{fig:case2}
\end{figure}

Fig.~\ref{fig:case2} displays the evolution of the level of attractiveness and density for case 2 given by the PDE solver. Indeed, we observe that the system tends towards a steady state 
with large hotspots. In the top row of 
Fig.~\ref{fig:case2}, we see large spots of high attractiveness surrounded by areas
of low attractiveness. Variable $\rho$ mirrors the pattern in $A$: in the bottom row of 
Fig.~\ref{fig:case2} we see the emergence of large
spots with high concentrations of criminals
surrounded by areas with much lower concentrations.

\begin{figure}[htb!]
     \centering
         \begin{overpic}[percent,width=0.19\textwidth]{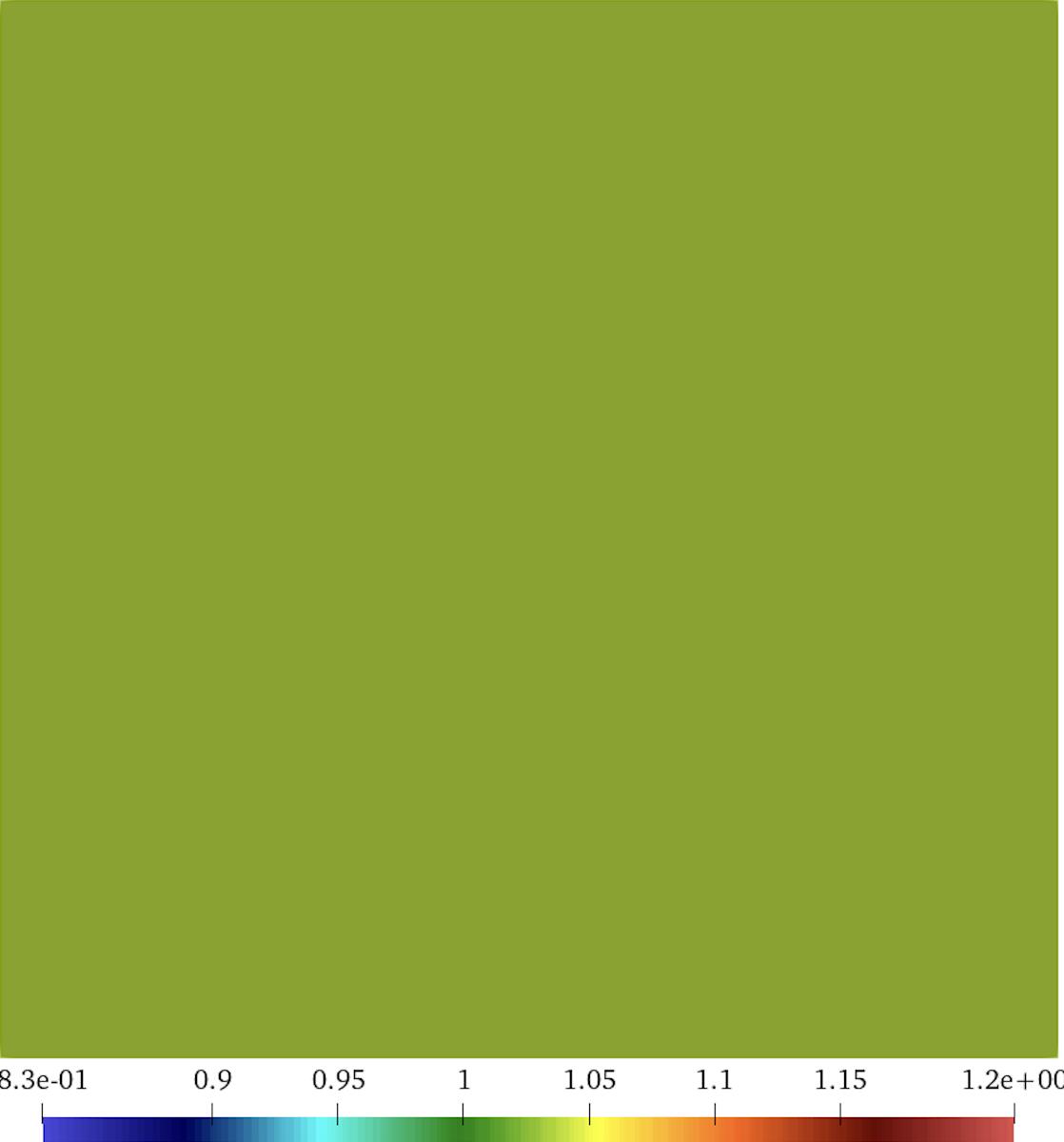}
    \put(30,103){\footnotesize{$t = 0$}}
    \put(-10,55){\footnotesize{$A$}}
    \end{overpic} 
     \begin{overpic}[percent,width=0.19\textwidth, grid=false]{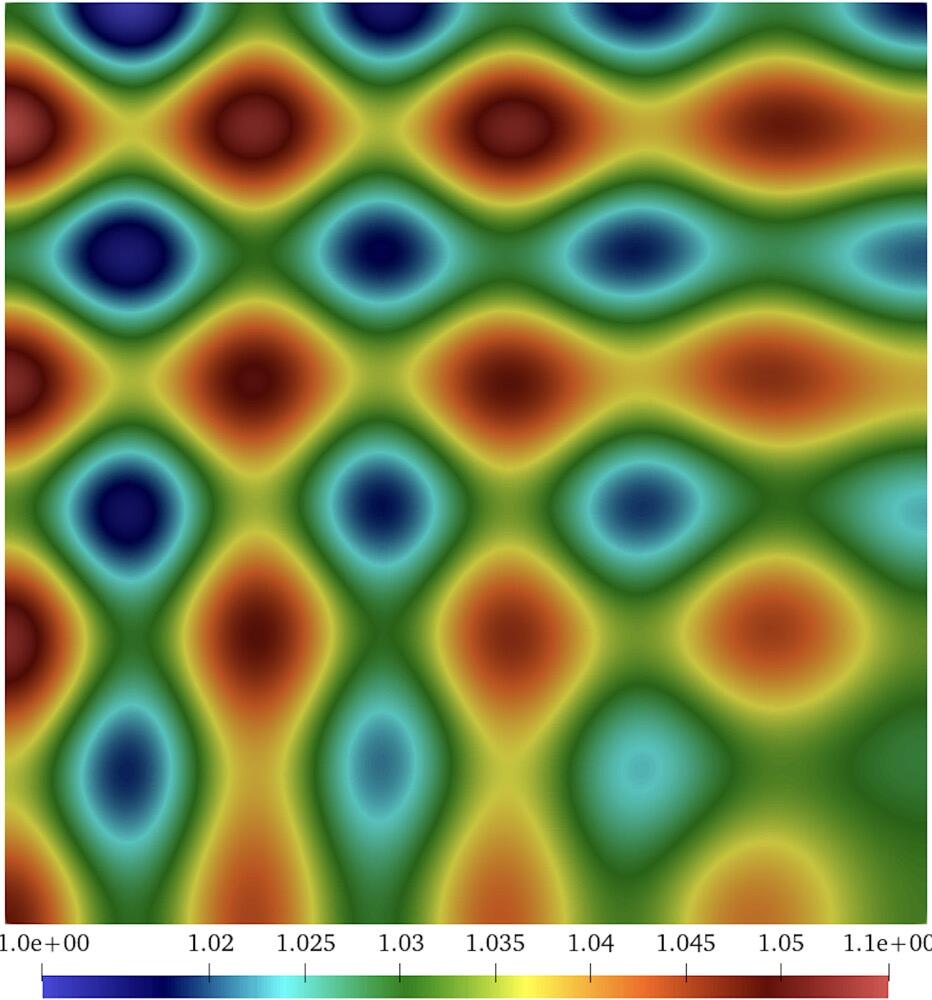}
    \put(35,103){\footnotesize{$t = 110$}}
    \end{overpic}
    \begin{overpic}[percent,width=0.19\textwidth]{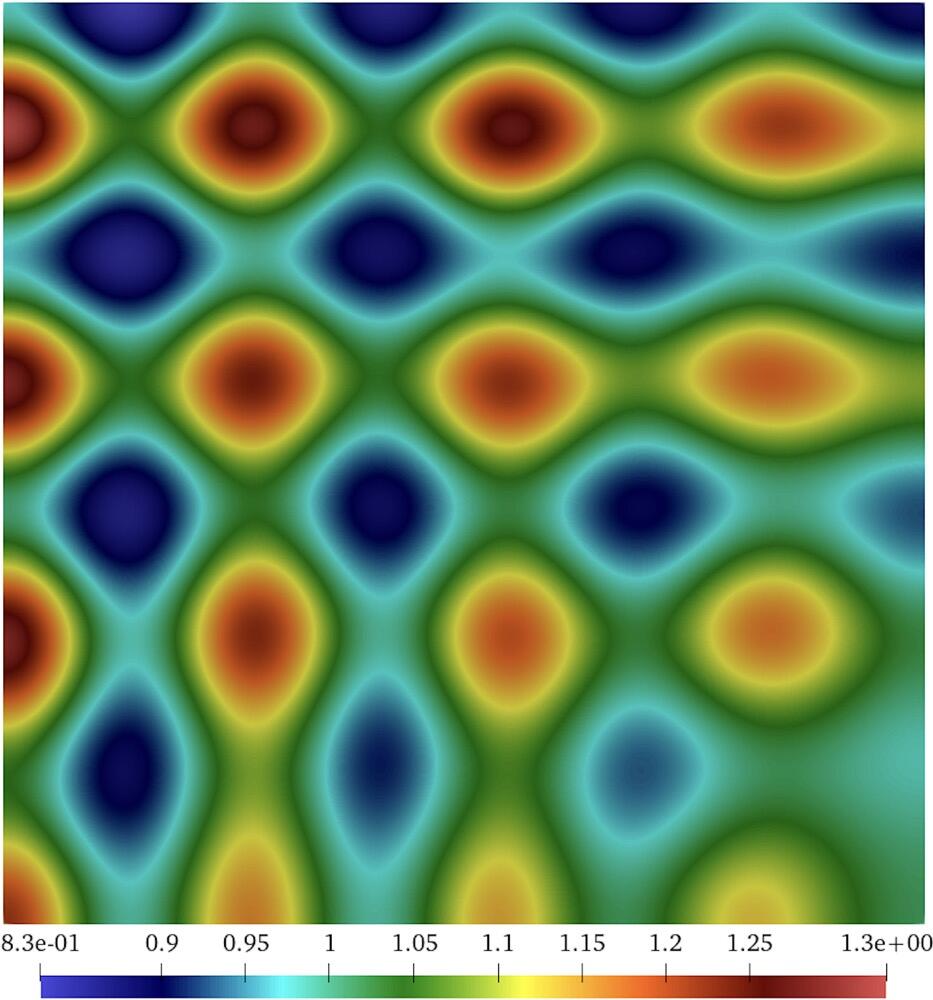}
    \put(35,103){\footnotesize{$t = 120$}}
    \end{overpic} 
        \begin{overpic}[percent,width=0.19\textwidth]{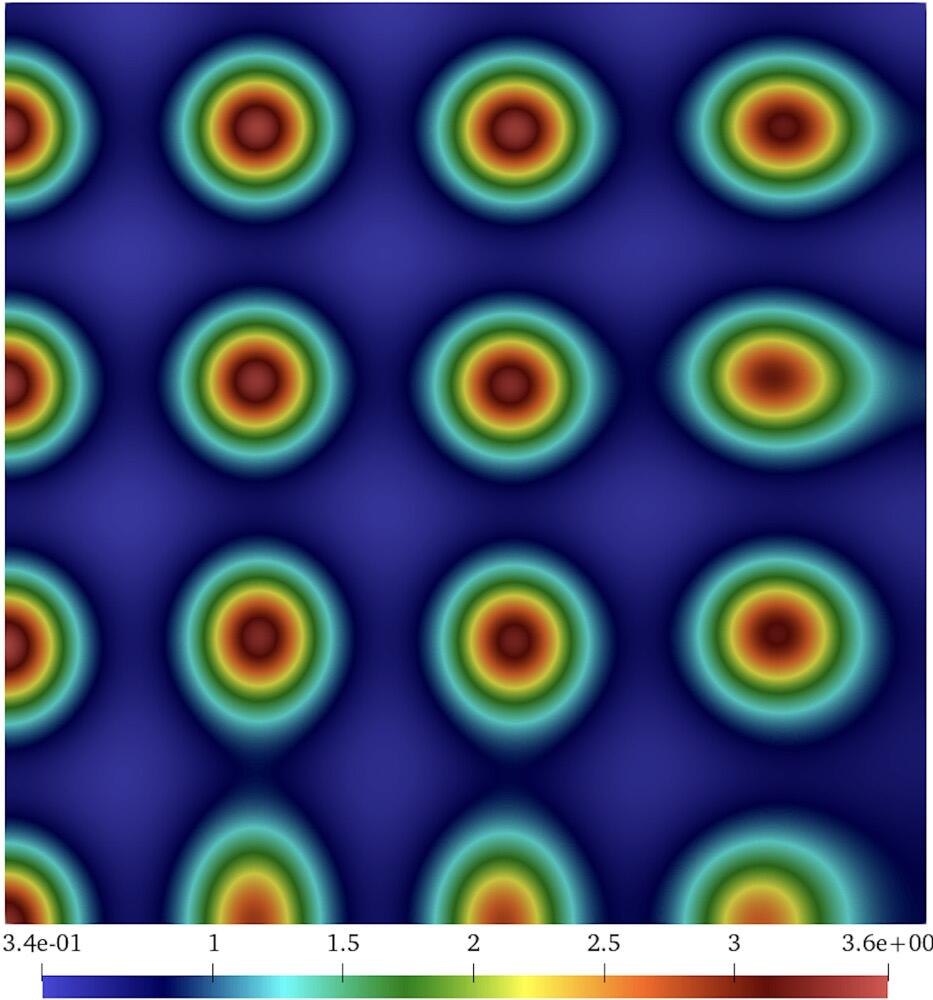}
    \put(32,103){\footnotesize{$t = 130$}}
    \end{overpic} 
    \begin{overpic}[percent,width=0.19\textwidth]{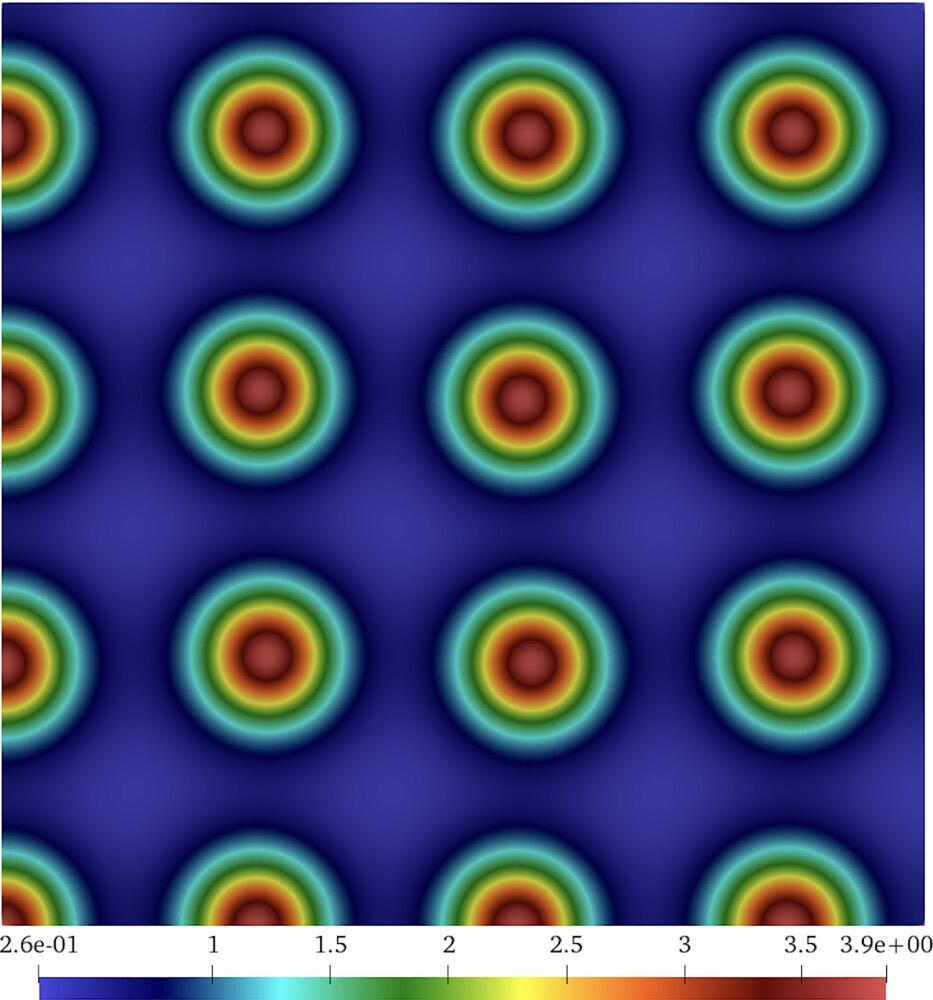}
    \put(30,103){\footnotesize{$t = 200$}}
    \end{overpic}
    \\
         \begin{overpic}[percent,width=0.19\textwidth]{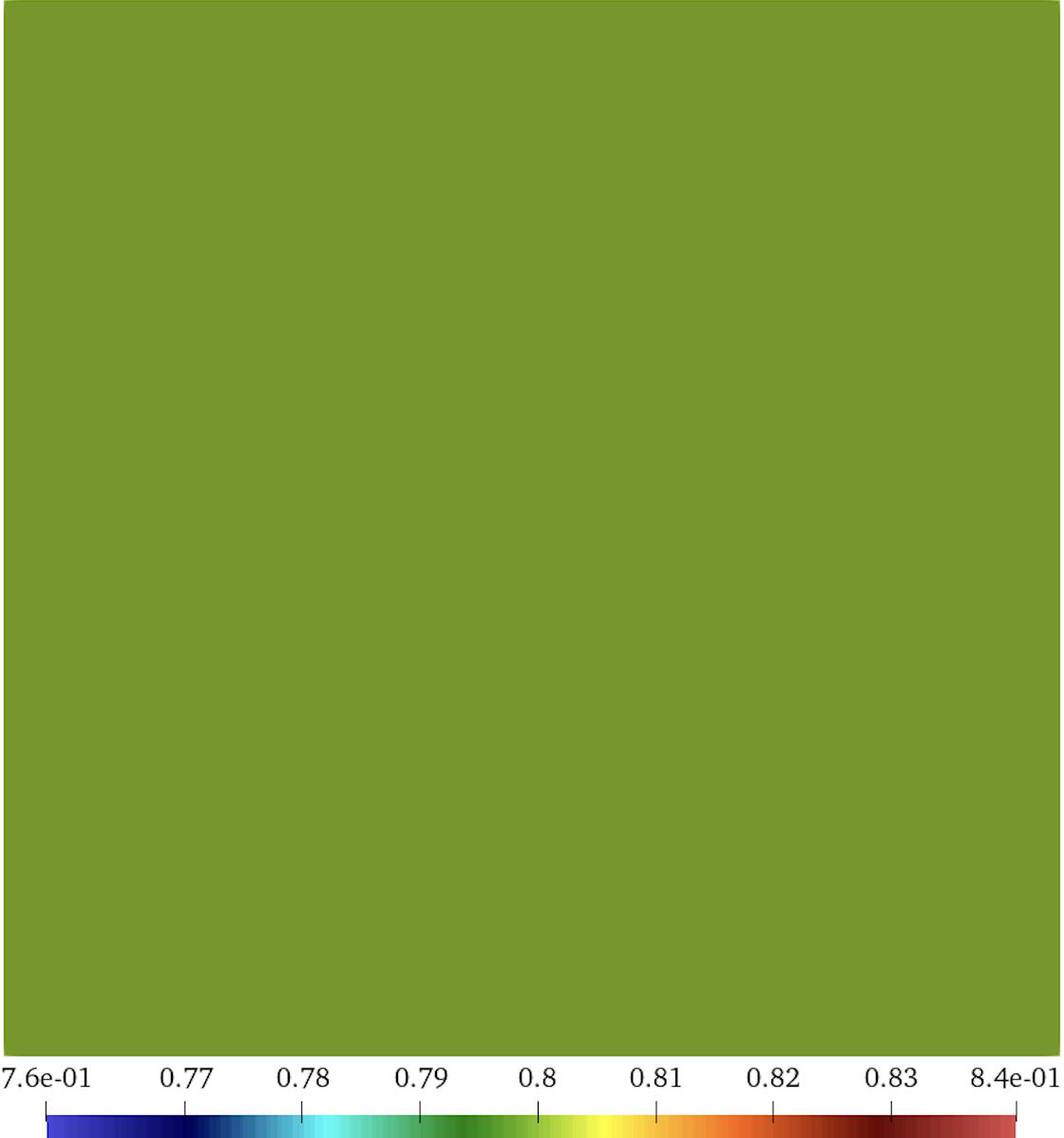}
    \put(-10,55){\footnotesize{$\rho$}}
    \end{overpic} 
     \begin{overpic}[percent,width=0.19\textwidth, grid=false]{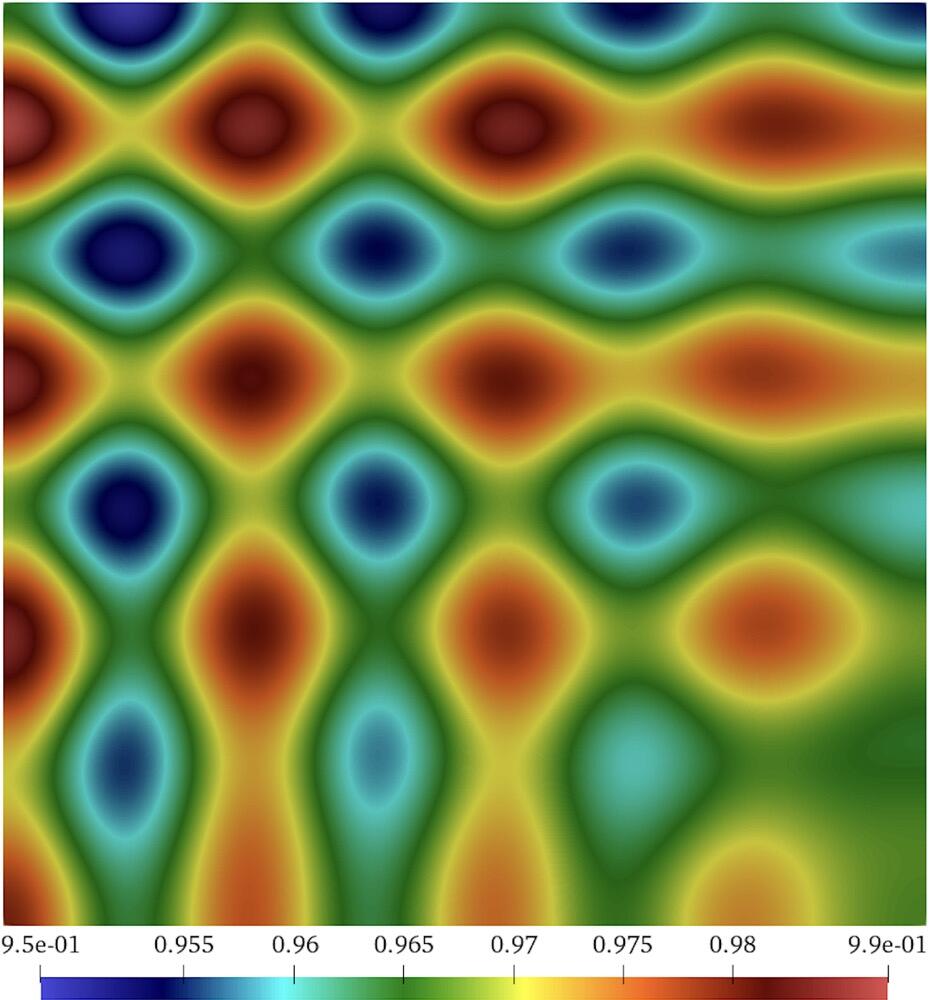}
    \end{overpic}
    \begin{overpic}[percent,width=0.19\textwidth]{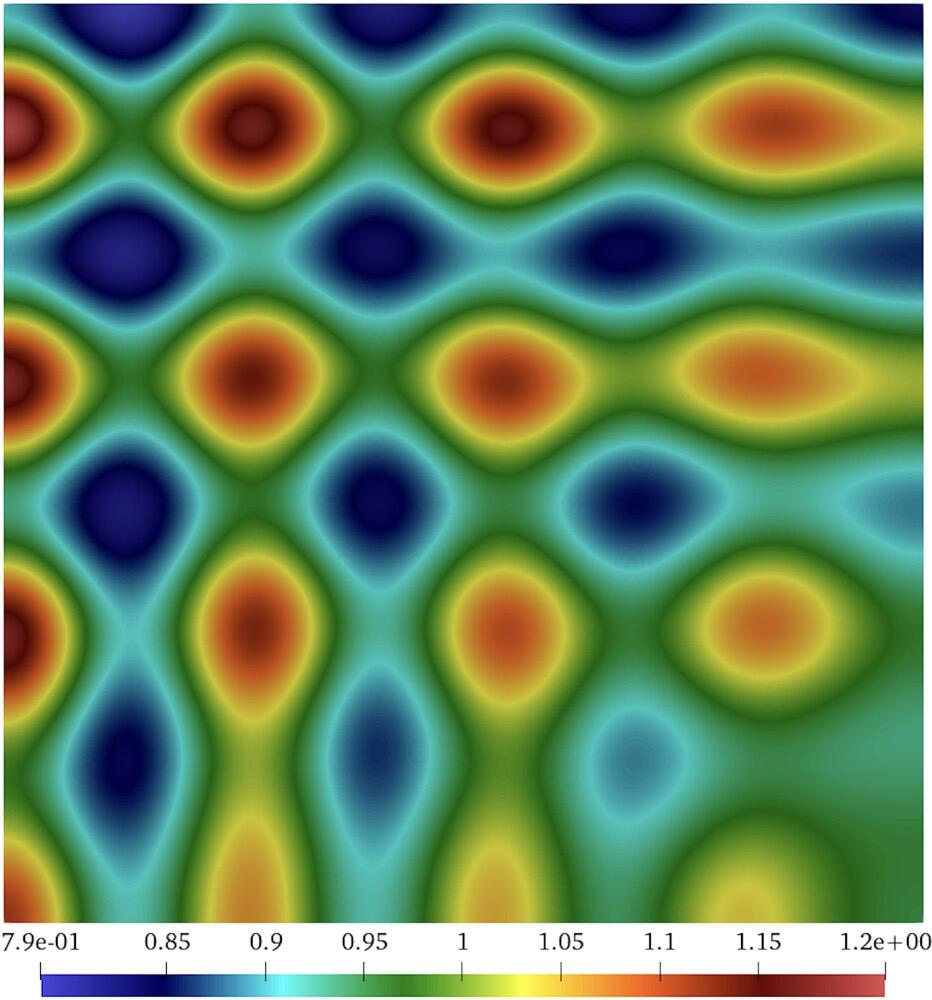}
    \end{overpic} 
        \begin{overpic}[percent,width=0.19\textwidth]{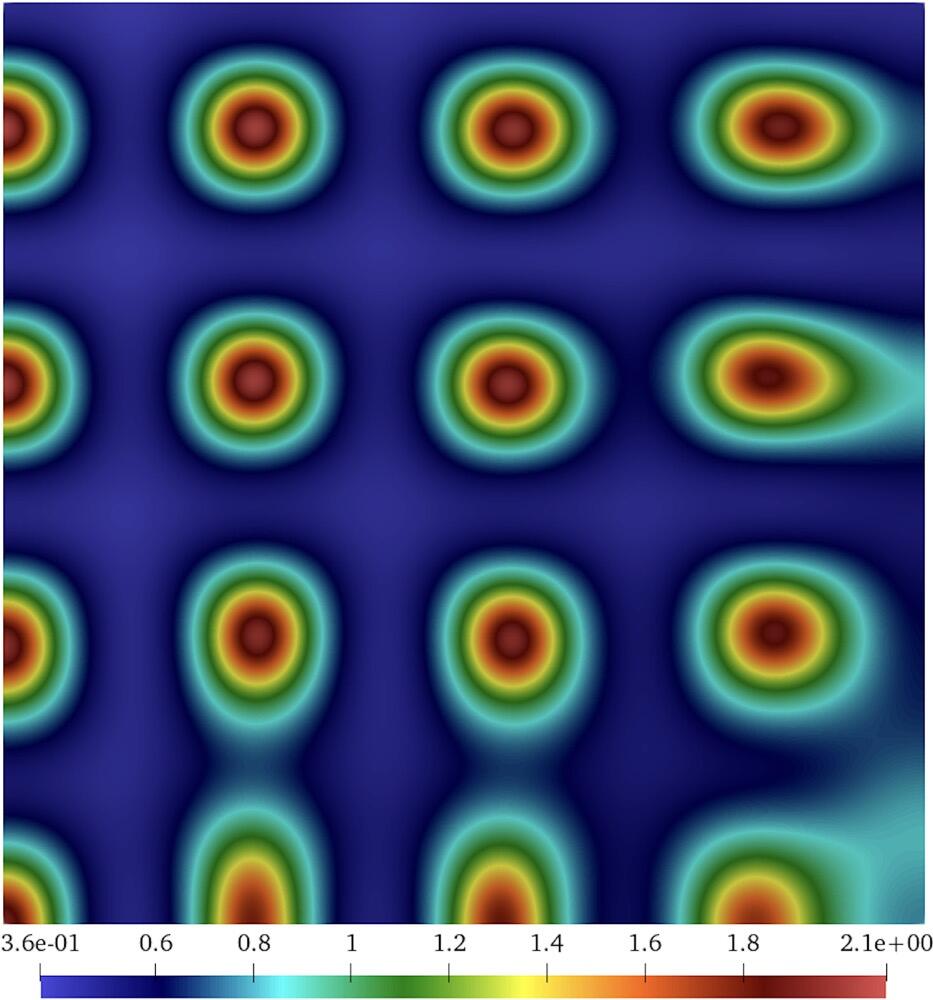}
    \end{overpic} 
    \begin{overpic}[percent,width=0.19\textwidth]{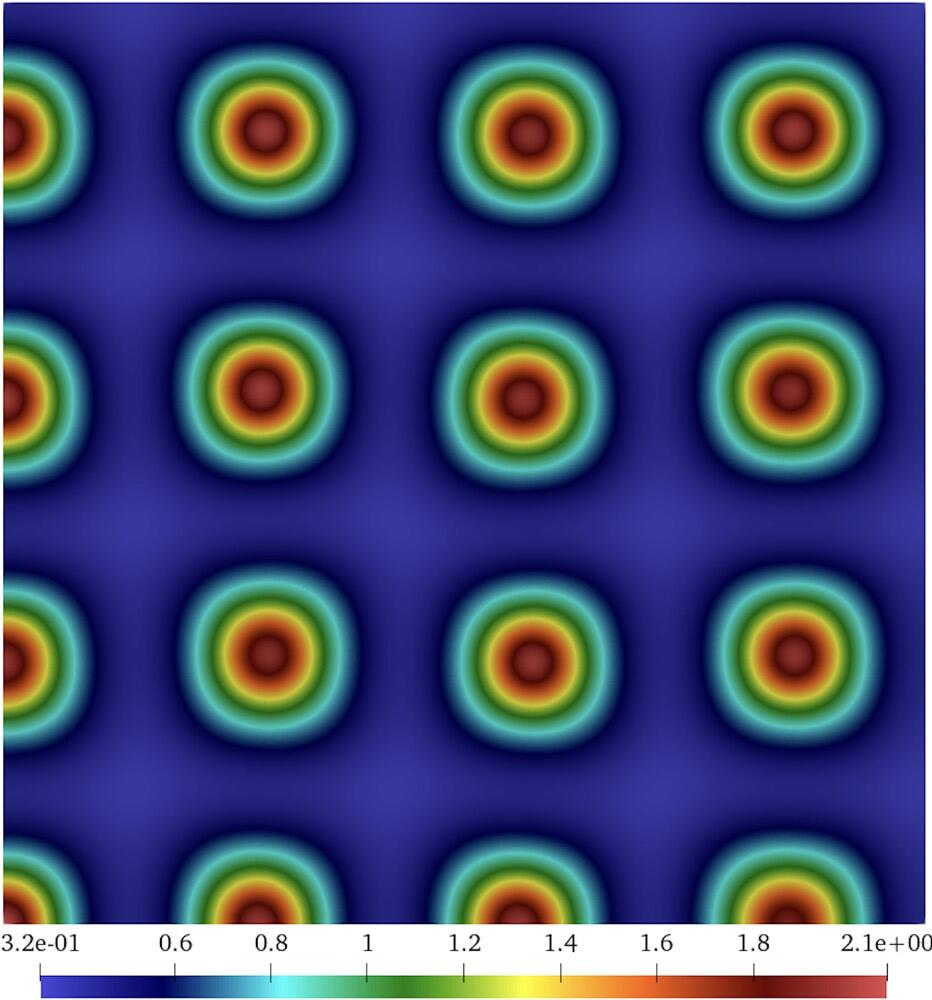}
    \end{overpic}
    \caption{Case 2 without noise - PDE model: Evolution of the level of attractiveness (top) and density (bottom).}
    \label{fig:case2-nonoise}
\end{figure}

The noise in the initial data of case 2 affects the position of the hotspots, not on their existence or rough size, and the time needed for the emergence of the pattern. See Fig.~\ref{fig:case2-nonoise}, which shows $A$ and $\rho$ for case 2 with no noise computed by the PDE solver. 
The hotspots have roughly the same size as those in Fig.~\ref{fig:case2}, but 
their shape is more regular and the pattern is more structured. Additionally, we see that while at
$t = 30$ a hotspot pattern has almost emerged in case 2 with noise, it takes
around $130$ time units for that to be the case when the noise is absent.

\begin{figure}[htb!]
     \centering
         \begin{overpic}[percent,width=0.19\textwidth]{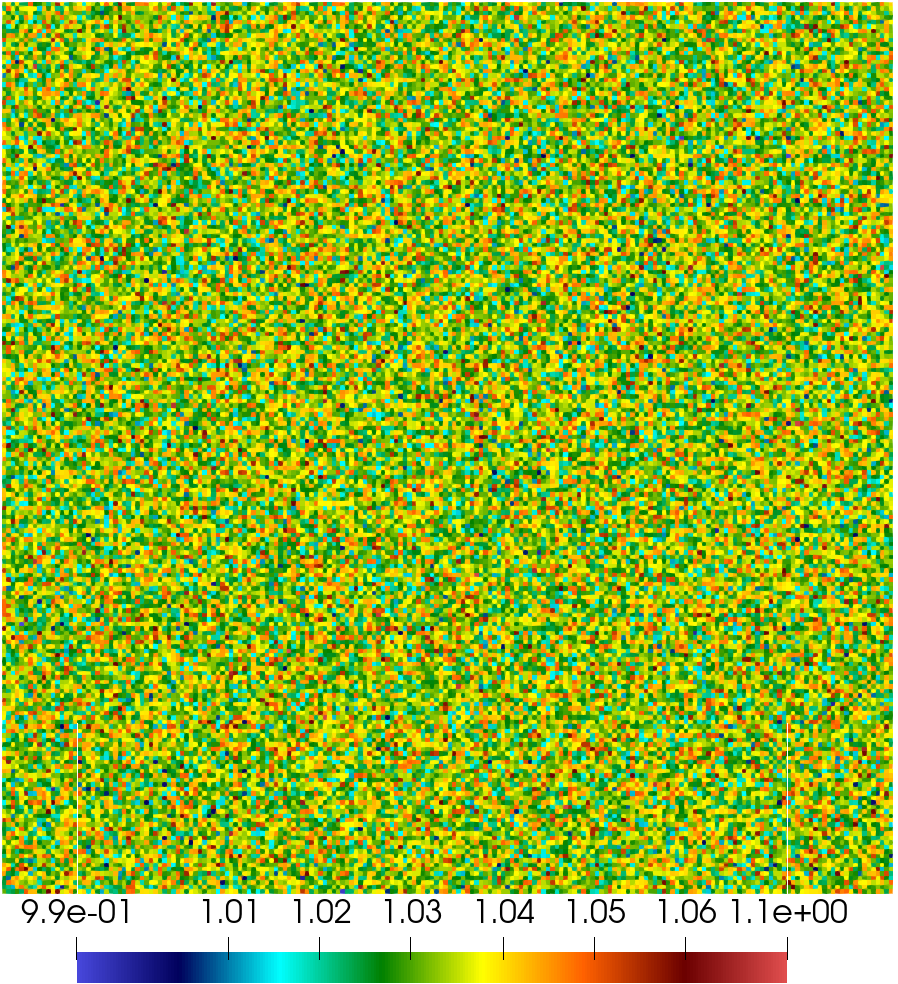}
    \put(30,103){\footnotesize{$t = 0$}}
    \put(-10,55){\footnotesize{$A$}}
    \end{overpic} 
     \begin{overpic}[percent,width=0.19\textwidth, grid=false]{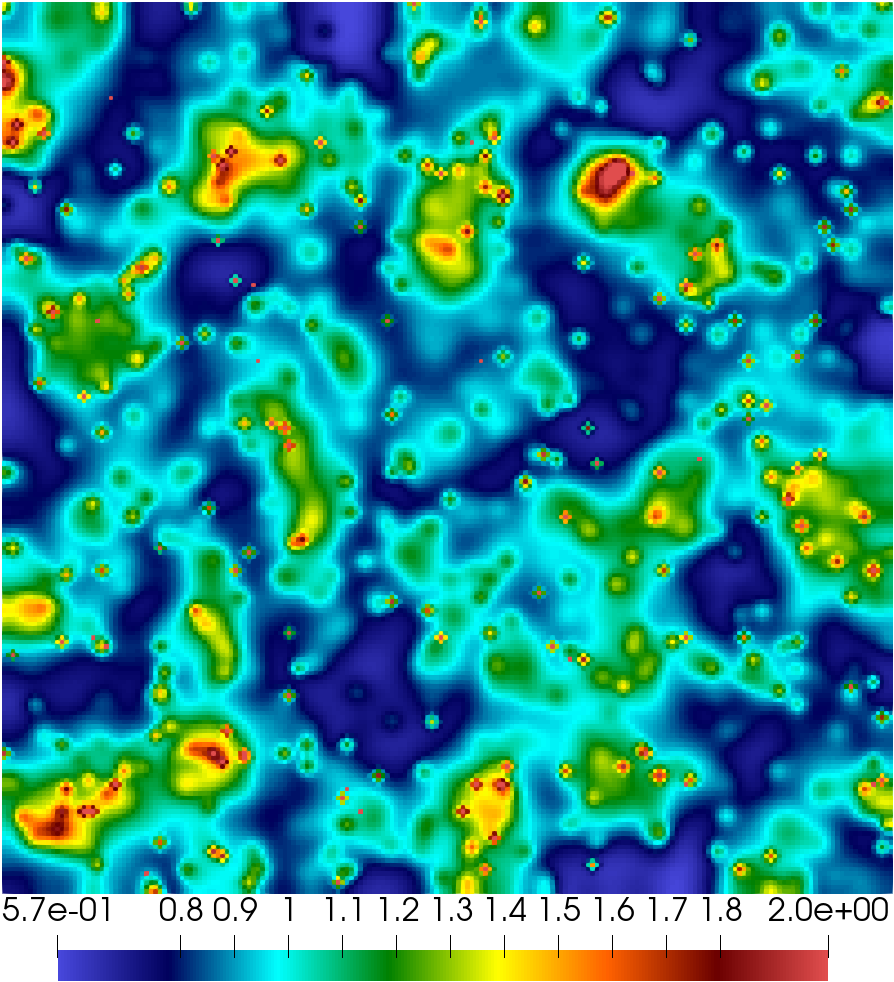}
    \put(35,103){\footnotesize{$t = 5$}}
    \end{overpic}
    \begin{overpic}[percent,width=0.19\textwidth]{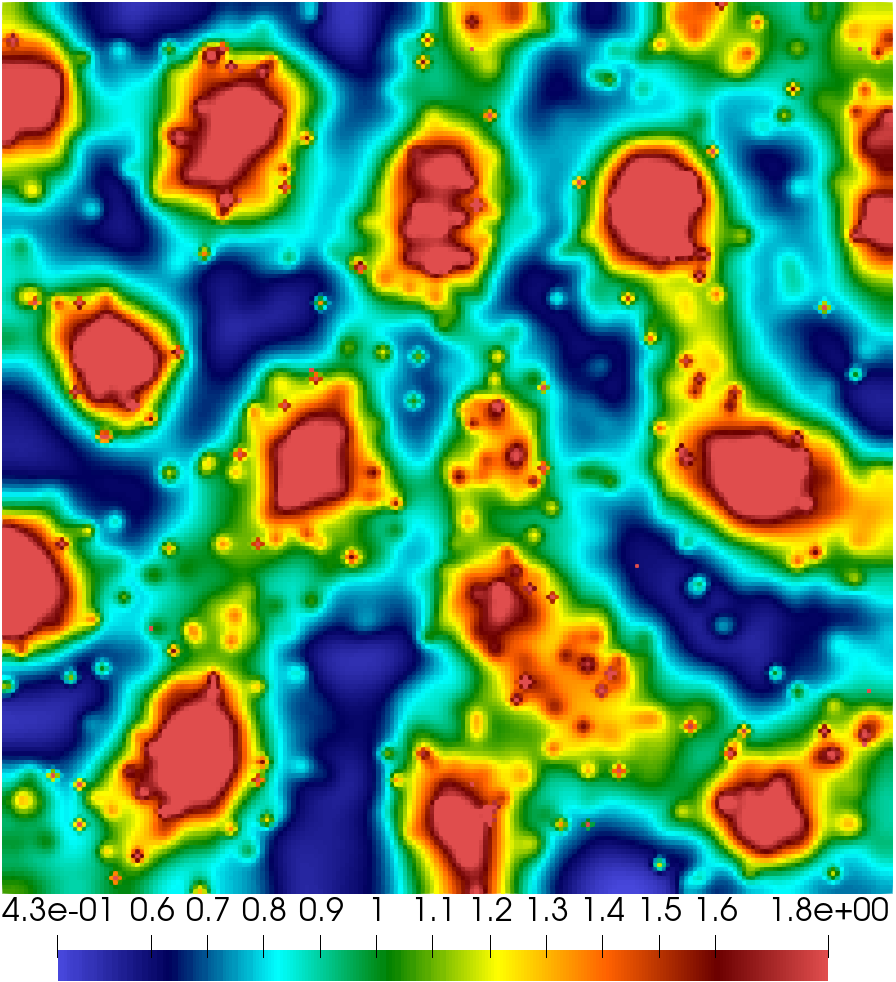}
    \put(35,103){\footnotesize{$t = 10$}}
    \end{overpic} 
        \begin{overpic}[percent,width=0.19\textwidth]{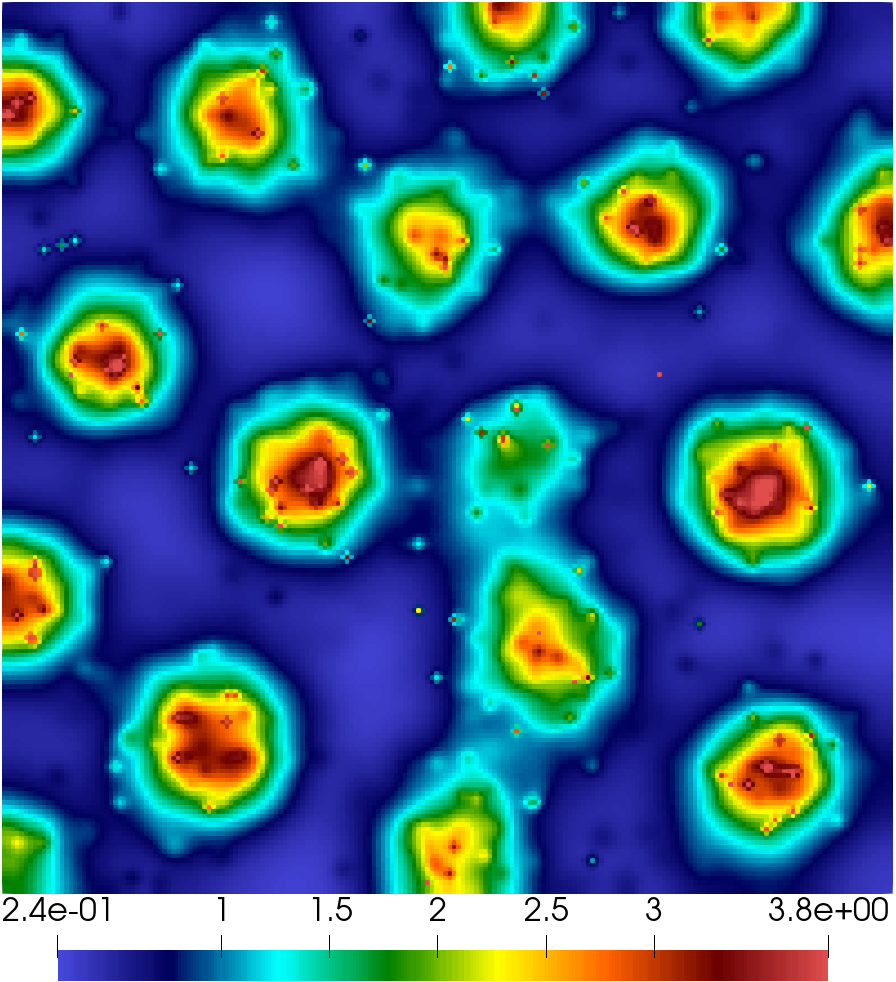}
    \put(32,103){\footnotesize{$t = 20$}}
    \end{overpic} 
    \begin{overpic}[percent,width=0.19\textwidth]{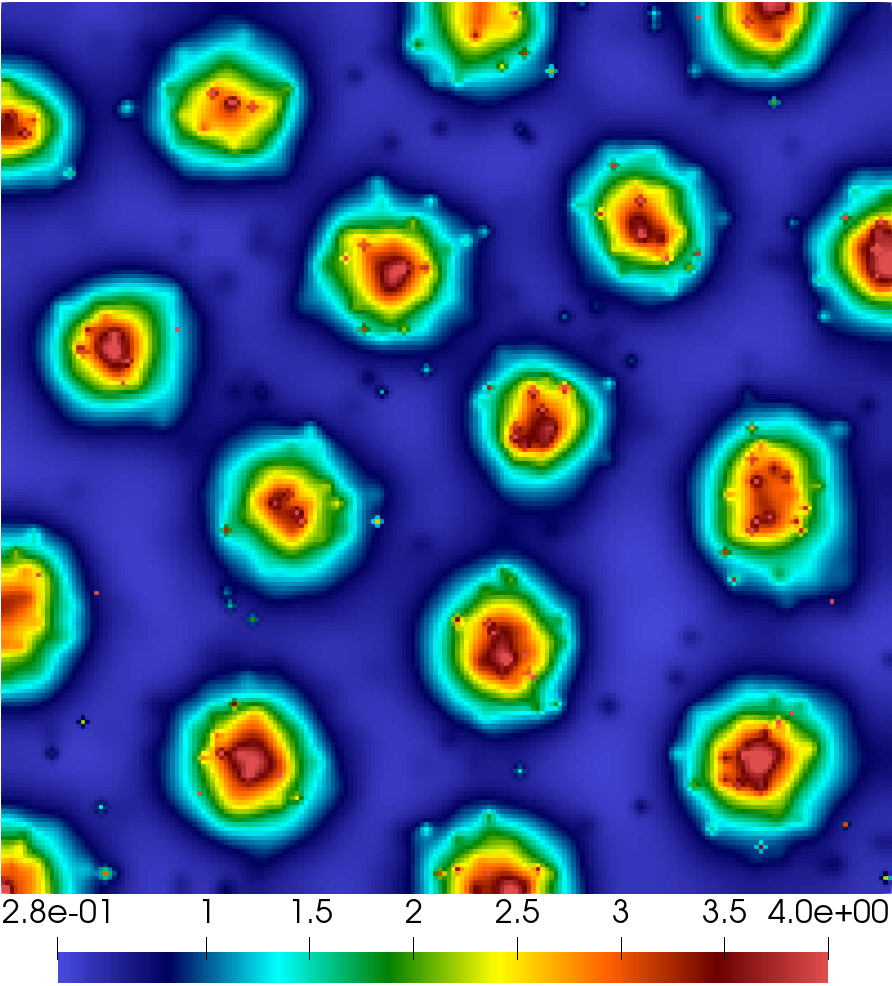}
    \put(30,103){\footnotesize{$t = 200$}}
    \end{overpic}
    \\
         \begin{overpic}[percent,width=0.19\textwidth]{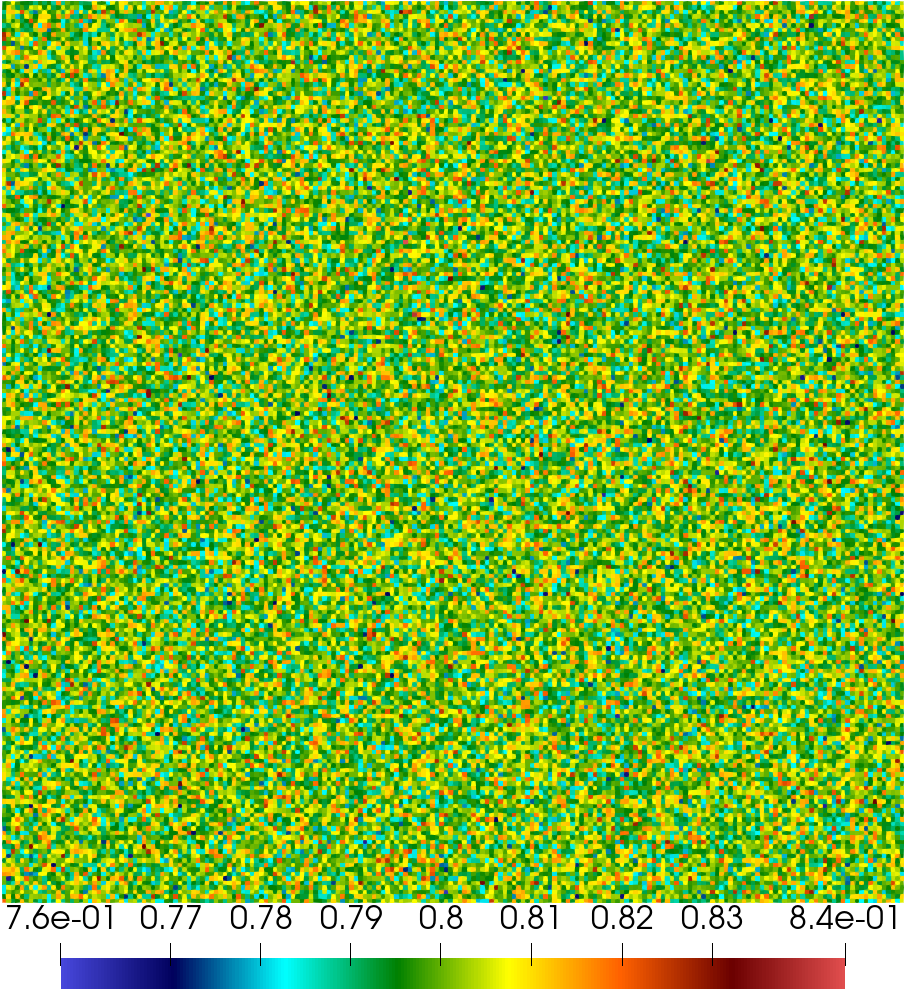}
    \put(-10,55){\footnotesize{$\rho$}}
    \end{overpic} 
     \begin{overpic}[percent,width=0.19\textwidth, grid=false]{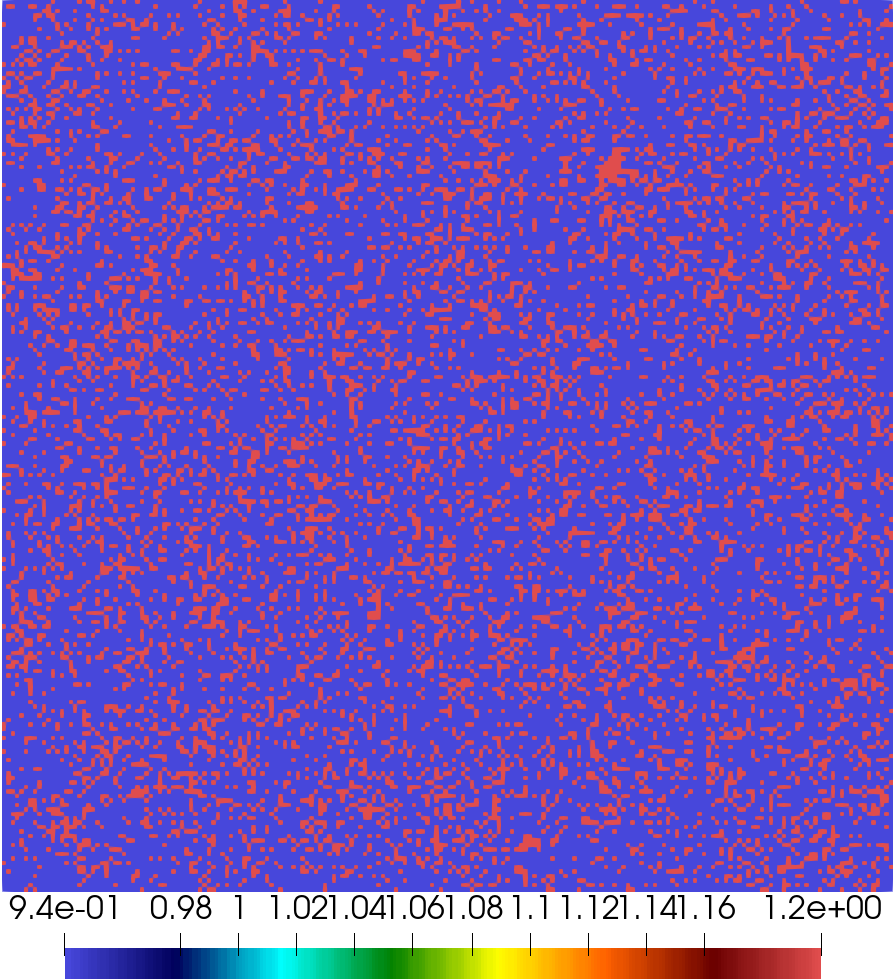}
    \end{overpic}
    \begin{overpic}[percent,width=0.19\textwidth]{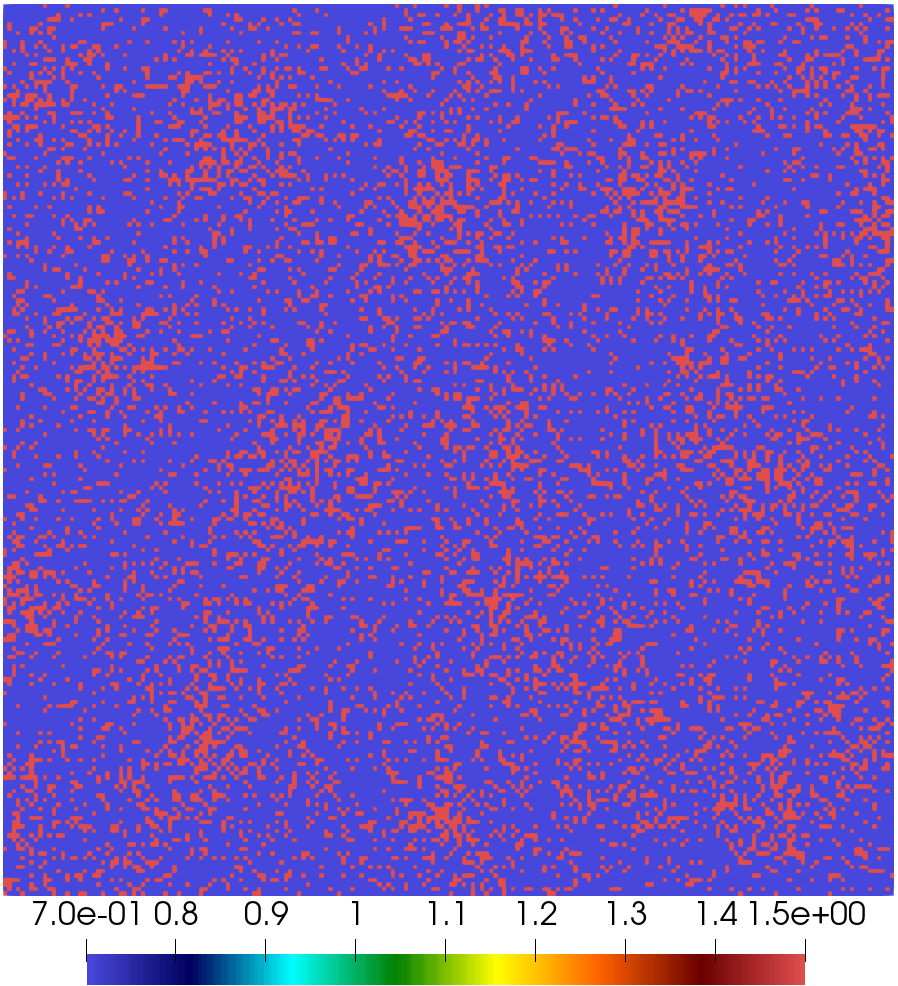}
    \end{overpic} 
        \begin{overpic}[percent,width=0.19\textwidth]{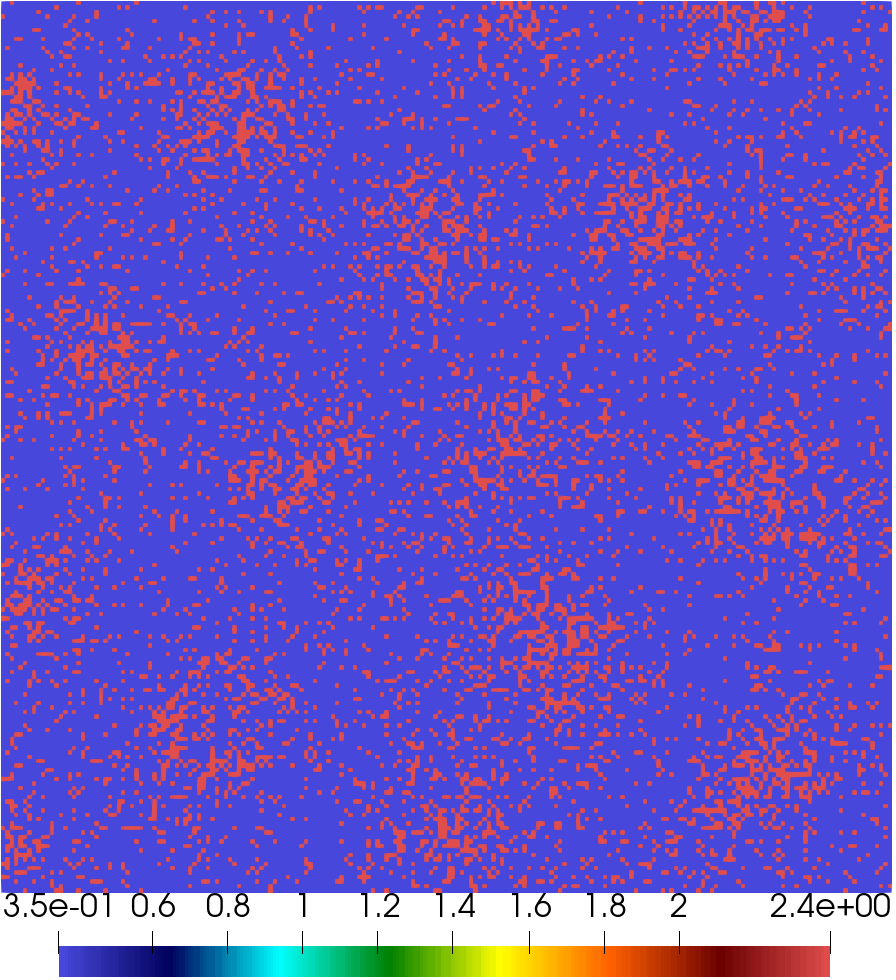}
    \end{overpic} 
    \begin{overpic}[percent,width=0.19\textwidth]{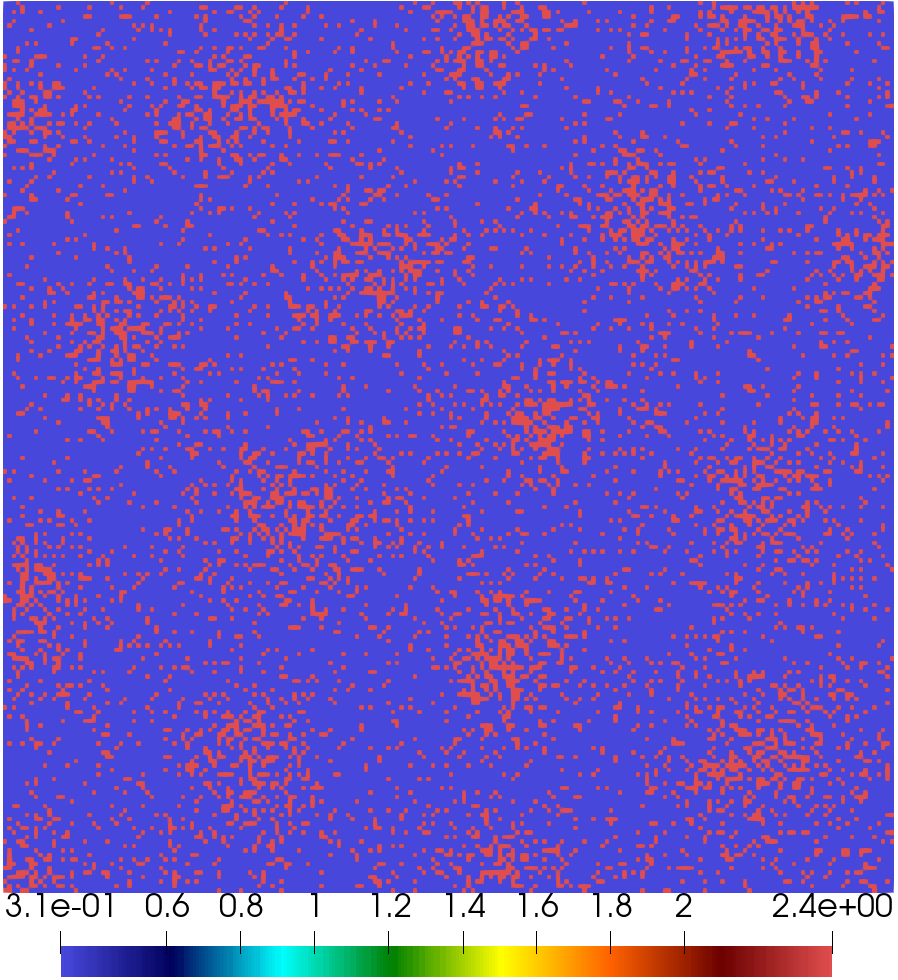}
    \end{overpic}
    \caption{Case 2 - agent-based model: Evolution of the level of attractiveness (top) and criminals (bottom).}
    \label{fig:case2DM}
\end{figure}

For comparison with Fig.~\ref{fig:case2}, 
Fig.~\ref{fig:case2DM} shows
the evolution of $A$ and $\rho$ for case 2 given by the agent-based solver after non-dimensionalization. For theis case, we set 
$\Gamma = 0.019$, $\theta = 0.2339$, and $\omega = 1/15$ to obtain $\Gamma \theta /\omega^2 \approx 1$.
We see that the magnitudes of 
$A$ and $\rho$ in Fig.~\ref{fig:case2}
and \ref{fig:case2DM} match
and the sizes for the hotspots
are comparable. The reason why we changed the values of $\Gamma$ and $\theta$ with respect to case 1 is to make the hotspots in $\rho$ more visibile. For $\Gamma = 0.0077$ and $\theta = 0.58$, the hotspots in 
$\rho$ are still there but look fainter in the plots. 

Further analysis from 
\cite{short2008statistical} indicates that
hotspot separation increases with $\eta$. The hotspot separation is affected by their size, which means that 
the hotspot size decreases with $\eta$.
For this reason, $\eta$ is further decreased in case 3, while all the other parameters remain the same as in case 1 and 2.

\begin{figure}[htb!]
     \centering
         \begin{overpic}[percent,width=0.19\textwidth]{Figs_ContinummModel/S2_A_T0.jpg}
    \put(30,103){\footnotesize{$t = 0$}}
    \put(-10,55){\footnotesize{$A$}}
    \end{overpic} 
     \begin{overpic}[percent,width=0.19\textwidth, grid=false]{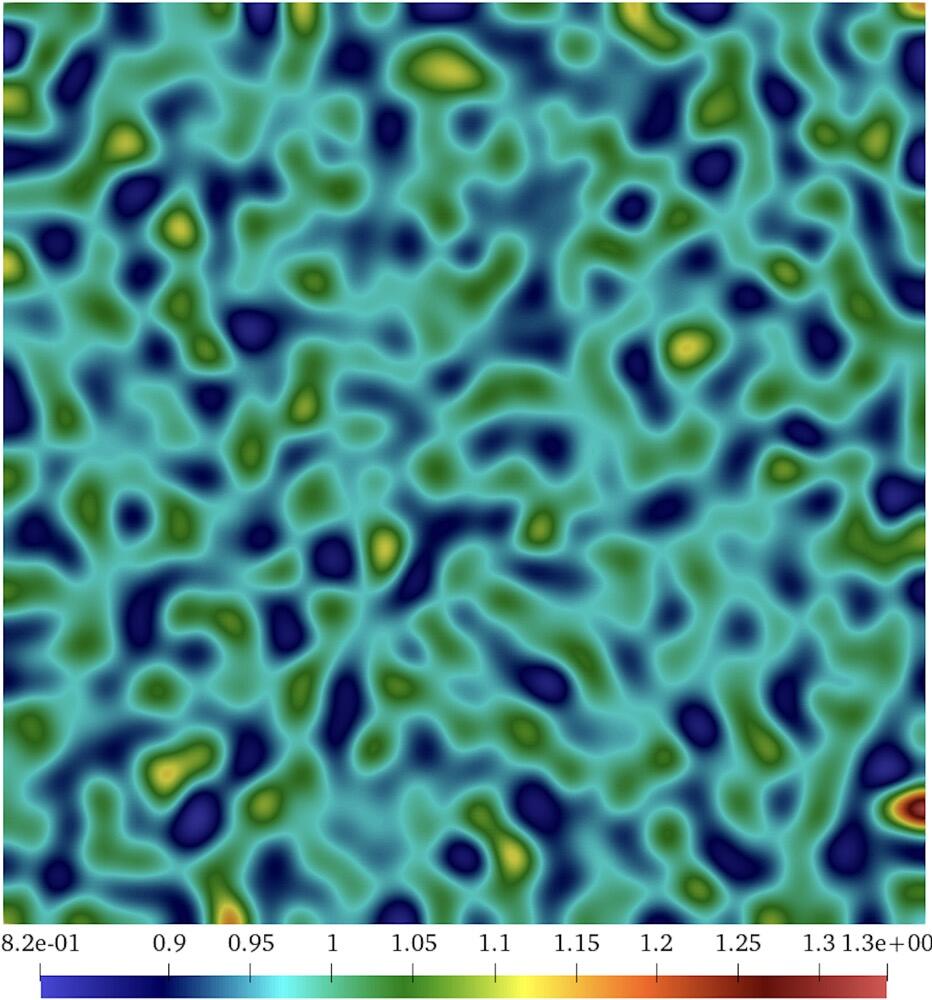}
    \put(35,103){\footnotesize{$t = 2$}}
    \end{overpic}
    \begin{overpic}[percent,width=0.19\textwidth]{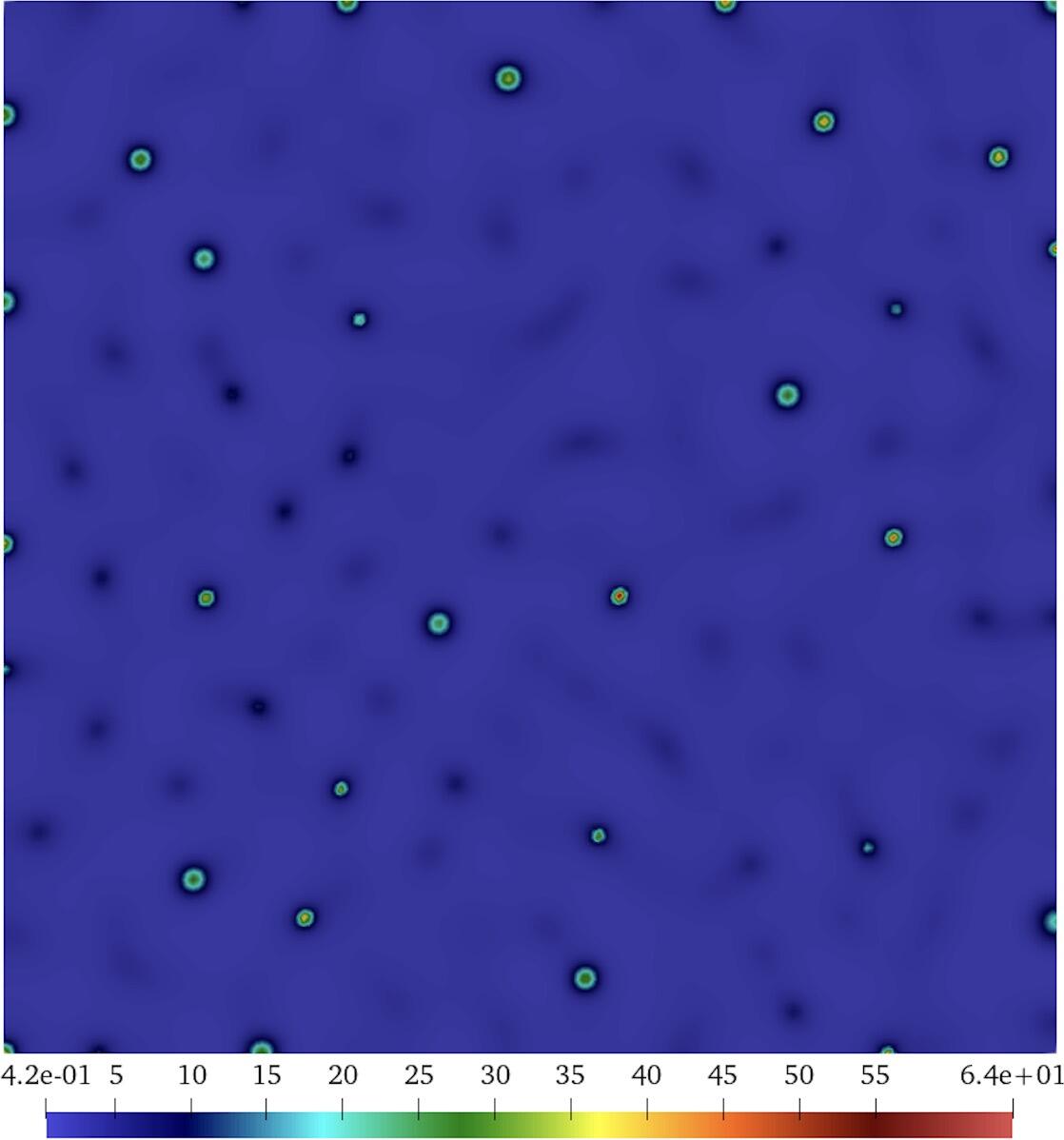}
    \put(35,103){\footnotesize{$t = 4$}}
    \end{overpic} 
        \begin{overpic}[percent,width=0.19\textwidth]{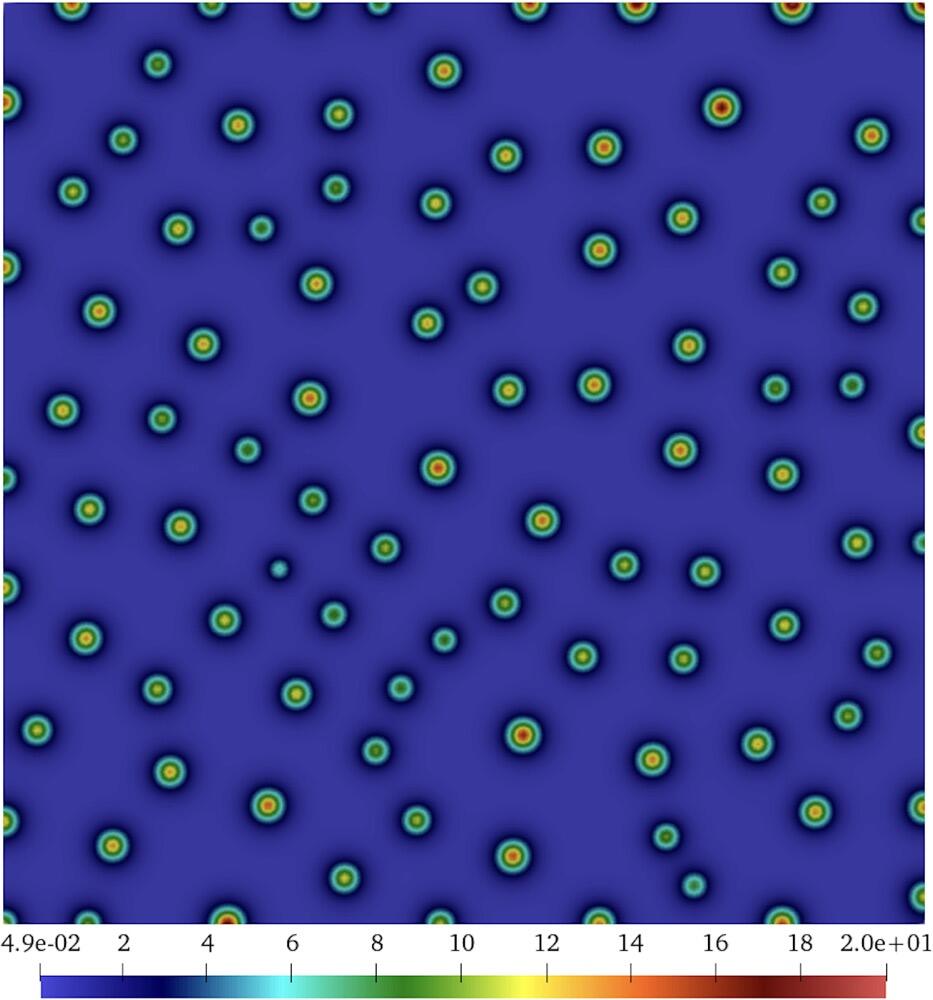}
    \put(32,103){\footnotesize{$t = 10$}}
    \end{overpic} 
    \begin{overpic}[percent,width=0.19\textwidth]{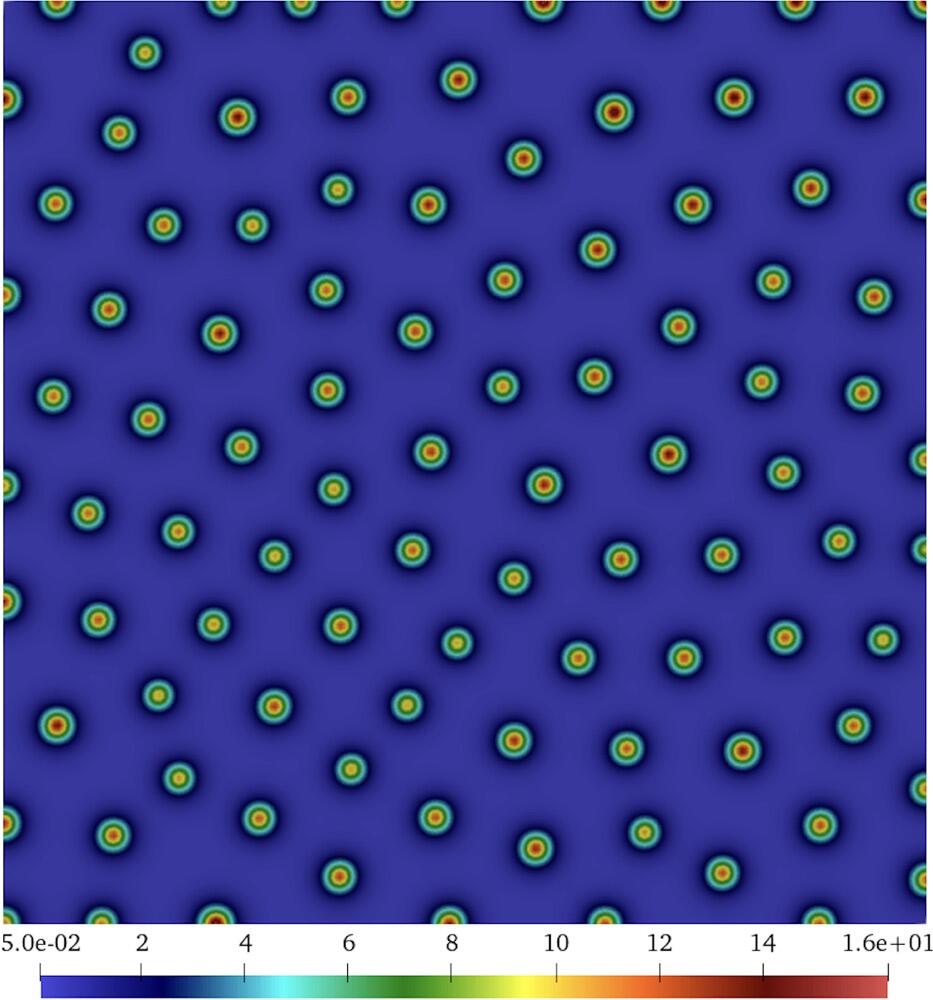}
    \put(30,103){\footnotesize{$t = 200$}}
    \end{overpic}
    \\
         \begin{overpic}[percent,width=0.19\textwidth]{Figs_ContinummModel/S2_rho_T0.jpg}
    \put(-10,55){\footnotesize{$\rho$}}
    \end{overpic} 
     \begin{overpic}[percent,width=0.19\textwidth, grid=false]{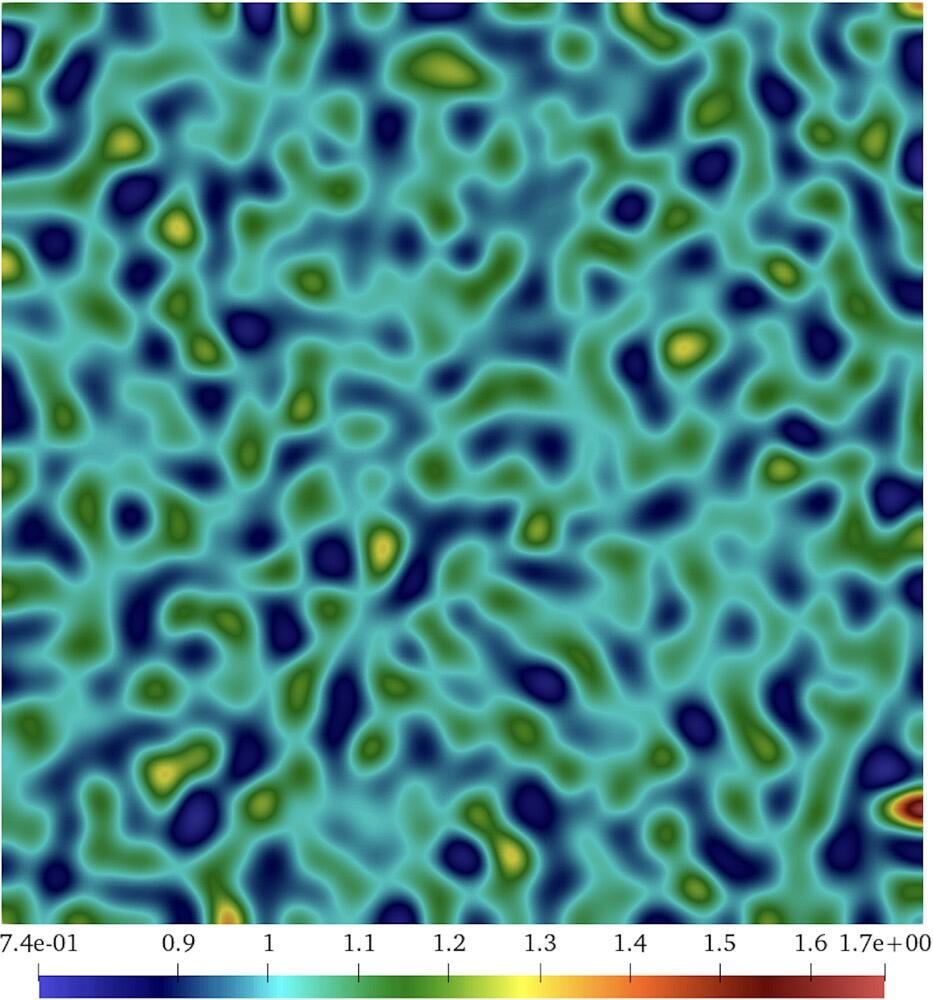}
    \end{overpic}
    \begin{overpic}[percent,width=0.19\textwidth]{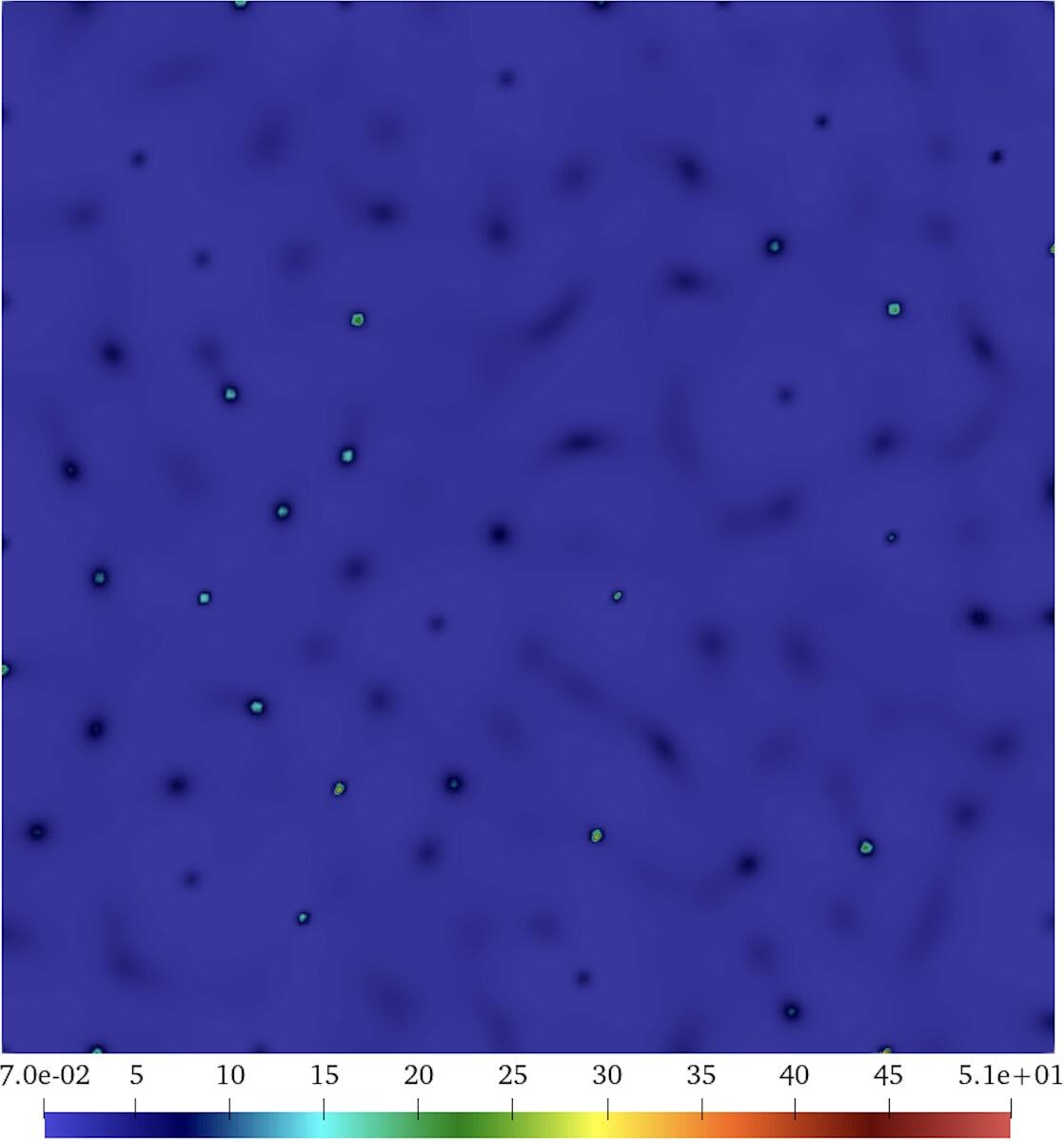}
    \end{overpic} 
        \begin{overpic}[percent,width=0.19\textwidth]{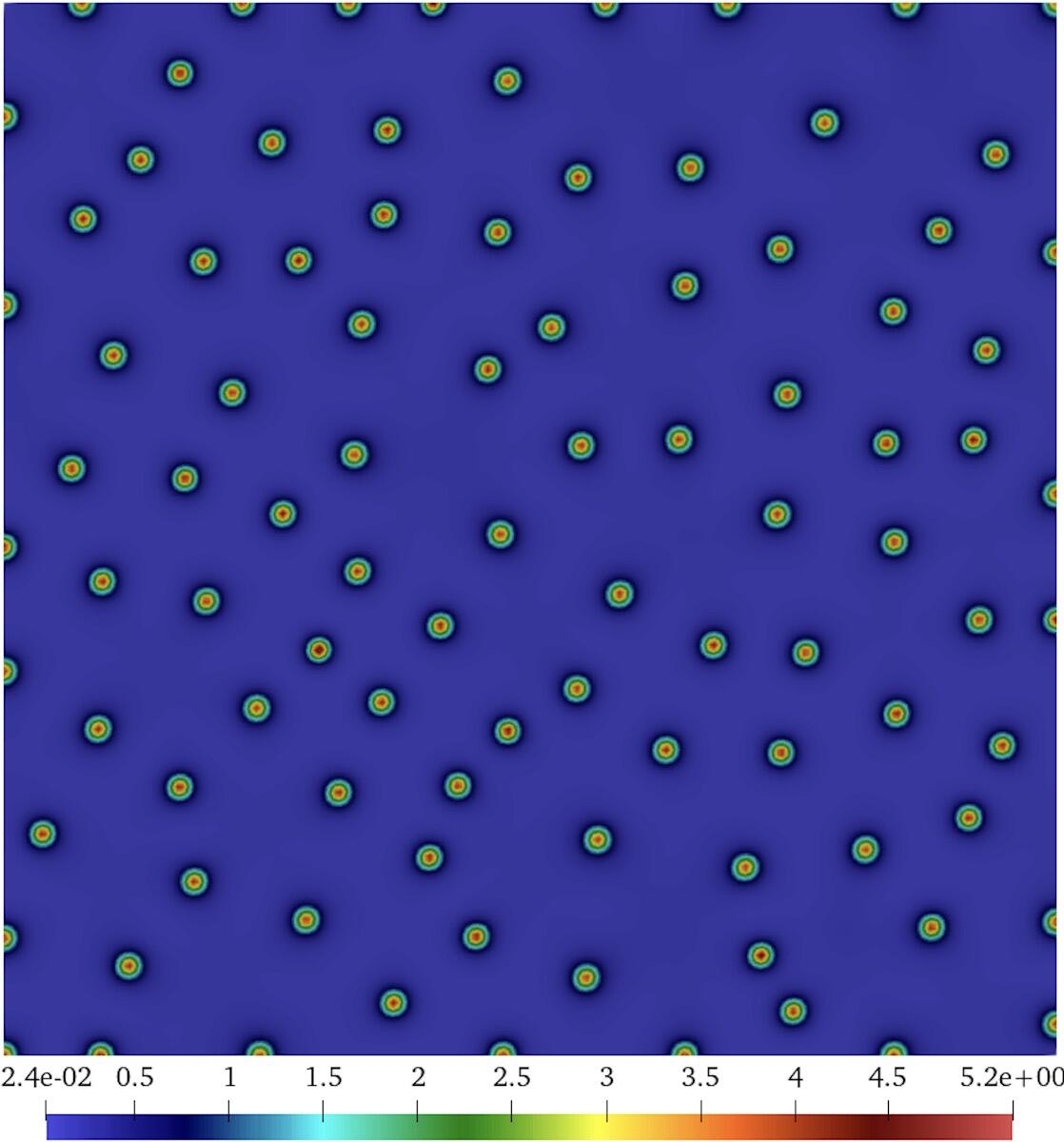}
    \end{overpic} 
    \begin{overpic}[percent,width=0.19\textwidth]{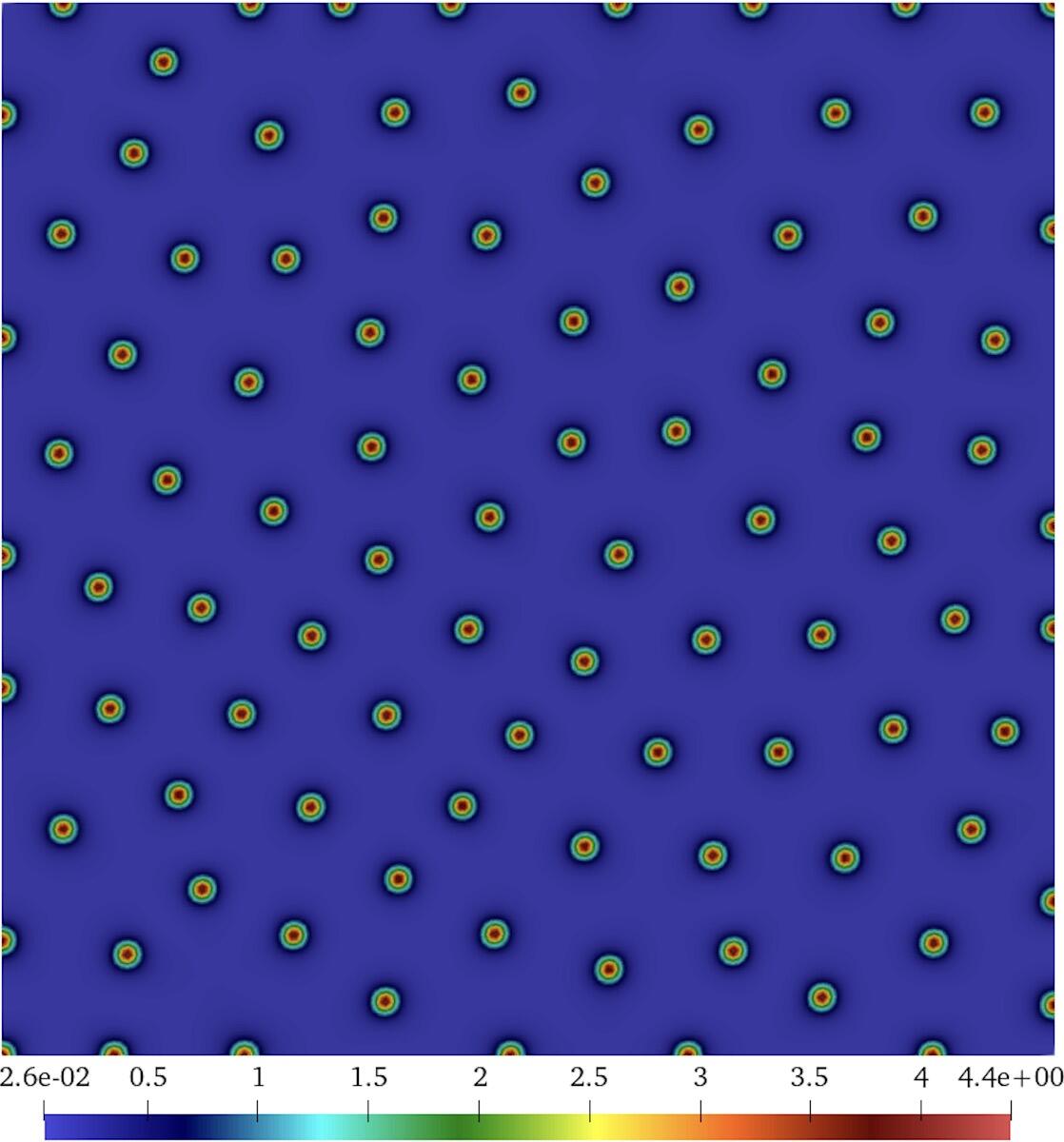}
    \end{overpic}
    \caption{Case 3 - PDE model: Evolution of the level of attractiveness (top) and density (bottom).}
    \label{fig:case3}
\end{figure}

\begin{figure}[htb!]
     \centering
         \begin{overpic}[percent,width=0.19\textwidth]{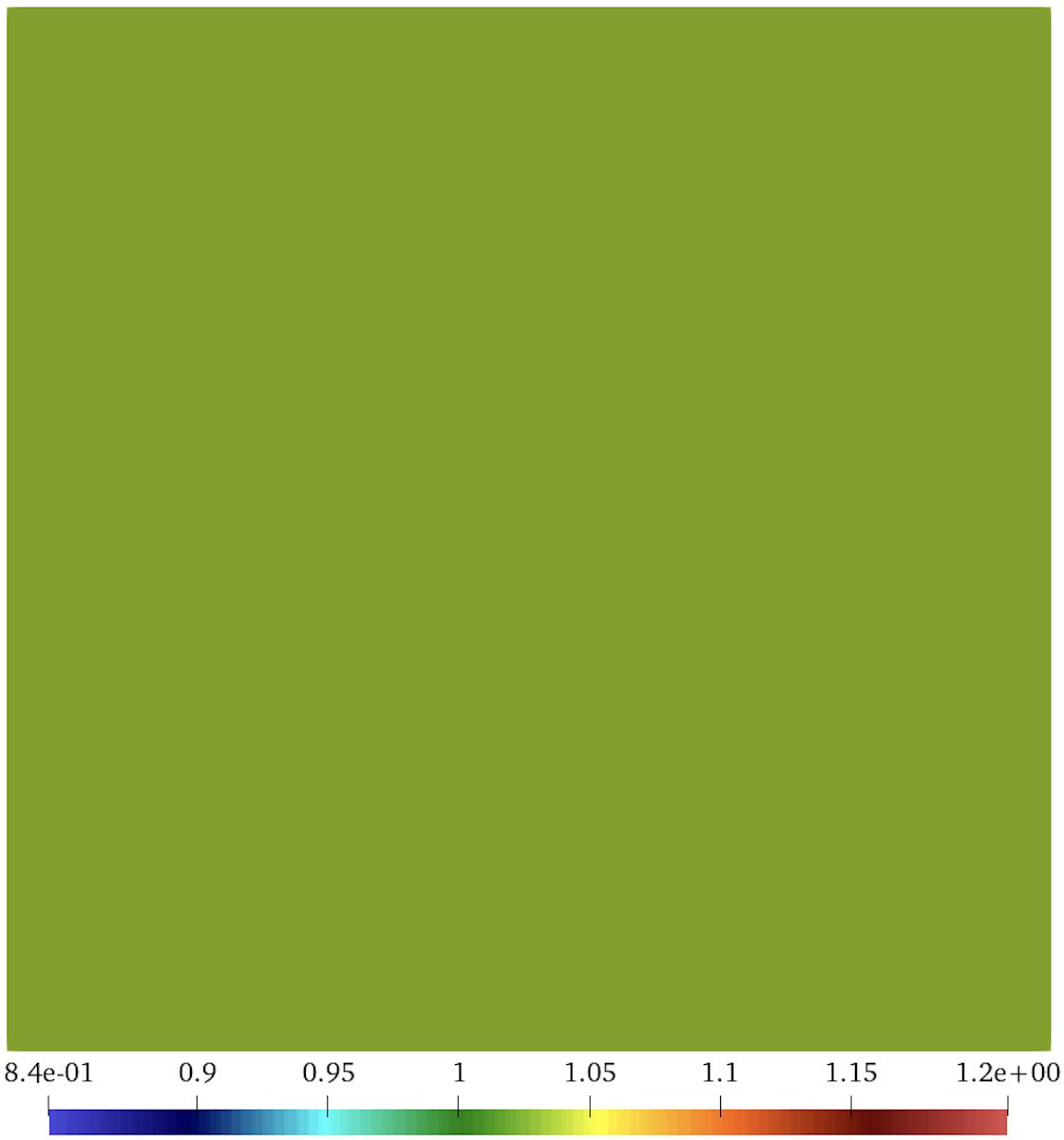}
    \put(30,103){\footnotesize{$t = 0$}}
    \put(-10,55){\footnotesize{$A$}}
    \end{overpic} 
     \begin{overpic}[percent,width=0.19\textwidth, grid=false]{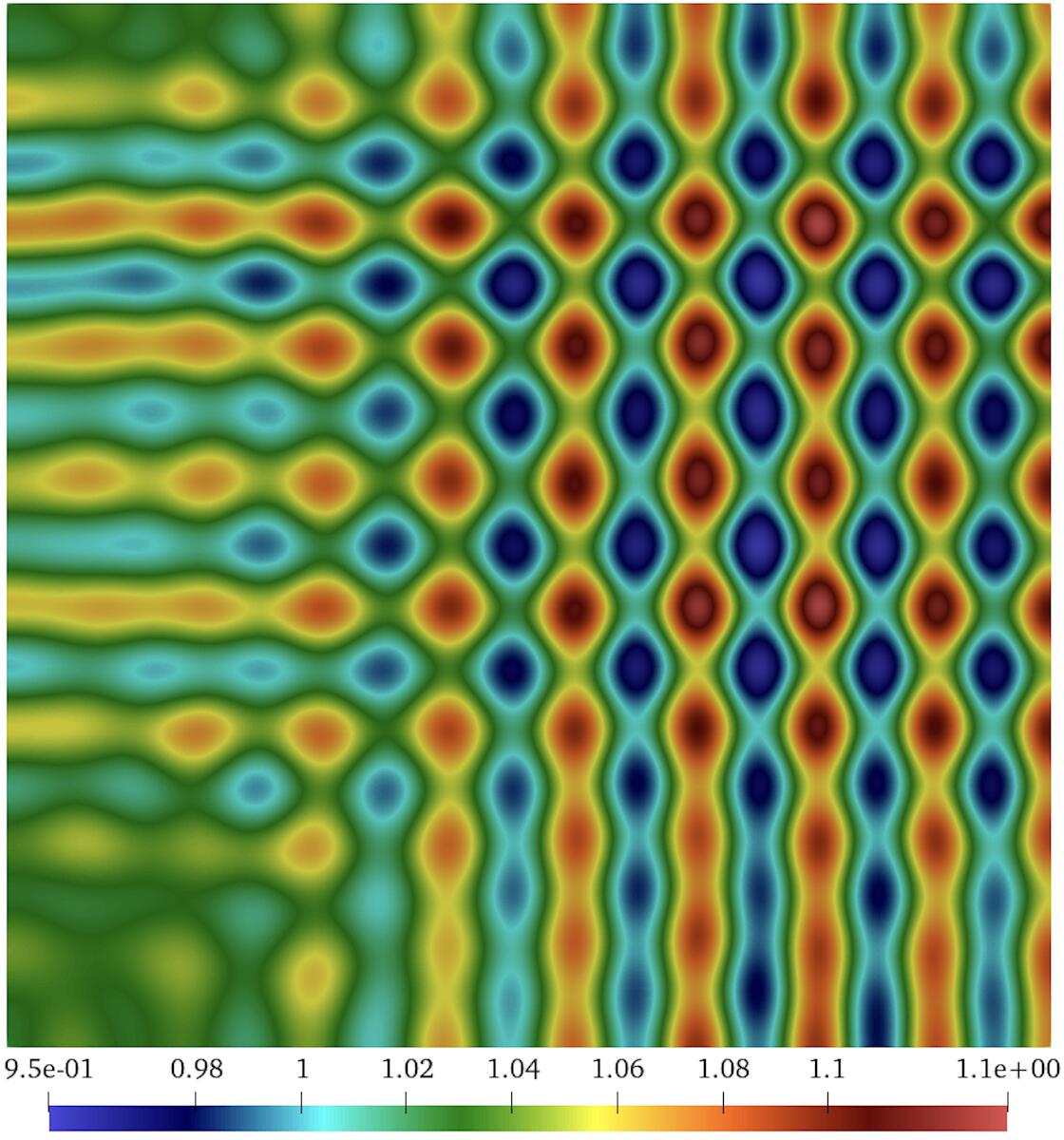}
    \put(35,103){\footnotesize{$t = 25$}}
    \end{overpic}
    \begin{overpic}[percent,width=0.19\textwidth]{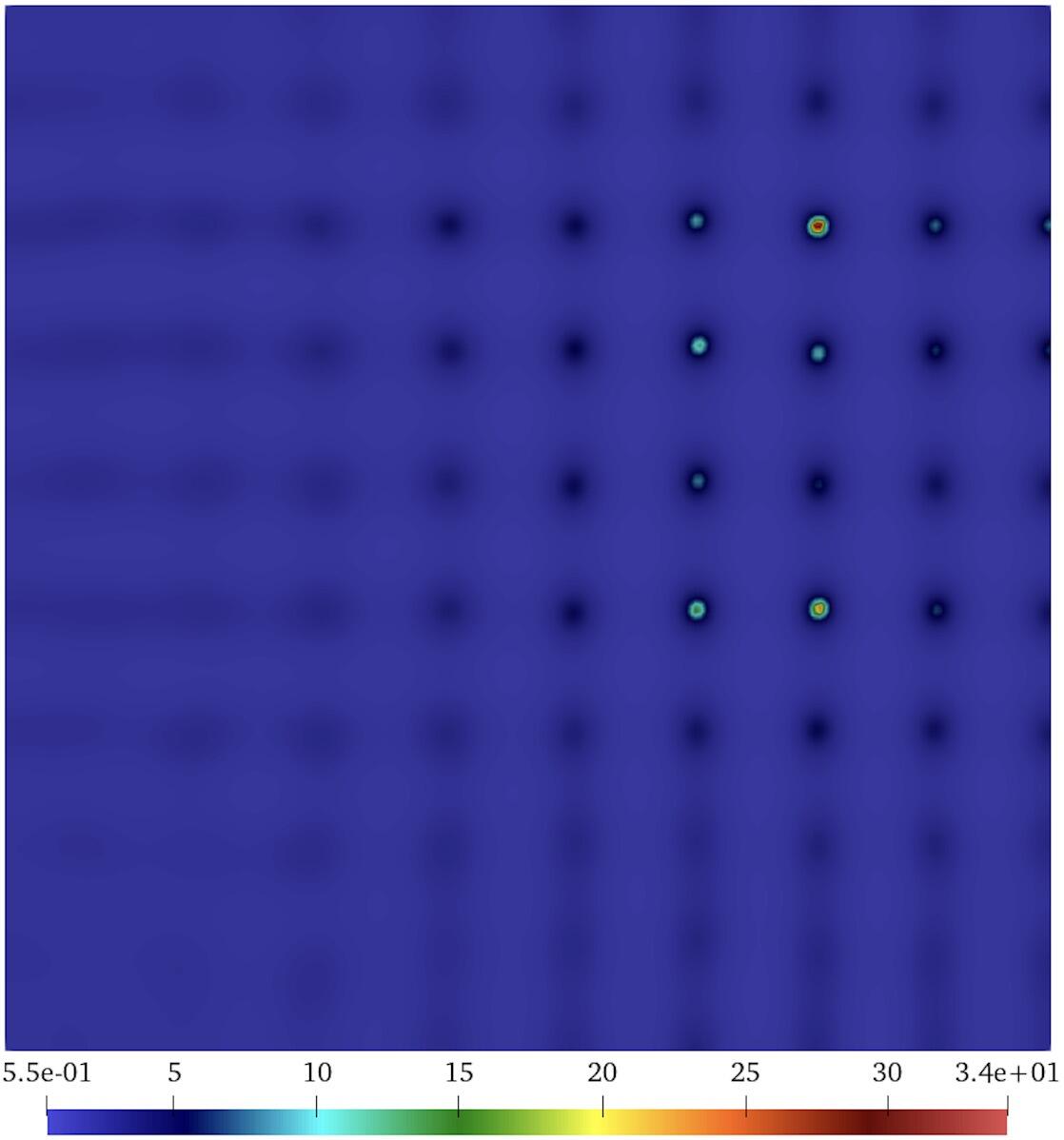}
    \put(35,103){\footnotesize{$t = 27$}}
    \end{overpic} 
        \begin{overpic}[percent,width=0.19\textwidth]{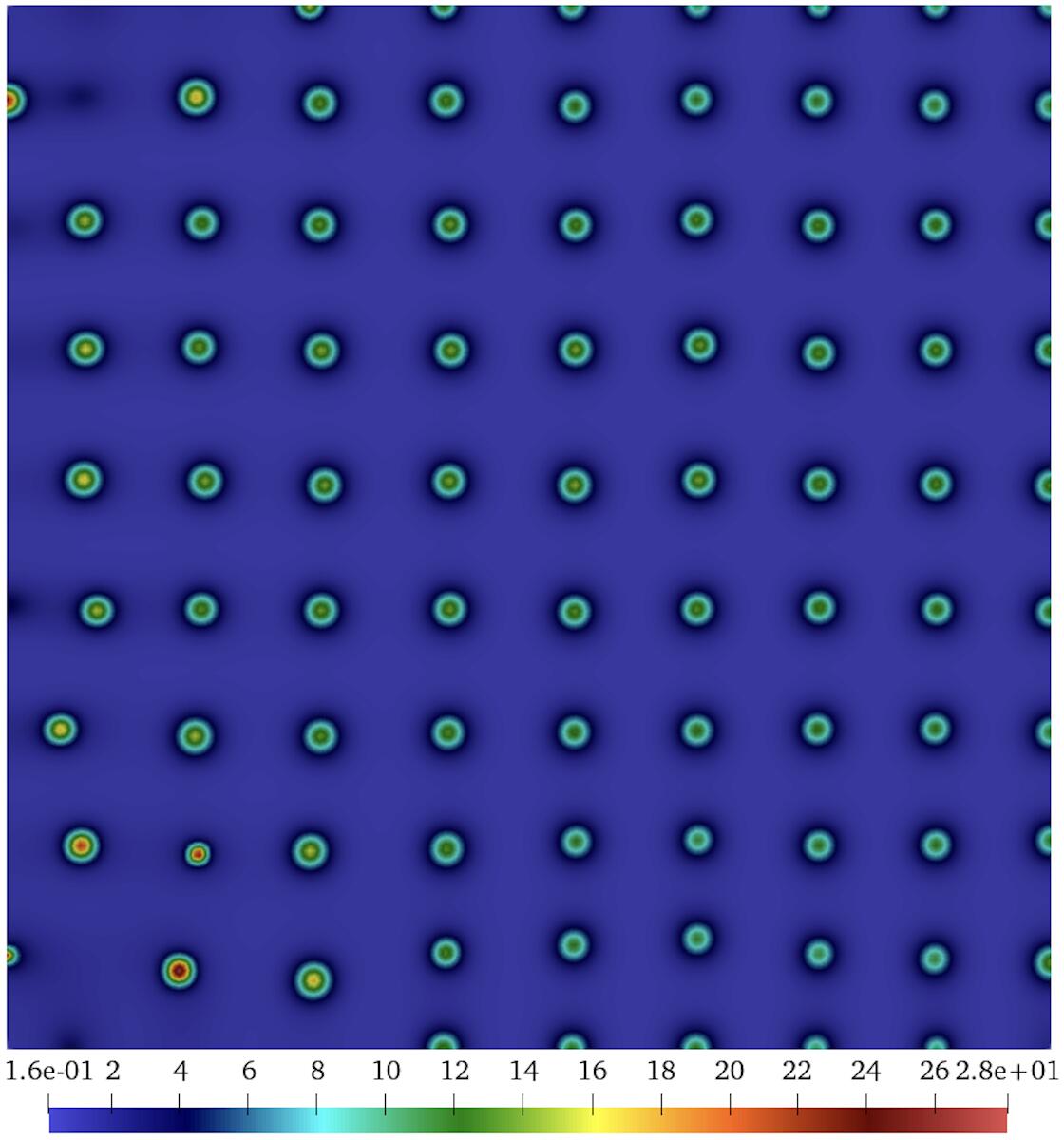}
    \put(32,103){\footnotesize{$t = 29$}}
    \end{overpic} 
    \begin{overpic}[percent,width=0.19\textwidth]{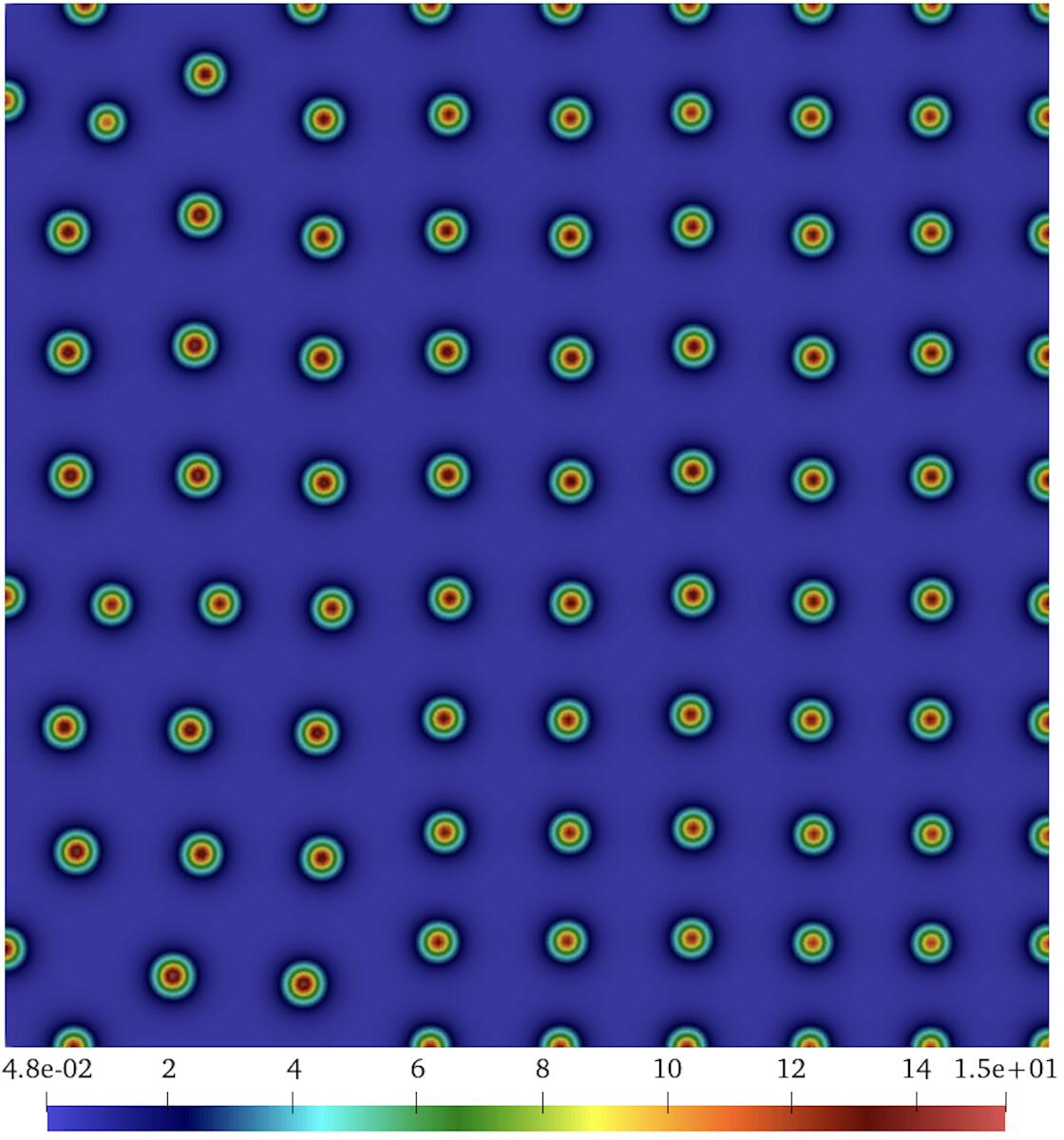}
    \put(30,103){\footnotesize{$t = 200$}}
    \end{overpic}
    \\
         \begin{overpic}[percent,width=0.19\textwidth]{Figs_ContinummModel/S2_rho_T0_nonoise.jpg}
    \put(-10,55){\footnotesize{$\rho$}}
    \end{overpic} 
     \begin{overpic}[percent,width=0.19\textwidth, grid=false]{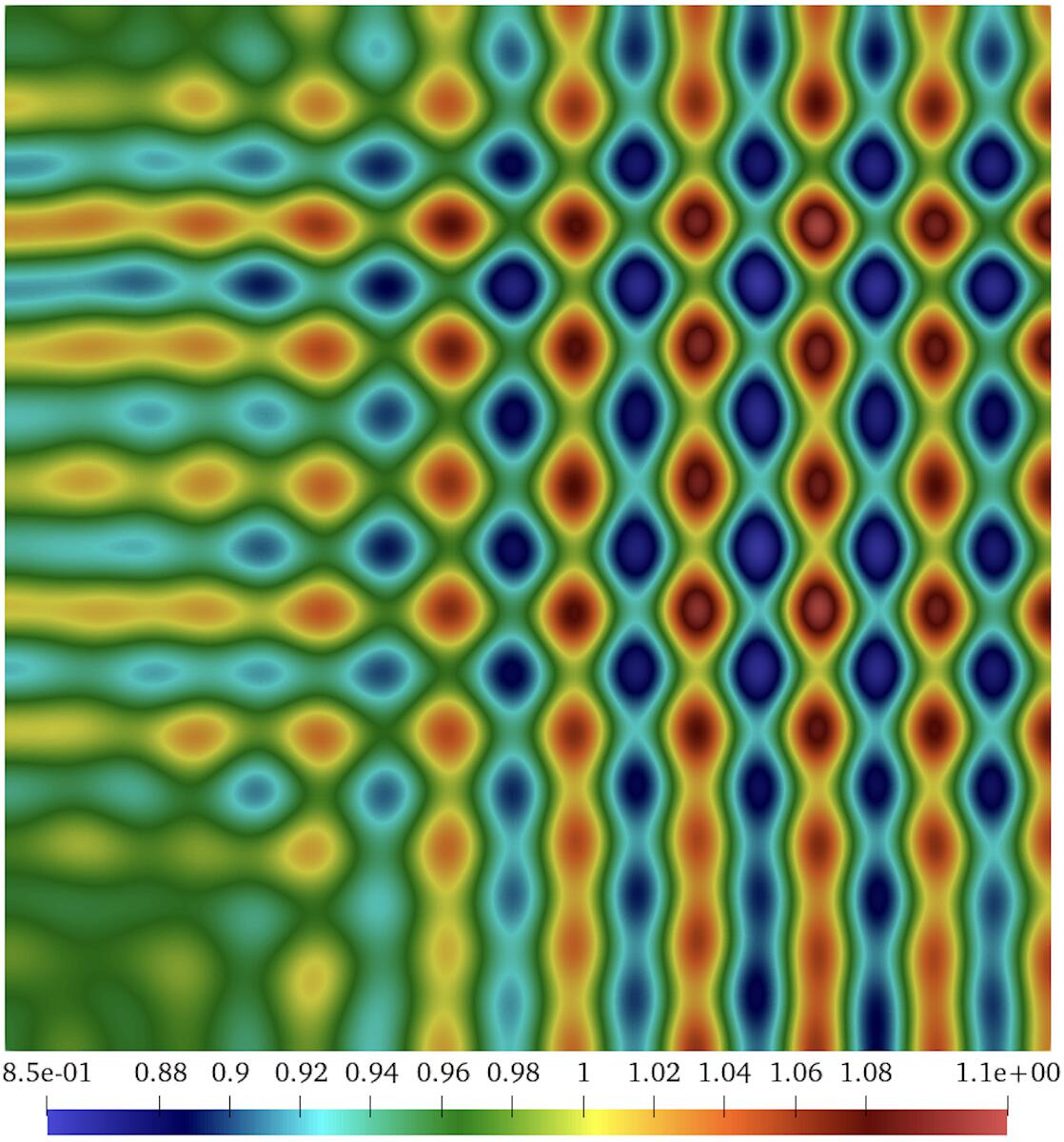}
    \end{overpic}
    \begin{overpic}[percent,width=0.19\textwidth]{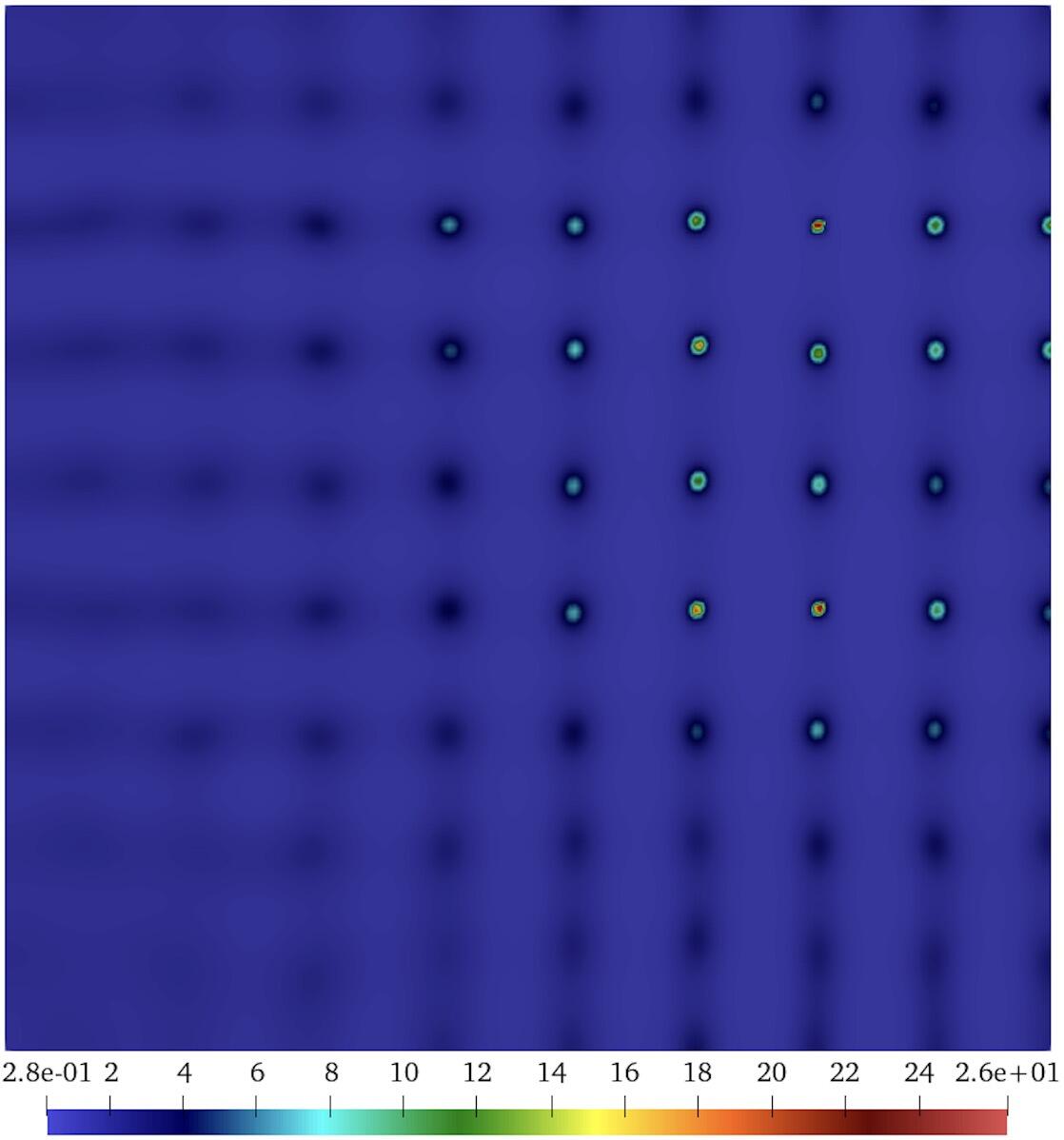}
    \end{overpic} 
        \begin{overpic}[percent,width=0.19\textwidth]{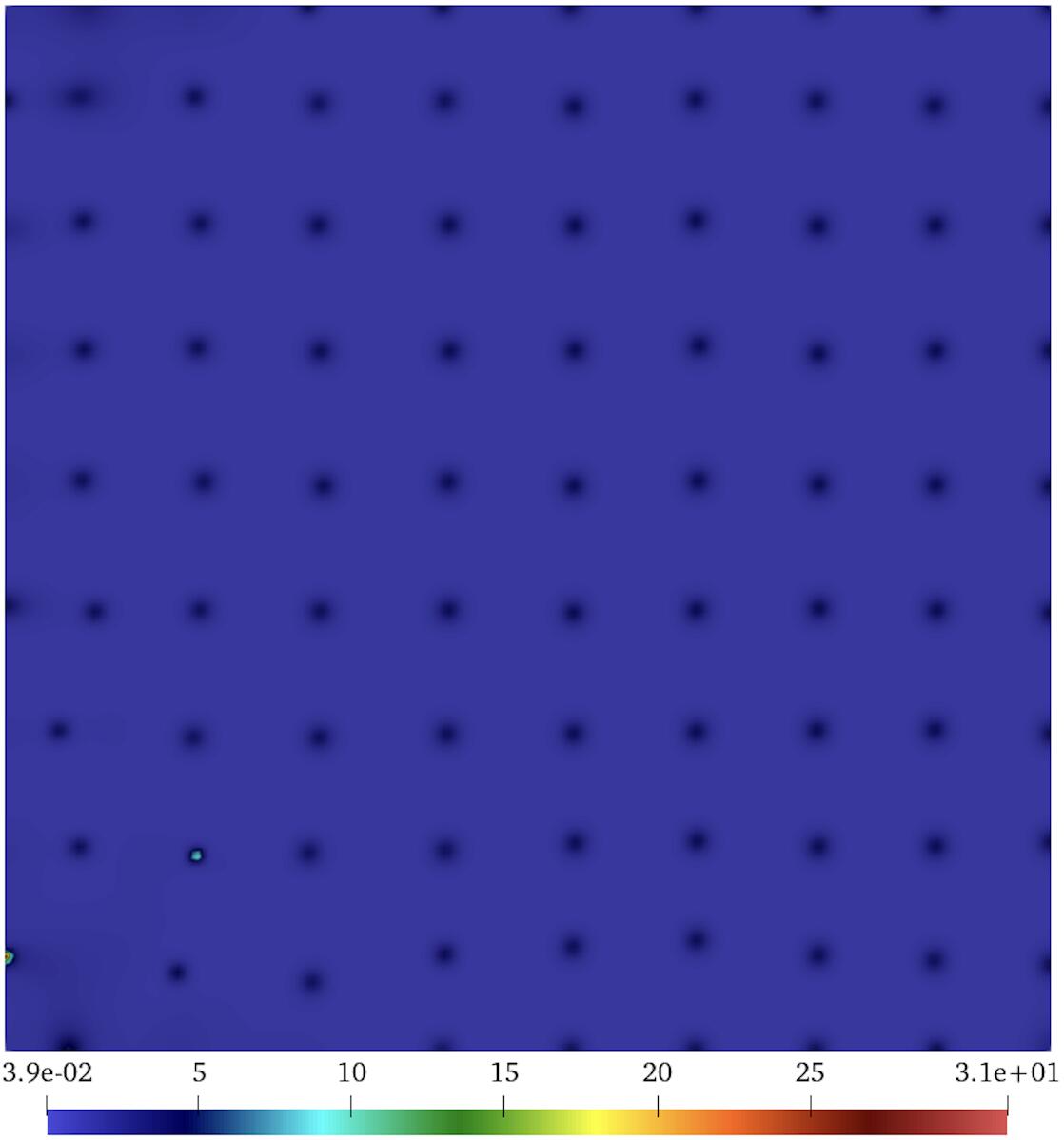}
    \end{overpic} 
    \begin{overpic}[percent,width=0.19\textwidth]{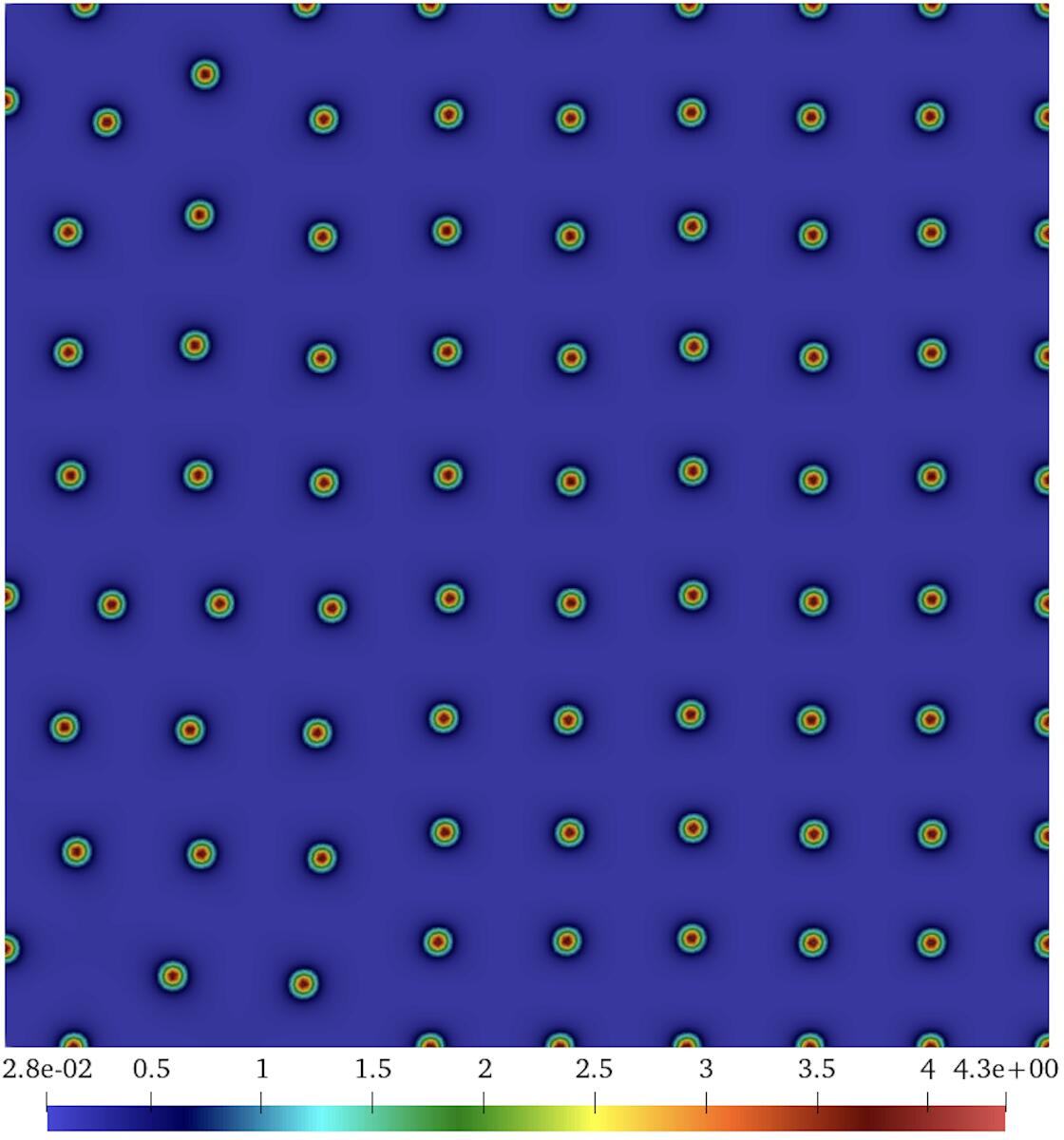}
 
    \end{overpic}
    \caption{Case 3 without noise - PDE model: Evolution of the level of attractiveness (top) and density (bottom).}
    \label{fig:case3-nonoise}
\end{figure}

Fig.~\ref{fig:case3} illustrates the evolution of attractiveness and density for case 3 given by the PDE solver. We observe the early emergence of small hotspots. Indeed, by $t=10$, distinct small hotspots are already visible in Fig.~\ref{fig:case3}, whereas in Fig.~\ref{fig:case2}, the large-scale hotspot structures have not yet emerged. This comparison highlights that a smaller value of $\eta$ leads to a faster onset of 
the hotspot pattern. For completeness, 
we report in Fig.~\ref{fig:case3-nonoise}
the evolution of $A$ and $\rho$ for case 3 with no noise. Again, we see that, in the absence of noise, the hotspot pattern is structured
and the pattern emergence is delayed.


\begin{figure}[htb!]
     \centering
         \begin{overpic}[percent,width=0.19\textwidth]{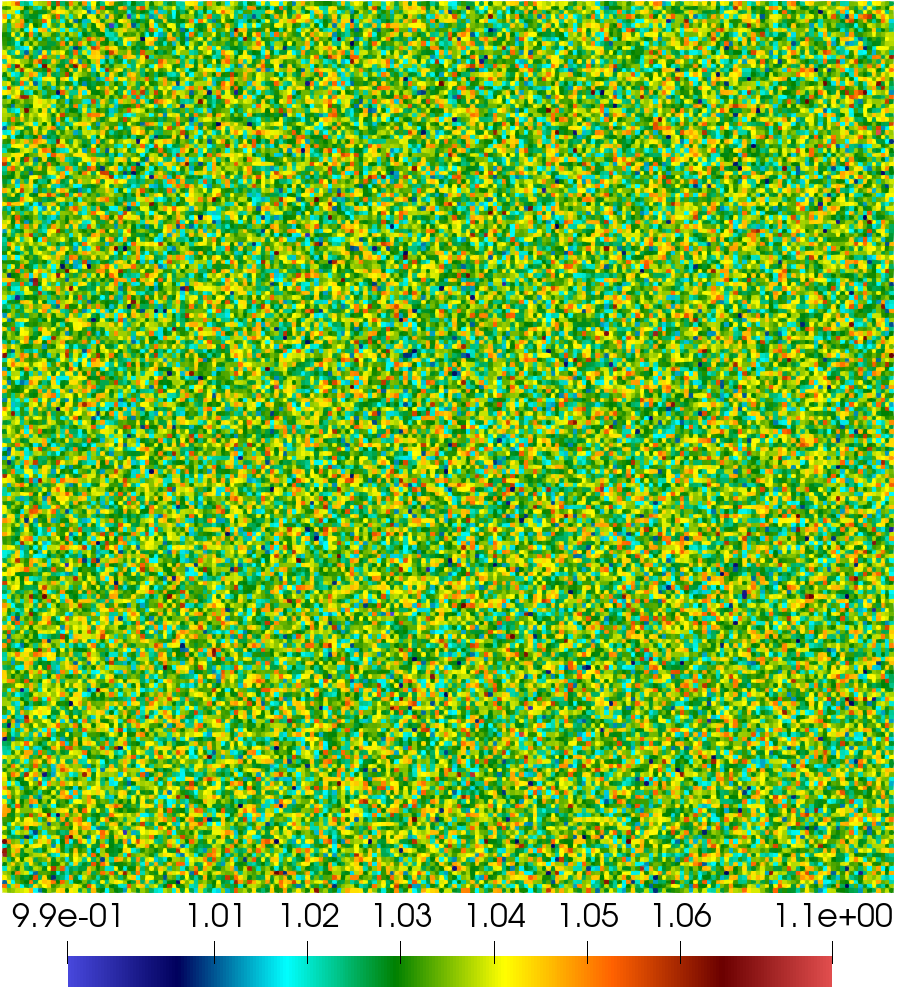}
    \put(30,103){\footnotesize{$t = 0$}}
    \put(-10,55){\footnotesize{$A$}}
    \end{overpic} 
     \begin{overpic}[percent,width=0.19\textwidth, grid=false]{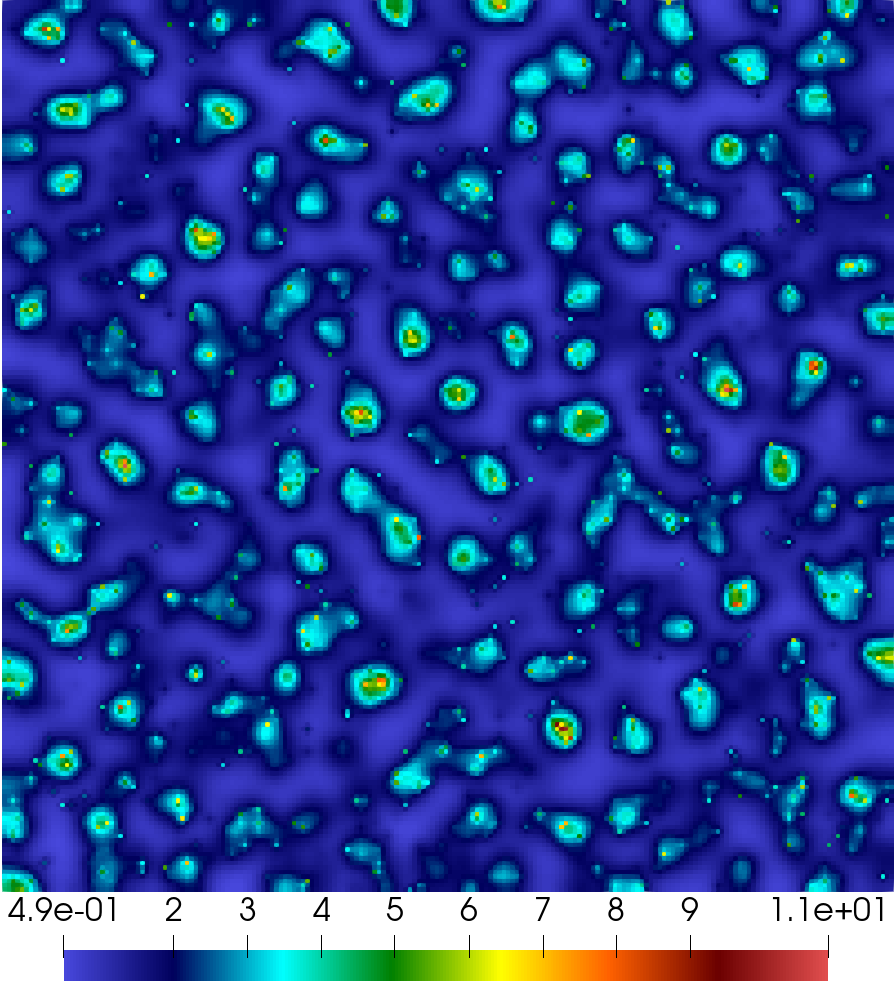}
    \put(35,103){\footnotesize{$t = 1.6$}}
    \end{overpic}
    \begin{overpic}[percent,width=0.19\textwidth]{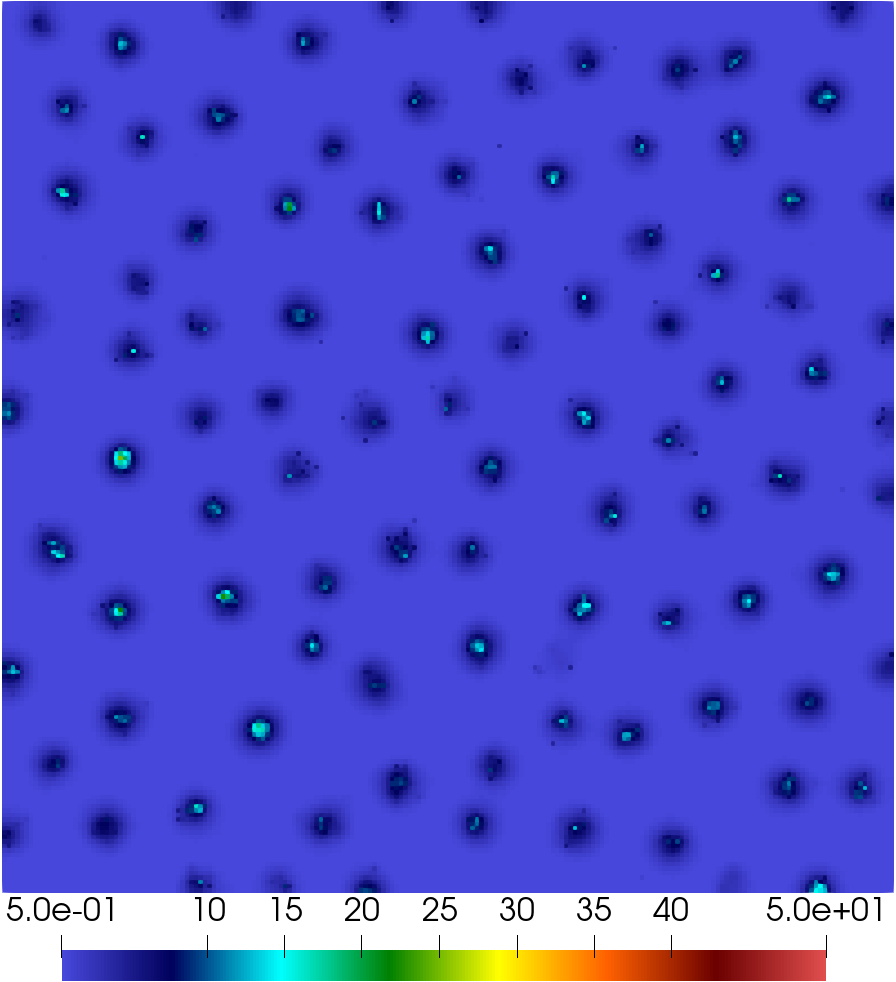}
    \put(35,103){\footnotesize{$t = 10$}}
    \end{overpic} 
        \begin{overpic}[percent,width=0.19\textwidth]{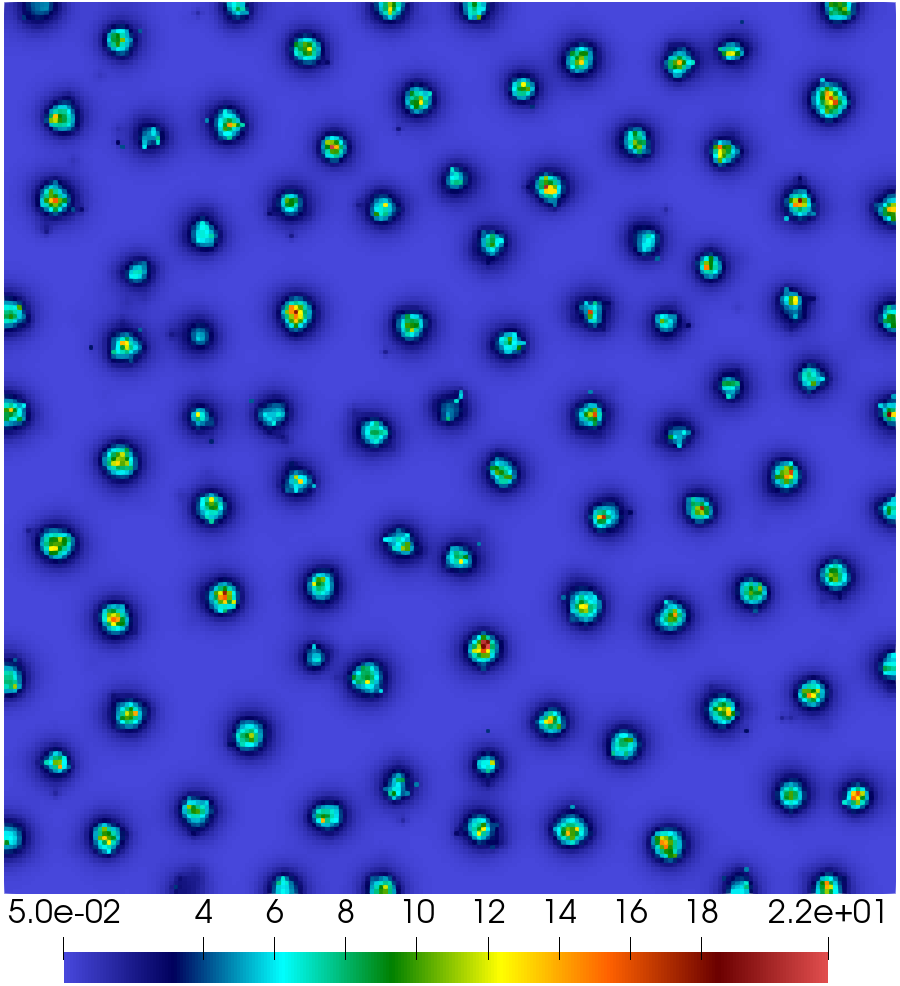}
    \put(32,103){\footnotesize{$t = 20$}}
    \end{overpic} 
    \begin{overpic}[percent,width=0.19\textwidth]{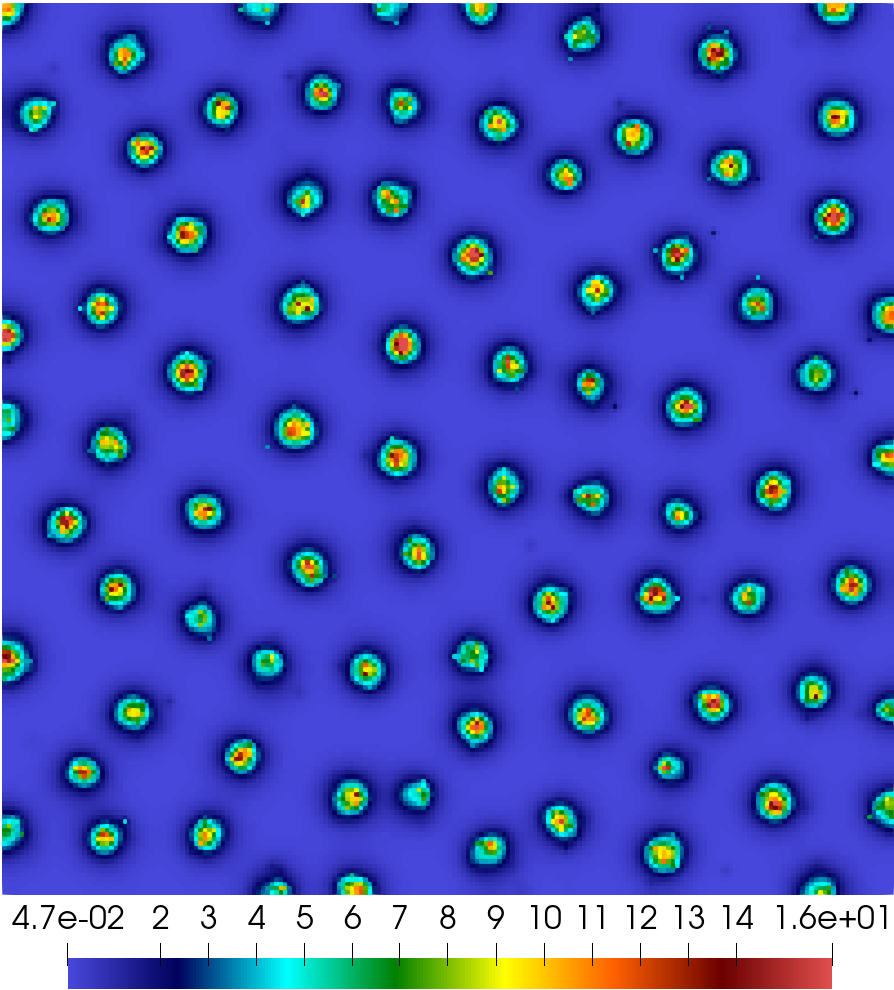}
    \put(30,103){\footnotesize{$t = 200$}}
    \end{overpic}
    \\
         \begin{overpic}[percent,width=0.19\textwidth]{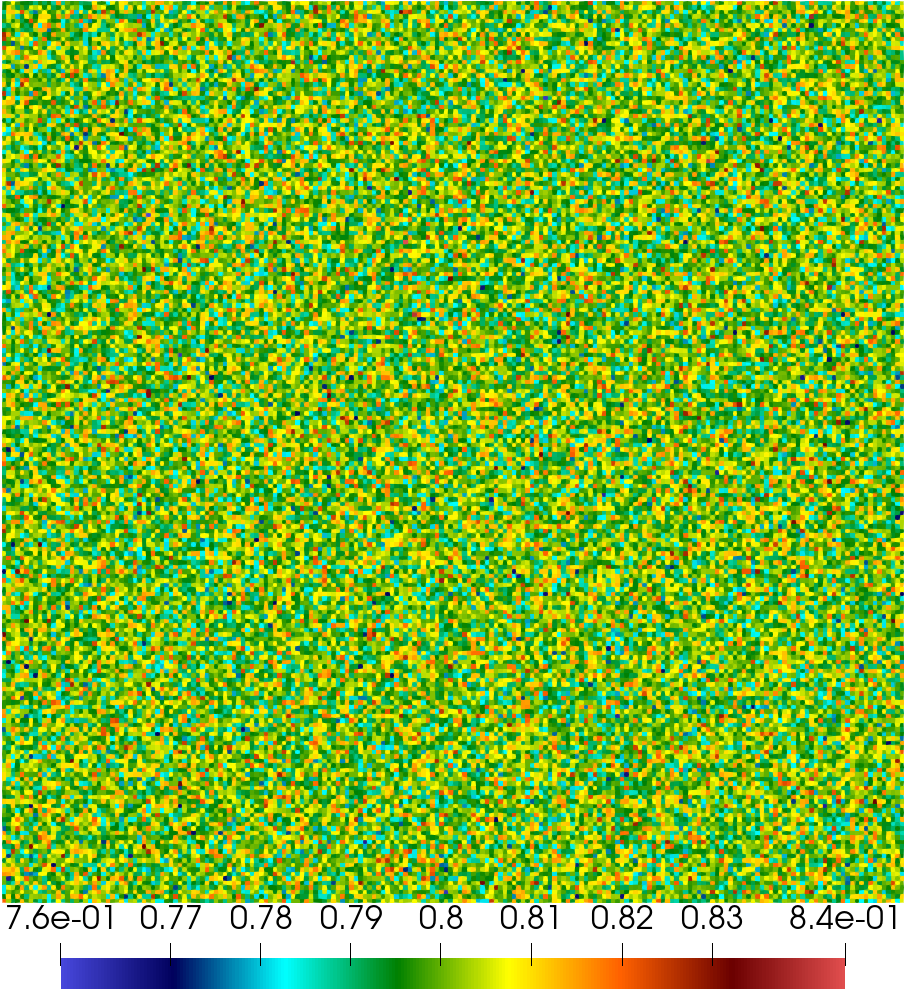}
    \put(-10,55){\footnotesize{$\rho$}}
    \end{overpic} 
     \begin{overpic}[percent,width=0.19\textwidth, grid=false]{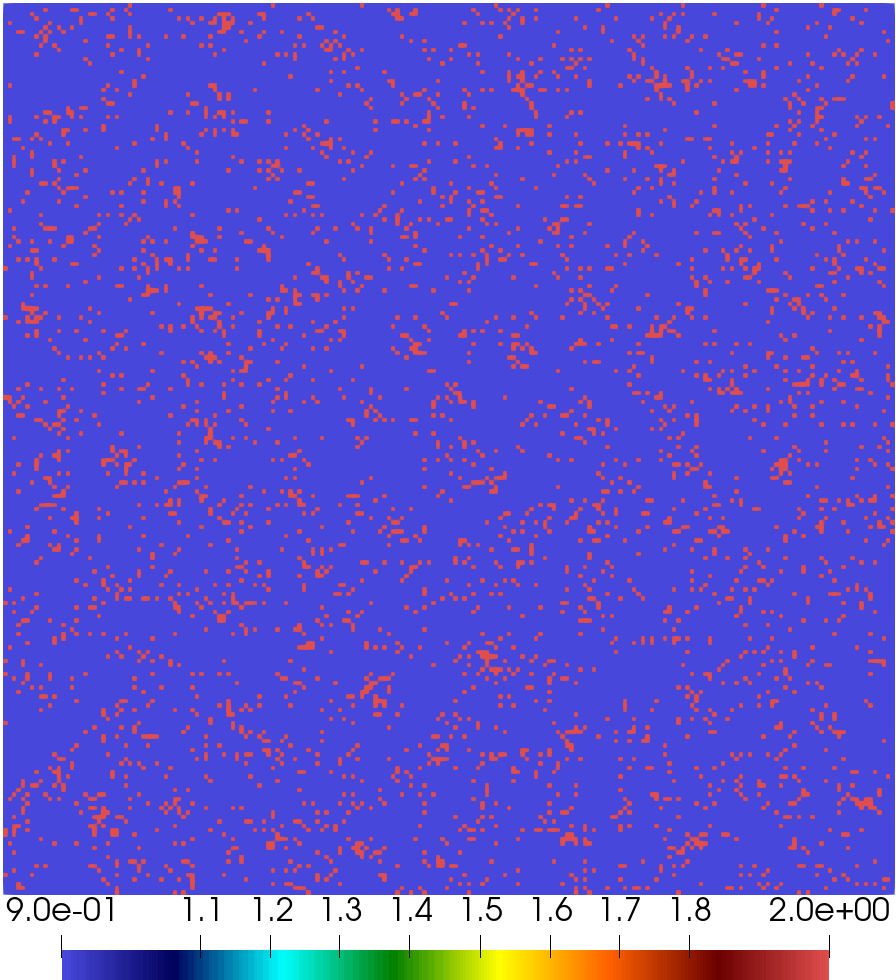}
    \end{overpic}
    \begin{overpic}[percent,width=0.19\textwidth]{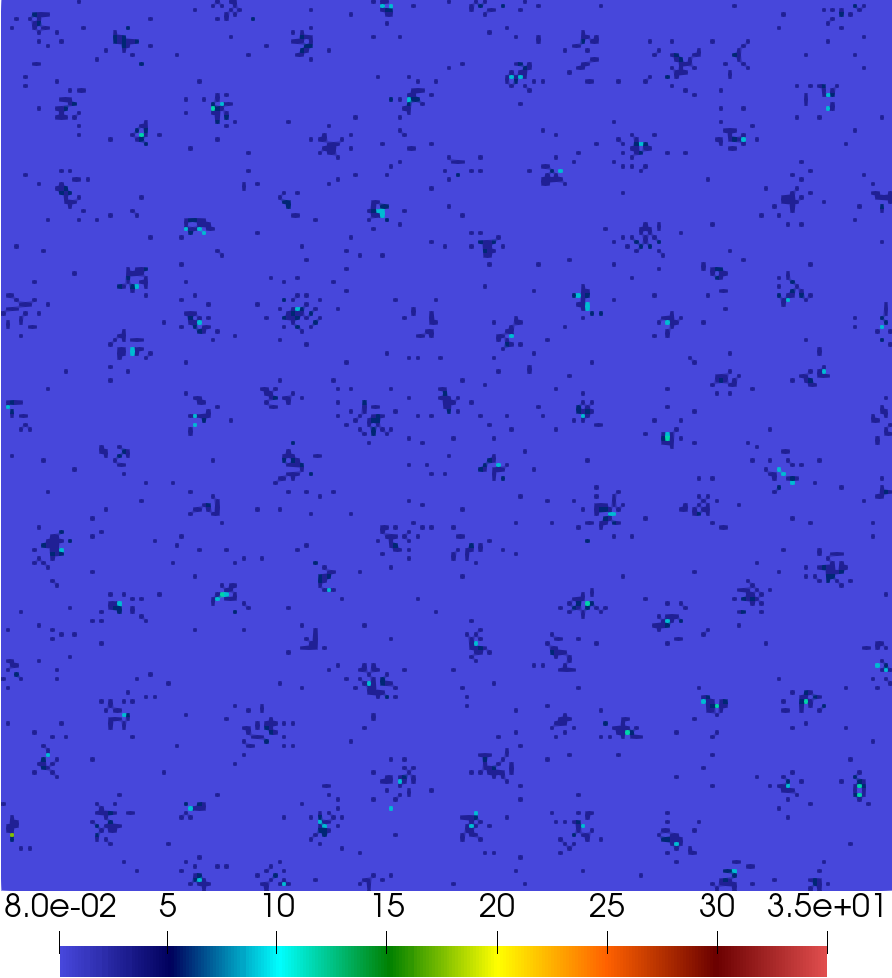}
    \end{overpic} 
        \begin{overpic}[percent,width=0.19\textwidth]{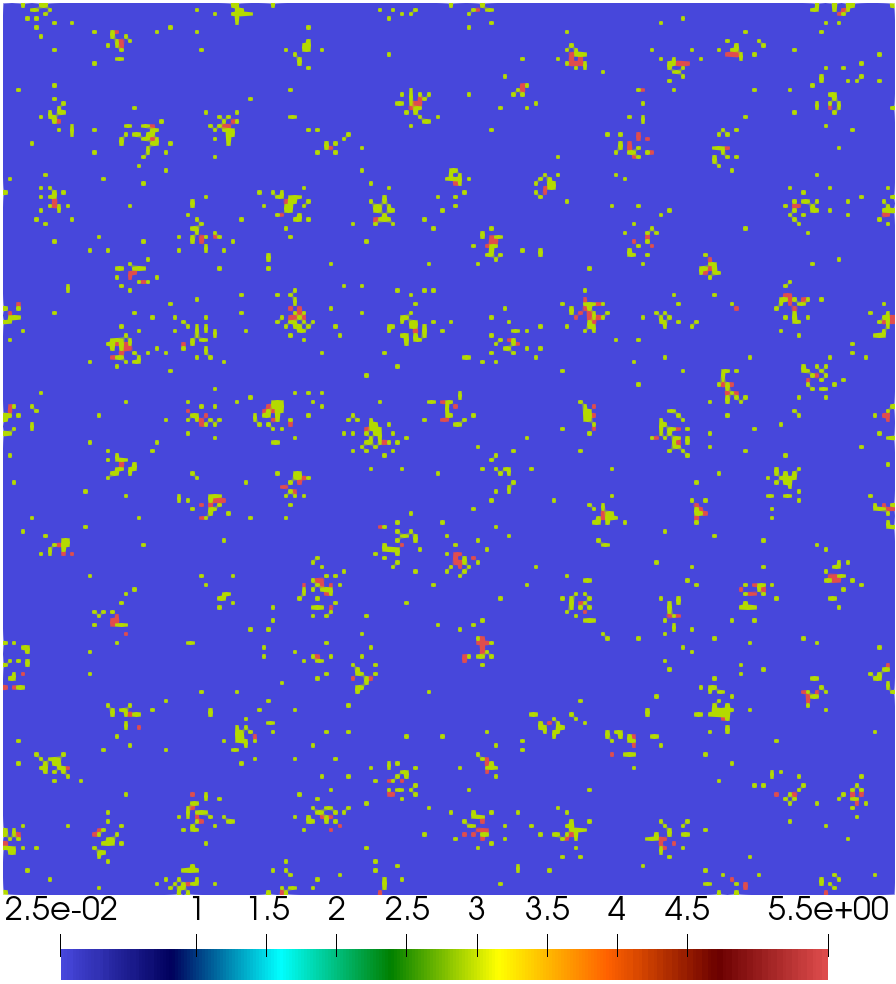}
    \end{overpic} 
    \begin{overpic}[percent,width=0.19\textwidth]{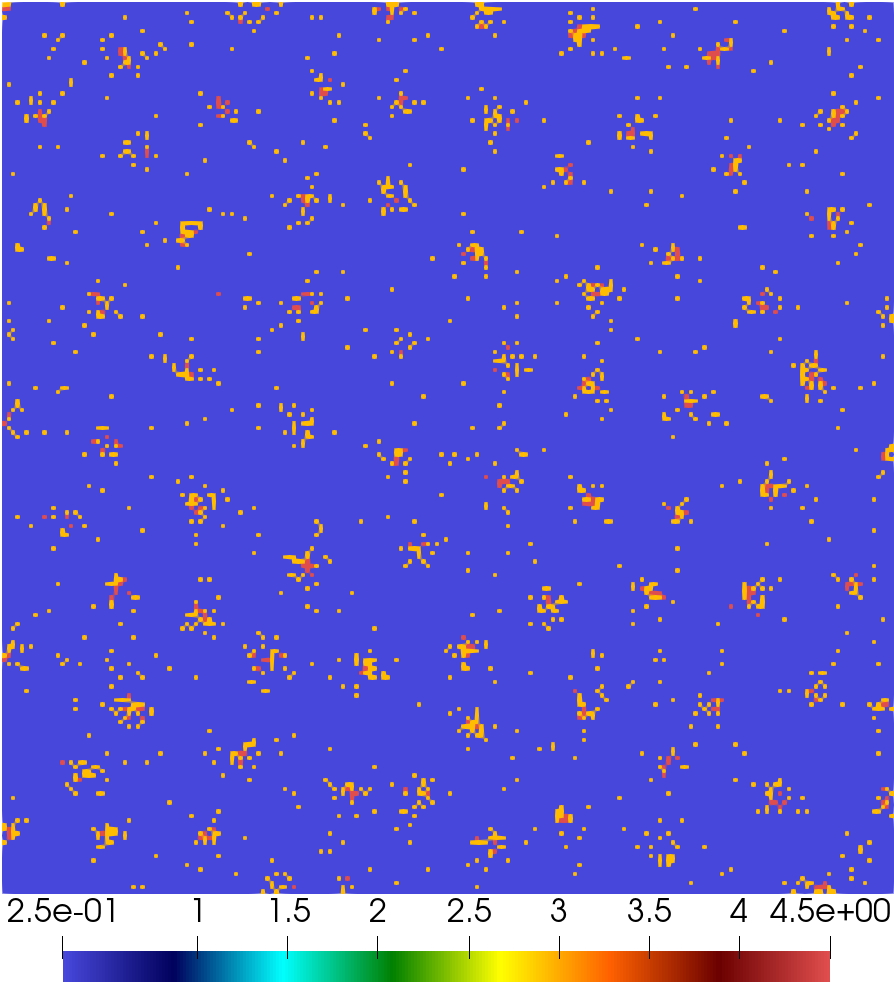}
    \end{overpic}
    \caption{Case 3 - agent-based model: Evolution of the level of attractiveness (top) and criminals (bottom).}
    \label{fig:case3DM}
\end{figure}

Finally, Fig.~\ref{fig:case3DM} shows
the evolution of $A$ and $\rho$ for case 3 given by the agent-based solver after non-dimensionalization. For these results,
we set again $\Gamma = 0.019$, $\theta = 0.2339$, and $\omega = 1/15$, which give
$\Gamma \theta /\omega^2 \approx 1$.
We observe that the hotspots 
in Fig.~\ref{fig:case3}
and \ref{fig:case3DM} are comparable in size and number. Additionally, 
the magnitudes of 
$A$ and $\rho$  match. 


From the results presented in this section, we can conclude that the PDE solver output
matches well the output
from the agent-based simulations when 
the density of criminals is large.

\subsection{Robustness and computational efficiency of the iterative partitioned solver}\label{sec:res_solver}

Now that the PDE solver has been validated, we use it to further investigate the relationship between $\eta$ and the size of the hotspots.
For this purpose, we vary the parameter $\eta$ over the interval $[0.025, 0.45]$.
Next, we examine the influence of discretization parameters and stopping tolerances for the three cases studied in Sec.~\ref{sec:validation}, in order to better understand the robustness and computational efficiency of the iterative partitioned solver.

\begin{figure}[htb!]
    \centering
    \begin{overpic}[width=0.7\textwidth,grid=false]{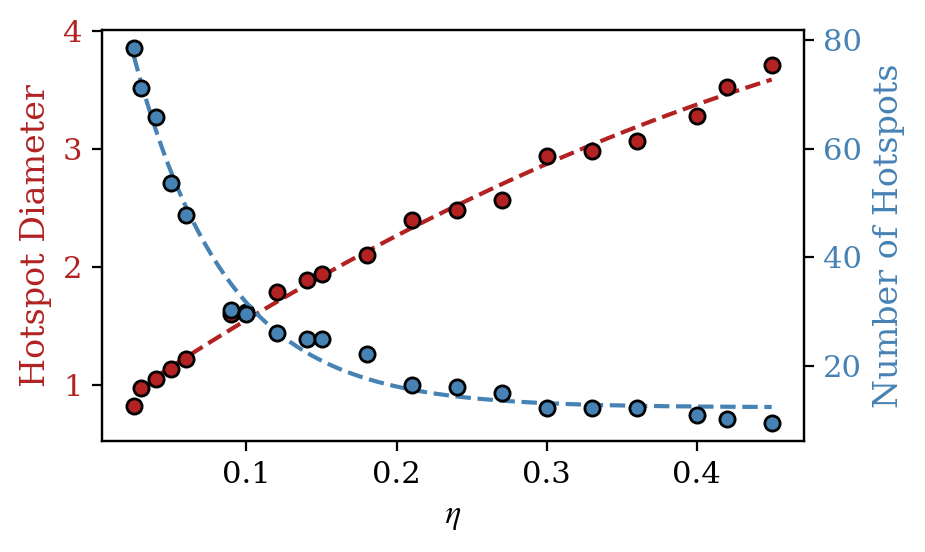}
    \put(55,27){\footnotesize{Case 2}}
    \put(59.05,26){\vector(0,-1){9}}
    \put(59.05,29){\vector(0,1){12}}
    \put(12,30){\footnotesize{Case 3}}
    \put(15,29){\vector(0,-1){10}}
    \put(15,32){\vector(0,1){16}}
    \end{overpic} 
    \caption{Hotspot diameter (left axis, red markers) and number of hotspots (right axis, blue markers) as the strength of the neighborhood effects $\eta$ is varied in interval [0.025,0.45]. The red and blue dashed lines fit the date with the corresponding color.}
    \label{fig:count-diam}
\end{figure}

  
Fig.~\ref{fig:count-diam} reports
the hotspot diameter and number of 
hotspots in domain $\Omega$ for $\eta \in [0.025,0.45]$. The two quantities are
measured from the computed results once the hotspot pattern has fully emerged
and shows no visual changes over several time units. 
To reduce 
the computational time need for this study, we considered a coarser mesh ($h = 16/100$)
and a larger time step ($\Delta t = 1/25$) than the ones used for the figures in the previous section
and we reduced the stopping tolerances ($\texttt{tol}_1 = \texttt{tol}_2 = 10^{-6}$).
Fig.~\ref{fig:count-diam} indicates a nontrivial dependence of both 
diameter and number of hotspots
on $\eta$. A quantitative characterization of the hotspot number and diameter is achieved via empirical curve fitting. Specifically, the number of hotspots is well approximated by the exponential decay  $N(\eta) \approx 96.82 * \exp(-16.14 \eta)+12.33$, while the hotspot diameter follows a quadratic law  $D(\eta) \approx -5.24 \eta^2+8.70\eta+0.73$. We emphasize that these functional forms are obtained via numerical regression and currently lack a rigorous analytical derivation. Developing such theoretical justification would involve intricate analytical work, which lies beyond the scope of the present work. We refer the reader to related analytical studies in similar contexts, such as \cite{alikakos2004ostwald,bray2002theory} for the Cahn–Hilliard equation.

\begin{figure}[htb!]
     \centering
         \begin{overpic}[percent,width=0.19\textwidth]{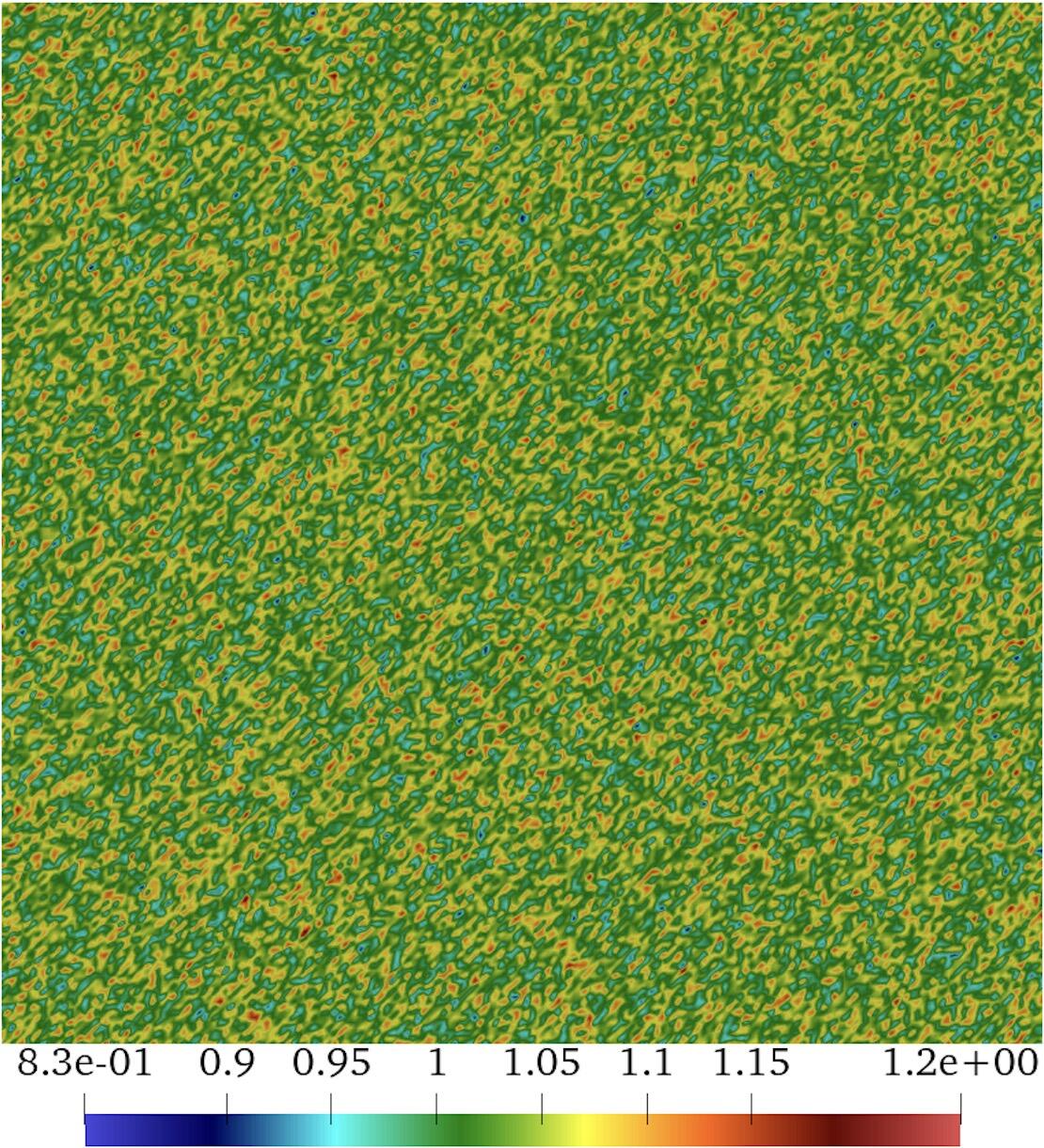}
    \put(30,103){\footnotesize{$t = 0$}}
    \put(-10,55){\footnotesize{$A$}}
    \end{overpic} 
     \begin{overpic}[percent,width=0.19\textwidth, grid=false]{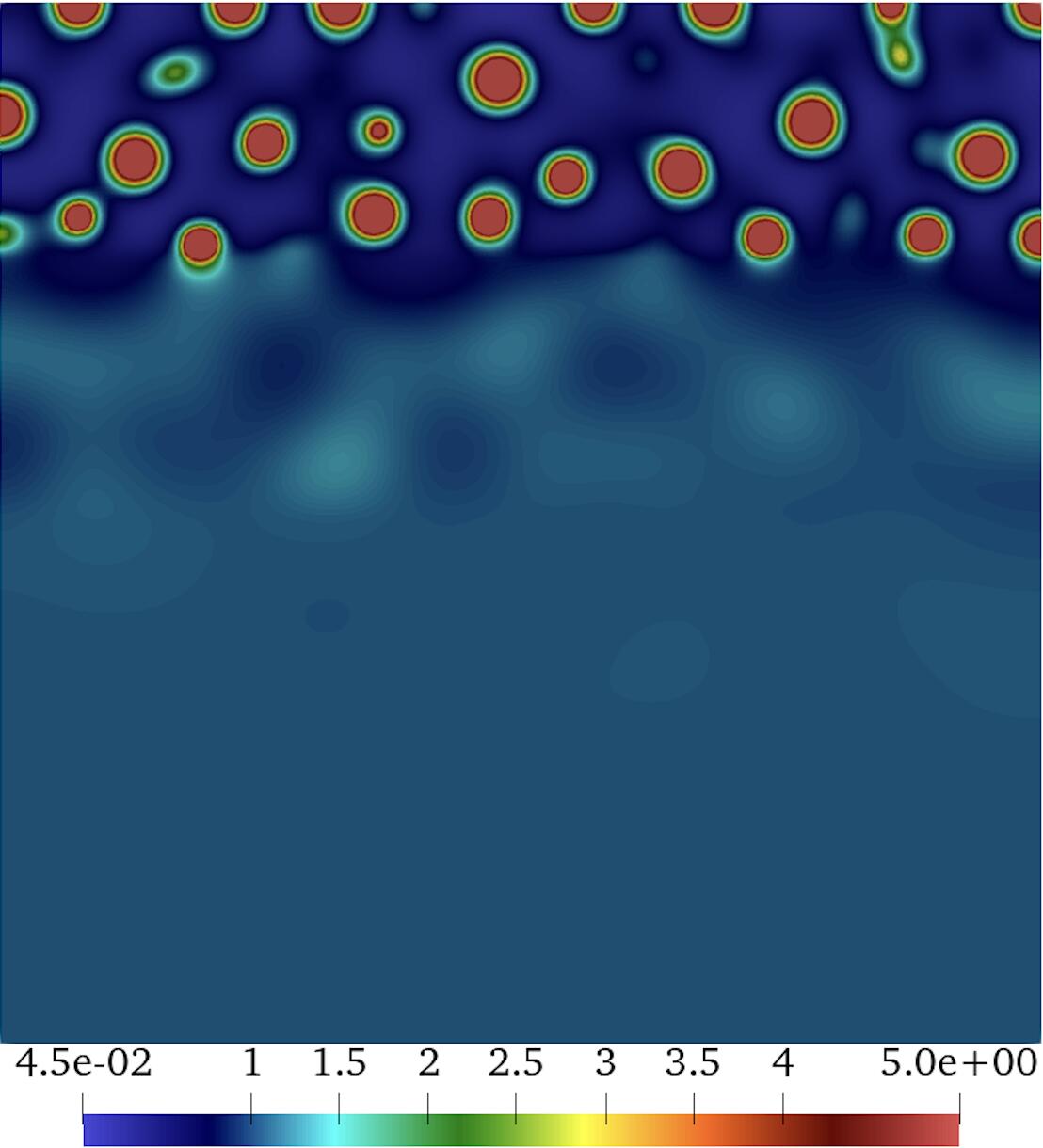}
    \put(35,103){\footnotesize{$t = 5$}}
    \end{overpic}
    \begin{overpic}[percent,width=0.19\textwidth]{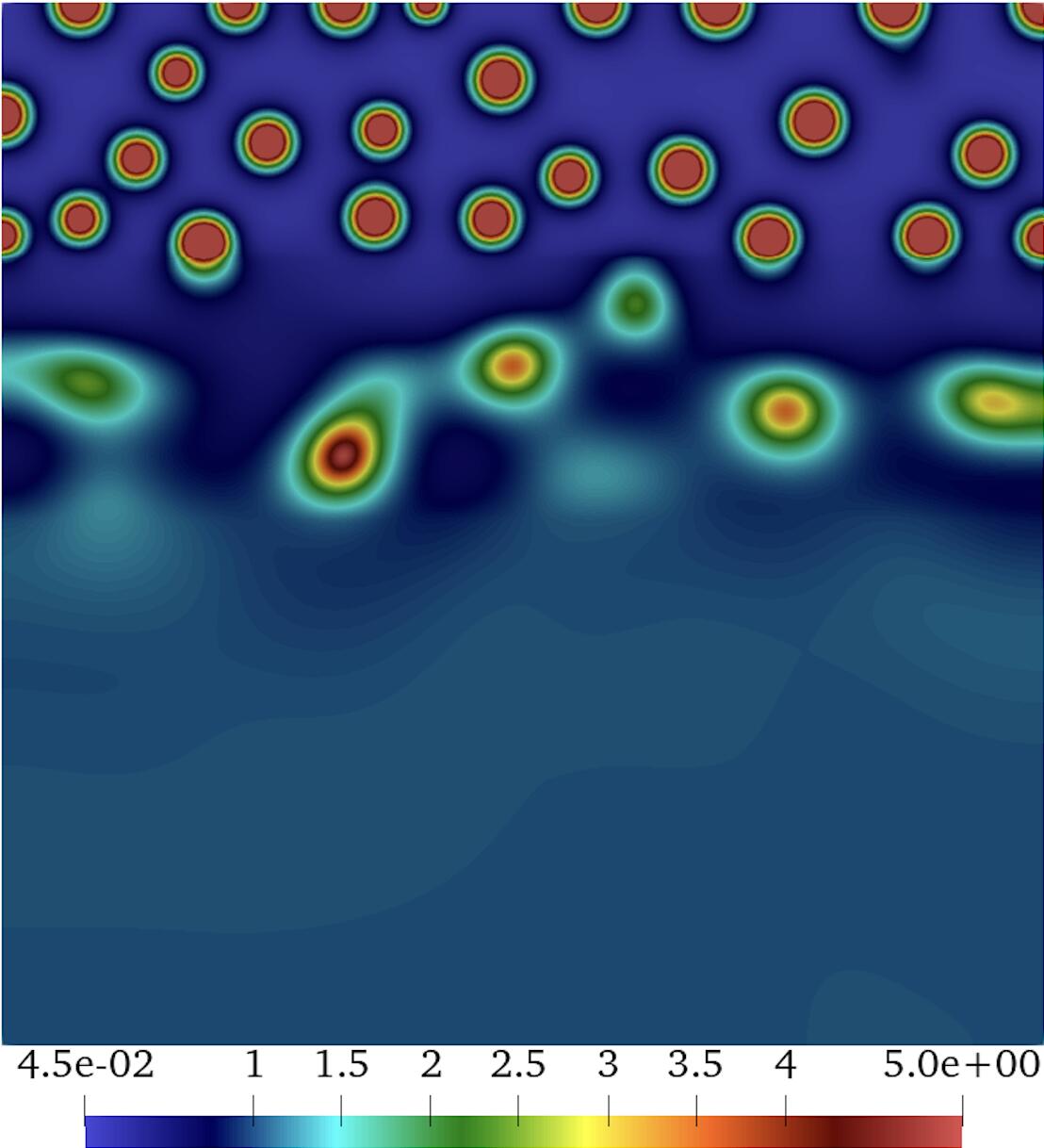}
    \put(35,103){\footnotesize{$t = 8$}}
    \end{overpic} 
        \begin{overpic}[percent,width=0.19\textwidth]{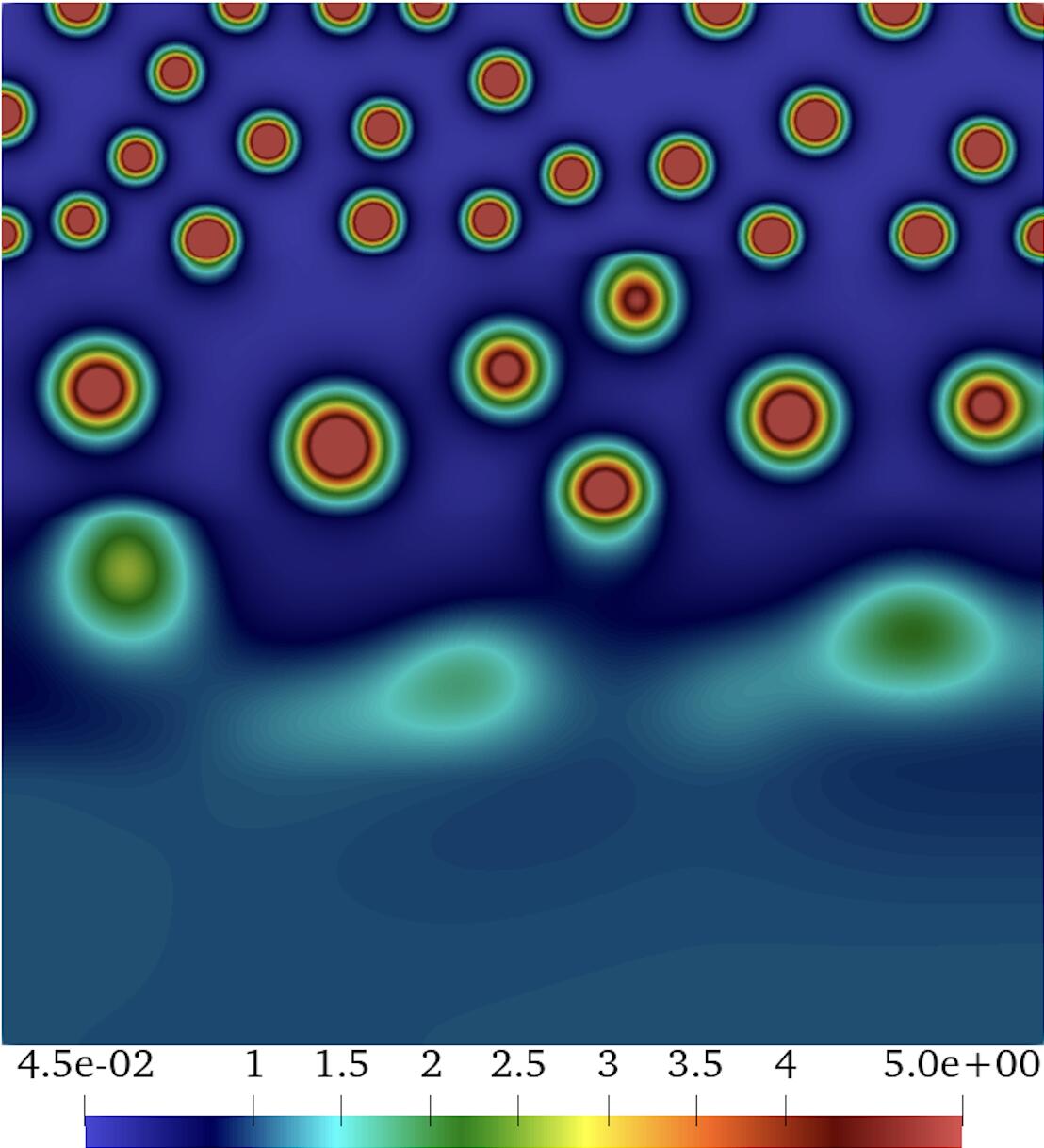}
    \put(32,103){\footnotesize{$t = 15$}}
    \end{overpic} 
    \begin{overpic}[percent,width=0.19\textwidth]{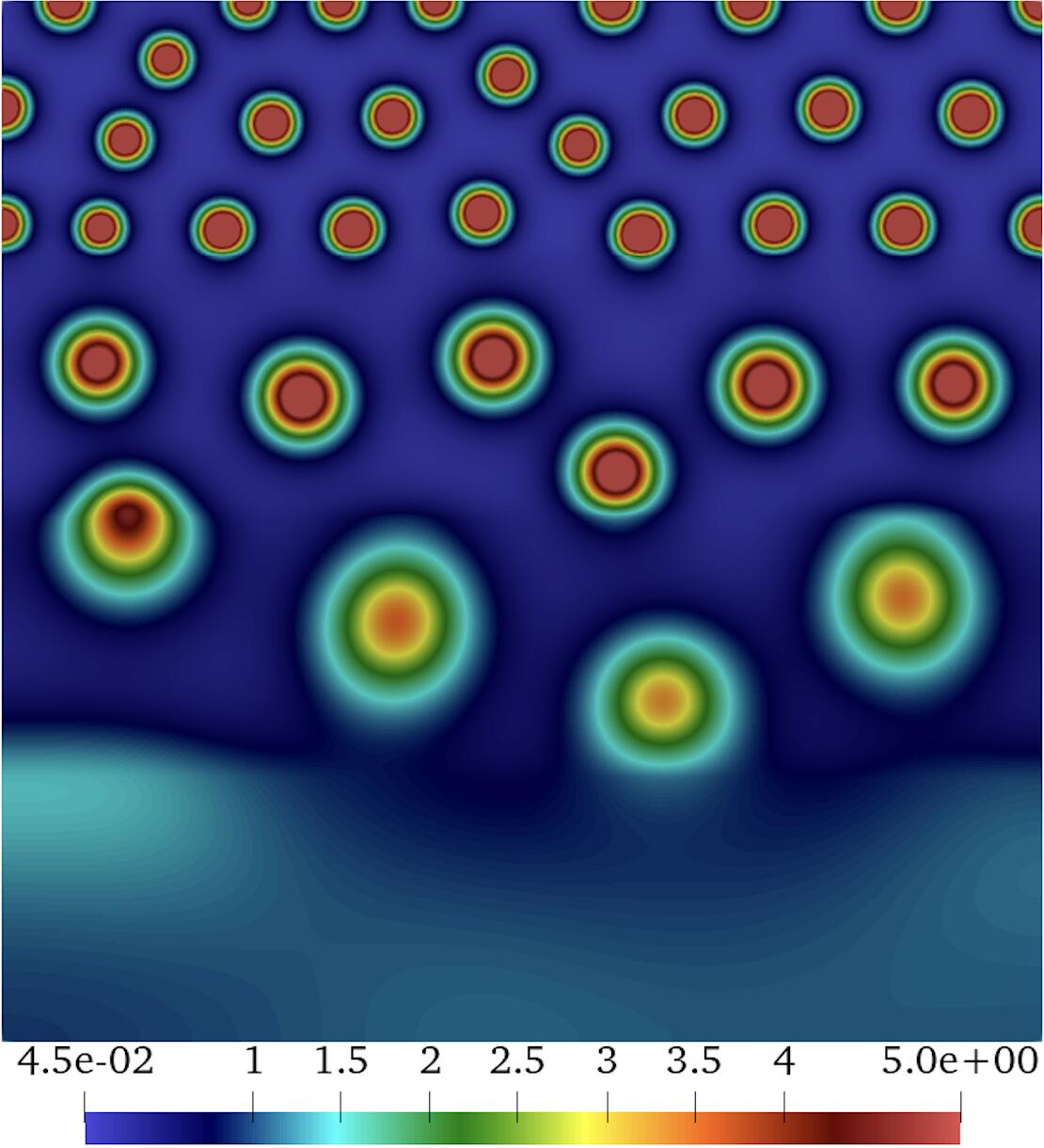}
    \put(30,103){\footnotesize{$t = 200$}}
    \end{overpic}
    \\
         \begin{overpic}[percent,width=0.19\textwidth]{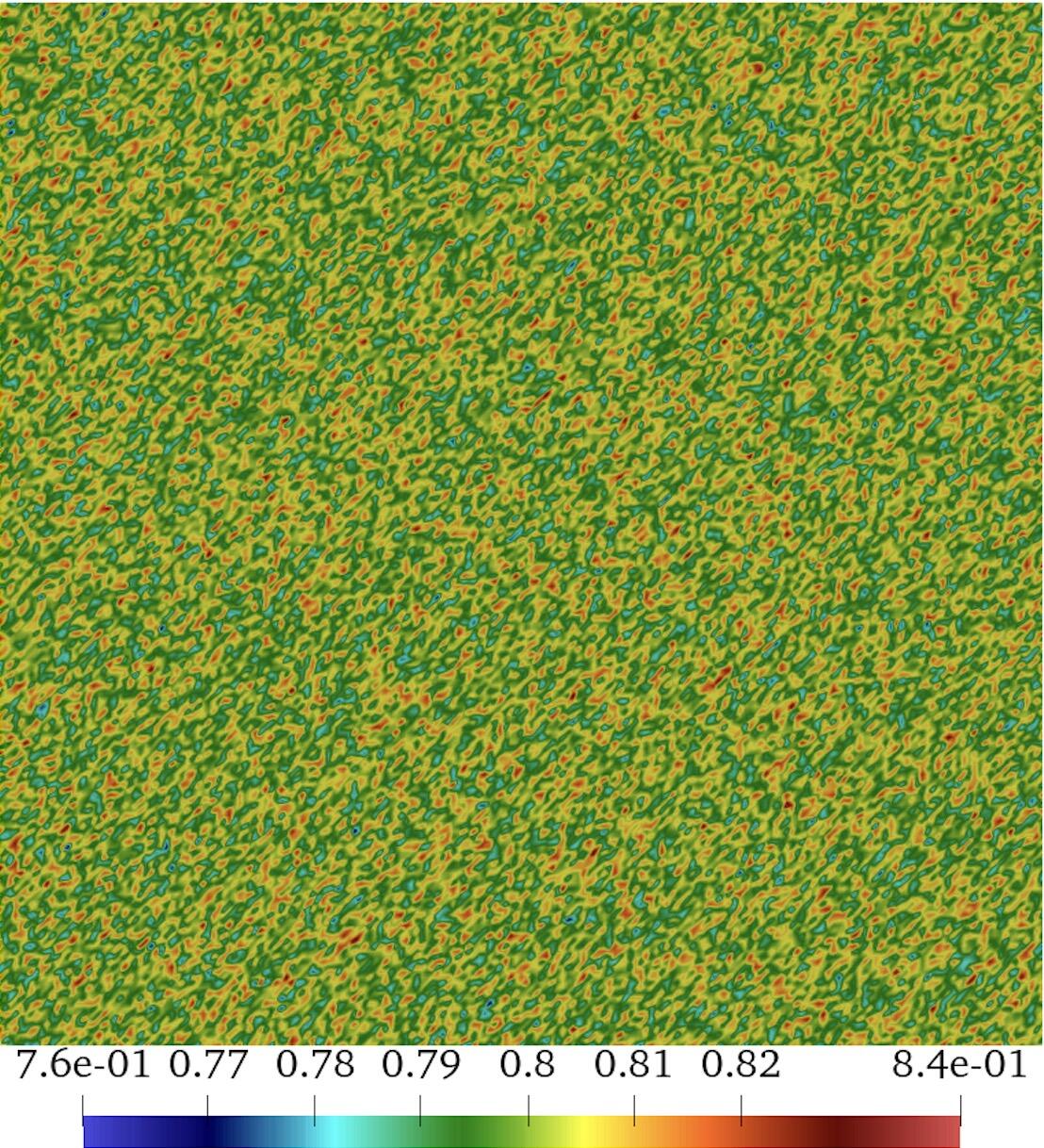}
    \put(-10,55){\footnotesize{$\rho$}}
    \end{overpic} 
     \begin{overpic}[percent,width=0.19\textwidth, grid=false]{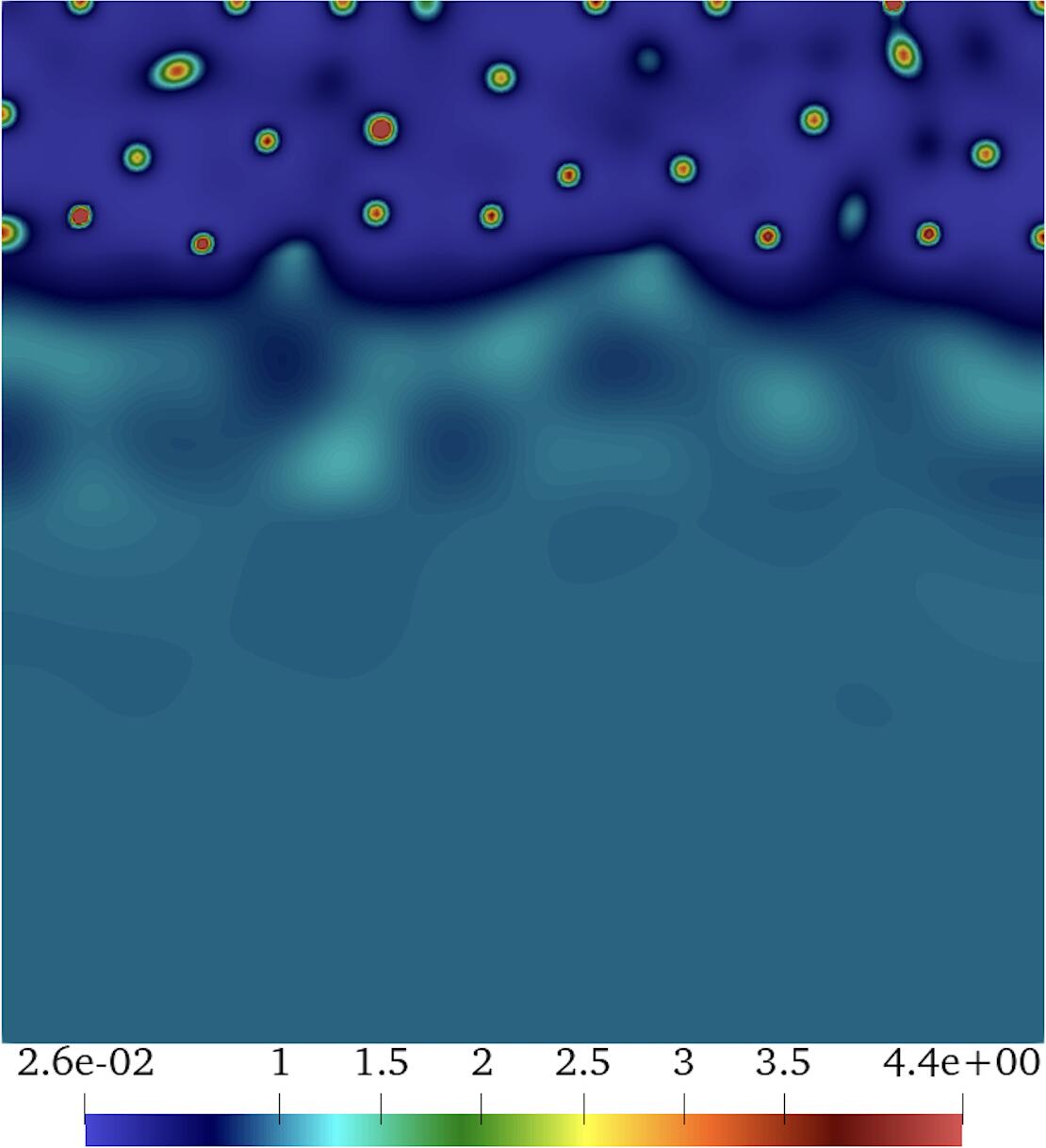}
    \end{overpic}
    \begin{overpic}[percent,width=0.19\textwidth]{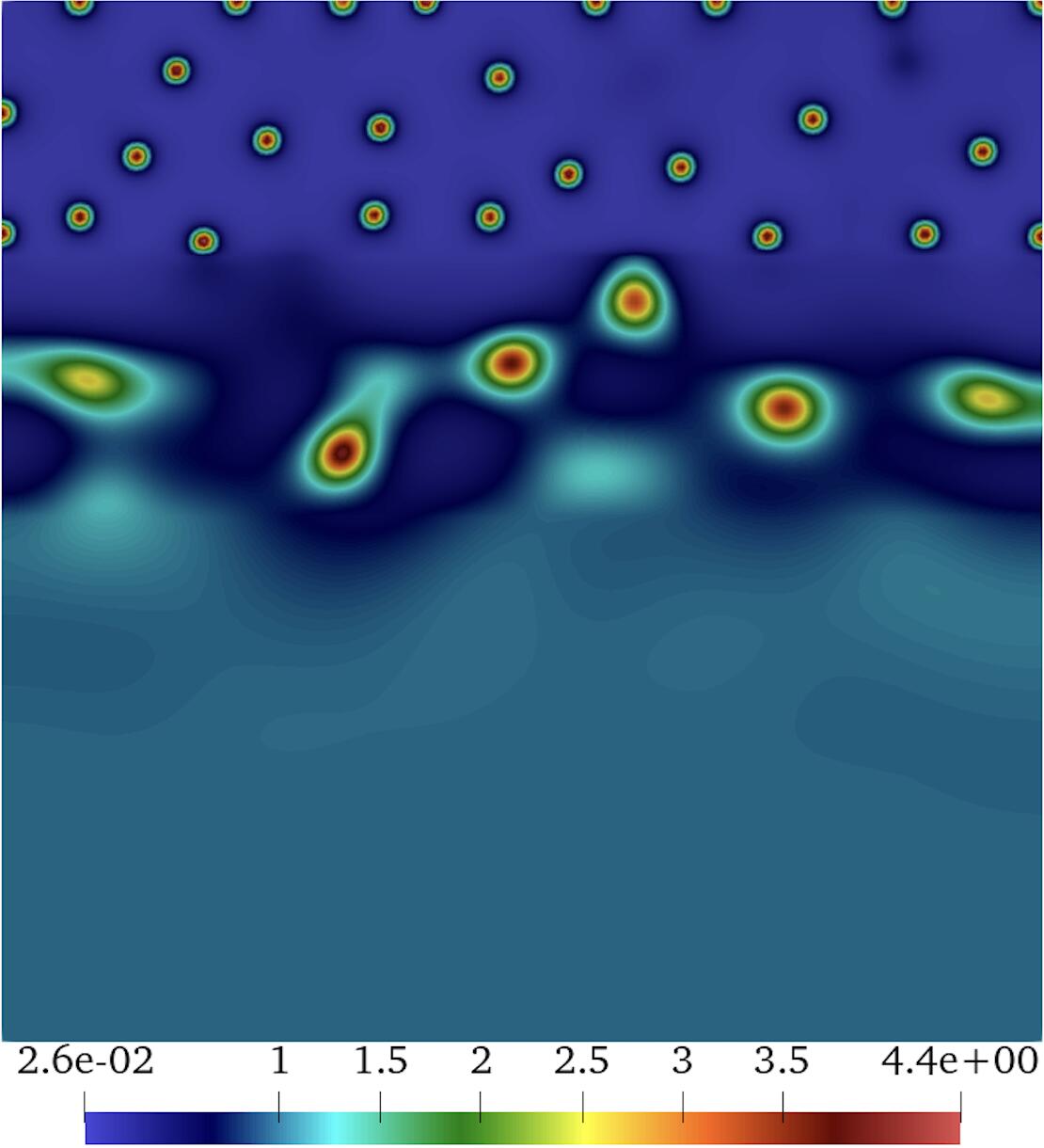}
    \end{overpic} 
        \begin{overpic}[percent,width=0.19\textwidth]{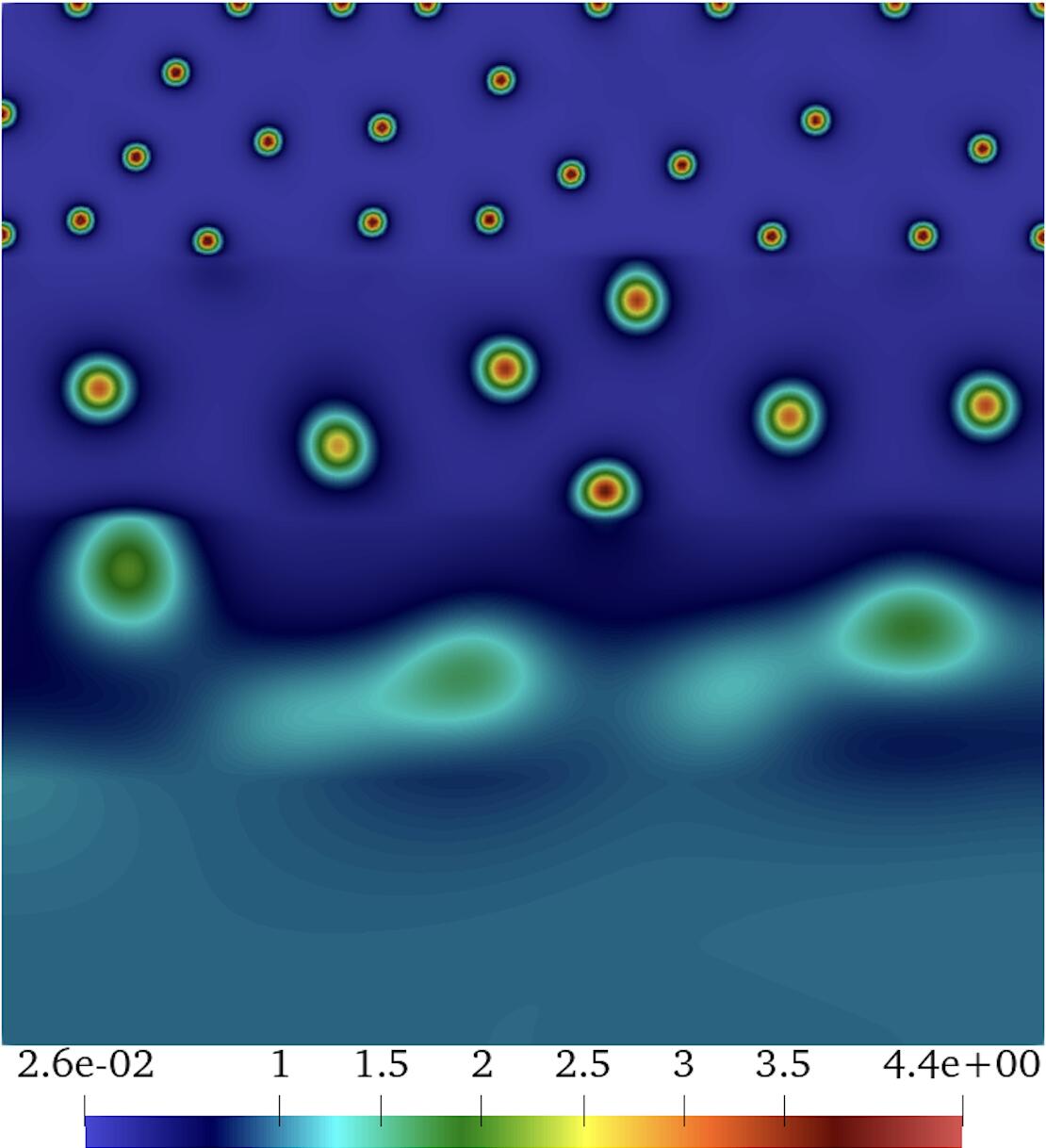}
    \end{overpic} 
    \begin{overpic}[percent,width=0.19\textwidth]{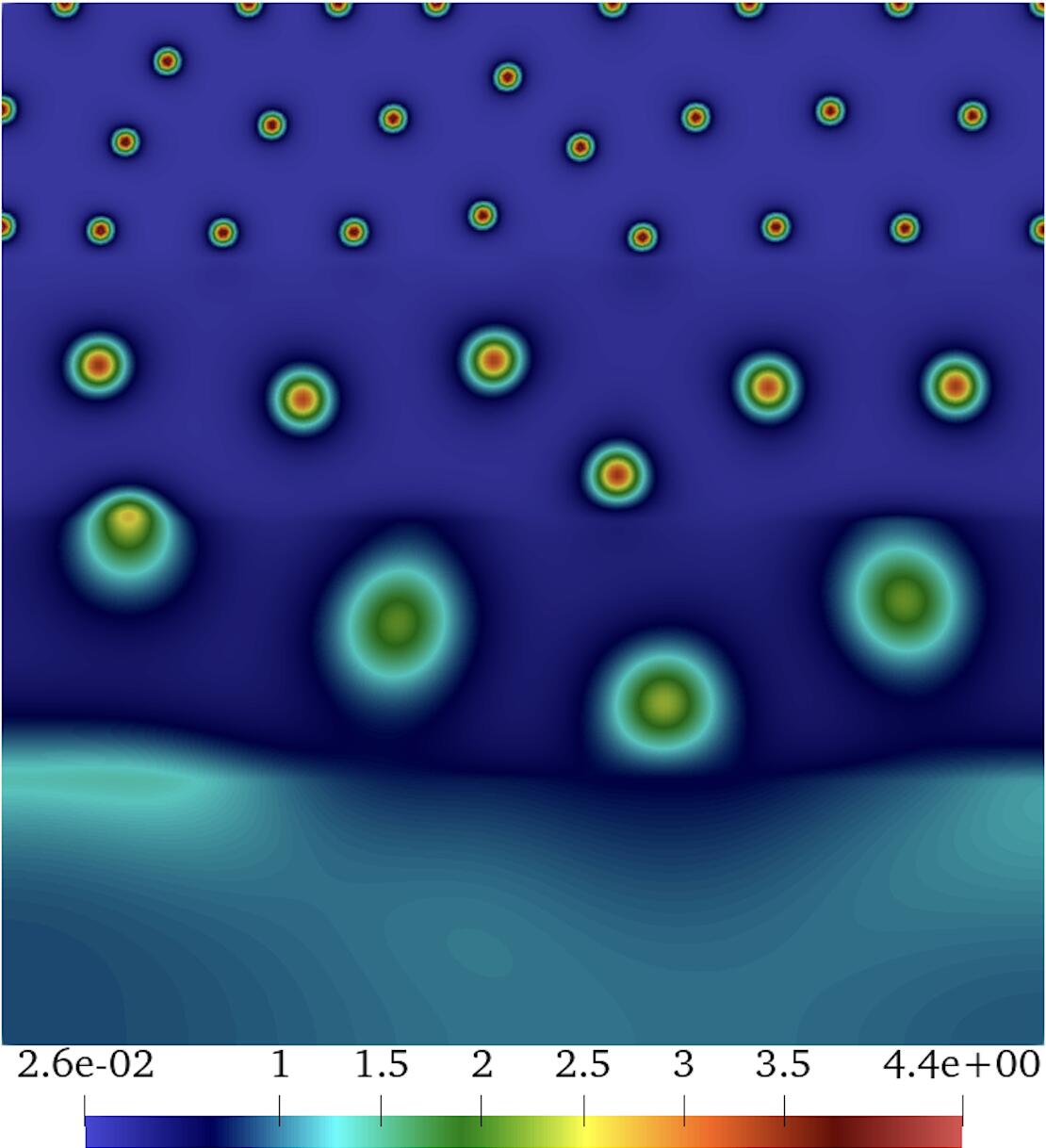}
    \end{overpic}
    \caption{Case with $\eta$ defined piecewise \eqref{eq:eta_piece}: Evolution of the level of attractiveness (top) and density (bottom).}
    \label{fig:eta-varying}
\end{figure}

Let us now define $\eta$ piecewise as follows:
\begin{equation}\label{eq:eta_piece}
    \eta = \begin{cases}
0.03 \text{ (as in case 3)}, & \text{if } 12 \le x_1 \le 16, \\
0.12, & \text{if } 8 \le x_1 < 12,\\
0.3 \text{ (as in case 2)}, & \text{if } 4 \le x_1 < 8,\\
0.9 \text{ (as in case 1)}, & \text{if } 0 \le x_1 < 4,
\end{cases}
\end{equation}
i.e., $\eta$ decreases with the vertical coordinate. We set all the other parameters
and noise as in case 2 and 3.
Fig.~\ref{fig:eta-varying} displays the evolution of $A$ and $\rho$ computed with 
mesh $h = 16/200$ and
time step $\Delta t = 1/50$. For a better visualization of the hotspots, we have manually set the maximum and minimum of $A$
and $\rho$ in the color bars.
From Fig.~\ref{fig:eta-varying}, we see that, by the end of the simulation, the solution transitions from uniform in space (with $A = \bar{A}$ and $\rho = \bar{\rho}$, roughly in the bottom quarter) to hotspots with decreasing size as $x_1$ increases. This is expected from \eqref{eq:eta_piece}. Notice also
that the smaller the hotsopts the faster they emerge.

Next, we vary the discretization 
parameters and stopping tolerances.
Tab.~\ref{tab:it} 
reports the average number of iterations over time required by the iterative partitioned algorithm to meet the stopping criterion \eqref{eq:crit} with $\texttt{tol}_1 = \texttt{tol}_2 = 10^{-3}, 10^{-6}, 10^{-9}$ for the three cases from Sec.~\ref{sec:validation} 
with different mesh sizes and time steps. We consider the combinations of $(h, \Delta t)$ used for the
results presented thus far, i.e., 
$(16/200, 1/50)$ and $(16/100, 1/25)$, 
and added one finer level, i.e., $(16/400, 1/100)$.
Additionally, we considered combination  $(16/100, 2/25)$, which keeps the same mesh
size as in one of the other combinations but doubles the time step. For case 1, the average number of iterations is around 1 for all $(h, \Delta t)$ combinations and tolerance values. 
This reflects the simplicity of the solution, which is spatially homogeneous and approaches
a steady state rather quickly.
Additionally, it indicates that even the algorithm in its non-iterative form, i.e., \eqref{eq:explicit-1}-\eqref{eq:explicit-2}, could be stable and accurate for case 1. For case 2 and 3, 
the average number of iterations is around 1
only for $\texttt{tol}_1 = \texttt{tol}_2 = 10^{-3}$. In case 2, it less than doubles each time the tolerance values are reduced by 3 orders of magnitude, while in case 3 it almost doubles. This suggests that case 3 is more challenging than case 2 for the partitioned algorithm. We speculate that the
the average number of iterations over time required by that algorithm grows with number of hotspots. 
Finally, we note that in case 2
the stopping criterion is always satisfied for every $(h, \Delta t)$ combination, while
the partitioned algorithm does not reach
convergence for case 3 when $h = 16/100$ and
$\Delta t = 2/25$ are selected. 

\begin{table}[htb!]
\centering
\renewcommand{\arraystretch}{1.2} 
\setlength{\tabcolsep}{6pt}       
\begin{tabular}{cccc ccc ccc}
\toprule
\addlinespace[0.3ex]
\multirow{2}{*}{\makecell[c]{($h$, $\Delta t$)}} 
& \multicolumn{3}{c}{\textbf{Case 1}} 
& \multicolumn{3}{c}{\textbf{Case 2}} 
& \multicolumn{3}{c}{\textbf{Case 3}} \\
\cmidrule(lr){2-4} \cmidrule(lr){5-7} \cmidrule(lr){8-10}
& $10^{-3}$ & $10^{-6}$ & $10^{-9}$ 
& $10^{-3}$ & $10^{-6}$ & $10^{-9}$ 
& $10^{-3}$ & $10^{-6}$ & $10^{-9}$ \\
\midrule
$(\frac{16}{400},\frac{1}{100})$ 
& 1.01 & 1.10 & 1.14
& 1.02 & 1.67 & 2.65
& 1.04 & 2.73 & 4.88 \\
$(\frac{16}{200}, \frac{1}{50})$ 
& 1.01 & 1.11 & 1.26 
& 1.04 & 1.89 & 2.96 
& 1.08 & 3.06 & 5.93 \\
$(\frac{16}{100}, \frac{1}{25})$ 
& 1.02 & 1.13 & 1.31
& 1.07 & 2.11 & 3.39
& 1.17 & 3.16 & 6.04 \\
$(\frac{16}{100}, \frac{2}{25})$ 
& 1.03 & 1.15 & 1.38 
& 1.11 & 2.34 & 4.17 
& - & - & - \\
\bottomrule
\end{tabular}
\caption{Average number of iterations over time required by the iterative partitioned algorithm to meet the stopping crieterion \eqref{eq:crit} with $\texttt{tol}_1 = \texttt{tol}_2 = 10^{-3}, 10^{-6}, 10^{-9}$ for the three cases from Sec.~\ref{sec:validation} 
with different values of mesh size $h$ and time step $\Delta t$.}
\label{tab:it}
\end{table}

As mentioned in Sec.~\ref{sec:num_met}, the iterative
partitioned algorithm \eqref{eq:iterative-1}-\eqref{eq:crit}
is a fixed point scheme that, at time $t^{n+1}$, 
starts the iterations from the solution at time $t^{n}$. Hence, 
if $\Delta t$ is too large, the solution at time 
$t^{n}$ is not a good guess for the solution at time $t^{n+1}$ and the algorithm fails to converge. Tab.~\ref{tab:it} shows that if one keeps the same mesh ($h = 16/100$) and halves the time step to 
$\Delta t = 1/25$, then algorithm \eqref{eq:iterative-1}-\eqref{eq:crit} converges even in the more challenging case. We conclude that 
Tab.~\ref{tab:it} indicates algorithm \eqref{eq:iterative-1}-\eqref{eq:crit} is robust 
to (reasonable) variations of time step, mesh size, stopping tolerance and parameter
$\eta$. The low numbers of iterations in 
Tab.~\ref{tab:it} also indicate numerical efficiency, since only two linear systems have to be solved per iteration.

\begin{figure}[htb!]
    \centering
    \includegraphics[width=0.7\linewidth]{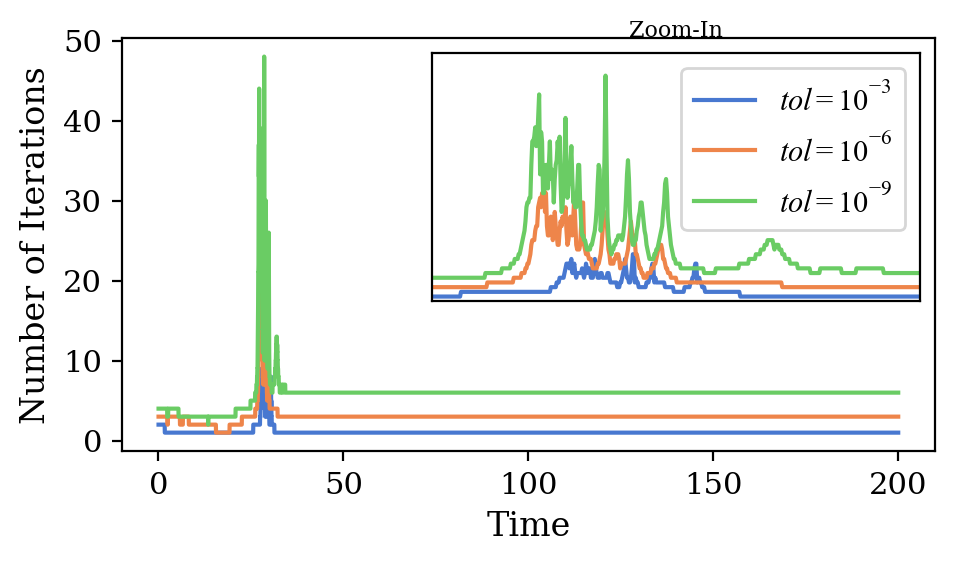}
    \caption{Case 3 without noise:
    number of iterations needed for the
    iterative partitioned algorithm \eqref{eq:iterative-1}-\eqref{eq:crit} to satisfy the stopping
    criterion with $\texttt{tol}_1 = \texttt{tol}_2 = 10^{-3}, 10^{-6}, 10^{-9}$ over time. The zoomed-in view
    highlights the time interval $t \in [25, 35]$, where the peaks in
    iteration numbers occur. }
    \label{fig:iters}
\end{figure}

To complement the information reported in 
Tab.~\ref{tab:it}, Fig.~\ref{fig:iters}
shows the number of iterations needed for the iterative partitioned algorithm \eqref{eq:iterative-1}-\eqref{eq:crit} to satisfy the stopping criterion with $\texttt{tol}_1 = \texttt{tol}_2 = 10^{-3}, 10^{-6}, 10^{-9}$ over time
for case 3 with no noise (see results in 
Fig.~\ref{fig:case3-nonoise}).
We see that for most of the time interval under consideration, the numbers of iterations are close to the averages in 
Fig.~\ref{fig:iters}. However, sharp peaks occur when the hotspost pattern
is forming. In case 3 with no noise, 
that happens for $t \in [25, 35]$, as confirmed by Fig.~\ref{fig:case3-nonoise}.

\subsection{Results for spatially heterogeneous parameters and a nontrivial geometry}\label{sec:space_dep}

All the results presented so far are associated
to constant model parameters. 
When that is the case, 
even with the addition of noise, the hotspot distributions are fairly regular. 
See Fig.~\ref{fig:case2}
and \ref{fig:case3}.
The hotspots observed from visualization of real data (see Fig.~\ref{fig:Houston}) lack such regularity. From  
Fig.~\ref{fig:Houston}, it seems that persistent hotspots emerge along highways. This is corroborated by the analysis
in \cite{IZAreport} based on a
dataset spanning all US counties between 1960-1993. It is found that a highway opening led to a 5\% rise in the local index crime, driven by property crime. Below, we present cases where persistent hotspots form along hypothetical highways in the simple square geometry considered so far and in a more realistic geometry.  

First, we stick to the domain $\Omega$ considered so far. The 
value of $\eta$ and the initial density ($\rho_0$ and associated noise) are the same  as in case 2.
The presence of a highway can be modeled through static component $A^{st}$ and
``burglar generation'' $\Gamma$ (hence, parameter combination $\Gamma \theta /\omega^2 = B_0$) heterogeneous in space $\bm{x} = (x_0,x_1)$. 
For example, we define them as: 
\begin{align}
 & A^{st}(\bm{x}) = \frac{1}{30} \cdot 
 \exp\{-20*(x_0+x_1-l)^2\} \cdot \left( x_0 > \frac{l}{2} ~ \& ~ x_1 < \frac{l}{2} \right)
 + \frac{1}{30}\cdot \chi_{1} (\bm{x}), \cl
 & \frac{\Gamma(\bm{x}) \theta}{\omega^2} = B_0(\bm{x}) = 
 \exp\{-20*(x_0+x_1-l)^2\} \cdot \left( x_0 > \frac{l}{2} ~ \& ~ x_1 < \frac{l}{2} \right)
 + 1 \cdot \chi_{0.1} (\bm{x}), \el
\end{align}
where $l = 16$ is the size of the computational domain. Note that the above $A^{st}$ and $B_0$
take the values from case 2, introduce randomness through the $\chi_\delta$, and 
add the highway, which
is only in the bottom right quadrant, i.e., where $x_0 > l/2$ and $x_1 < l/2$. For initial condition \eqref{eq:ic1}, we take
$\delta = 0.1$ and $\xi_B(\bm{x})$ as in case 2.

\begin{figure}[htb!]
     \centering
         \begin{overpic}[percent,width=0.19\textwidth]{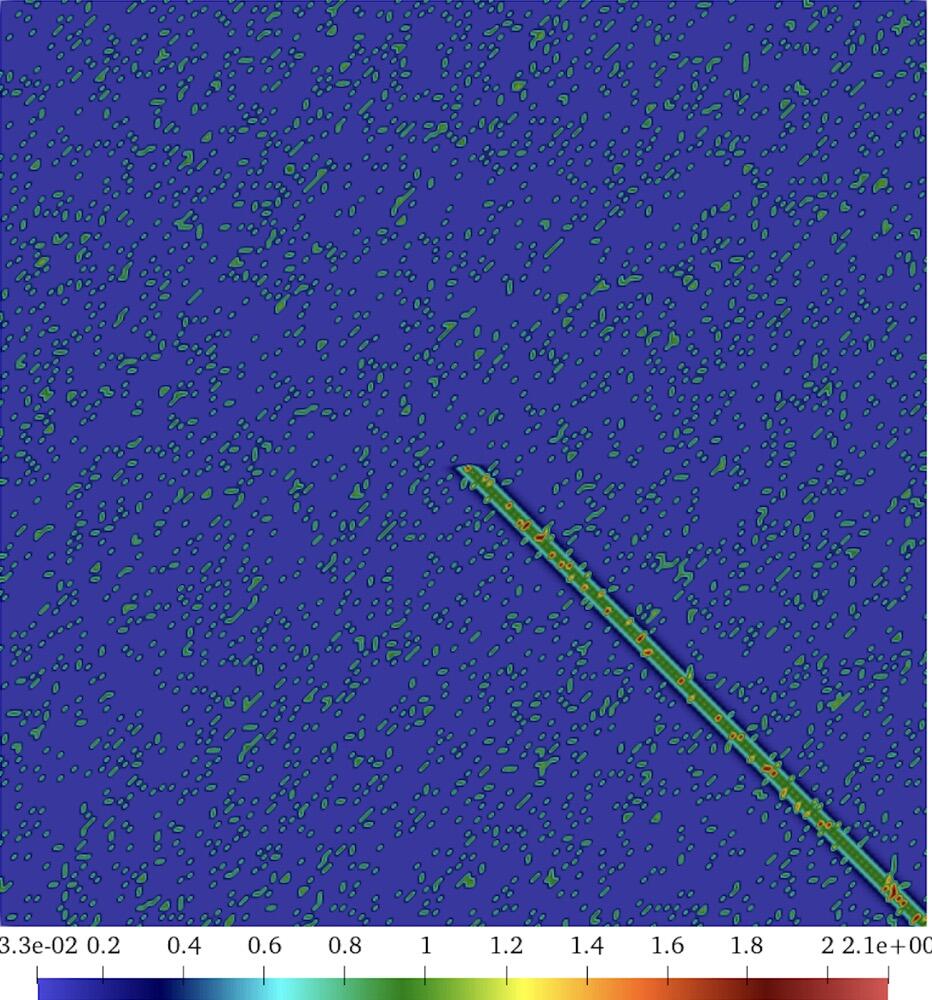}
    \put(30,103){\footnotesize{$t = 0$}}
    \put(-10,55){\footnotesize{$A$}}
    \end{overpic} 
     \begin{overpic}[percent,width=0.19\textwidth, grid=false]{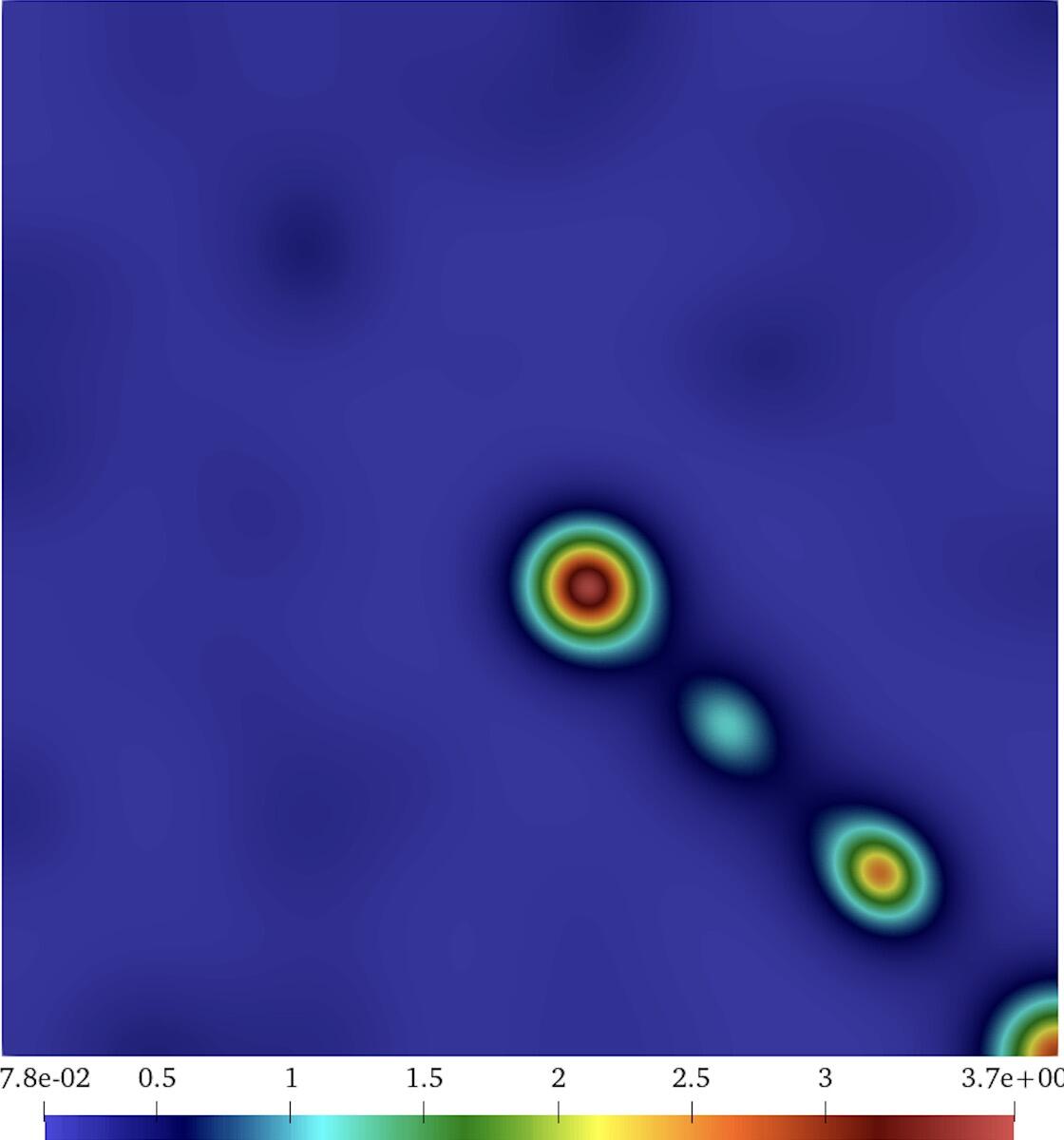}
    \put(35,103){\footnotesize{$t = 3$}}
    \end{overpic}
    \begin{overpic}[percent,width=0.19\textwidth]{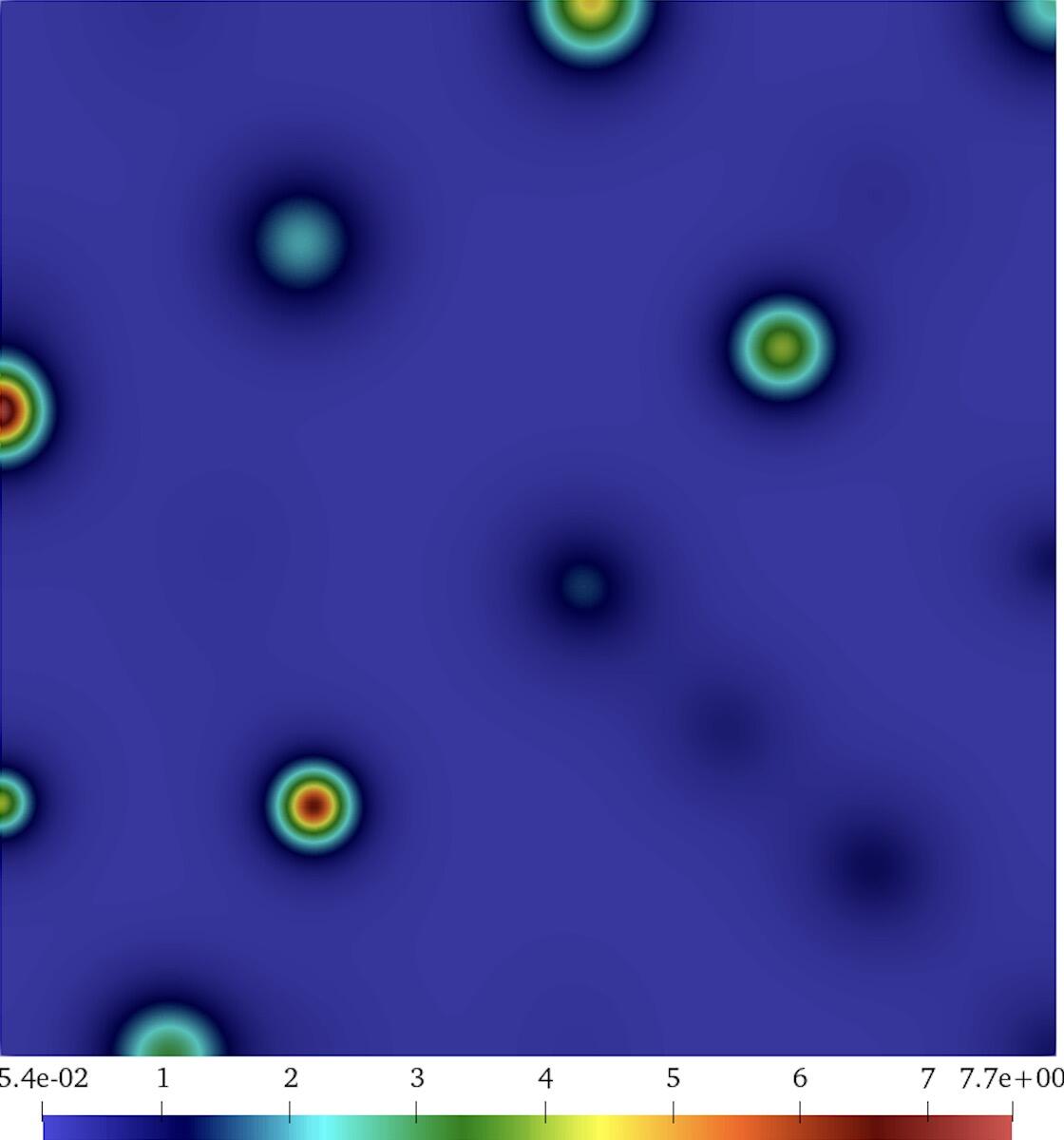}
    \put(35,103){\footnotesize{$t = 6$}}
    \end{overpic} 
        \begin{overpic}[percent,width=0.19\textwidth]{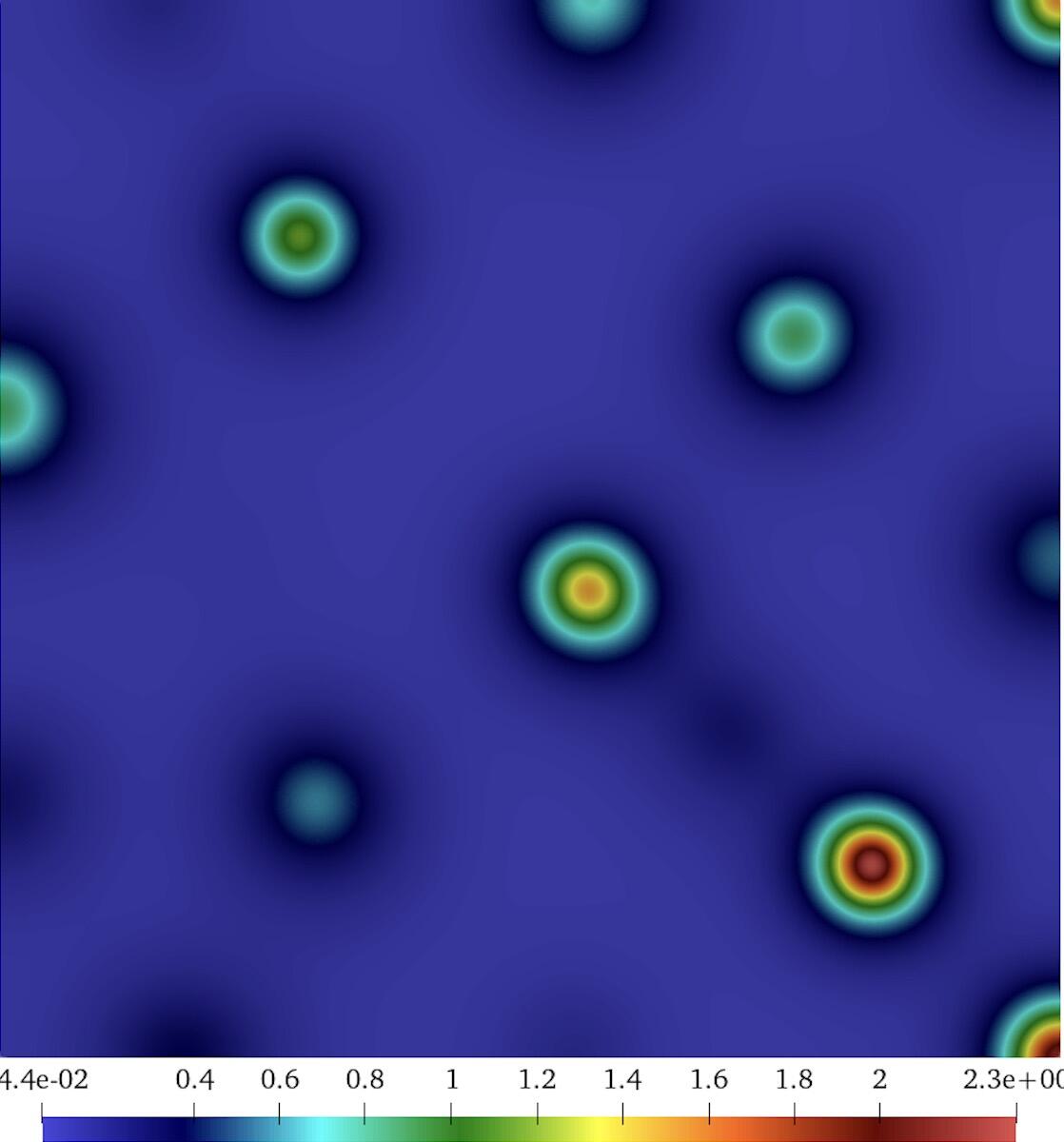}
    \put(32,103){\footnotesize{$t = 9$}}
    \end{overpic} 
    \begin{overpic}[percent,width=0.19\textwidth]{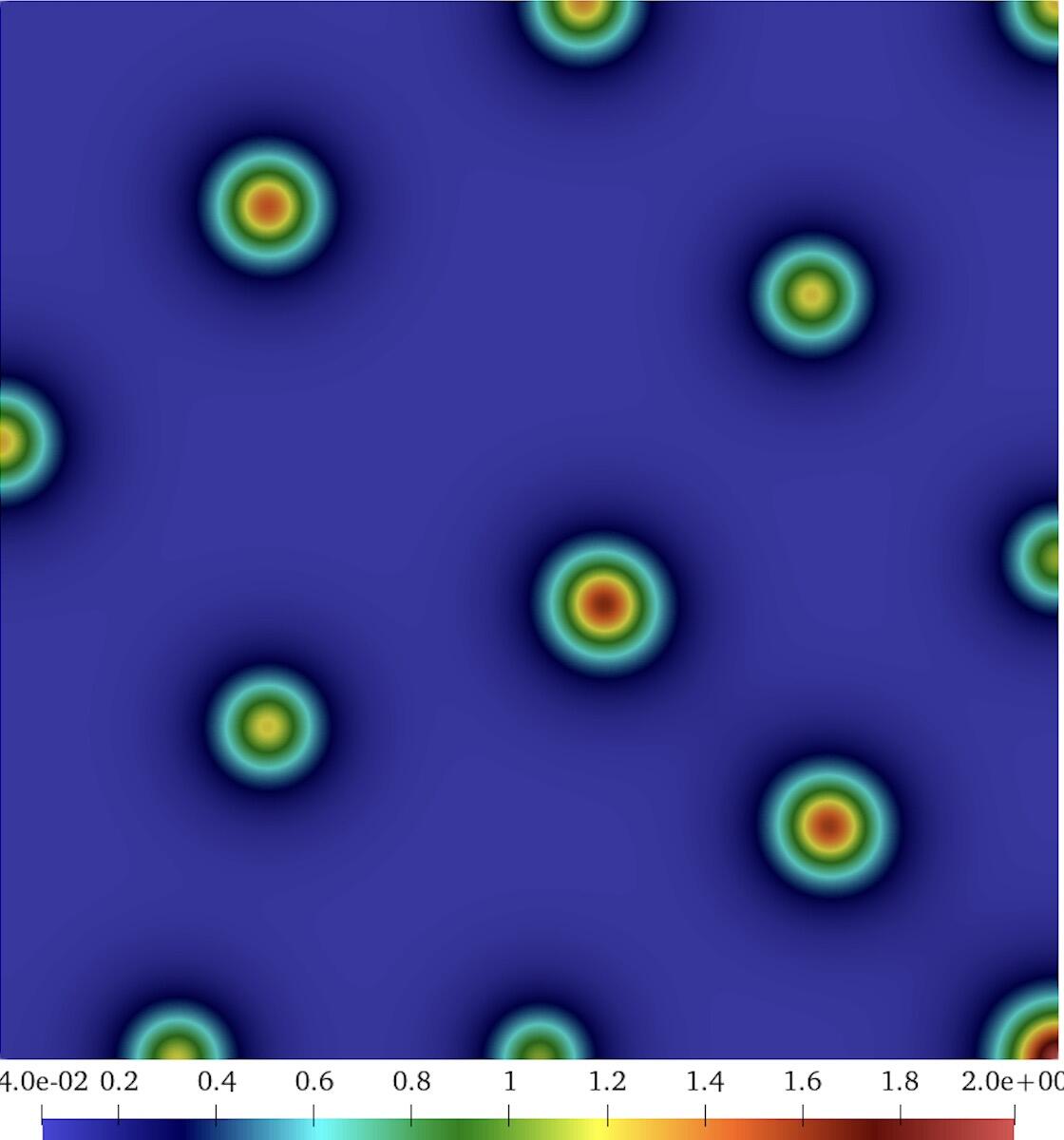}
    \put(30,103){\footnotesize{$t = 200$}}
    \end{overpic}
    \\
         \begin{overpic}[percent,width=0.19\textwidth]{Figs_ContinummModel/S2_rho_T0.jpg}
    \put(-10,55){\footnotesize{$\rho$}}
    \end{overpic} 
     \begin{overpic}[percent,width=0.19\textwidth, grid=false]{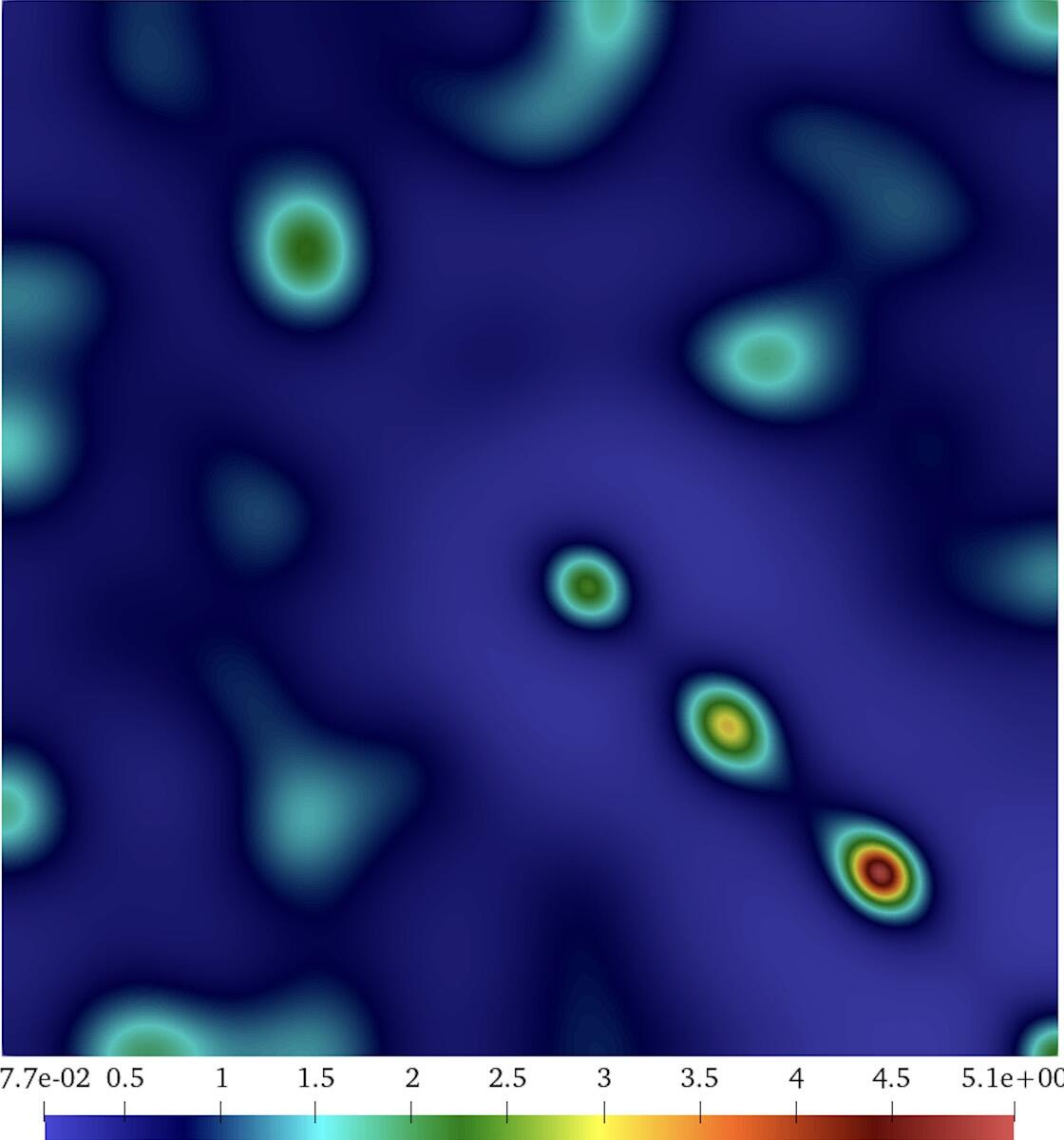}
    \end{overpic}
    \begin{overpic}[percent,width=0.19\textwidth]{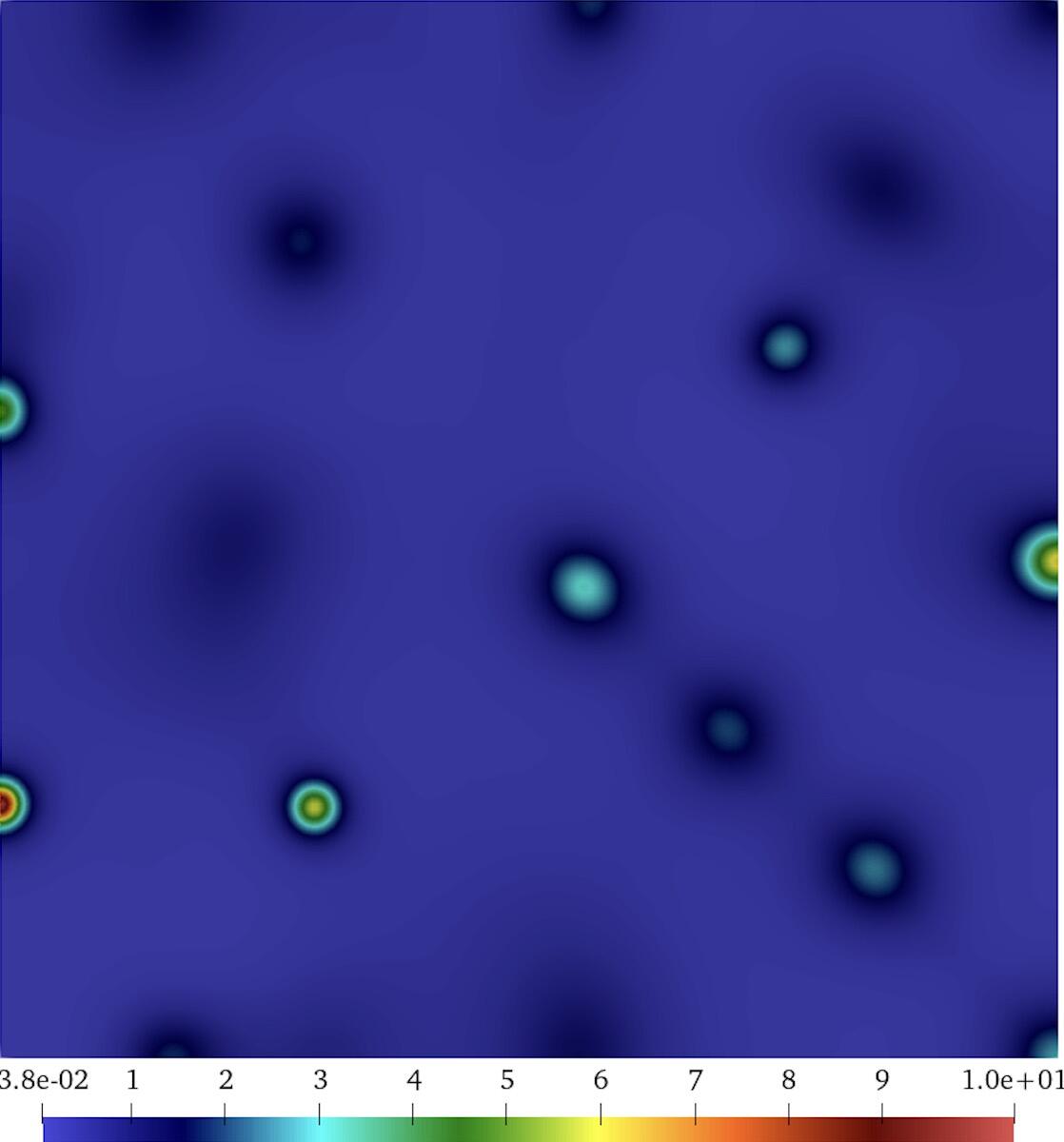}
    \end{overpic} 
        \begin{overpic}[percent,width=0.19\textwidth]{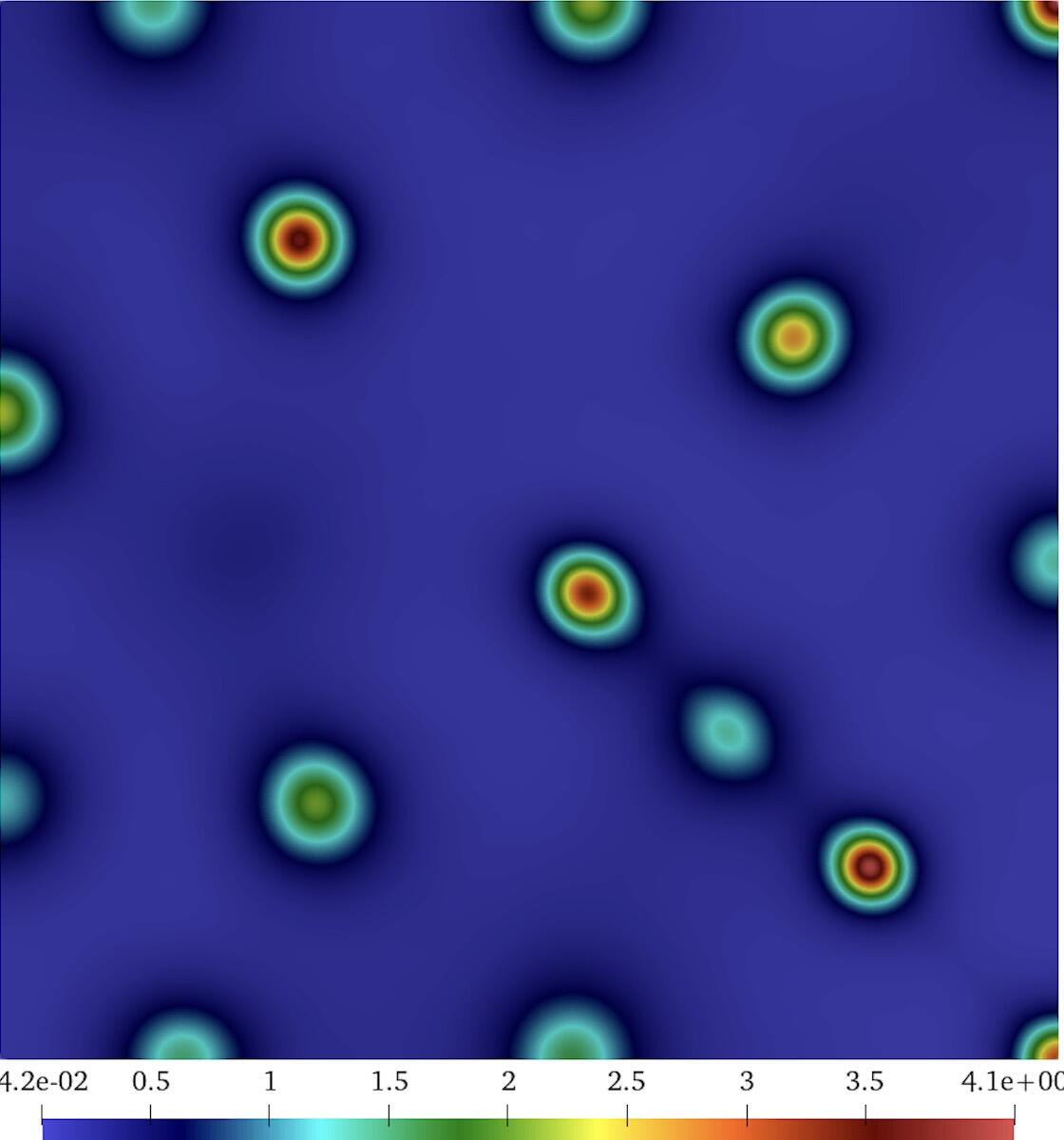}
    \end{overpic} 
    \begin{overpic}[percent,width=0.19\textwidth]{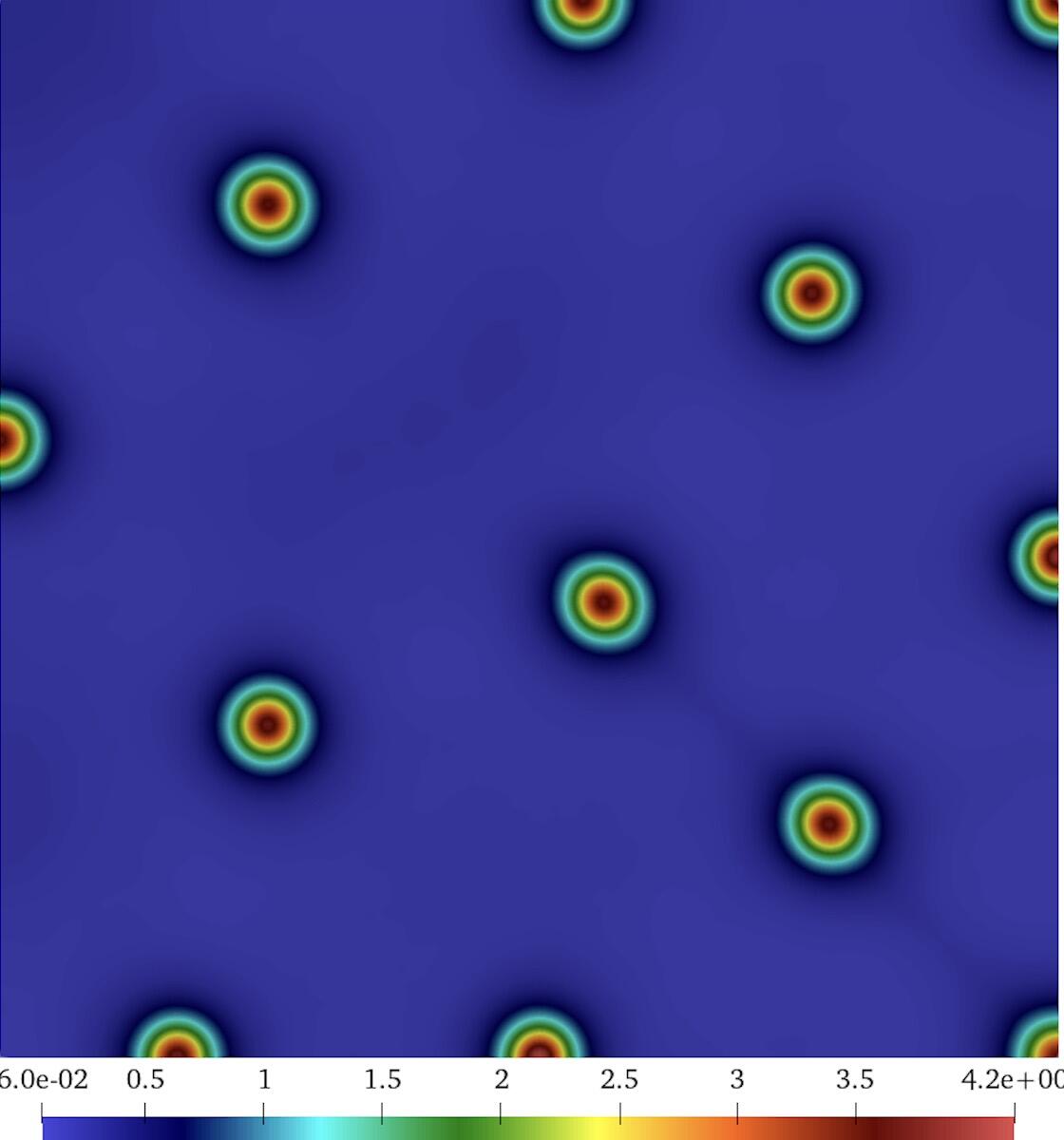}
    \end{overpic}
    \caption{Simple highway case: Evolution of the level of attractiveness (top) and density (bottom).}
    \label{fig:case_hw}
\end{figure}

We consider again mesh size $h = 16/200$ and time step to $\Delta t = 1/50$. 
Fig.~\ref{fig:case_hw}
shows the evolution of $A$ and $\rho$
computed by the iterative partitioned algorithm \eqref{eq:iterative-1}-\eqref{eq:crit}, with stopping
tolerances $\texttt{tol}_1 = \texttt{tol}_2 = 10^{-9}.$ We observe that by $t = 9$ large hotspots
have emerged along the highway, and smaller transient hotspots form in directions parallel to the highway. These hotspots remain transient for the entire duration of the time interval under consideration, while those along the highway persist.

\begin{figure}[htb!]
     \centering
         \begin{overpic}[percent,width=0.19\textwidth]{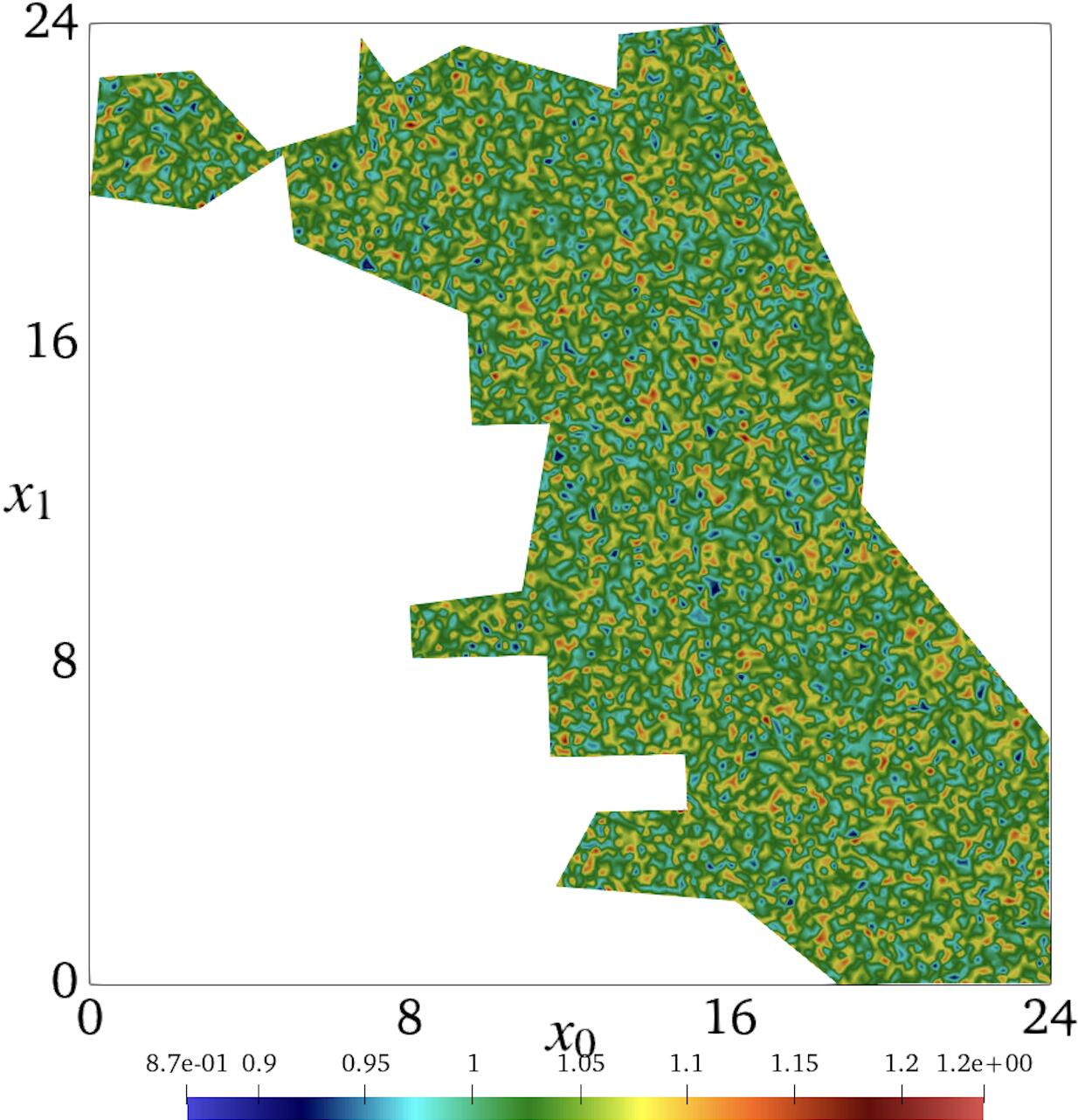}
    \put(30,103){\footnotesize{$t = 0$}}
    \put(12,35){\footnotesize{$A$}}
    \end{overpic} 
     \begin{overpic}[percent,width=0.19\textwidth, grid=false]{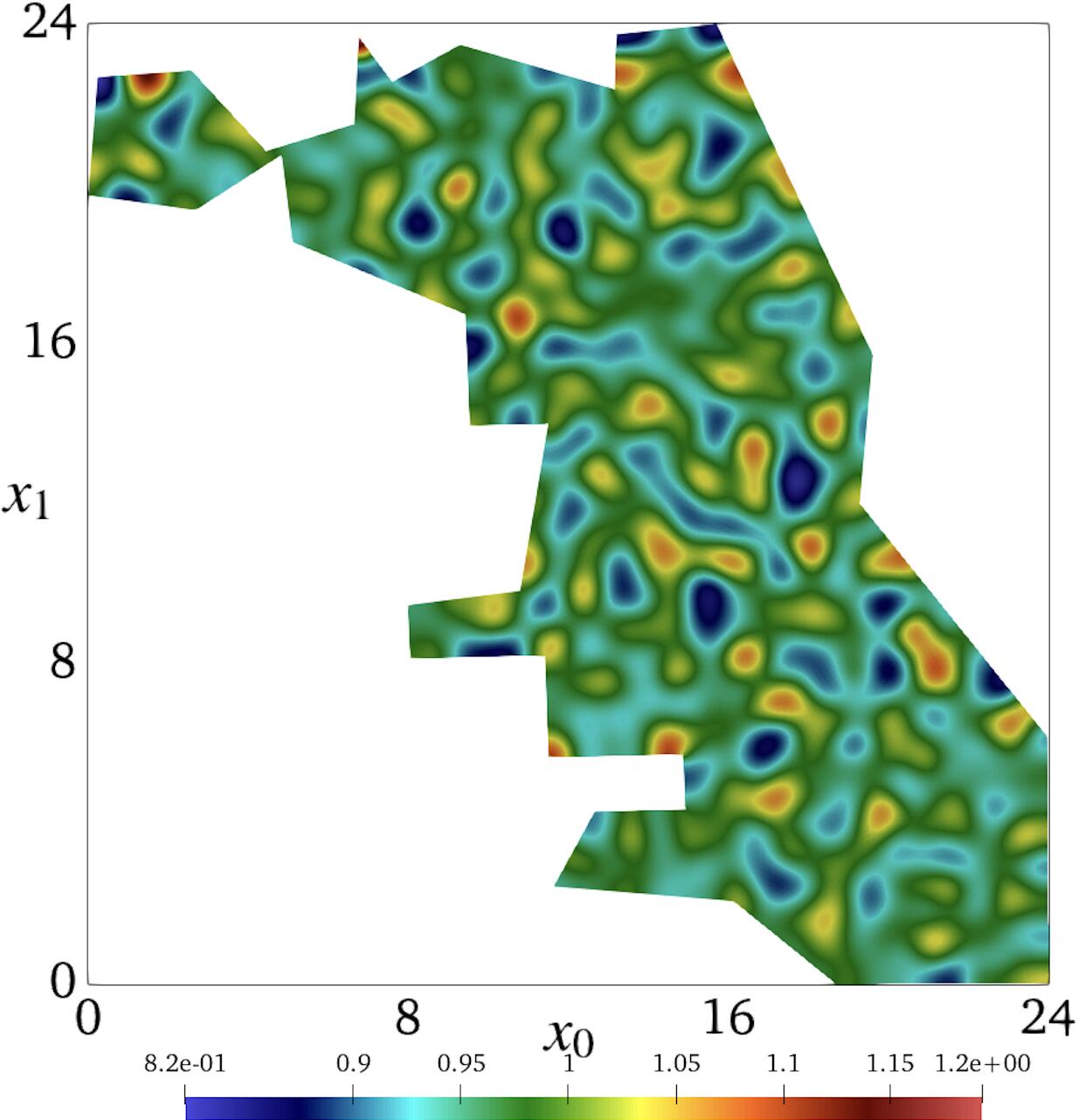}
    \put(35,103){\footnotesize{$t = 2$}}
    \end{overpic}
    \begin{overpic}[percent,width=0.19\textwidth]{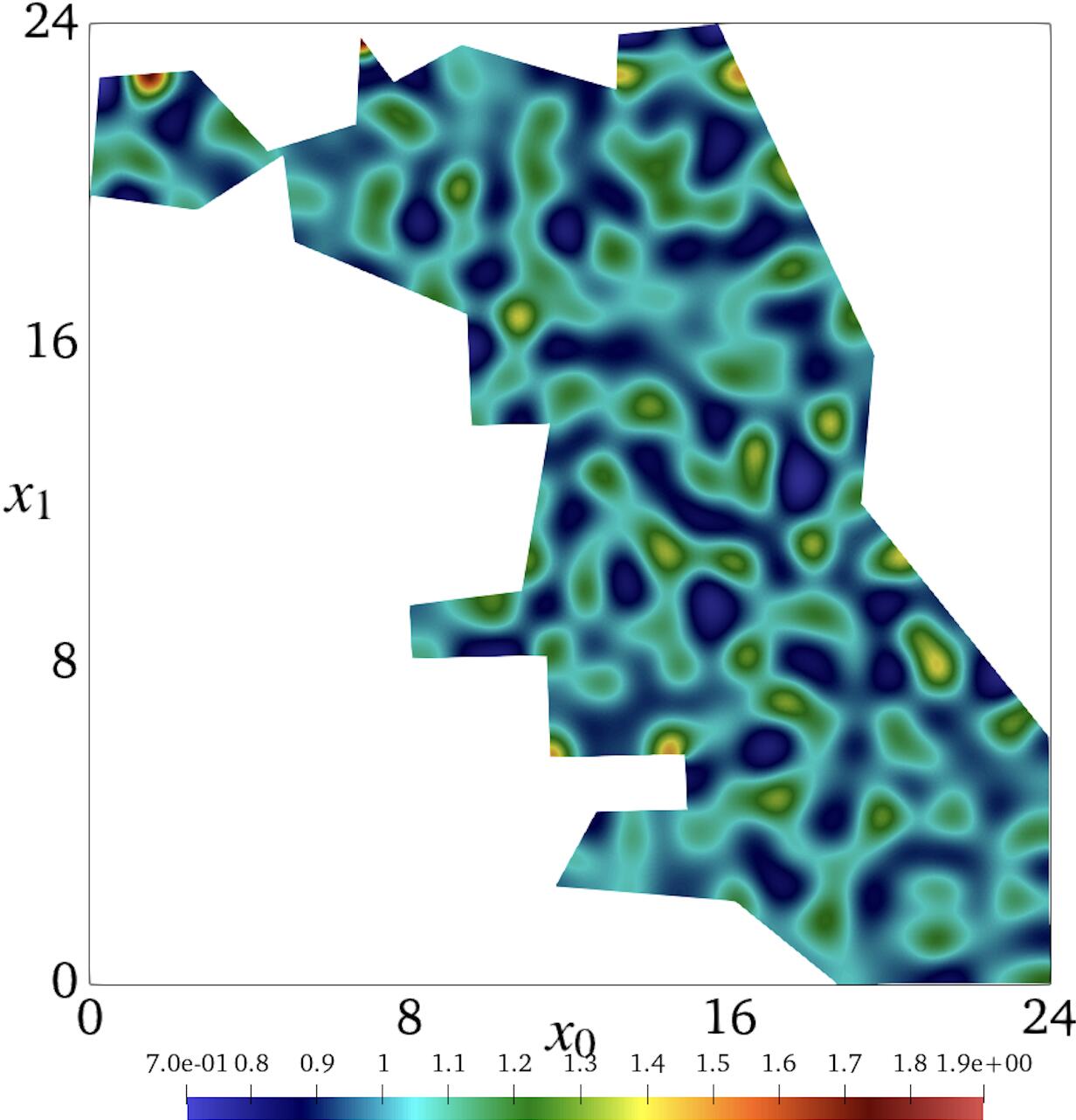}
    \put(35,103){\footnotesize{$t = 3$}}
    \end{overpic} 
        \begin{overpic}[percent,width=0.19\textwidth]{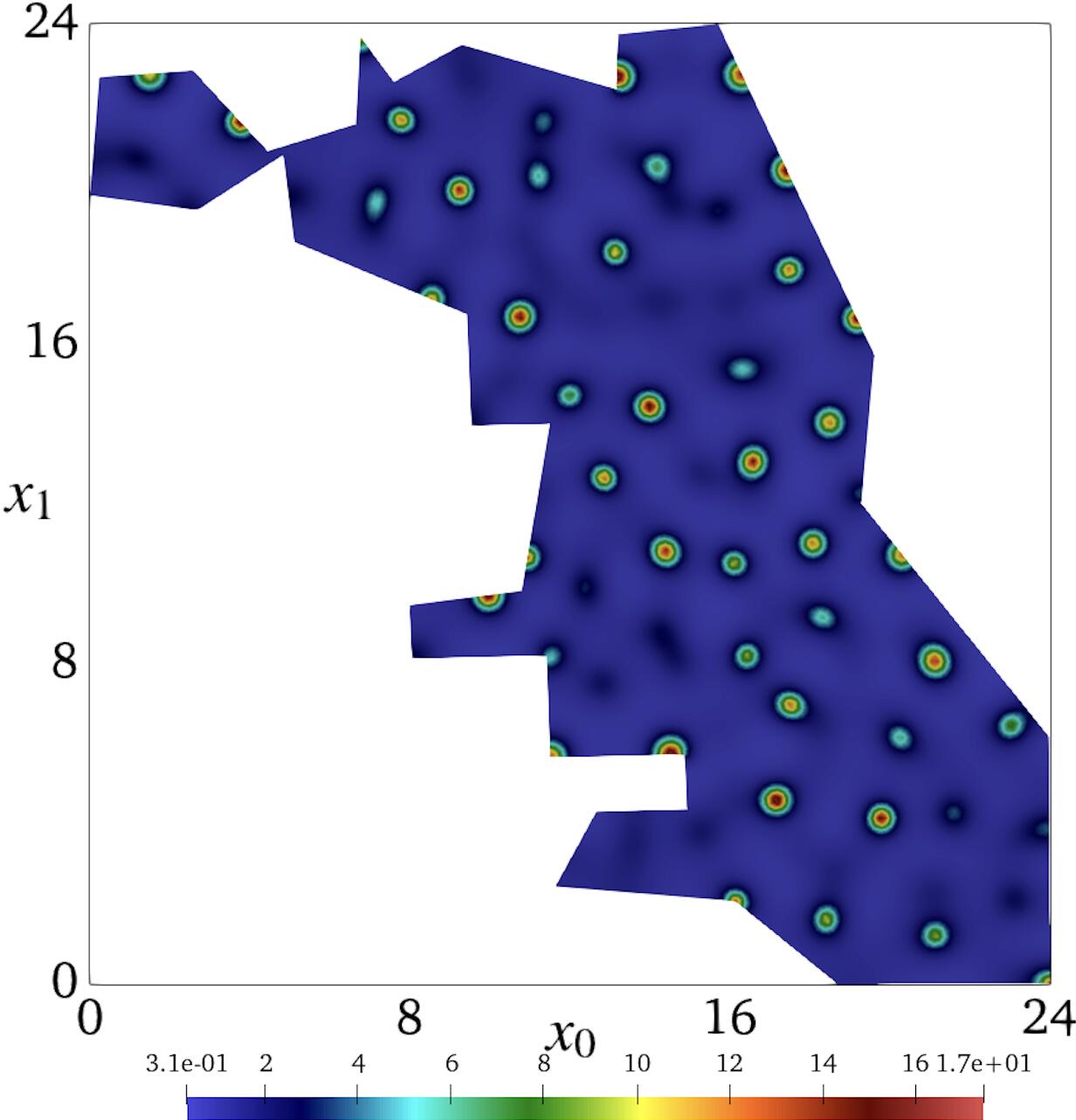}
    \put(32,103){\footnotesize{$t = 5$}}
    \end{overpic} 
    \begin{overpic}[percent,width=0.19\textwidth]{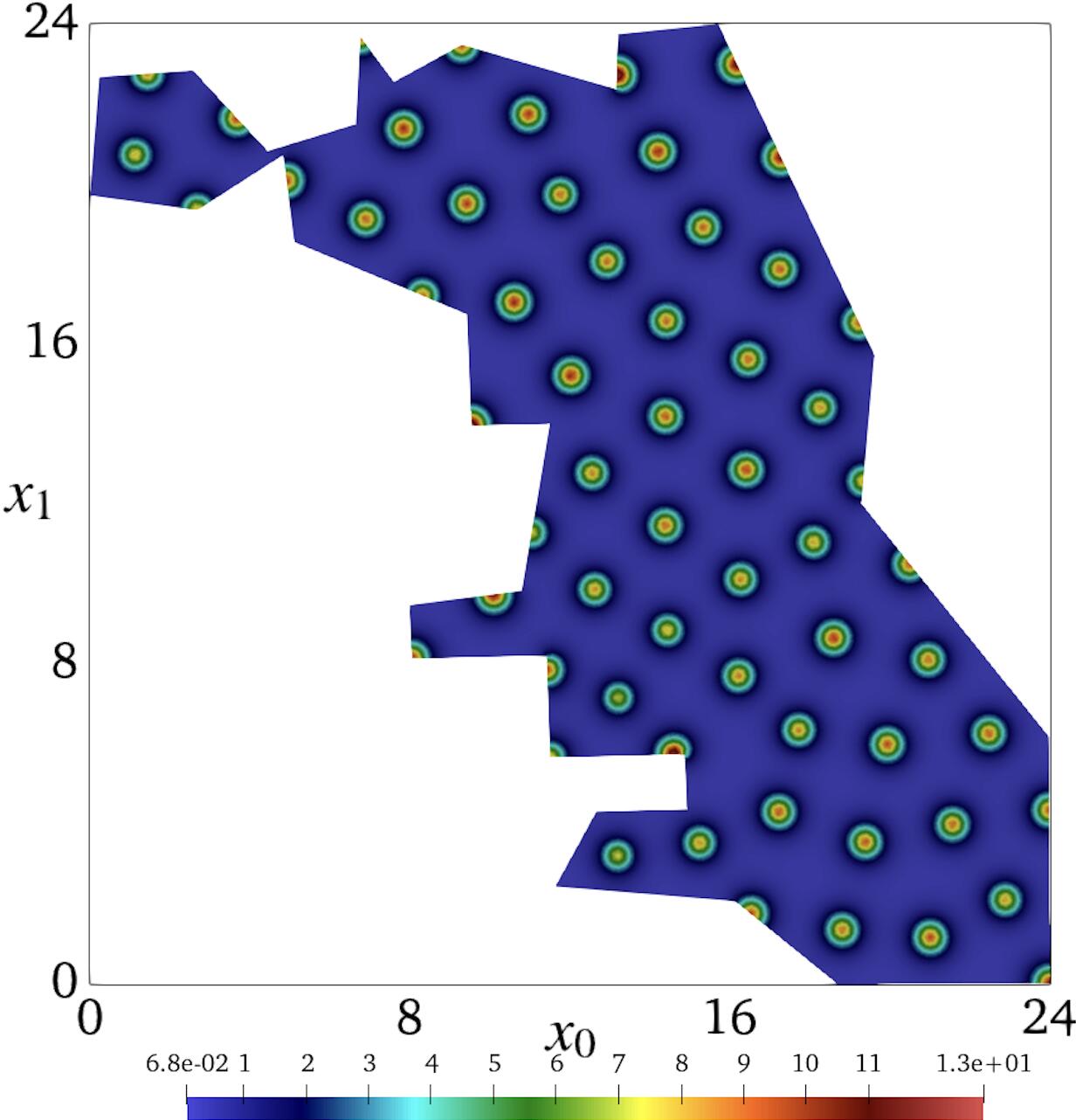}
    \put(30,103){\footnotesize{$t = 200$}}
    \end{overpic}
    \\
         \begin{overpic}[percent,width=0.19\textwidth]{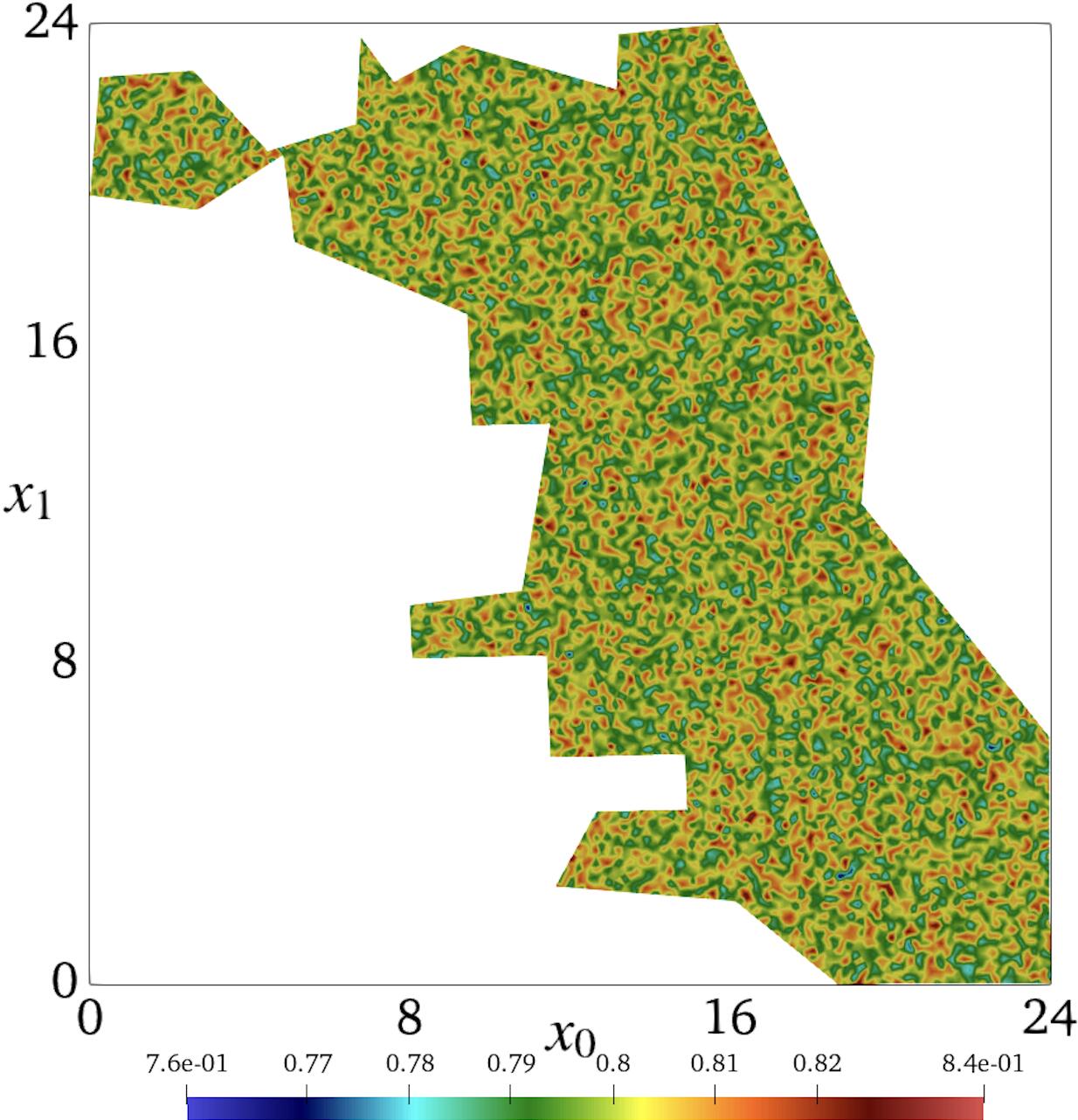}
    \put(12,35){\footnotesize{$\rho$}}
    \end{overpic} 
     \begin{overpic}[percent,width=0.19\textwidth, grid=false]{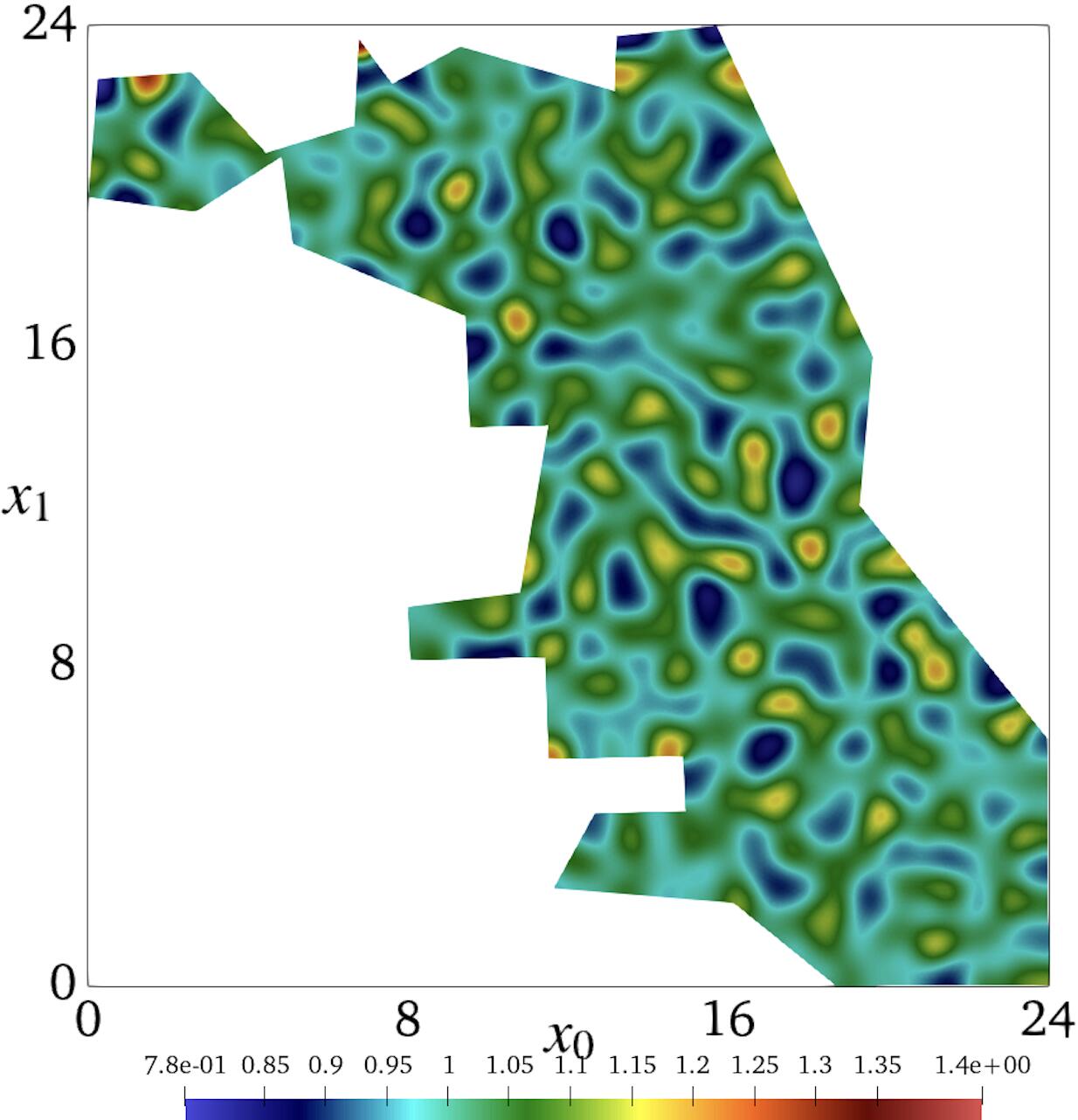}
    \end{overpic}
    \begin{overpic}[percent,width=0.19\textwidth]{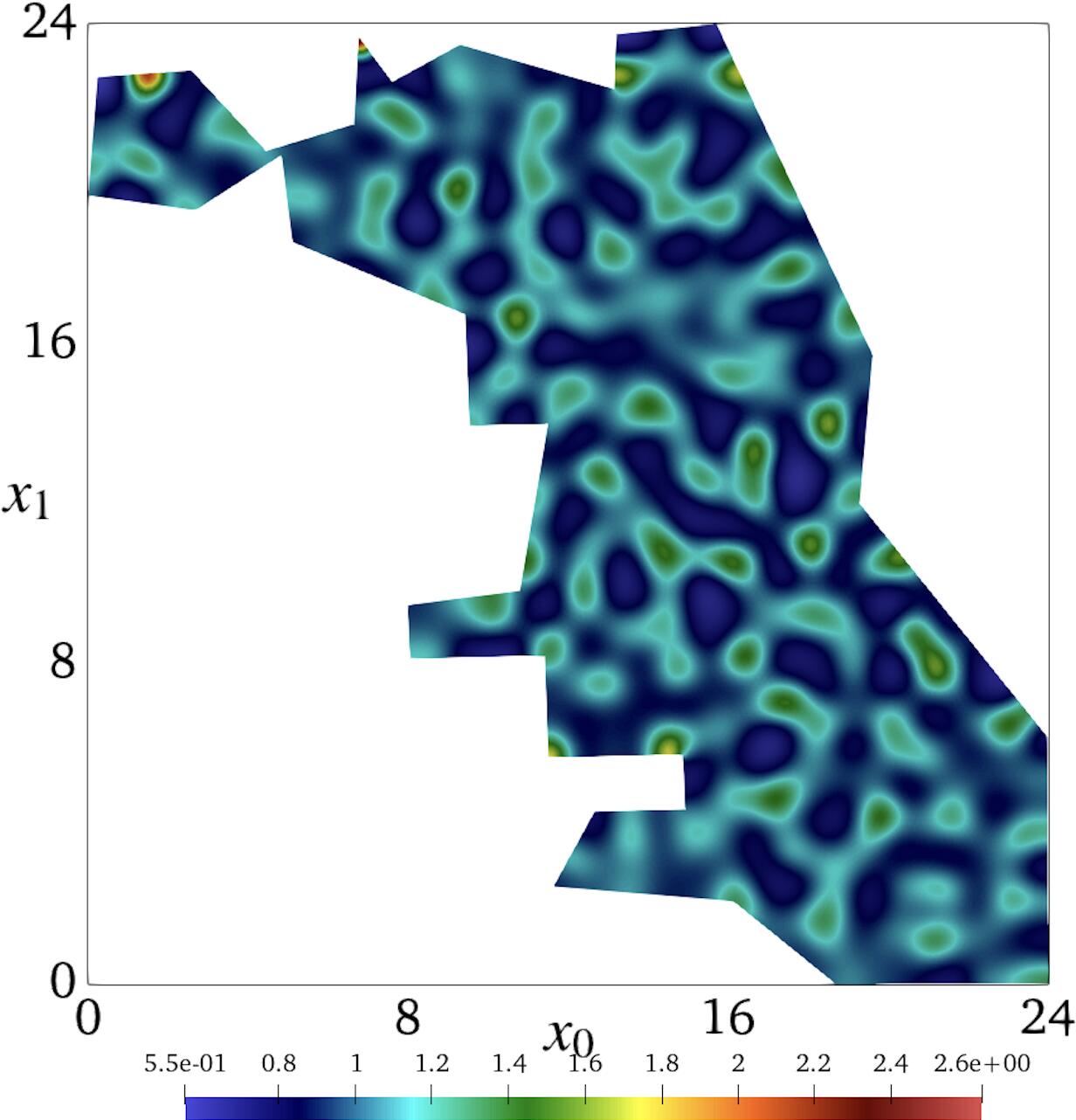}
    \end{overpic} 
        \begin{overpic}[percent,width=0.19\textwidth]{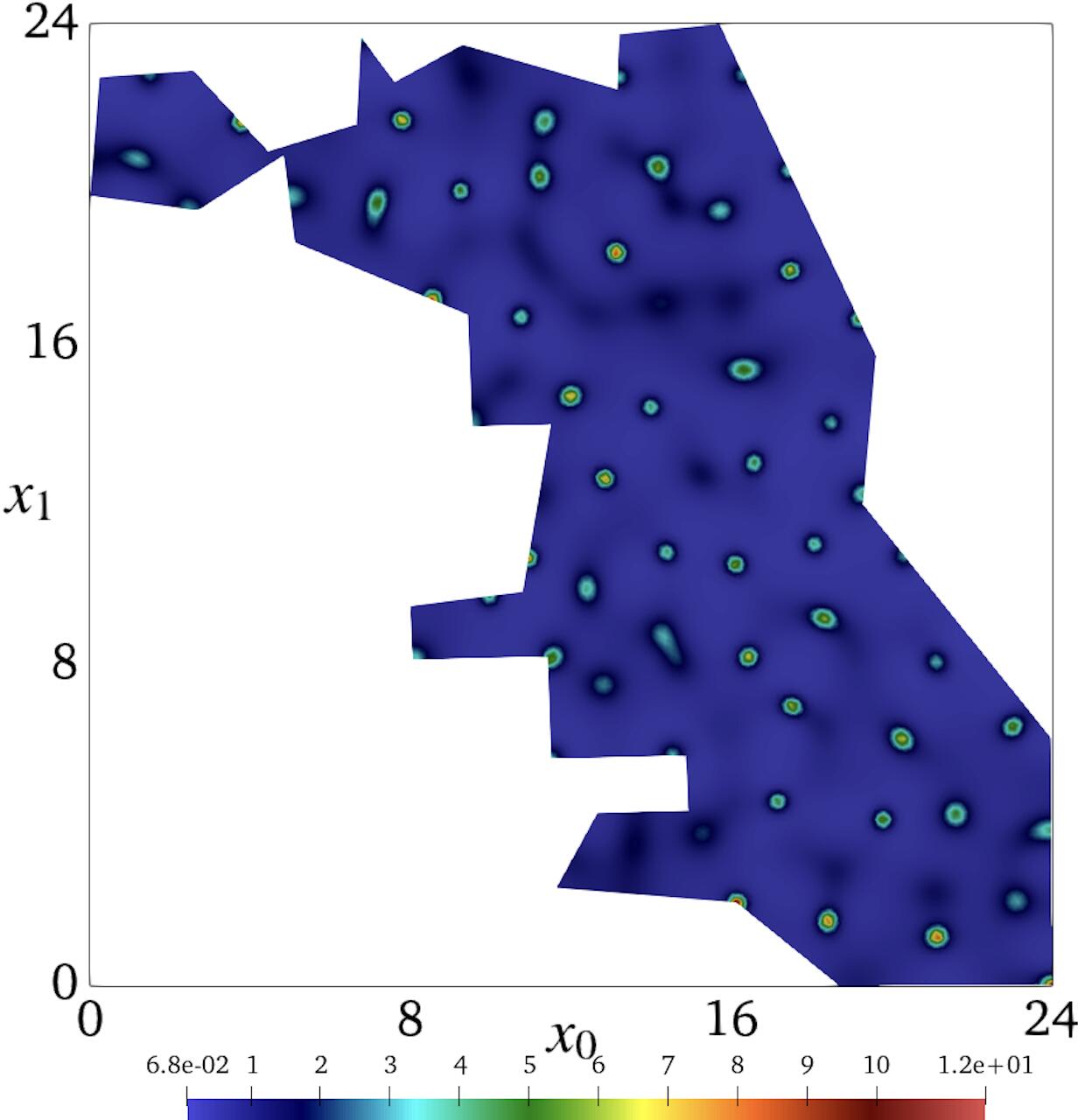}
    \end{overpic} 
    \begin{overpic}[percent,width=0.19\textwidth]{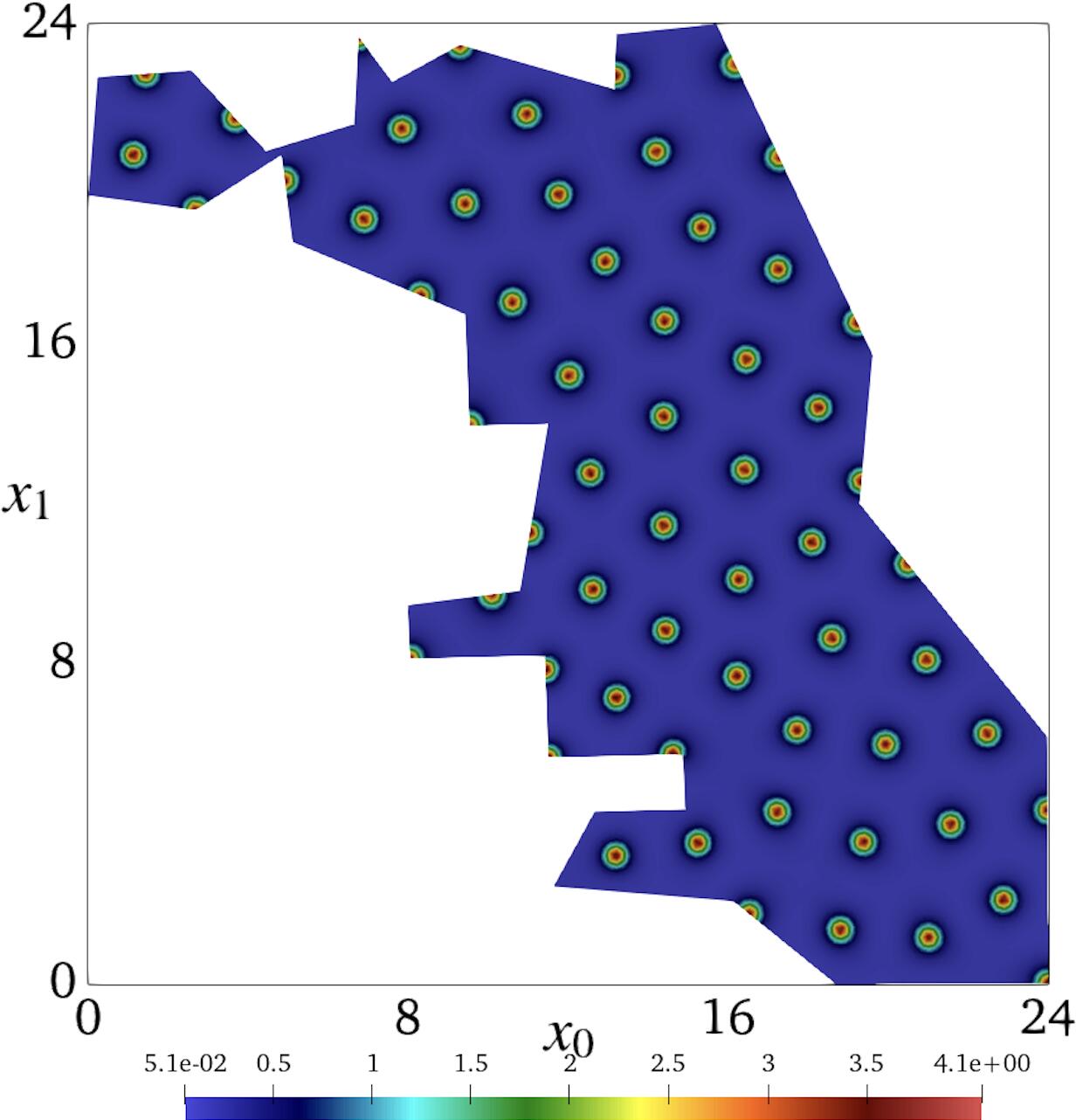}
    \end{overpic}
    \caption{Chicago without highway I-90: Evolution of the level of attractiveness (top) and density (bottom).}
    \label{fig:chi}
\end{figure}

Next, we consider a more realistic domain, with the same time interval of interest (i.e., $t \in (0,200]$). Since the greater Houston are is mostly a circle, we take the geometry of Chicago \cite{Chicago}, which is more interesting. 
This domain is meshed with 11218 nodes and 22466 elements. The time step is set to $\Delta t = 1/50$. Before adding highway I-90, we run a simulation with spatially homogeneous parameters: $\eta = 0.06$, $A^{st} = 1/30$, $\Gamma \theta /\omega^2 = {B}_0 = 1$, and   ${\rho}_0 = 0.8$. Uniform Gaussian noise is added  to the
initial data, with $\chi_\delta(\bm{x}) = 1$, $\xi_B(\bm{x}) \sim \mathcal{N}(0, 0.05^2)$ and $\xi_\rho(\bm{x}) \sim \mathcal{N}(0, 0.01^2)$ in \eqref{eq:ic1}–\eqref{eq:ic2}, sampled independently at each spatial point $\bm{x}$. 
Fig.~\ref{fig:chi} displays the evolution of $A$ and $\rho$  given by the PDE solver with $\texttt{tol}_1 = \texttt{tol}_2 = 10^{-9}.$ As expected from the low value of $\eta$, we see the formation of persistent small hotspots. 
The averaged number of iterations  over time required by the iterative partitioned algorithm is 5.46, which is in line with the number reported in Tab.~\ref{tab:it}
for a similar case (case 3, $(16/200,1/50)$, $10^{-9}$).

\begin{figure}[htb!]
     \centering
         \begin{overpic}[percent,width=0.19\textwidth]{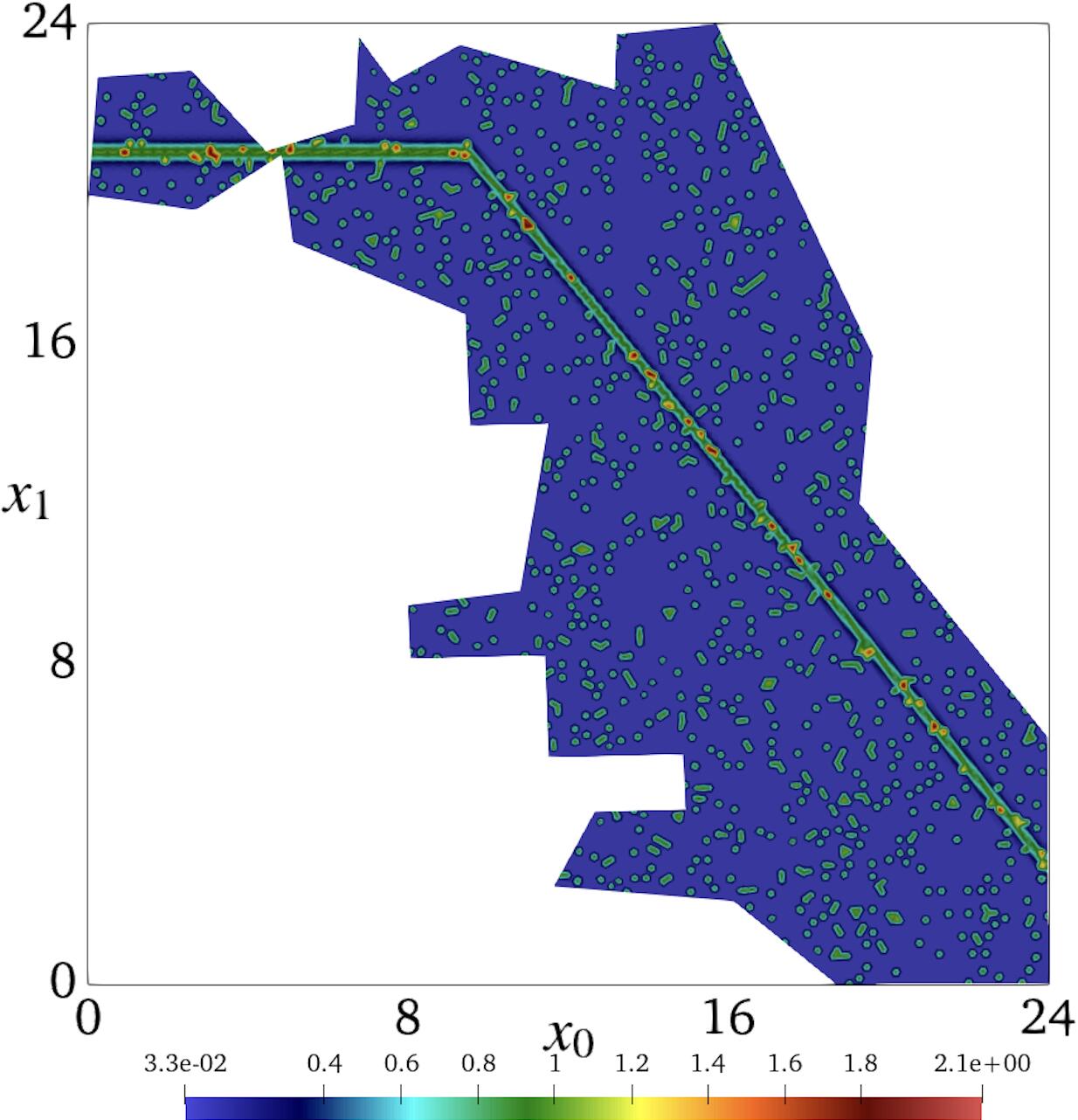}
    \put(30,103){\footnotesize{$t = 0$}}
    \put(12,35){\footnotesize{$A$}}
    \end{overpic} 
     \begin{overpic}[percent,width=0.19\textwidth, grid=false]{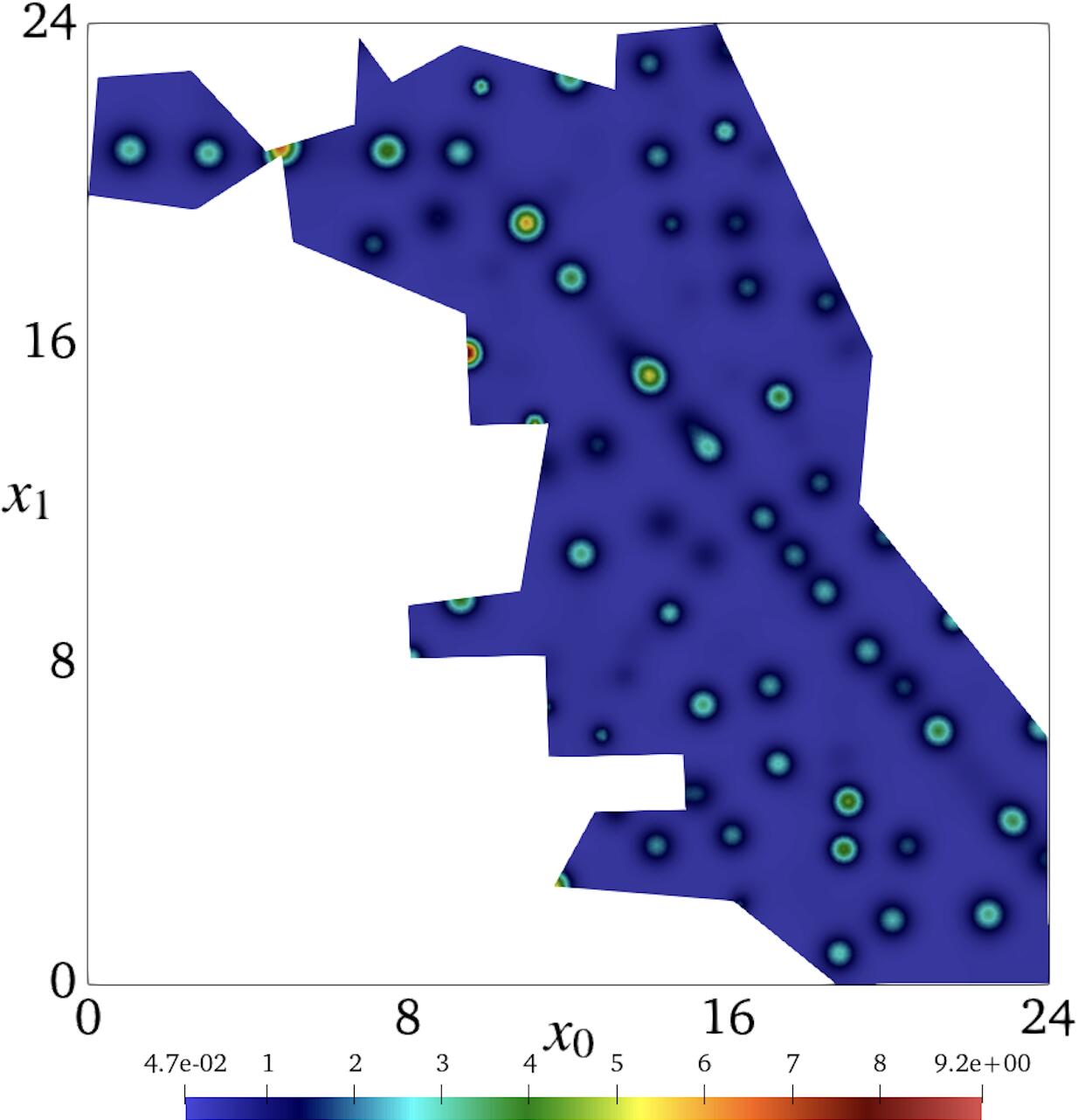}
    \put(35,103){\footnotesize{$t = 2$}}
    \end{overpic}
    \begin{overpic}[percent,width=0.19\textwidth]{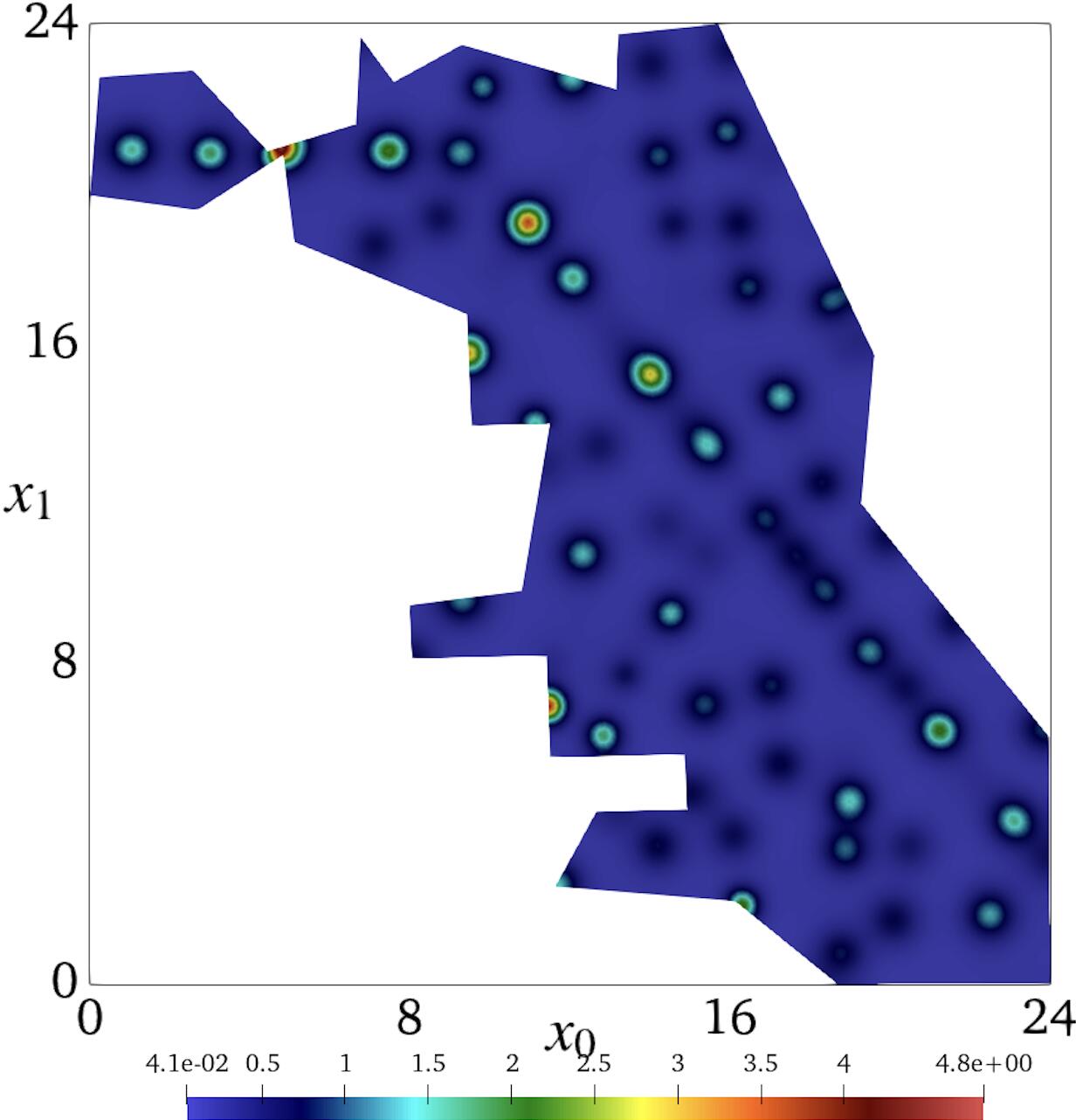}
    \put(35,103){\footnotesize{$t = 3$}}
    \end{overpic} 
        \begin{overpic}[percent,width=0.19\textwidth]{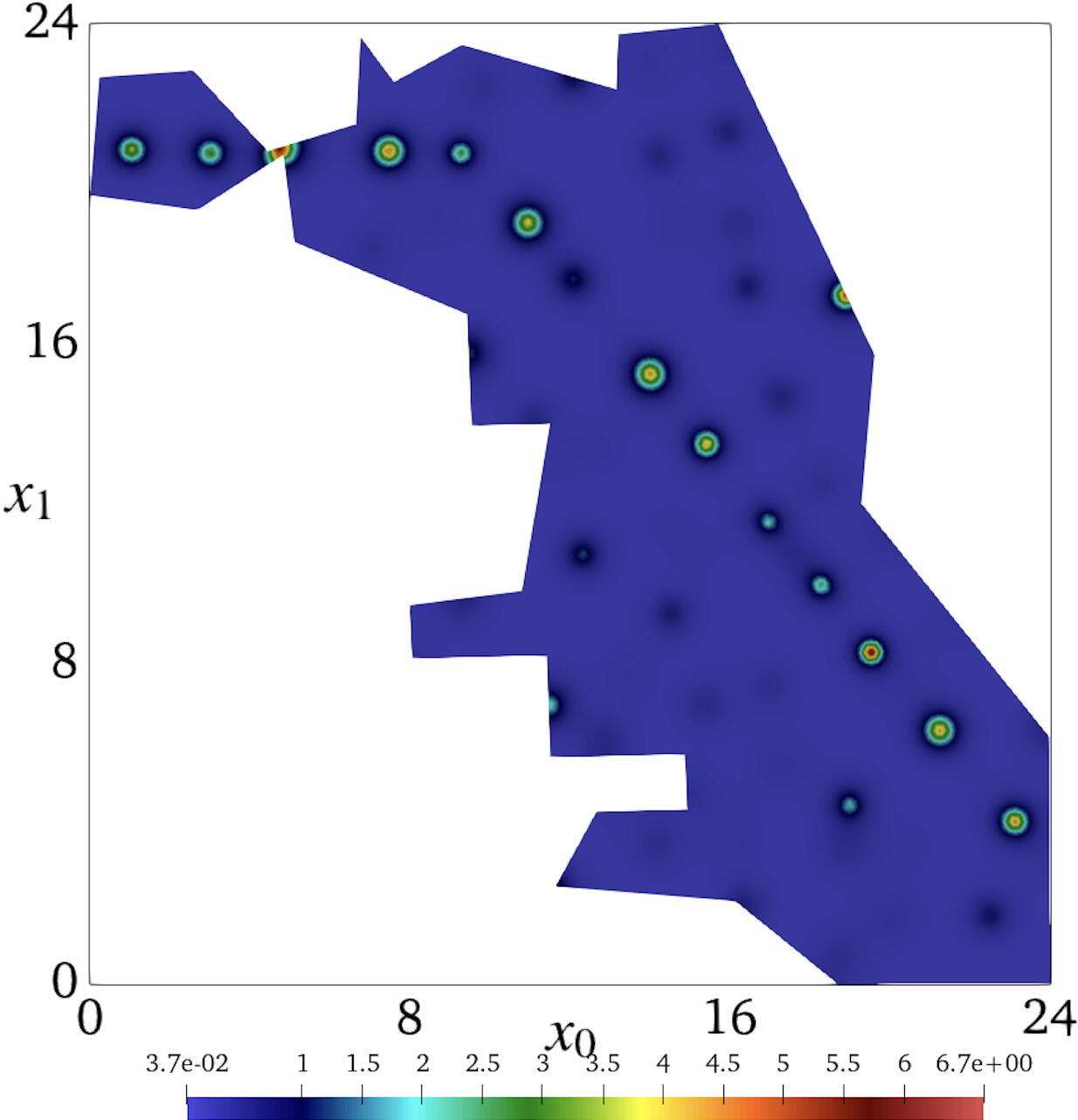}
    \put(32,103){\footnotesize{$t = 5$}}
    \end{overpic} 
    \begin{overpic}[percent,width=0.19\textwidth]{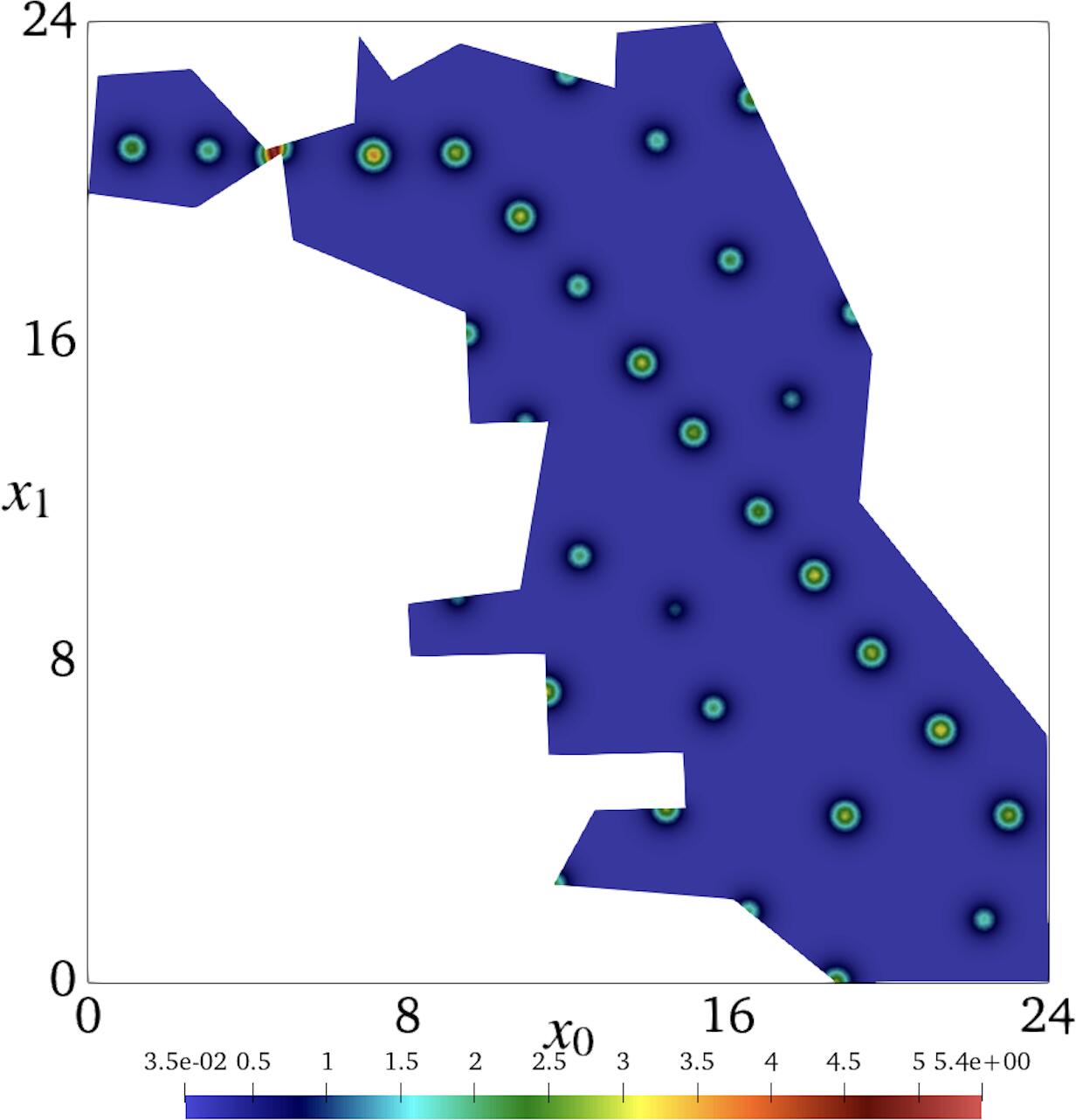}
    \put(30,103){\footnotesize{$t = 200$}}
    \end{overpic}
    \\
         \begin{overpic}[percent,width=0.19\textwidth]{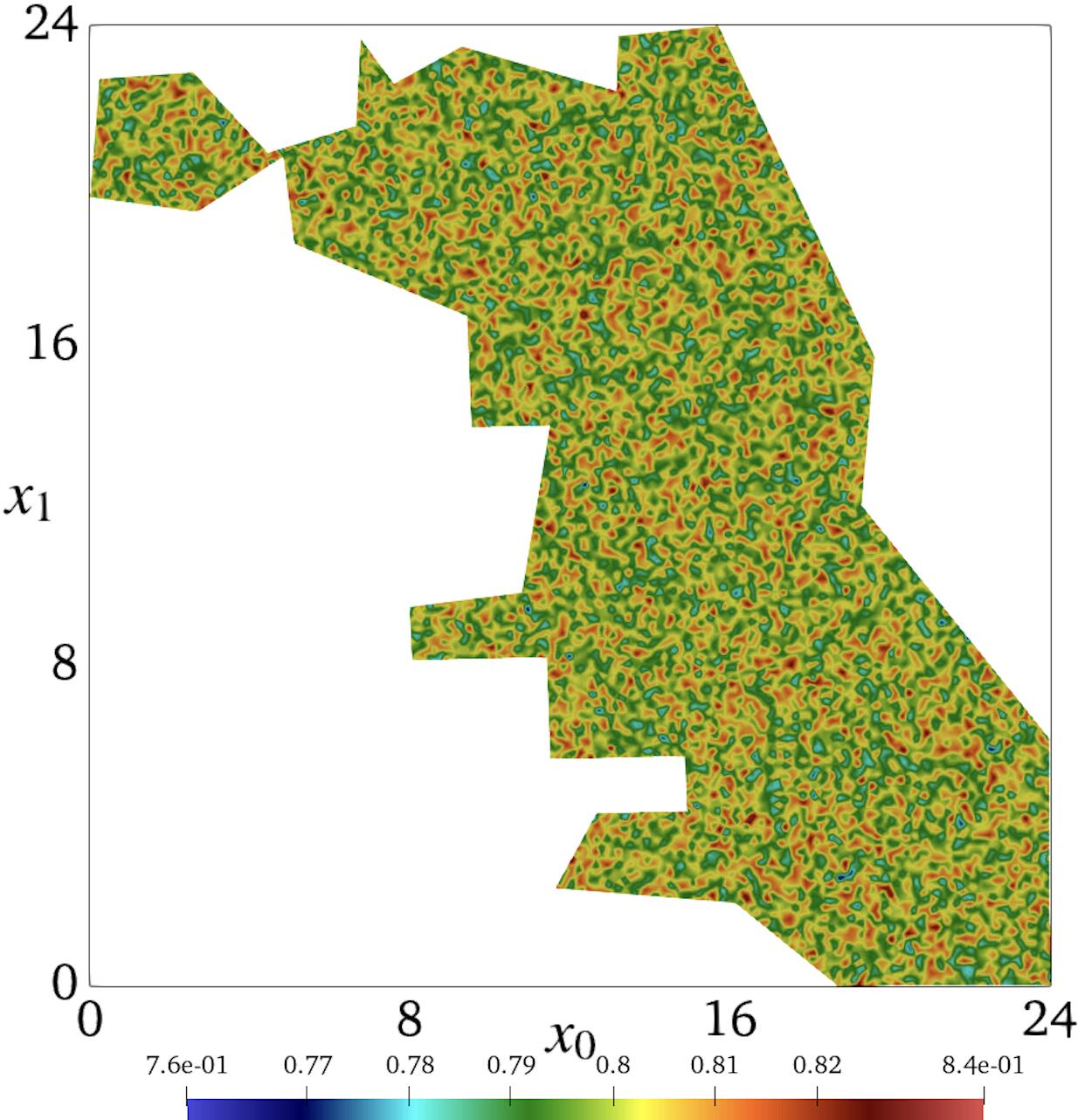}
    \put(12,35){\footnotesize{$\rho$}}
    \end{overpic} 
     \begin{overpic}[percent,width=0.19\textwidth, grid=false]{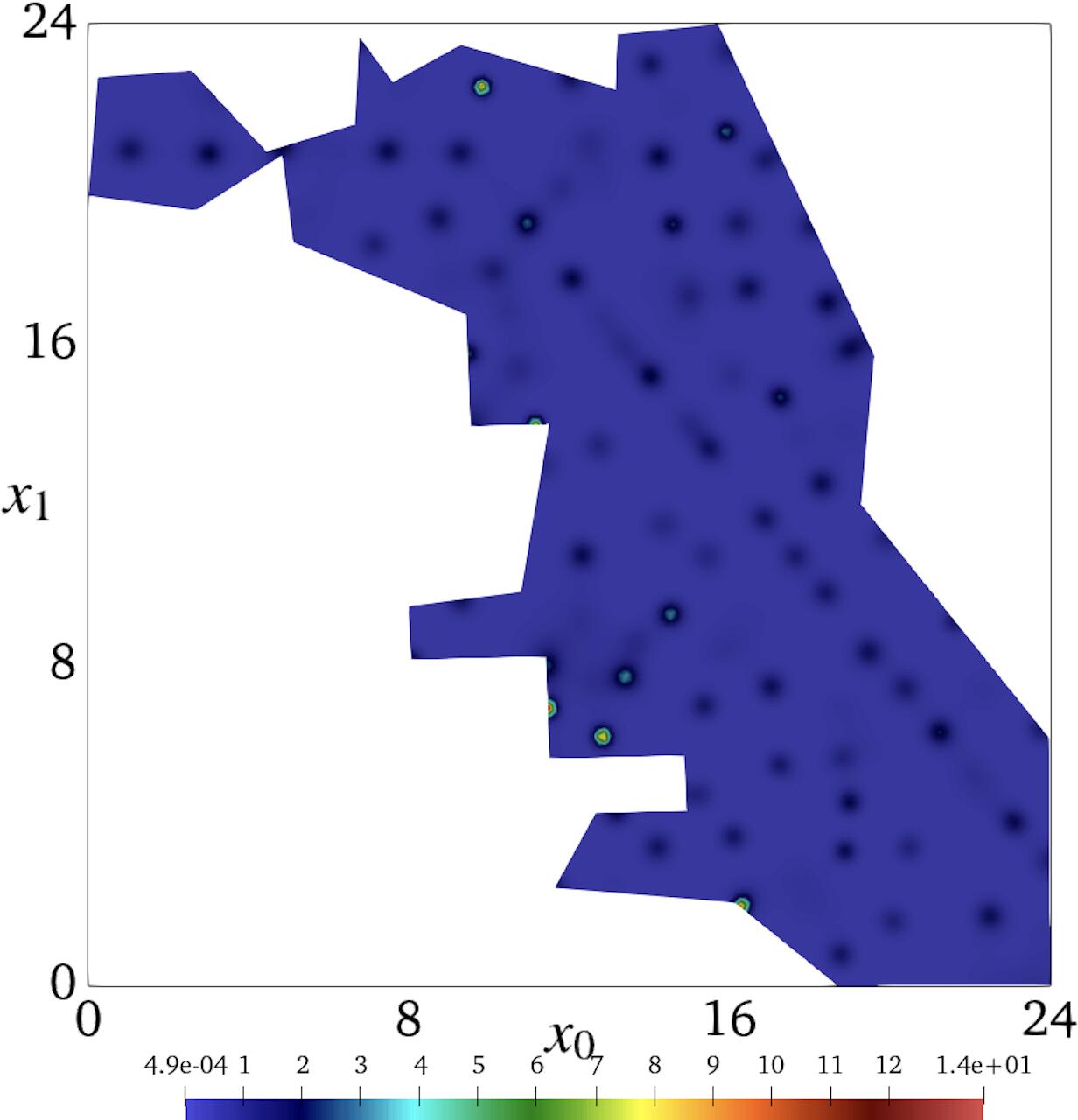}
    \end{overpic}
    \begin{overpic}[percent, width=0.19\textwidth]{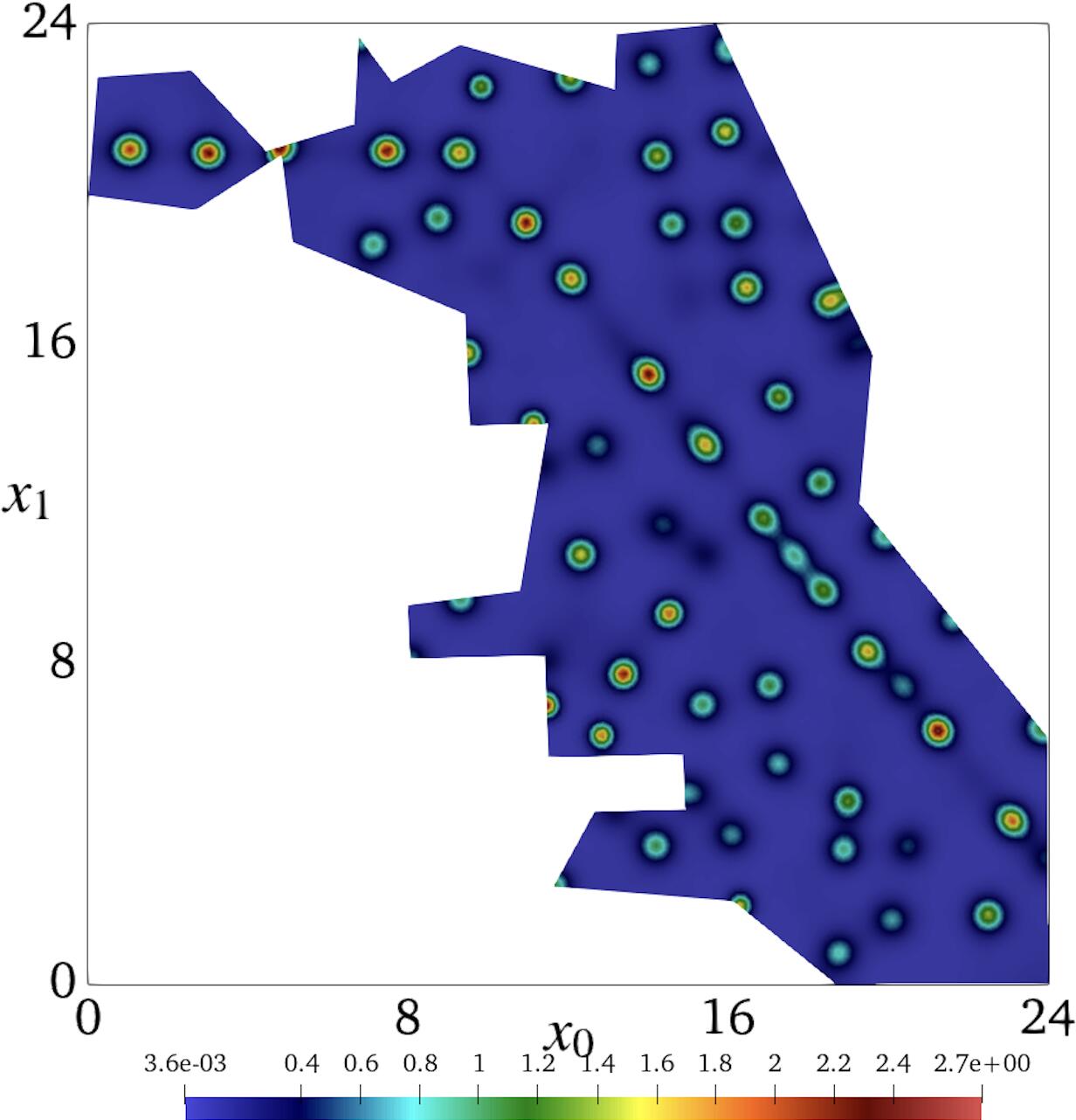}
    \end{overpic} 
        \begin{overpic}[percent,width=0.19\textwidth]{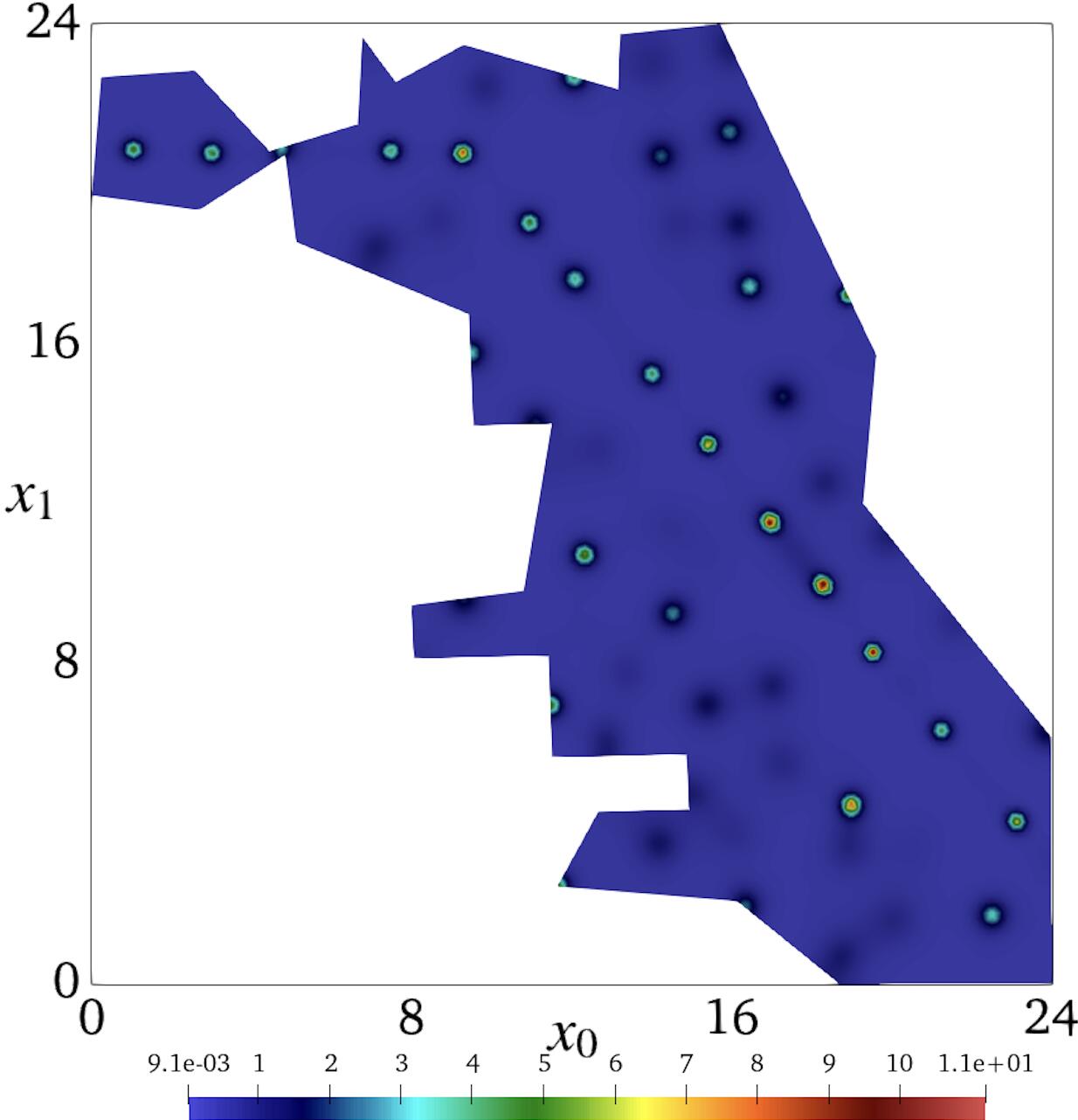}
    \end{overpic} 
    \begin{overpic}[percent,width=0.19\textwidth]{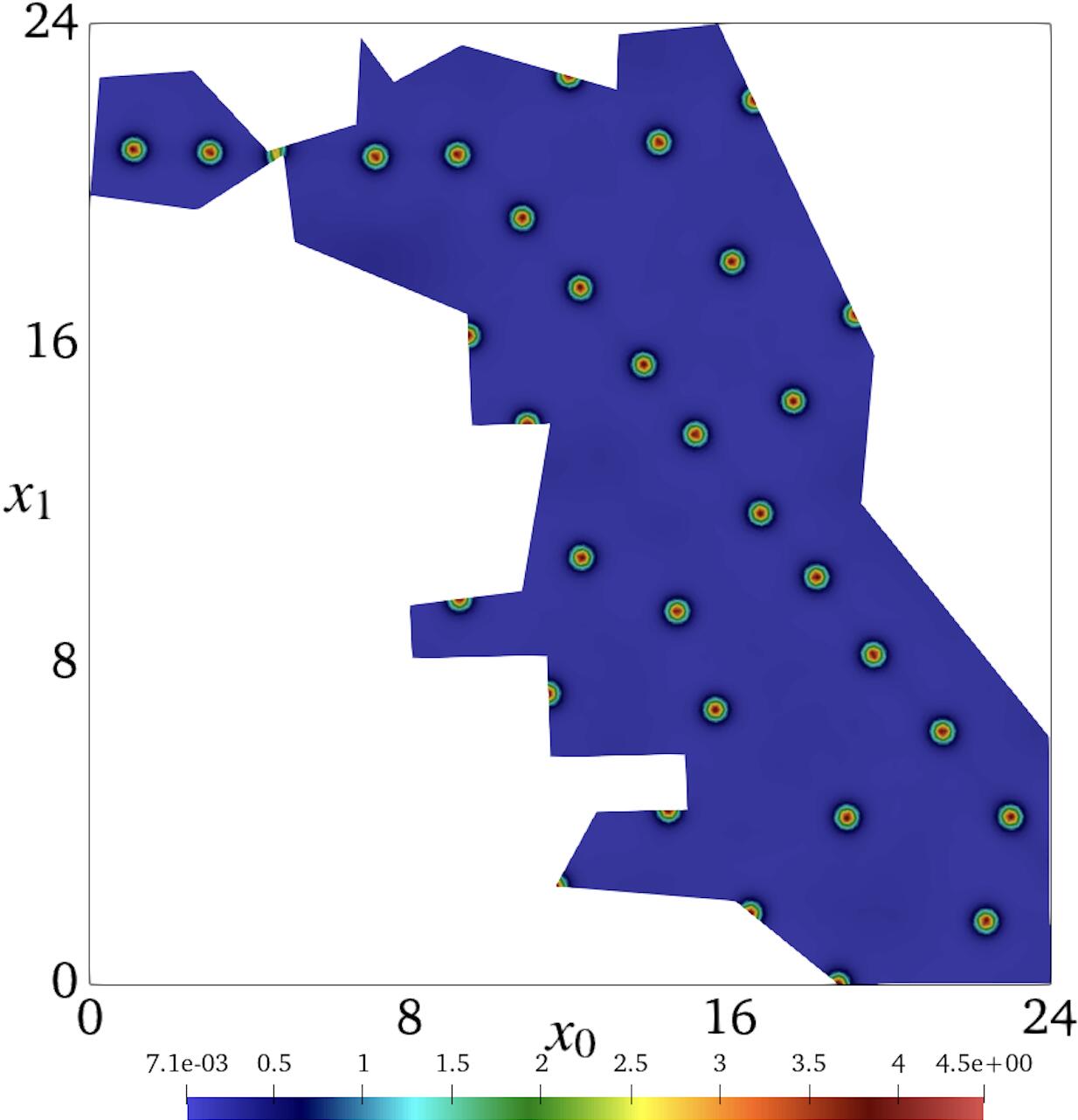}
 
    \end{overpic}
    \caption{Chicago with highway I-90: Evolution of the level of attractiveness (top) and density (bottom).}
    \label{fig:chi-hw}
\end{figure}

In order to add highway I-90, we embed the geometry of Chicago in $[0,24] \times [0,24]$. See Fig.~\ref{fig:chi-hw}.
We keep the same $\eta$ and $\rho_0$ as in the case with no highway, but make $A^{st}$ and $\Gamma$ dependent on $\bm{x}$ as follows:
\begin{align}
 & A^{st}(\bm{x}) = \begin{cases}
\frac{1}{30} \cdot \chi_{1} (\bm{x}) + \frac{1}{30}\cdot\exp\left\{ -20* \left( 20.8 - x_1 \right)^2 \right\}, & \text{if } x_0 < 9.5, \\
\frac{1}{30} \cdot \chi_{1} (\bm{x}) + \frac{1}{30}\cdot\exp\left\{ -20 *\left( 32.4923 - \dfrac{16}{13} x_0 - x_1 \right)^2 \right\}, & \text{if } x_0 \ge 9.5.
\end{cases} \cl
 & \frac{\Gamma(\bm{x}) \theta}{\omega^2} = B_0(\bm{x}) = \begin{cases}
1 \cdot \chi_{0.1} (\bm{x}) + \exp\left\{ -20* \left( 20.8 - x_1 \right)^2 \right\}, & \text{if } x_0 < 9.5, \\
1 \cdot \chi_{0.1} (\bm{x}) + \exp\left\{ -20 * \left( 32.4923 - \dfrac{16}{13} x_0 - x_1 \right)^2 \right\}, & \text{if } x_0 \ge 9.5.
\end{cases} \el
\end{align}
Note that to create the highway we have followed the same strategy used in the simple case (square geometry). Gaussian noise is added  to the
initial data, with $\chi_{\delta_B}(\bm{x}) = 0.1$, $\chi_{\delta_\rho}(\bm{x}) = 1$, $\xi_B(\bm{x}) \sim \mathcal{N}(0, 0.05^2)$ and $\xi_\rho(\bm{x}) \sim \mathcal{N}(0, 0.01^2)$ in \eqref{eq:ic1}–\eqref{eq:ic2}, sampled independently at each spatial point $\bm{x}$. 
Fig.~\ref{fig:chi-hw} shows the evolution of $A$ and $\rho$ given by the PDE solver with $\texttt{tol}_1 = \texttt{tol}_2 = 10^{-9}.$ We observe the formation of small persistent hotspots aligned with the highway, with smaller transient hotspots in directions parallel to
the highway.
The averaged number of iterations  over time required by the iterative partitioned algorithm is 6.39, which is slightly higher than the case with no highway.

In this paper, we have shown the evolution of the variables in the continuum model. 
However, from a practical point of view, one is interested in the number of burglaries over a certain period of time, which is the quantity shown 
in Fig.~\ref{fig:Houston}. Such number of b
comes from the integral of 
$E_s(t)$ in \eqref{eq:Es} over 
the given time interval and it 
requires dimensional variables. This means that specific values of parameters $\theta$ and $\omega$
need to be assigned, as explained below, while the PDE
model \eqref{eq:continous_bertozzi-1}-\eqref{eq:continous_bertozzi-2} 
needs only the value of $\Gamma \theta/ \omega^2$. 

Let us assume that we are interested in the number of burglaries $n_B$ occurred over each element $Q_k \in \mathcal{T}_h $
during the entire interval
$(0,T)$:
\begin{equation*}
    n_B |_{Q_k} = \int_0^T E_s(t) |_{Q_k} dt \approx \sum_{n = 0}^{N_t-1} E_s(t^n) |_{Q_k} \Delta t
    = \sum_{n = 0}^{N_t-1} n_s(t^n) |_{Q_k} p_s(t^n) |_{Q_k} \Delta t
\end{equation*}
Quantities $ n_s(t^n) |_{Q_k}$ and $p_s(t^n) |_{Q_k}$
can be obtained from the finite element approximation of the continnum model:
\begin{equation*}
    n_s(t^n) |_{Q_k} = \int_{Q_k} \epsilon L^2 \rho^{n}_h d\bx \approx \frac{\theta h^2}{\omega} \int_{Q_k} \rho^{n}_h d\bx, \quad p_s(t^n) |_{Q_k} = 1-e^{-A_h^{n}|_{Q_k} \Delta t}.
\end{equation*}
Thus, $ n_B |_{Q_k}$ could be obtained from post-processing the finite element solution of the 
continnum model as follows:
\begin{equation*}
    n_B |_{Q_k}  \approx \frac{\theta h^2 \Delta t}{\omega}  \sum_{n = 0}^{N-1} \left( 1-e^{-A_h^{n}|_{Q_k} \Delta t} \right)  \int_{Q_k} \rho^{n}_h d\bx. 
\end{equation*}
This shows that ratio $\Gamma \theta/ \omega^2$ is not enough to obtain $n_B$ from the computed 
level of attractiveness and density
of criminals. The ratio $\theta/\omega$ needs to be provided as well.

\section{Future perspectives}\label{sec:future}

\anna{We believe the work presented in this paper is a stepping stone towards the application of the PDE model to predict
residential burglary in real cities. 
Imminent future work will be focused on making the simulations more and more realistic by, e.g., considering more realistic geometries, 
incorporating real data on property (see \cite{Woodworth2024}), adding to the model police intervention and community investment, and
developing data-driven methods for parameter estimation.
{As shown in Sec.~\ref{sec:num_res}, solutions are
sensitive to key parameters $\eta$, $ A^{st}$, and $\Gamma \theta /\omega^2$.}
However, as mentioned in Sec.~\ref{sec:intro}, the scale of the PDE model (i.e., the mascroscopic scale) is only one of three possible scales of observation and complex self-organizing, collective dynamics, such crime,
are best represented by a mathematical description covering all scales.}

Systems like crime are made of living entities (as opposed to inert matter) that possess the behavioral ability to
develop and adapt moving strategies in response to the context, which significantly complicates their mathematical modeling.
A general mathematical structure able to capture the specific complexity features of crime as a living system at all scales is still missing. While there exists a seizable body of literature on microscopic and macroscopic models (see Sec.~\ref{sec:intro}), the mesoscopic scale has received much less attention. 
Taking inspiration from the kinetic theory of gases,
mesoscopic (or kinetic) models derive a Boltzmann-type evolution equation for a statistical distribution function. In these models, people
are seen as {active} particles, as opposed
to classical particles in gases. This leads
to a key difference: the interactions in 
kinetic models for gas dynamics are
conservative and reversible, while the interactions in the kinetic models for living systems are irreversible, 
non-conservative and, in some cases, nonlocal and nonlinearly additive. 
The reader interested in the kinetic theory
of active particles is referred to \cite{BellomoBurini2021,BELLOMO2026}. 
{In general, Ref.~\cite{BELLOMO2026} provides the mathematical tools necessary for the approaches outlined in this section.}
To enable a \emph{multiscale framework}, models must be derived at 
each scale using the same principles and similar parameters.
It is worth noticing that substantial model development at the three different scales, even
with the same principles and analogous parameters
for an integrated multiscale vision, 
has been done for a different kind of collective social system: crowd dynamics \cite{Bellomo2022,Bellomo2023}. So, existing theory
for crowd dynamics could be extended to crime dynamics.

\anna{
The concept of \emph{activity variable} mentioned in Sec.~\ref{sec:intro} 
belongs to the kinetic theory of active particles. It is a 
vector variable at the microscopic scale that collects a number of behavioral variables, such as emotional state, and changes through interactions with others.
In the context of crowd dynamics, several works include the activity variable and some model its 
dynamics through a \emph{transport equation} \cite{Bellomo2022}.
With the exception of \cite{Bellomo2015}, the concept of 
activity variable remains unexplored for models of crime dynamics. 
We believe the activity variable, and corresponding transport equation, is 
the key to modeling behavioral heterogeneity in criminals and policemen. Ideally, following H.~A.~Simon's theories, models 
should account for rules of interactions that depend on the environment where interactions occur, with the 
environment in turn evolving as a result of human intervention \cite{Simon2019,BellomoEgidi2024,BELLOMO2026}.
}

\anna{Finally, we would like to mention the importance of developing ``\emph{multi-physics}'' models that couple
crime dynamics to other factors correlated with crime, such as
wealth inequalities \cite{ITSKOVICH2023102065} and
level of education \cite{Baron2024}. 
Ideas on how to do this at the macroscopic level are discussed in \cite{BellomoWinkler2022}.
It is also imperative to think about safety checks for these models, and constant improvement through fact-based evidence, to avoid turning them into 
weapons of math destruction \cite{ONeil2016}. }

\section{Conclusions}\label{sec:concl}

We introduced a finite element approximation of a well-known PDE model to predict residential burglary. For the numerical solution of the discretized PDE problem, we proposed an algorithm that decouples the computation of the model variables
(attractiveness  of residential sites and density of criminals) at each time step, resulting in the solution of two linear algebraic systems per iteration.
We validated the numerical results
obtained with our approach against the results provided by the agent-based model from which the PDE model was derived. We showed that the PDE solver output
matches well the output
from the agent-based simulations when the density of criminals is large, i.e., when the criminal population can be looked at as a continuum versus separate individuals. Through numerical experimentation, 
we have characterized the dependence of both diameter and number of crime hotspots
on the strength of the neighborhood effects. 
Concerning the solution algorithm,
we showed that the average number of iterations over time required to meet the stopping criterion ranges between 1 and 6 when model parameters, discretization and algorithmic parameters are varied by orders of magnitude. 
These low numbers of iterations indicate robustness to parameter variations and numerical efficiency, since only two linear systems have to be solved per iteration. Thanks to 
the flexibility allowed by the finite element method, we simulated crime dynamics with spatially heterogeneous model parameters and a realistic geometry.

\anna{Finally, while the robust and efficient code developed for this paper has a value per se, it is mostly intended as the solid base for the broader research program described in Sec.~\ref{sec:future}.}


\section*{Acknowledgments}
We thank Maria Lozano for the visualization of the 2024 Houston police data in Fig.~\ref{fig:Houston}. Zhong is supported by NSF-AoF grant $\#2225507$. Zhong also gratefully acknowledges funding provided by the Oak Ridge Associated Universities Ralph E. Powe Junior Faculty Enhancement Award for FY$2024$.

\end{document}